\documentclass[aos]{imsart}
\usepackage{amsthm,amsmath,natbib,amssymb,graphicx,color,marvosym}
\bibpunct{(}{)}{;}{a}{,}{,}

\startlocaldefs
\newtheorem{theorem}{Theorem}
\newtheorem{corollary}{Corollary}

\newtheorem{lemma}{Lemma}
\let\hat\widehat
\let\tilde\widetilde
\def\E{\mathbb{E}}
\def\P{\mathbb{P}}

\def\given{\,|\,}
\def\H{\mathcal{H}}
\def\Var{\text{Var}}
\def\reals{\mathbb{R}}

\def\H{\mathcal{H}}
\def\cH{{\cal H}}

\def\qed{\hskip1pt $\;\;\scriptstyle\Box$}
\def\c3po{$\text{C3PO}$}

\def\P{{\mathbb P}}
\def\E{{\mathbb E}}
\def\reals{{\mathbb R}}
\def\ip#1#2{\left\langle #1 ,#2\right\rangle}

\def\H{{\mathcal H}}
\def\M{{\mathcal M}}
\def\S{{\mathcal S}}

\def\ds{\displaystyle}
\def\argmin{\mathop{\text{arg\,min}}}

\def\Var{\textrm{Var}}
\def\mytop{{\scriptscriptstyle\top}}

\def\hs{}
\def\Res{R}
\def\Proj{{P}}

\def\supp{{\text{supp}}}
\def\ip#1#2{\left\langle #1, #2\right\rangle}

\def\Sc{S^c}
\def\df{\mathop{\text{df}}}
\def\GCV{\mathop{\text{GCV}}}

\def\bm#1{{#1}}
\def\onen{{\textstyle{\frac{1}{n}}}}
\def\L{{\mathcal L}}
\def\comment#1{}
\def\iff{\Leftrightarrow}
\endlocaldefs

\definecolor{Gray}{rgb}{.6,.6,.6}
\begin{document}

\begin{frontmatter}

\title{Sparse Additive Models}
\runtitle{Sparse Additive Models}

\author{Pradeep~Ravikumar, John Lafferty, Han Liu and~Larry~Wasserman}\ead[label=e1]{pradeepr@cs.cmu.edu}
\address{Carnegie Mellon University\\5000 Forbes Ave\\ Pittsburgh, PA 15213} %\printead{e1}}
\affiliation{Carnegie Mellon University}
\begin{center}
\today
\end{center}

\runauthor{Ravikumar, Lafferty, Liu, and Wasserman}

\begin{abstract}
  We present a new class of methods for high-dimensional nonparametric
  regression and classification called sparse additive models (SpAM).
  Our methods combine ideas from sparse linear modeling and additive
  nonparametric regression.  We derive an algorithm for fitting the
  models that is practical and effective even when the number of
  covariates is larger than the sample size.  SpAM is essentially a
  functional version of the grouped lasso of Yuan and Lin (2006). SpAM
  is also closely related to the COSSO model of \cite{LZ:Cosso}, but
  decouples smoothing and sparsity, enabling the use of arbitrary
  nonparametric smoothers. We give an analysis of the theoretical
  properties of sparse additive models, and present empirical results
  on synthetic and real data, showing that SpAM can be effective in
  fitting sparse nonparametric models in high dimensional data.
\end{abstract}

\begin{keyword}
\kwd{lasso}
\kwd{nonparametric regression}
\kwd{additive model}
\kwd{sparsity}
\end{keyword}
\end{frontmatter}

\section{Introduction}
\label{sec:intro}

Substantial progress has been made recently on 
the problem of fitting high dimensional
linear regression models of the form
$Y_i = X_i^T\beta + \epsilon_i$, for $i=1,\ldots, n$.
Here $Y_i$ is a real-valued response,
$X_i$ is a predictor
and $\epsilon_i$ is a mean zero error term.
Finding an estimate of $\beta$ when $p > n$
that is both statistically well-behaved and computationally efficient 
has proved challenging; however, under the assumption that the vector
$\beta$ is sparse, the lasso estimator (\cite{Tibshirani:96}) has been
remarkably successful.
The lasso estimator $\hat\beta$ minimizes the $\ell_1$-penalized sum of squares
%\begin{equation}
%\label{eq:lasso}
$\sum_i(Y_i - X_i^T\beta)^2 + \lambda\sum_{j=1}^p |\beta_j|$
%\end{equation}
with the $\ell_1$ penalty $\|\beta\|_1$ encouraging sparse solutions,
where many components $\hat\beta_j$ are zero.
The good empirical success of this estimator has been recently
backed up by results confirming that it has strong theoretical
properties; see \citep{Bunea2:07,GR:04,ZY07,MY06,Wai06}.

The nonparametric regression model
$Y_i = m(X_i) + \epsilon_i$, 
where $m$ is a general smooth function, relaxes the strong
assumptions made by a linear model, but 
is much more challenging in high dimensions.
\cite{Hast:Tibs:1999} introduced
the class of additive models
of the form
\begin{equation}
\label{eq:additive}
Y_i = \sum_{j=1}^p f_j(X_{ij})+\epsilon_i .
\end{equation}
This additive combination of univariate functions---one for each 
covariate $X_j$---is less general than joint multivariate nonparametric models, but can be 
more interpretable and easier to fit;
in particular, an additive model can be estimated using a coordinate descent 
Gauss-Seidel procedure, called backfitting.  
Unfortunately, additive models only have good statistical
and computational behavior when the number of variables $p$ is not
large relative to the sample size $n$, so their usefulness is
limited in the high dimensional setting.

In this paper we investigate sparse additive models (SpAM), which
extend the advantages of sparse linear models to the additive, nonparametric 
setting.  The underlying model is the same as in \eqref{eq:additive}, 
but we impose a sparsity constraint
on the index set $\{j: f_{j} \not\equiv 0\}$ of 
functions $f_j$ that are not identically zero. \cite{LZ:Cosso} have proposed COSSO, an extension 
of lasso to this setting, for the case where the component functions $f_j$ belong to 
a reproducing kernel Hilbert space (RKHS). 
They penalize the sum of the RKHS norms of the component functions.
\citet{Yuan:07} proposed an extension of  
the non-negative garrote to this setting.  As with the parametric 
non-negative garrote, the success of this method depends on the initial estimates of 
component functions $f_j$.

In Section~\ref{sec:backfitting},
we formulate an optimization problem in the population  
setting that induces sparsity.
Then we derive a sample version of the solution.
The SpAM estimation procedure we introduce allows the use of arbitrary nonparametric smoothing techniques,
effectively resulting in a combination of the lasso and backfitting.
The algorithm extends
to classification problems using generalized additive models.
As we explain later, SpAM can also be thought of as a functional version
of the grouped lasso \citep{Yuan:Lin:06}.

\comment{
let $(Y_1,X_1),\ldots, (Y_n,X_n)$ be $n$ data points
where $X_i=(X_{i1},\ldots, X_{ij}, \ldots, X_{ip})^T\in\mathbb{R}^p$ and $Y_i\in \mathbb{R}$.
We form an additive model 
\begin{equation}
\label{eq:am}
Y_i = \alpha + \sum_{j=1}^p \beta_j g_j(X_{ij})+\epsilon_i
\end{equation}
with the identifiability conditions that the component functions
have mean zero and norm one: 
$\int g_j(x_j) d P(x_j) = 0$ and 
$\int g_j^2(x_j) d P(x_j) = 1$.
Further, we impose the sparsity condition
$\sum_{j=1}^p | \beta_j| \leq L_n$
and the smoothness condition
$g_j\in {\cal T}_j$
where
${\cal T}_j$ is some class of smooth functions.  
While this problem
is not convex, it makes clear the way in which sparsity is encouraged,
through the $\ell_1$ penalty
$\sum_{j=1}^p | \beta_j| \leq L_n$ which induces sparsity.  Below, we 
derive an alternative formulation that is convex.
This approach is closely related to the 
COSSO, introduced by \cite{LZ:Cosso}, in which 
a regression function $m(x)$
is assumed to be a sparse linear combination of smooth functions.
However, \cite{LZ:Cosso} put a sparsity constraint on
the second derivatives of the $g_j$.  Our formulation of sparse
additive models allows the use of general smoothing operators,
not only smoothing splines.
As we explain later, SpAM can also be though of as a functional version
of the grouped lasso \citep{Yuan:Lin:06}.
}

The main results of this paper include the formulation of a convex optimization
problem for estimating a sparse additive model, an efficient backfitting
algorithm for constructing the estimator, and theoretical results that analyze
the effectiveness of the estimator in the high dimensional setting.  Our
theoretical results are of several different types.  First, we show that, under
suitable choices of the design parameters, the SpAM backfitting algorithm
recovers the correct sparsity pattern asymptotically; this is a property we call
\textit{sparsistence}, as a shorthand for ``sparsity pattern consistency.''
Second, we show that that the estimator is \textit{persistent}, in the sense of
\cite{GR:04}, which is a form of risk consistency.  
Specifically, we show:

\def\ds{\displaystyle}
\begin{center}
\begin{tabular}{lll}
&&\\
Sparsistence: & $\ds \P\left(\hat{S} = S\right)\rightarrow 1$ &
if $p_n =
O(e^{n^\xi})$, for $\xi < \frac{3}{5}$ \\[15pt]
Persistence: & $\ds R(\hat m_n) - \inf_{h\in \M_n} R(h) \stackrel{P}{\rightarrow} 0$ &
if $p_n = O(e^{n^\xi})$, for $\xi < 1.$ \\[15pt]
\end{tabular}
\end{center}

Here $S = \{j\,:\, f_j \neq 0\}$ is the index set for the nonzero
components, 
$\hat{S} = \{j\,:\, \hat{f}_j \neq 0\}$
and $\M_n$ is a class of functions defined by the level of
regularization.

In the following section we establish notation and assumptions.  In
Section~\ref{sec:backfitting} we formulate SpAM as an optimization
problem and derive a scalable backfitting algorithm.  An extension to
sparse nonparametric logistic regression is presented in
Section~\ref{sec:logistic}.  Examples showing the use of our sparse
backfitting estimator on high dimensional data are included in
Section~\ref{sec:examples}.  In Section~\ref{sec:sparsistency} we
formulate the sparsistency result, when orthogonal function regression
is used for smoothing.
In Section~\ref{sec:persistence} we give the persistence result.  
Section~\ref{sec:discussion} contains
a discussion of the results and possible extensions.
Proofs are contained in Section~\ref{sec:proofs}.

\section{Notation and Assumptions}
\label{sec:assumptions}

We assume that we are given data
$(X_1,Y_1), \ldots, (X_n,Y_n)$ 
where $X_i = (X_{i1}, \ldots, X_{ij},\ldots, X_{ip})^T \in [0,1]^p$
and
\begin{equation}
Y_i = m(X_i) + \epsilon_i
\end{equation}
with $\epsilon_i \sim N(0,\sigma^2)$
and
\begin{equation}
m(x) =\sum_{j=1}^p f_j(x_j).
\end{equation}
Denote the joint distribution of $(X_i,Y_i)$ by $P$. For a function 
$f$ on $[0,1]$ denote its $L_{2}(P)$ norm by
\begin{equation}
\|f\| = \sqrt{\int_0^1 f^2(x) dP(x)}  = \sqrt{\E (f)^{2}} .
\end{equation}

For $j \in \{1,\ldots, p\}$, let $\H_{j}$ denote the Hilbert subspace of
$L_{2}(P)$, of $P-$measurable functions $f_j(x_j)$ of the single scalar variable $X_j$
with zero mean, $\E(f_j(X_j)) = 0$. Thus, $\H_j$ has the inner product
\begin{equation}
\ip{f_j}{f'_j} = \E\left(f_j(X_j) f'_{j}(X_j)\right)
\end{equation}
and $\|f_j\| = \sqrt{\E(f_j(X_j)^2)} < \infty$.
Let $\H  = \H_{1} \oplus \H_{2} \oplus \ldots \oplus \H_{p}$ denote
the Hilbert space of functions of $(x_1,\ldots, x_p)$ that have
the additive form: $m(x) = \sum_{j} f_j(x_j)$, with $f_j \in \H_j, j = 1,\hdots, p$.

Let $\{\psi_{jk}, k=0,1,\ldots \}$ denote a uniformly
bounded, orthonormal basis with respect to Lebesgue measure on $[0,1]$.
Unless stated otherwise,
we assume that $f_j\in {\cal T}_j$ where
\begin{equation}
{\cal T}_j = \Biggl\{ f_j\in\H_j:\ 
f_j(x_j) =\sum_{k=0}^\infty \beta_{jk}\psi_{jk}(x_j),\ \ \ 
\sum_{k=0}^\infty \beta_{jk}^2 j^{2\nu_j}\leq C^2 \Biggr\}
\end{equation}
for some $0 < C < \infty$.
We shall take $\nu_j=2$ although the extension to other levels of
smoothness is straightforward.
It is also possible to adapt to $\nu_j$ although we do not
pursue that direction here.
 
Let $\Lambda_{\rm min}(A)$ and
$\Lambda_{\rm max}(A)$
denote the minimum and maximum eigenvalues of a square matrix $A$.
If $v=(v_1,\ldots, v_k)^T$ is a vector, we use the norms
\begin{equation}
\|v\| = \sqrt{\sum_{j=1}^k v_j^2},\ \ \ 
\|v\|_1 = \sum_{j=1}^k |v_j|,\ \ \ 
\|v\|_\infty = \max_j |v_j|.
\end{equation}

\section{Sparse Backfitting}
\label{sec:backfitting}

The outline of the derivation of our algorithm is as follows.  We
first formulate a population level optimization problem, and show that
the minimizing functions can be obtained by iterating through a series
of soft-thresholded univariate conditional expectations. We then plug
in smoothed estimates of these univariate conditional expectations, to
derive our sparse backfitting algorithm.

\vspace{10pt}

{\em Population SpAM.}
For simplicity, assume that $\E(Y_i)=0$.
The standard additive model optimization 
problem in $L_{2}(P)$ (the population setting) is
\begin{equation}
\min_{f_j\in\H_{j},\, 1\leq j\leq p} \; \E\left(Y - \textstyle\sum_{j=1}^p f_j(X_j)\right)^2 
\end{equation}
where the expectation is taken with respect to $X$ and
the noise~$\epsilon$.  Now consider the following modification of this
problem that introduces a scaling parameter for each function, and
that imposes additional constraints:
\begin{eqnarray}
\label{eq:nonconv}
\min_{\beta\in\reals^p, g_j\in\H_j}  && 
\E\left(Y - \textstyle\sum_{j=1}^p \beta_j g_j(X_j)\right)^2  \\ 
\text{subject to:} &&  \sum_{j=1}^p |\beta_j| \leq L,\\
&& \E\left(g_j^2\right) = 1,\; j=1,\ldots, p.
\end{eqnarray}
noting that $g_j$ is a function while 
$\beta = (\beta_1,\ldots, \beta_p)^T$ is a vector.
The constraint that $\beta$ lies in the $\ell_1$-ball
$\{\beta : \|\beta\|_1 \leq L\}$ encourages sparsity of the estimated $\beta$,
just as for the parametric lasso \citep{Tibshirani:96}.

%When $\beta$ is sparse, 
%the additive function 
%\begin{equation}
%m(x) = \sum_{j=1}^p f_j(x_j) = \sum_{j=1}^p \beta_j g_j(x_j)
%\end{equation}
%will also be sparse, meaning that
%many of the component functions $f_j(\cdot) = \beta_j g_j(\cdot)$ are
%identically zero.  
It is convenient
to re-express the minimization in the following equivalent form:
\begin{eqnarray}
\label{eq:optim}
\min_{f_j\in\H_j}  && \E\left(Y-\textstyle \sum_{j=1}^p f_j(X_j)\right)^2 \\
\text{subject to:} && \sum_{j=1}^p \sqrt{\E(f_j^2(X_j))} \;\leq\; L.
%% \E(f_j) = 0,\; j=1,\ldots, p.
\end{eqnarray}

The optimization problem in (\ref{eq:optim}) can also be written in the penalized
Lagrangian form,
\begin{equation}\label{eq::lag}
\L(f, \lambda) = \frac{1}{2} \E\left(Y-\textstyle \sum_{j=1}^p
  f_j(X_j)\right)^2  + \lambda \sum_{j=1}^p   \sqrt{\E(f_j^2(X_j))}.
\end{equation}

\begin{theorem}\label{thm::backfit}
The minimizers $f_j\in\H_j$ of
(\ref{eq::lag}) satisfy
\begin{equation}
\label{eq:soft}
f_j = \left[1 - \frac{\lambda}{\sqrt{\E(\Proj_j^2)}}\right]_+ \Proj_j\qquad \text{a.s.}
\end{equation}
where $[\cdot]_+$ denotes the positive part,
and $\Proj_j = \E[\Res_j\given X_j]$ denotes the projection of the residual
$\Res_j = Y - \sum_{k\neq j} f_k(X_k)$ 
onto $\H_j$.
\end{theorem}
The proof is given in Section~\ref{sec:proofs}.

At the population level, the $f_j$s can be found
by a coordinate descent procedure 
that fixes $(f_k:\ k\neq j)$ and fits $f_j$
by equation (\ref{eq:soft}),
then iterates over $j$.

\vskip20pt
{\em Data version of SpAM.}
To obtain a sample version
of the population solution,
we insert sample estimates into the
population algorithm,
as in standard backfitting \citep{Hast:Tibs:1999}.
Thus, we estimate the
projection $\Proj_j =\E(R_j\given X_j)$ by smoothing the residuals:
\begin{equation}
\hat \Proj_j = \S_j \Res_j
\end{equation}
where $\S_j$ is a linear smoother, such as a local linear or kernel smoother.
Let 
\begin{equation}
\label{eq:sjsimple}
\hat s_j = \frac{1}{\sqrt{n}} \|\hat P_j\| = \sqrt{\text{mean}(\hat P_j^2)} 
\end{equation}
be the estimate of $\sqrt{\mathbb{E}(P_j^2)}$.
Using these plug-in estimates
in the coordinate descent procedure
yields the SpAM backfitting algorithm given
in Figure~\ref{fig:backfitting:algo}.

This algorithm can be seen as a functional version of the coordinate
descent algorithm for solving the lasso.  In particular, if we solve
the lasso by iteratively minimizing with respect to a single
coordinate, each iteration is given by soft thresholding;
see Figure~\ref{fig:coordlasso}.  Convergence properties of variants
of this simple algorithm have been recently treated by
\cite{Daubechies:04,Daubechies:07}.  Our sparse backfitting algorithm
is a direct generalization of this algorithm, and it reduces to it in
case where the smoothers are local linear smoothers with large bandwidths.

\comment{As an alternative to estimating the conditional expectations in \eqref{eq::lag}
by smoothing, we can define estimators by minimizing
a sample version of the problem.
Thus, we would minimize
\begin{equation}
\frac{1}{n}\sum_{i=1}^n (Y_i - \sum_{j=1}^p f_j(X_{ij}))^2
\end{equation}
subject to $f_j \in {\cal T}_j$, and
\begin{equation}
\sum_{j=1}^p \sqrt{\frac{1}{n}\sum_{i=1}^n f_j^2(X_{ij})} \leq L,\ \ \ \ \ 
\frac{1}{n}\sum_{i=1}^n f_j(X_{ij}) =0,\ \ j=1, \ldots, p.
\end{equation}
Note that disregarding the functional constraints $f_j\in{\cal
  T}_j$, and optimizing only over the $np$ values
$f_j(X_i)$ leads to a finite dimensional convex optimization problem.
}

\begin{figure}[t]
{\sc SpAM Backfitting Algorithm \hfill}
\vskip5pt
\begin{center}
\hrule
\vskip7pt
\normalsize
\begin{enumerate}
\item[] \textit{Input}:  Data $(X_i,Y_i)$, regularization parameter $\lambda$.
\item[] \textit{Initialize} $\hat f_j  = 0$, for $j=1,\ldots, p$.
\item[] \textit{Iterate} until convergence:
\begin{enumerate}
   \item[] \textit{For each $j=1,\ldots, p$}:
   \begin{enumerate}
        \item[(1)] Compute the residual: $\Res_{j} = Y - \sum_{k\neq j} \hat f_k(X_{k})$;
        \item[(2)] Estimate $\Proj_j = \E[\Res_j\given
          X_j]$ by smoothing: $\hat \Proj_j = \S_j \Res_j$; \\[5pt]
        \item[(3)] Estimate norm: $ \hat s_j^2 = \frac{1}{n} \sum_{i=1}^n \hat P_j^2(i)$;
        \item[(4)] Soft-threshold:  $\hat f_j = \left[1 - {\lambda}/{\hat s_j}\right]_+ \hat \Proj_j$;
        \item[(5)] Center:  $\hat f_j \leftarrow \hat f_j - \text{mean}(\hat f_j)$.
   \end{enumerate}
\end{enumerate}
\item[] \textit{Output}: Component functions $\hat f_j$ and estimator $\hat
  m(X_i) = \sum_j \hat f_j(X_{ij})$.
\end{enumerate}
\vskip3pt
\hrule
\end{center}
\vskip0pt
\caption{The SpAM backfitting algorithm.  The first two steps in the iterative
  algorithm are the usual backfitting procedure; the remaining steps carry out
  functional soft thresholding.}
\label{fig:backfitting:algo}
\end{figure}

\begin{figure}[t]
{\sc Coordinate Descent Lasso \hfill}
\vskip5pt
\begin{center}
\hrule
\vskip7pt
\normalsize
\begin{enumerate}
\item[] \textit{Input}:  Data $(X_i,Y_i)$, regularization parameter $\lambda$.
\item[] \textit{Initialize} $\hat \beta_j  = 0$, for $j=1,\ldots, p$.
\item[] \textit{Iterate} until convergence:
\begin{enumerate}
   \item[] \textit{For each $j=1,\ldots, p$}:
   \begin{enumerate}
        \item[(1)] Compute the residual: $\Res_{j} = Y - \sum_{k\neq j} \hat \beta_k X_{k}$;
        \item[(2)] Project residual onto $X_j$: $\Proj_j = X_j^T R_j$ \\[5pt]
        \item[(3)] Soft-threshold:  $\hat \beta_j = \left[1 - {\lambda}/{|P_j|}\right]_+ P_j$;
   \end{enumerate}
\end{enumerate}
\item[] \textit{Output}: Estimator $\hat  m(X_i) = \sum_j \hat\beta_j X_{ij}$.
\end{enumerate}
\vskip3pt
\hrule
\end{center}
\vskip0pt
\caption{The SpAM backfitting algorithm is a functional version of the
  coordinate descent algorithm for the lasso, which computes 
  $\hat\beta = \arg\min \frac{1}{2} \|Y -  X\beta\|_2^2 + \lambda \|\beta\|_1$.}
\label{fig:coordlasso}
\end{figure}

\vspace{10pt}
{\em Basis Functions.}
It is useful to express the model in terms of basis functions.
Recall that $B_j =(\psi_{jk}:\ k=1,2,\ldots)$ is an orthonormal basis
for ${\cal T}_j$ and
that $\sup_x |\psi_{jk}(x)| \leq B$ for some $B$.
Then
\begin{equation}
f_j(x_j) = \sum_{k=1}^\infty \beta_{jk}\psi_{jk}(x_j)
\end{equation}
where $\beta_{jk} = \int f_j(x_j) \psi_{jk}(x_j) dx_j$.

Let us also define
\begin{equation}
\tilde{f}_j(x_j) = \sum_{k=1}^d \beta_{jk}\psi_{jk}(x_j)
\end{equation}
where $d=d_n$ is a truncation parameter.
For the Sobolev space ${\cal T}_j$ of order two
we have that
$\left\|f_j - \tilde{f}_j\right\|^2 = O(1/d^4)$.
Let $S = \{ j:\ f_j \neq 0\}$.
Assuming the sparsity condition $|S|=O(1)$ it follows that
$\left\|m - \tilde{m}\right\|^2 = O(1/d^4)$
where
$\tilde{m} = \sum_j \tilde{f}_j$.
The usual choice is $d \asymp n^{1/5}$ yielding truncation bias
$\left\|m - \tilde{m}\right\|^2 = O(n^{-4/5})$.

In this setting, the smoother can be taken to be the least squares
projection onto the truncated set of basis functions
$\{\psi_{j1},\ldots, \psi_{jd}\}$; this is also called orthogonal series
smoothing.
Let $\Psi_j$ denote the
$n\times d_n$ matrix given by
$\Psi_j(i,\ell) = \psi_\ell(X_{ij})$.
The smoothing matrix is the projection matrix
${\cal S}_j = \Psi_j (\Psi_j^T \Psi_j)^{-1} \Psi_j ^T$.
In this case, the backfitting algorithm 
in Figure \ref{fig:backfitting:algo} is 
exactly the
coordinate descent algorithm
for minimizing
\begin{equation}
\frac{1}{2n}\left\|Y - \sum_{j=1}^{p}\Psi_{j}\beta_{j}\right\|_{2}^{2} +
\lambda \sum_{j=1}^p \sqrt{\frac{1}{n} \beta_j^T \Psi_j^T \Psi_j \beta_j}
\end{equation}
which is the sample version of
(\ref{eq::lag}).
In Section~\ref{sec:sparsistency} we prove theoretical properties
assuming that this particular smoother is being used.

\vspace{1cm}
{\em Connection with the Grouped Lasso.}
The SpAM model can be thought of as
a functional version of the grouped lasso
\citep{Yuan:Lin:06} as we now explain.
Consider the following linear regression model with multiple factors, 
\begin{equation}
Y = \sum_{j=1}^{p_n} X_j \beta_j + \epsilon = X \beta + \epsilon,
\end{equation}
where $Y$ is an $n \times 1$ response vector, $\epsilon$ is an $n
\times 1$ vector of iid mean zero noise, $X_j$ is an $n \times d_j$
matrix corresponding to the $j$-th factor, and $\beta_j$ is the
corresponding $d_j \times 1$ coefficient vector.  Assume for
convenience (in this subsection only)
that each $X_j$ is orthogonal, so that $X_j^T X_j =
I_{d_j}$, where $I_{d_j}$ is the $d_j \times d_j$ identity matrix.
We use $X=(X_1, \ldots,
X_{p_n})$ to denote the full design matrix and use $\beta=(\beta_1^T,
\dots, \beta_{p_n}^T)^T$ to denote the parameter.

The \textit{grouped lasso} estimator
is defined as the solution of the following convex optimization problem:
\begin{eqnarray}\label{eqn.grouplasso}
\hat{\beta} (\lambda_n) = \mathop{\arg\min}_\beta \|Y -
X\beta\|^2_2 + \lambda_n
\sum_{j=1}^{p_n}\sqrt{d_j}\|\beta_j\|
\end{eqnarray}
where $\sqrt{d_j}$ scales the $j$th term to compensate for different group sizes.

It is obvious that when $d_j=1$ for $j=1,\ldots,p_n$, the grouped lasso
becomes the standard lasso.  From the KKT optimality conditions,
a necessary and sufficient condition 
for $\hat{\beta}=(\hat{\beta}_1^T, \ldots, \hat{\beta}_{p}^T)^T$ to be
the grouped lasso solution is 
\begin{eqnarray}\label{kkt}
- X_j^T (Y-X\hat{\beta}) + 
\frac{\lambda \sqrt{d_j} \hat{\beta}_j }{\| \hat{\beta}_j \|} \ \ \ = 
& \mathbf{0}, & \forall \hat{\beta}_j \ne \mathbf{0},\\
\| X_j^T (Y-X\hat{\beta})\| \ \ \ \le & \lambda\sqrt{d_j},
& \forall \hat{\beta}_j = \mathbf{0} . \nonumber
\end{eqnarray}
Based on this stationary condition, an 
iterative blockwise coordinate descent algorithm can be derived; as shown by 
\cite{Yuan:Lin:06}, a solution to (\ref{kkt}) satisfies
\begin{eqnarray}\label{fixedpoint}
\hat{\beta}_j = \left[ 1 - \frac{\lambda\sqrt{d_j}}{\|S_j\|}
\right]_+S_j
\end{eqnarray}
where $S_j= X^T_j(Y - X\beta_{\backslash j})$, with $\beta_{\backslash
  j} =
(\beta^T_1,\ldots,\beta^T_{j-1},\mathbf{0}^T,\beta^T_{j+1},\ldots,
\beta^T_{p_n})$. By iteratively applying (\ref{fixedpoint}), the
grouped lasso solution can be obtained.

As discussed in the introduction,
the COSSO model of \cite{LZ:Cosso}
replaces the lasso constraint on $\sum_j |\beta_j|$
with a RKHS constraint.
The advantage of our formulation is that it decouples 
smoothness
($g_j\in {\cal T}_j$) and sparsity
($\sum_j |\beta_j|\leq L$).
This leads to a simple algorithm that
can be carried out with any nonparametric smoother
and scales easily to high dimensions.

\section{Sparse Nonparametric Logistic Regression}
\label{sec:logistic}

The SpAM backfitting procedure can be extended to nonparametric
logistic regression for classification.  The additive logistic model is 
\begin{equation}
\P(Y=1\given X) \equiv p(X;f) = 
\frac{\exp\left(\sum_{j=1}^p f_j(X_j)\right)}
{1 + \exp\left(\sum_{j=1}^p f_j(X_j)\right)}
\end{equation}
where $Y\in\{0,1\}$, and the population log-likelihood is 
\begin{equation}
\ell(f) = \E\left[Yf(X) - \log\left(1+\exp f(X)\right)\right].
\end{equation}
Recall that in the local scoring algorithm for
generalized additive models \citep{Hast:Tibs:1999} in the logistic case, one 
runs the backfitting procedure within Newton's method.  Here one
iteratively computes the transformed response for the current estimate $f_0$
\begin{equation}
Z_i = f_0(X_i) + \frac{Y_i-p(X_i;f_0)}{p(X_i;f_0)(1-p(X_i;f_0))}
\end{equation}
and weights $w(X_i) = p(X_i;f_0)(1-p(X_i;f_0)$, and carries out a weighted
backfitting of $(Z,X)$ with weights~$w$. The weighted smooth is given
by 
\begin{equation}
\label{eq:unregscore}
\hat P_j = \frac{\S_j (w R_j)}{\S_j w}.
\end{equation}
To incorporate the sparsity penalty, we first note that the Lagrangian
is given by 
\begin{equation}
\L(f, \lambda) =
\E\left[\log\left(1+e^{f(X)}\right) - Yf(X)\right]
  + \lambda \left(\sum_{j=1}^p \sqrt{\E(f_j^2(X_j))} - L\right)
\end{equation}
and the stationary condition for component function $f_j$ is
$\E\left(p - Y\given X_j\right) + \lambda v_j= 0$
where $v_j$ is an element of the subgradient  $\partial {\sqrt{\E(f_j^2)}}$.
As in the unregularized case, this condition is nonlinear in $f$,
and so we linearize the gradient of the log-likelihood 
around $f_0$.  This yields the linearized
condition $\E\left[ w(X) (f(X) - Z) \given X_j\right] + \lambda v_j = 0$.
%\begin{eqnarray}
%0 &=& 
%\E\left(p(X;f_0) - Y + p(X;f_0)(1-p(X;f_0))(f(X) - f_0(X)) \given X_j\right)
%  + \lambda v_j\\
%&=& \E\left[ w(X) (f(X) - Z) \given X_j\right] + \lambda v_j
%\end{eqnarray}
When $\E(f_j^2) \neq 0$, this implies the condition
\begin{eqnarray}
\left(\E\left(w \given X_j\right) +
  \frac{\lambda}{\sqrt{\E(f_j^2)}}\right) f_j(X_j) = \E(w R_j\given X_j).
\end{eqnarray}
In the finite sample case, in terms of the smoothing matrix $\S_j$, this becomes
\begin{eqnarray}
f_j = \frac{\S_j (w R_j)}{\S_j w + {\lambda}\left/{\sqrt{\E(f_j^2)}}\right.} .
\end{eqnarray}
If $\|\S_j(w R_j)\| < \lambda$, then $f_j = 0$. Otherwise,
this implicit, nonlinear equation for $f_j$ cannot be solved explicitly, so
we propose to iterate until convergence:
\begin{equation}
f_j \;\leftarrow\; \frac{\S_j (w R_j)}{\S_j w + {\lambda
    \sqrt{n}}\left/{\|f_j\|}\right.} .
\end{equation}
When $\lambda=0$, this yields the standard local scoring
update~\eqref{eq:unregscore}.  An example of logistic SpAM is
given in Section~\ref{sec:examples}.

\section{Choosing the Regularization Parameter}
\label{sec:tuning}

We choose $\lambda$ by minimizing
an estimate of the risk.
Let $\nu_j$ be the effective degrees of freedom for
the smoother on the $j^{\rm th}$ variable, that is,
$\nu_j = {\rm trace}({\cal S}_j)$
where $\S_j$ is the smoothing matrix for the $j$-th dimension.
Also let $\hat{\sigma}^2$ be an estimate of the variance. 
Define the total
effective degrees of freedom as
\begin{equation}
\df(\lambda) = \sum_{j} \nu_j I\left(\left\|\hat{f}_j\right\|\neq 0\right).
\end{equation}
Two estimates of risk are
\begin{equation}
C_p  = \frac{1}{n} \sum_{i=1}^n\left( Y_i -  \sum_{j=1}^p\widehat{f}_j(X_j)\right)^2 + 
\frac{2\hat{\sigma}^2}{n} \df(\lambda)
\end{equation}
and
\begin{equation}
\GCV(\lambda) = 
\frac{\frac{1}{n}\sum_{i=1}^n (Y_i - \sum_j \hat{f}_j(X_{ij}))^2}
     {(1- \df(\lambda)/n)^2}.
\end{equation}
The first is $C_p$ and the second is
generalized cross validation but with
degrees of freedom defined by $\df(\lambda)$.
A proof that
these are valid
estimates of risk is not currently available;
thus,  these should be regarded as heuristics.

Based on the results in
\cite{WR}
about the lasso, it seems likely that
choosing $\lambda$ by risk estimation
can lead to overfitting.
One can further clean the estimate
by testing 
$H_0: f_j =0$ for all 
$j$ such that $\hat{f}_j\neq 0$.
For example, the tests in
\cite{Fan:Jiang:05} could be used.

\section{Examples}
\label{sec:examples}

To illustrate the method, we consider a few examples.

\vspace{.3cm}
{\em Synthetic Data.} \enspace
Our first example is from \citep{hardle:04}. 
We generated 
$n = 150$ observations
from the following 200-dimensional
additive model:
\begin{gather}
Y_i = f_1(x_{i1})+f_2(x_{i2}) + f_3(x_{i3}) + f_4(x_{i4}) + \epsilon_i \\
\nonumber
f_1(x) = -2\sin(2x), \; f_2(x) = x^2 - \frac{1}{3},\; f_3(x) =
x-\frac{1}{2},\; f_4(x) = e^{-x} +e^{-1} - 1
\end{gather}
and $f_j(x) = 0$ for $j\geq 5$ with noise $\epsilon_i~\sim \mathcal{N}(0,1)$.
These data therefore have 196 irrelevant dimensions.  

\begin{figure}[htp]
\begin{center}
\begin{tabular}{ccc}
\hskip-.0in
\includegraphics[width=.3\textwidth]{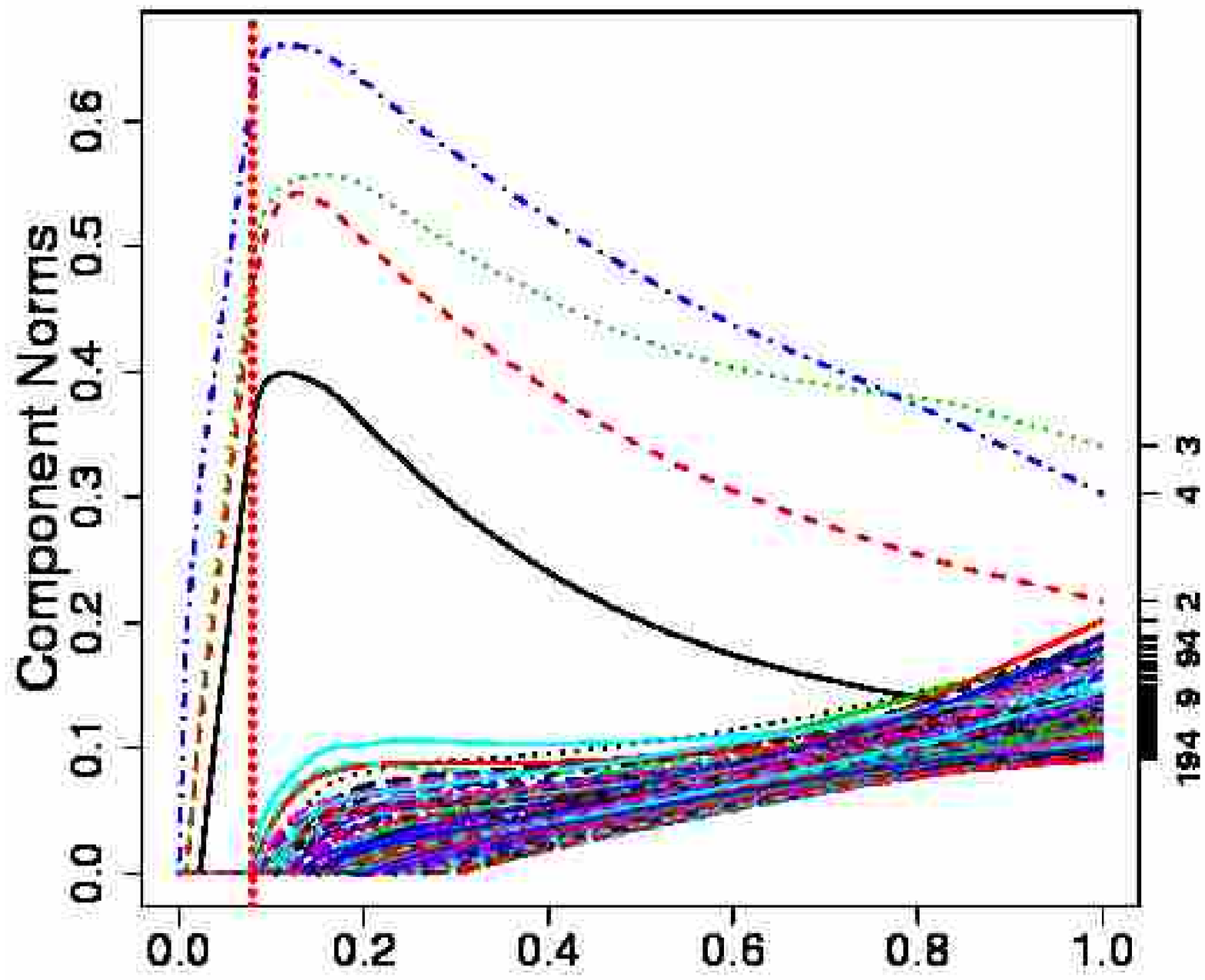} &
\includegraphics[width=.3\textwidth]{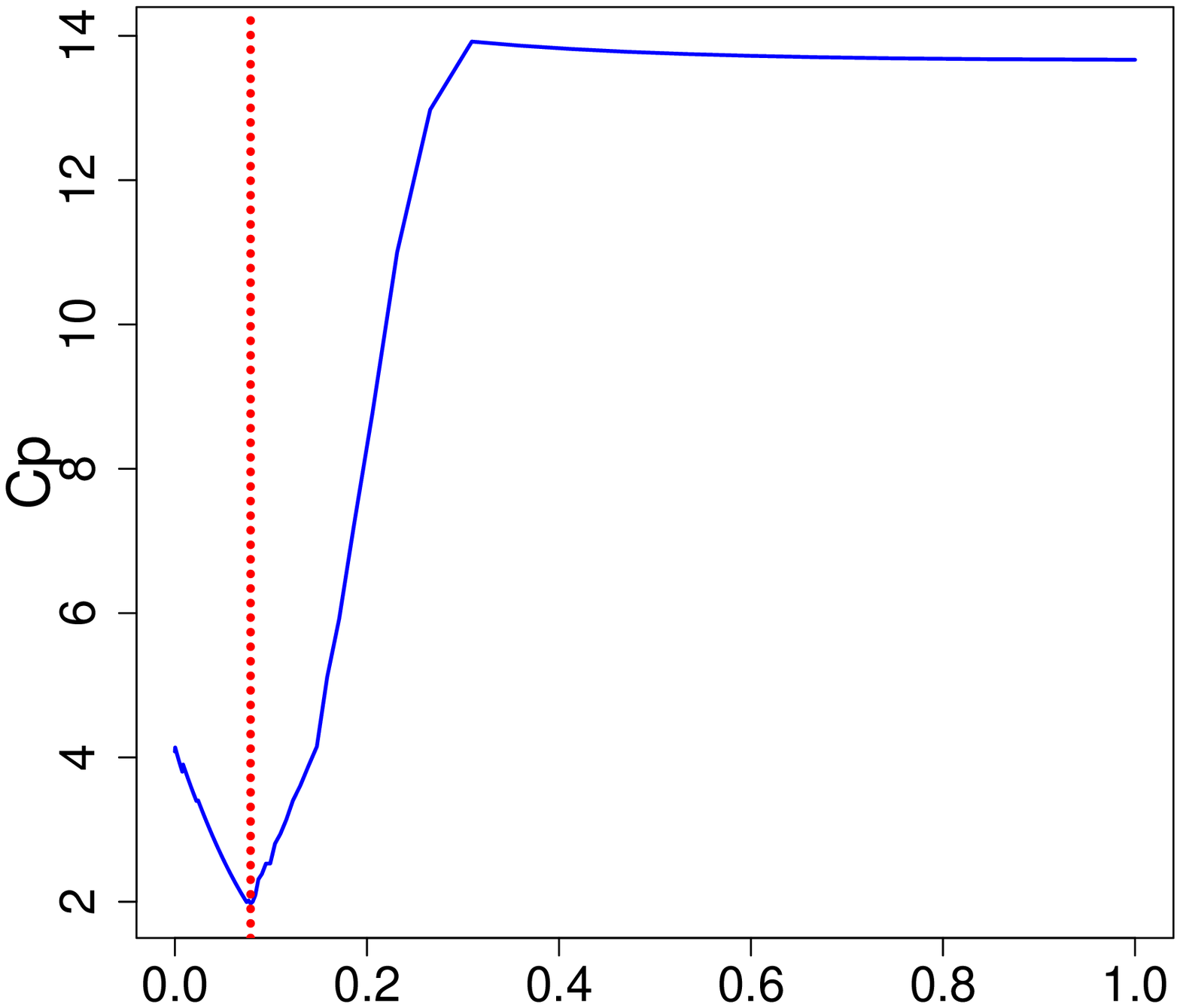}  &
\includegraphics[width=.3\textwidth]{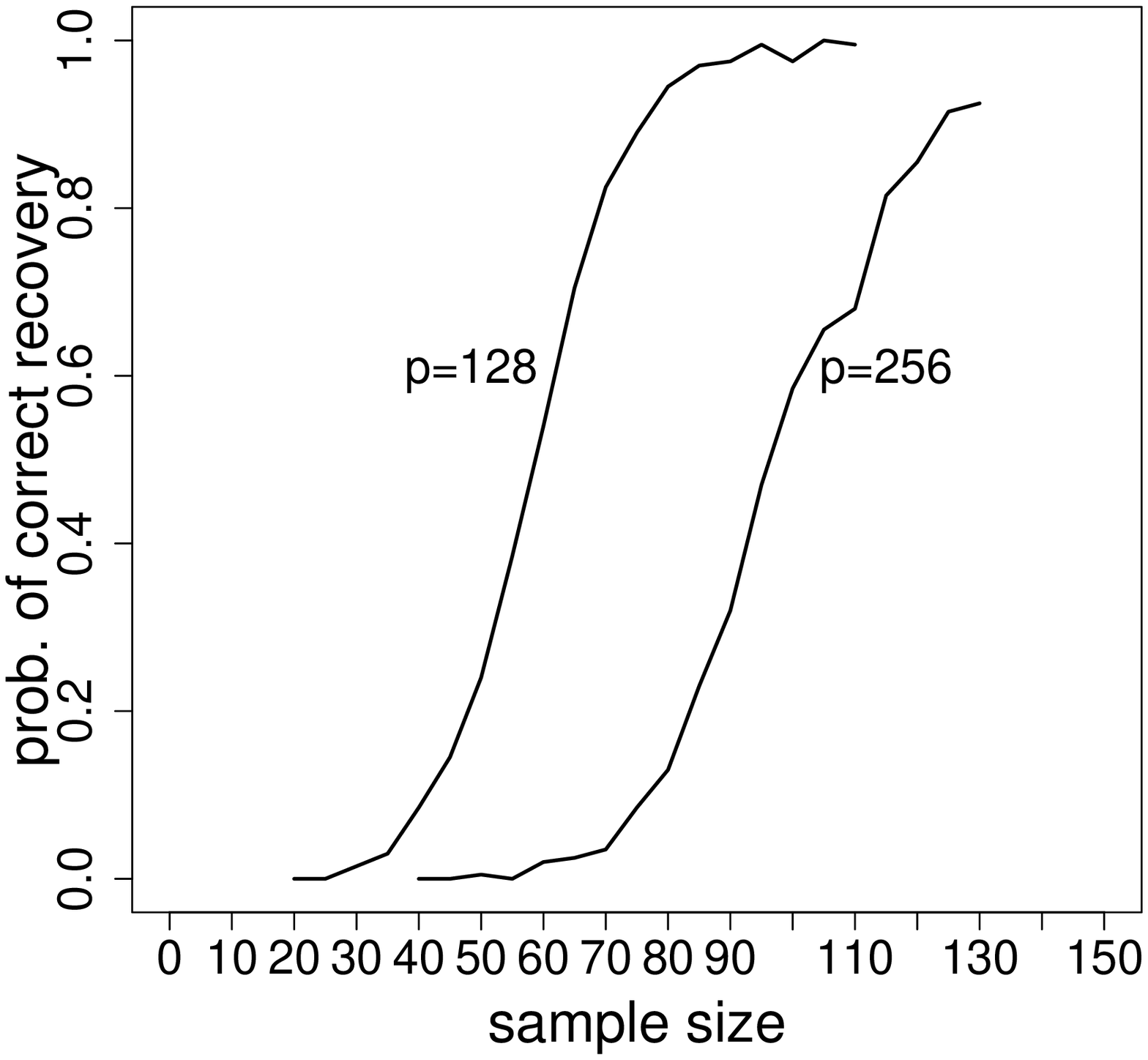}\\[5pt]
\end{tabular}
\begin{tabular}{ccc}
\includegraphics[width=.28\textwidth]{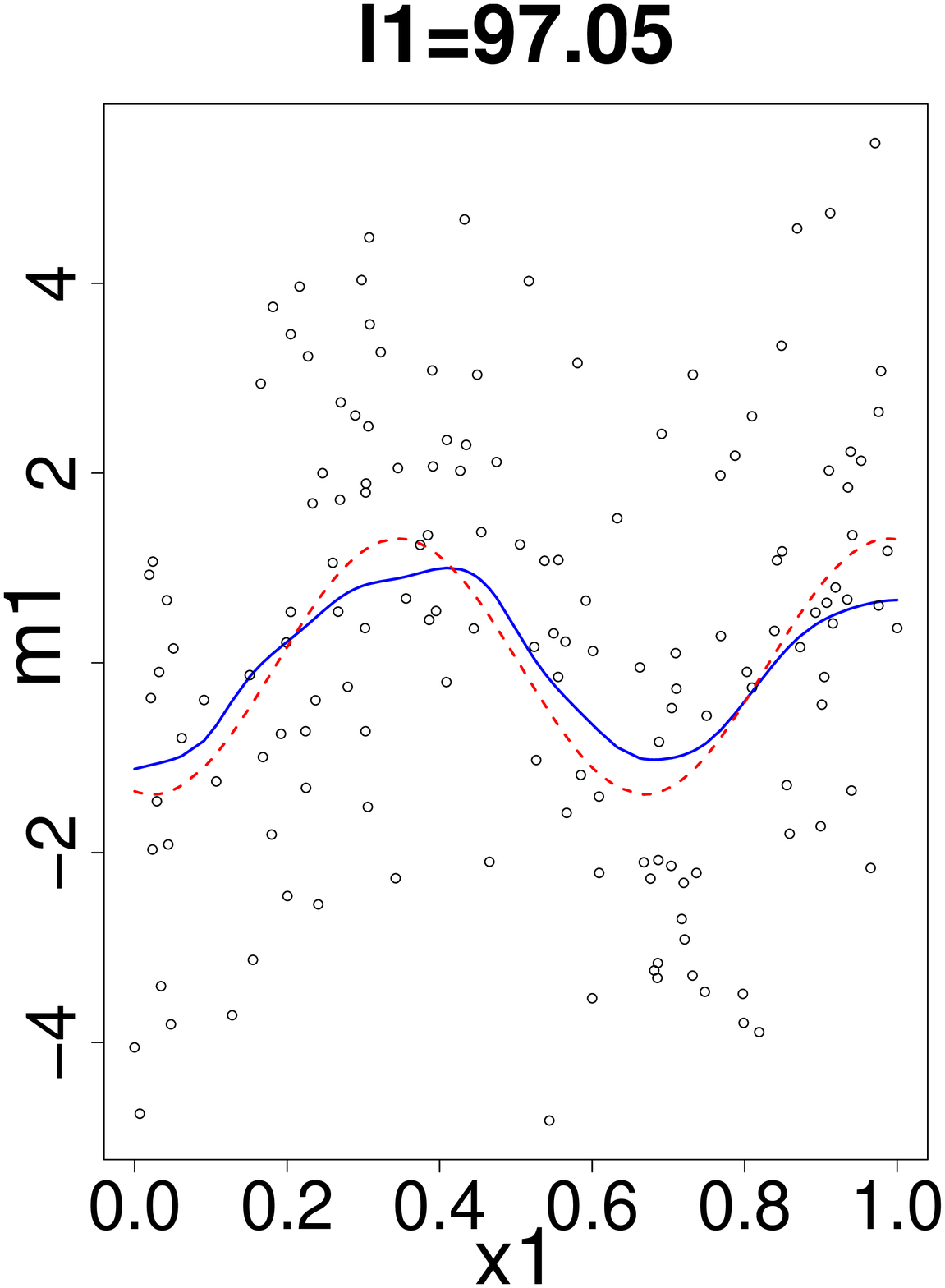}&
\includegraphics[width=.28\textwidth]{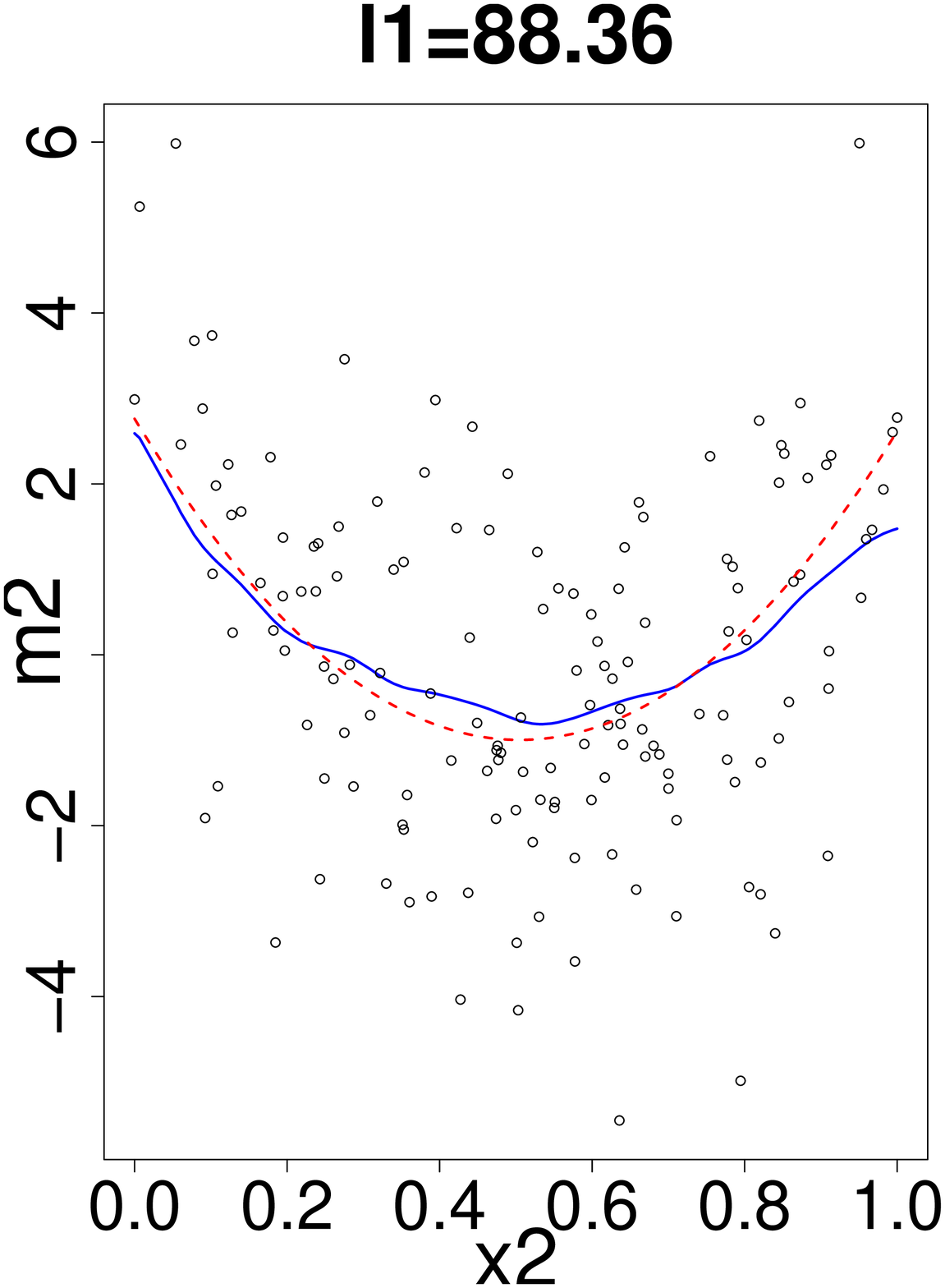}&
\includegraphics[width=.28\textwidth]{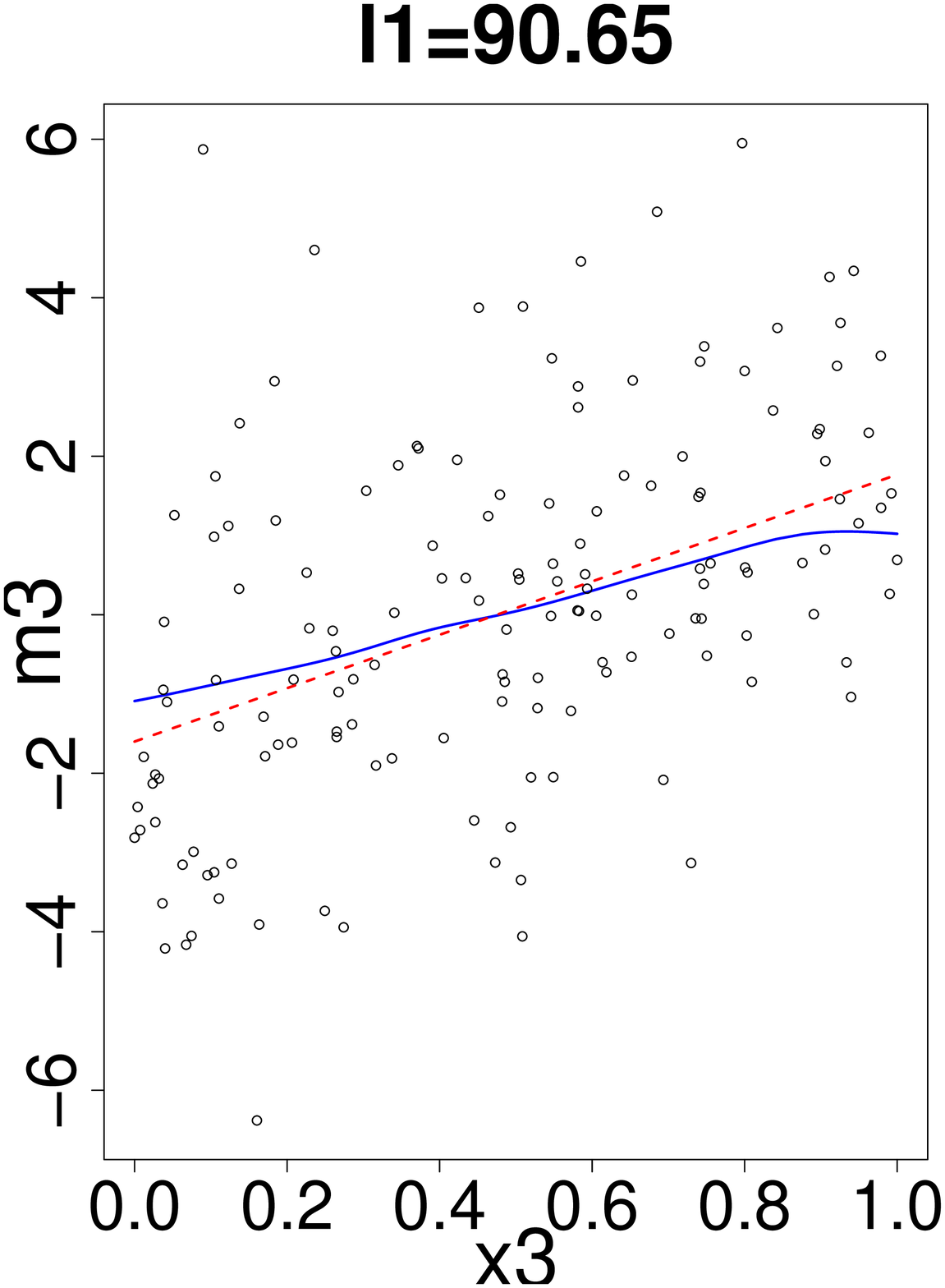}\\[10pt]
\includegraphics[width=.28\textwidth]{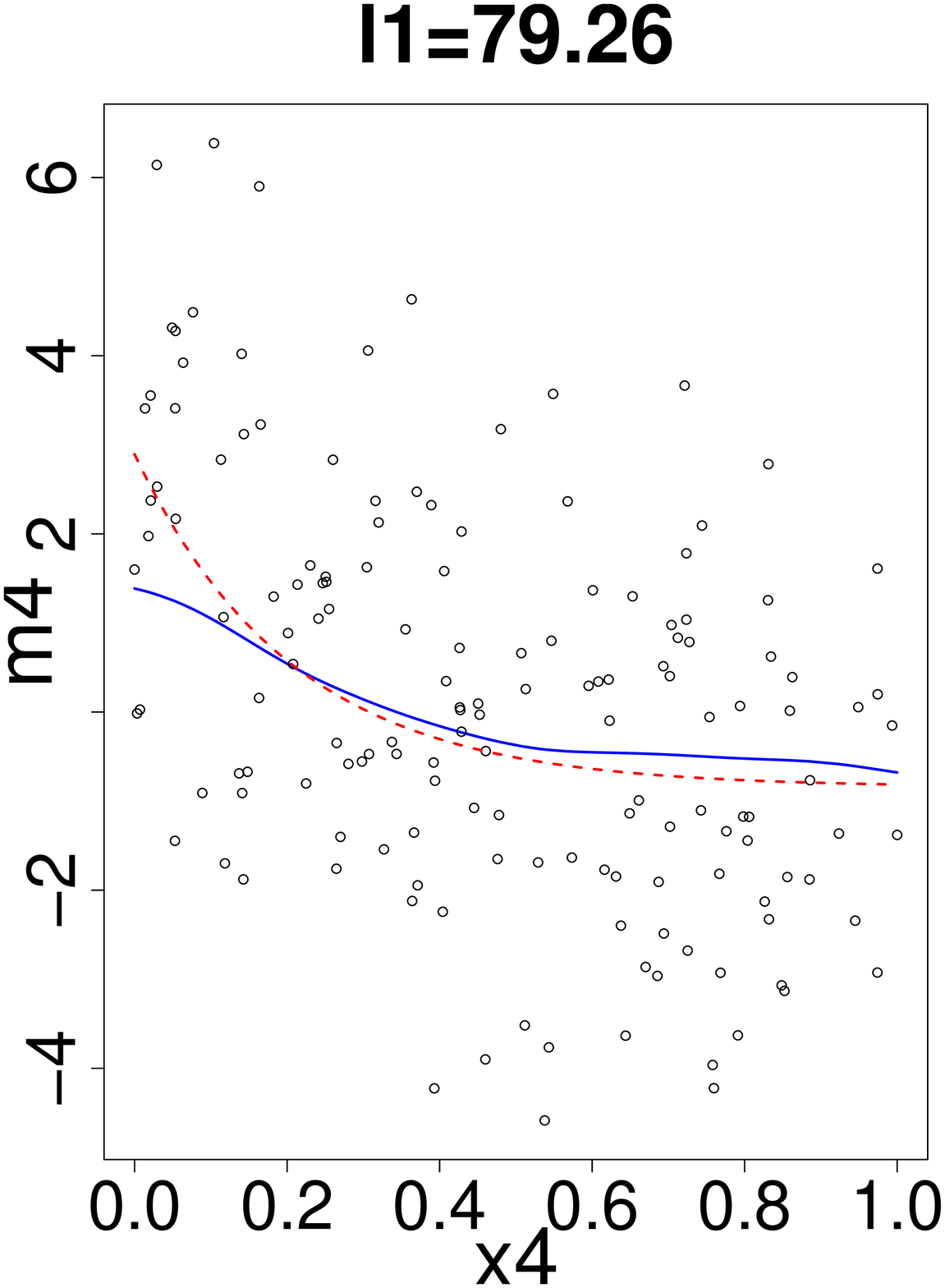}&
\includegraphics[width=.28\textwidth]{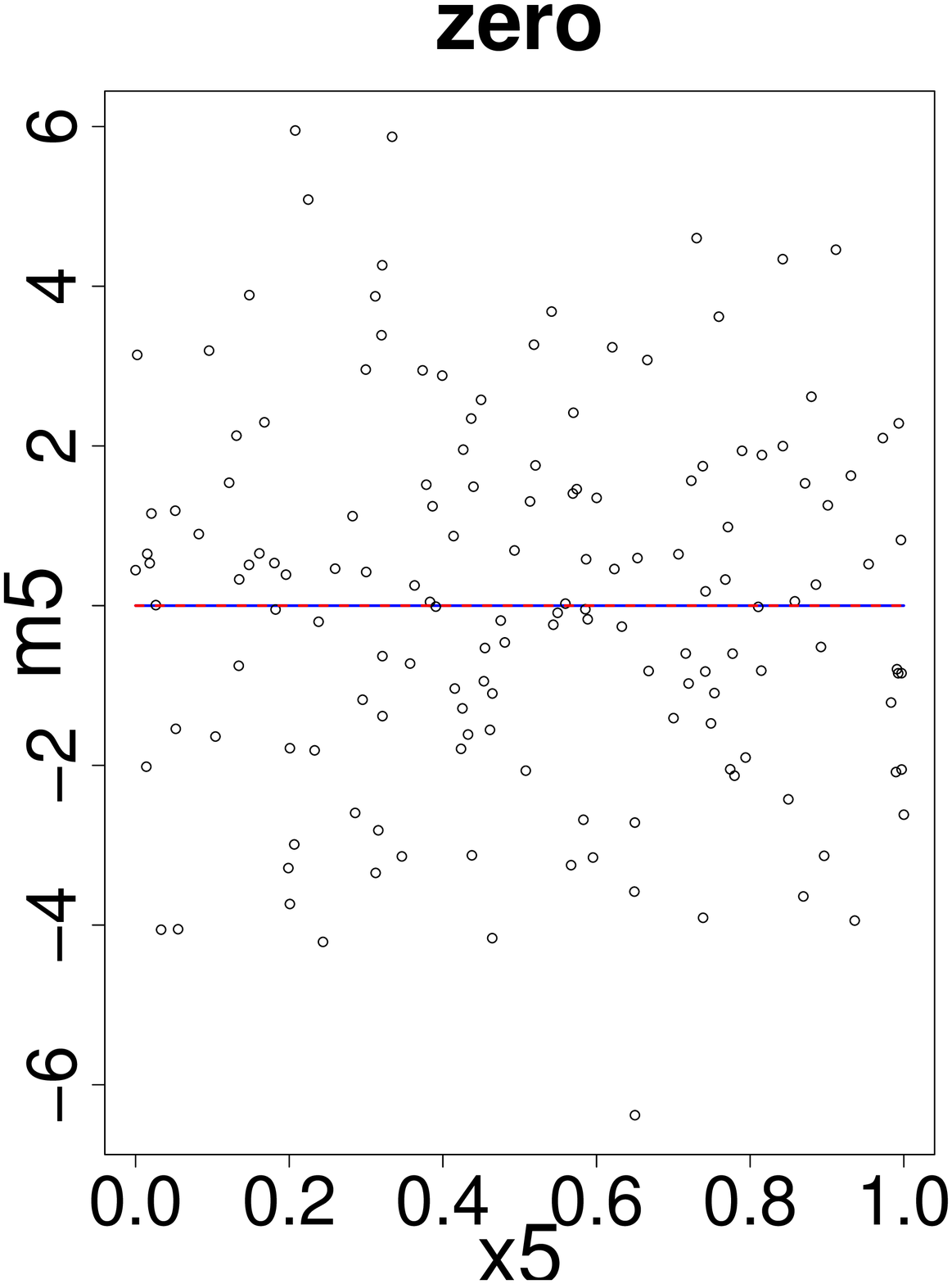}&
\includegraphics[width=.28\textwidth]{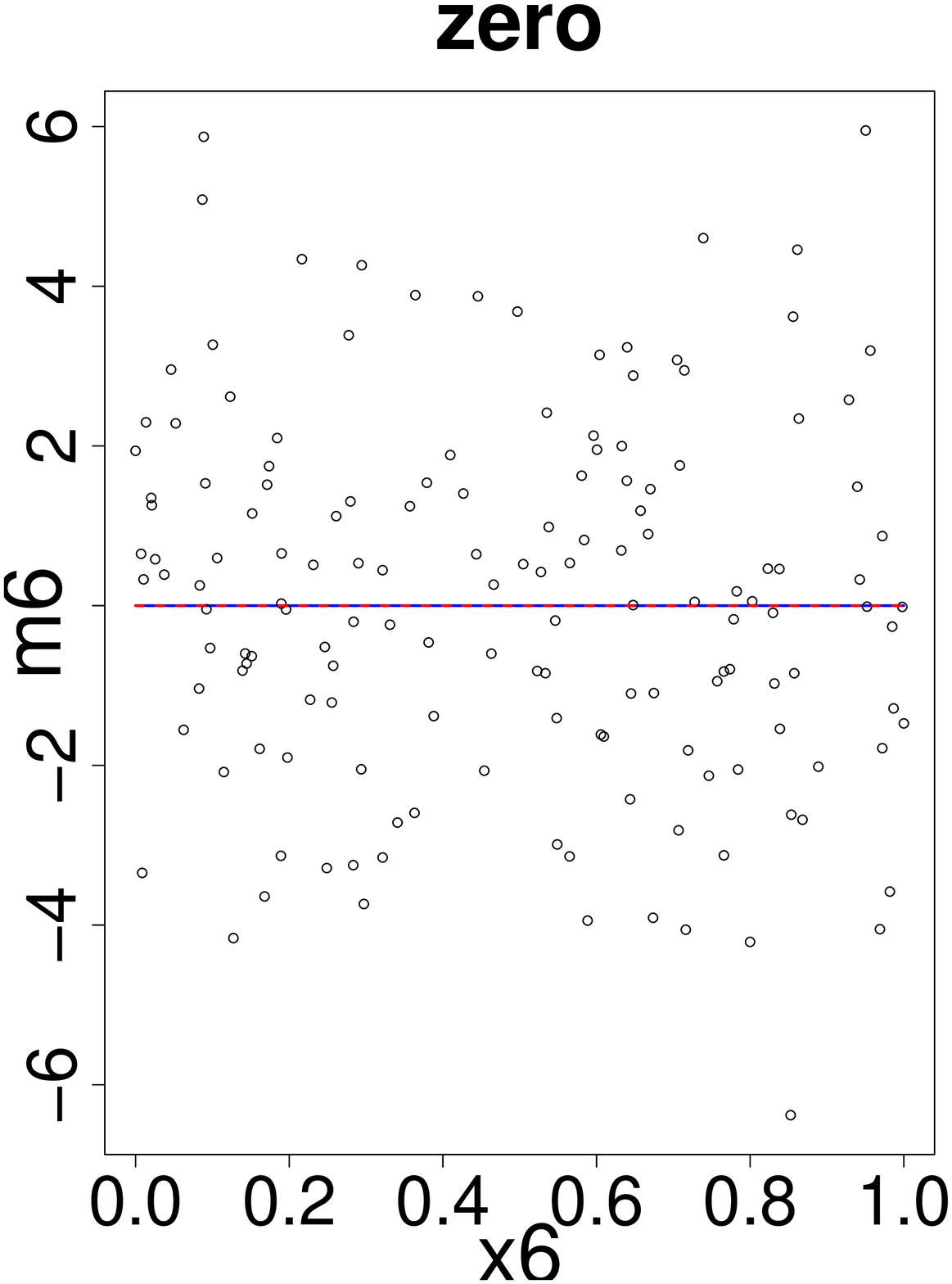}
\end{tabular}
\end{center}
\vskip-.1in
\caption{(Simulated data) Upper left: The empirical $\ell_2$ norm
of the estimated components as plotted against the regularization parameter
$\lambda$; the value on the $x$-axis is proportional to $\sum_j \|\hat f_j\|$.
Upper center: The $C_p$ scores against the regularization 
parameter $\lambda$; the dashed vertical line corresponds to the
value of $\lambda$ which has the smallest $C_p$ score. Upper right:
The proportion of 200 trials where the correct relevant 
variables are selected, as a function of sample size $n$.
Lower (from left to right): Estimated (solid lines) versus true additive component
functions (dashed lines) for the first 6 dimensions; 
the remaining components are zero.
 }\label{fig.sim}
\end{figure}

The results
of applying SpAM with the plug-in bandwidths are summarized
in Figure~\ref{fig.sim}.
The top-left plot in Figure~\ref{fig.sim} shows regularization paths 
as the parameter $\lambda$ varies; each curve is a plot
of $\|\hat f_j(\lambda)\|$ versus 
\begin{equation}
\frac{\sum_{k=1}^p \| \hat f_k(\lambda)\|}{\max_{\lambda}{\sum_{k=1}^p \| \hat f_k(\lambda)\|}}
\end{equation}
for a particular variable $X_j$.
The estimates are generated efficiently over a sequence of $\lambda$
values by ``warm starting'' $\hat f_j(\lambda_t)$ at the previous
value $\hat f_j(\lambda_{t-1})$.  The top-center plot shows the
$C_p$ statistic as a function of $\lambda$.  The top-right plot
compares the empirical probability of correctly selecting the 
true four variables as a function of sample size $n$, for
$p=128$ and $p=256$.  This behavior suggests the same threshold
phenomenon that was shown for the lasso by \cite{Wai06}.

\vspace{.5cm}
{\em Boston Housing.} \enspace
The Boston housing data were collected to study
house values in the suburbs of Boston.
There are 506
observations with 10 covariates. The dataset 
has been studied by many other authors \citep{hardle:04,LZ:Cosso},
with various transformations proposed for different covariates.  To
explore the sparsistency properties of our method, we add 20
irrelevant variables. Ten of them are randomly drawn from
$\text{Uniform}(0,1)$, the remaining ten are a random permutation of
the original ten covariates.
The model is
\begin{eqnarray}
 Y& =& \alpha + f_1({\tt crim}) + f_2({\tt indus})+ f_3({\tt
nox})+ f_4({\tt rm}) +  f_5({\tt age})  \nonumber \\
&& \; + \; f_6({\tt dis})+ f_7({\tt tax}) + f_8({\tt ptratio}) + f_9({\tt b})+ f_{10}({\tt lstat}) + \epsilon.
\end{eqnarray}

\def\hs{\hspace{-10pt}}
\begin{figure*}
\begin{center}
\begin{tabular}{cc}
\hskip-.2in
\includegraphics[width=.40\textwidth]{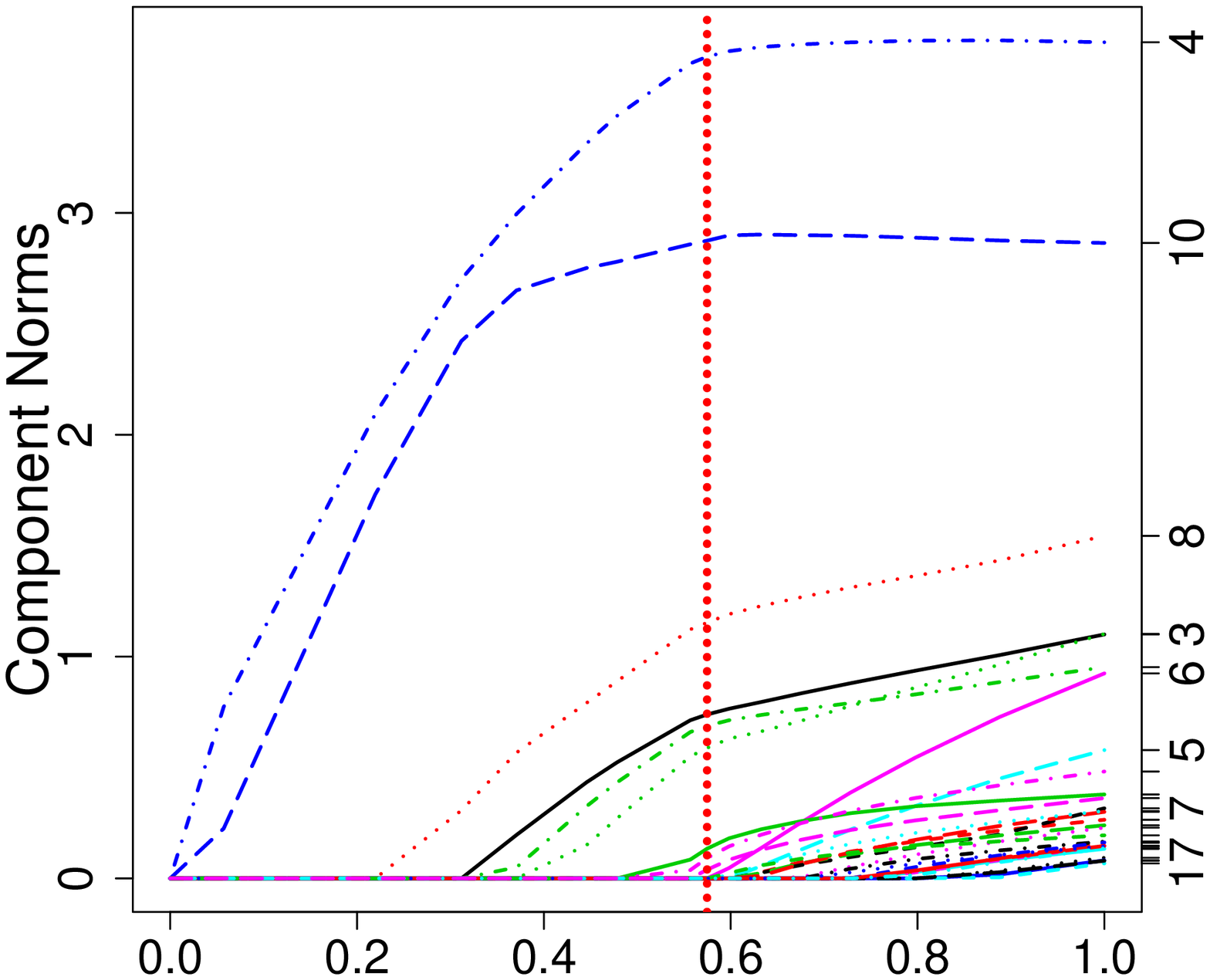} &
\hskip-.1in 
\includegraphics[width=.40\textwidth,
angle=0]{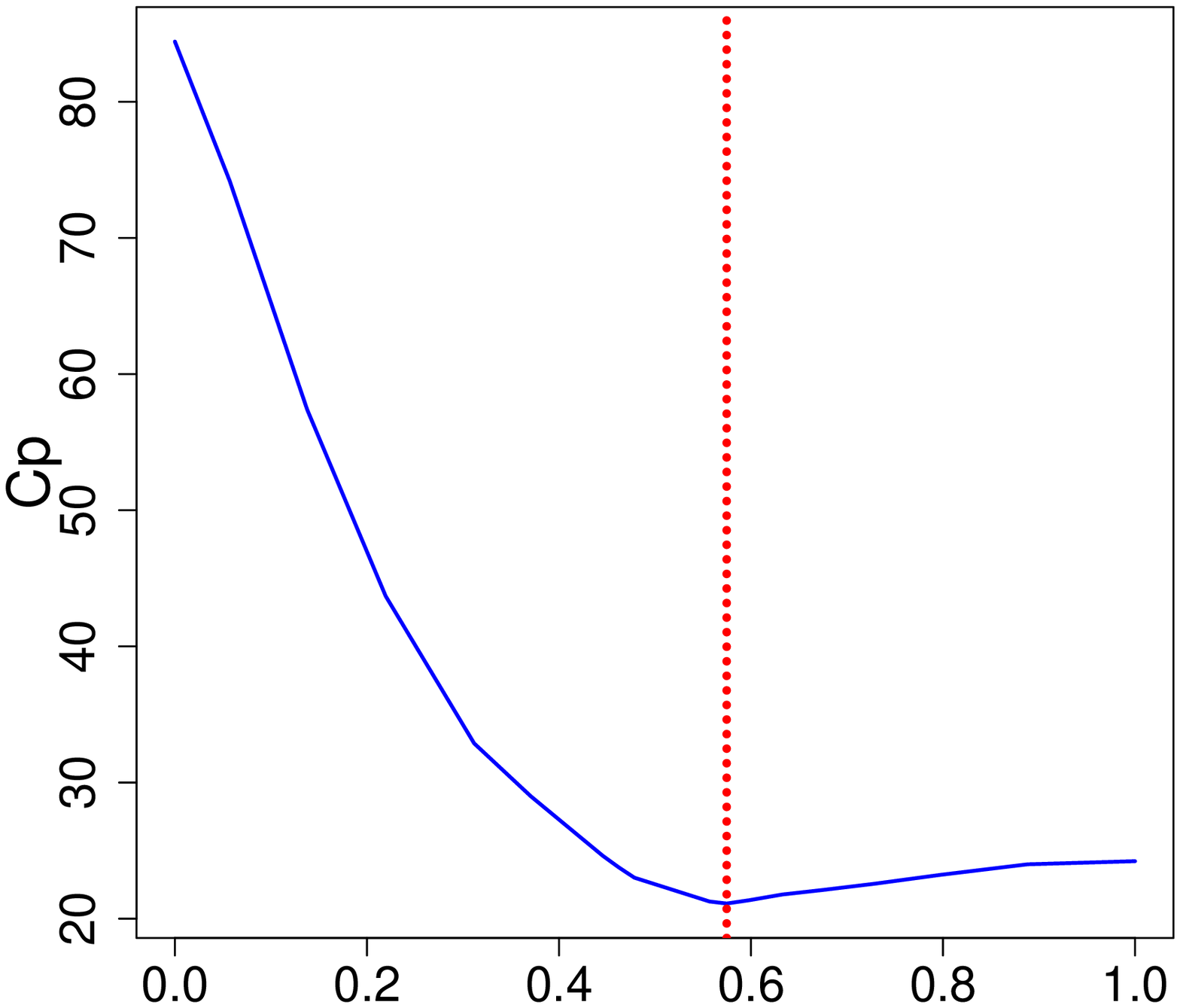} 
\end{tabular}
\begin{tabular}{cc}
\hs\includegraphics[width=.3\textwidth]{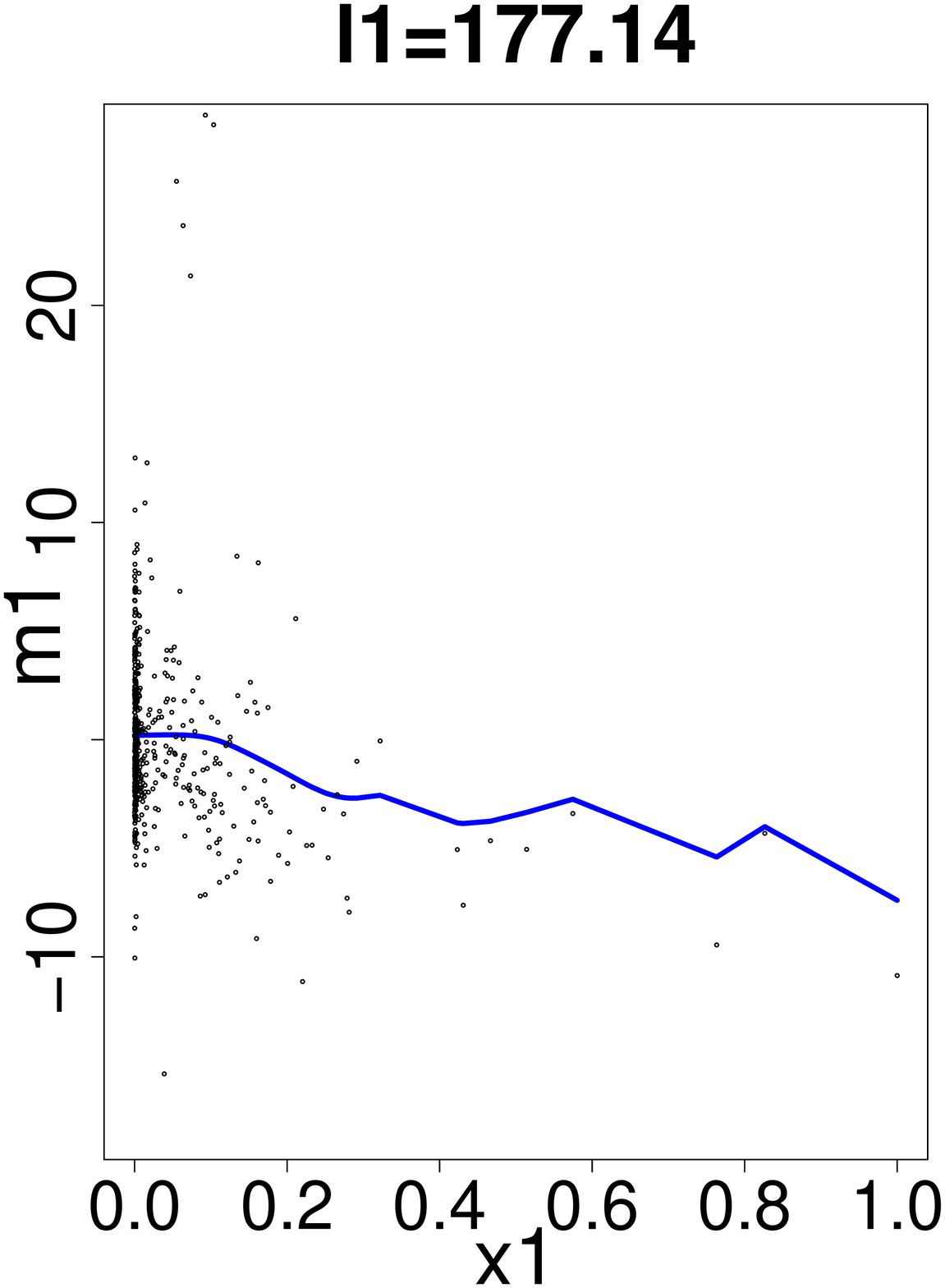} &
\hs\includegraphics[width=.3\textwidth]{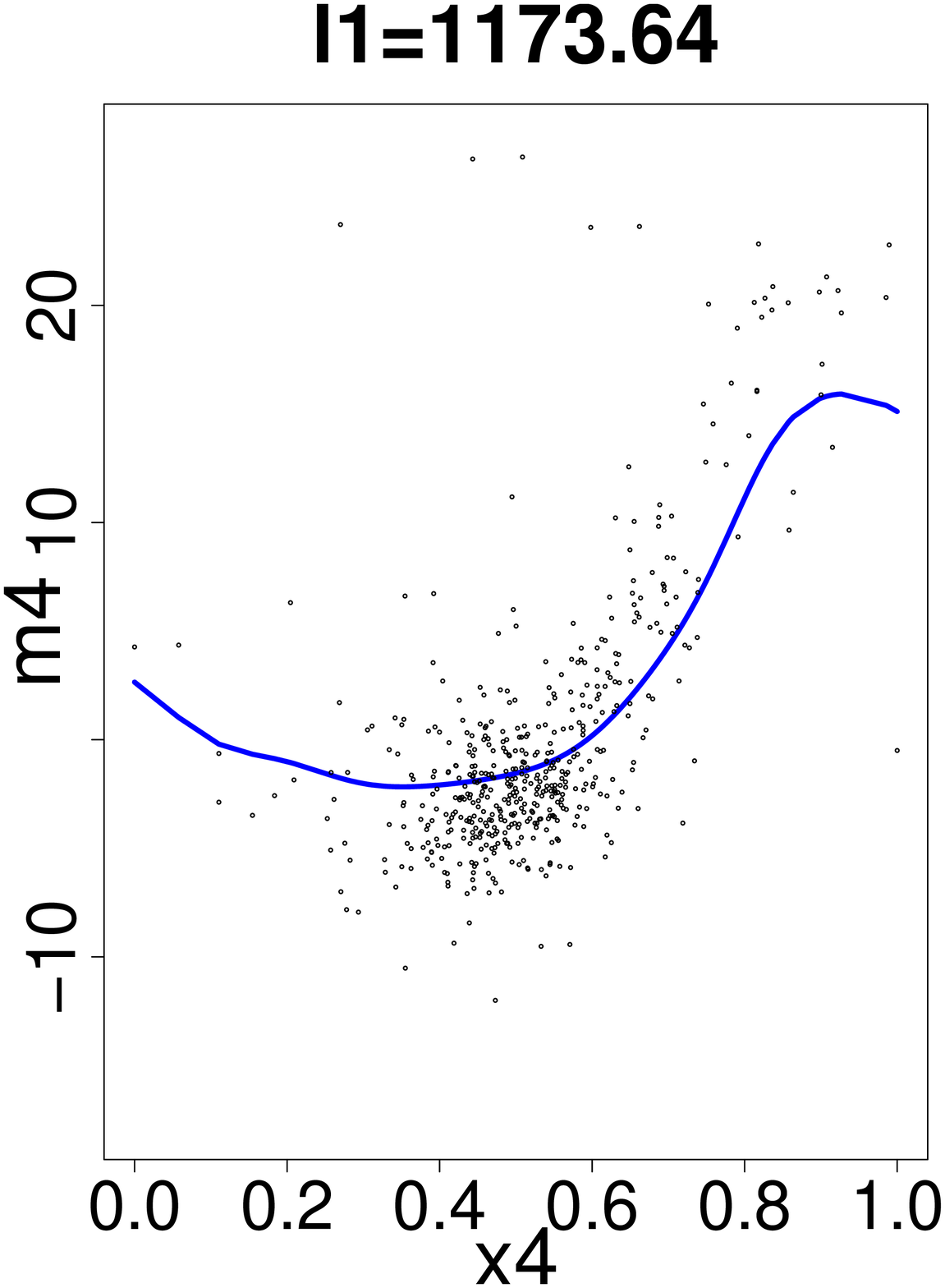} \\[10pt]
\hs\includegraphics[width=.3\textwidth]{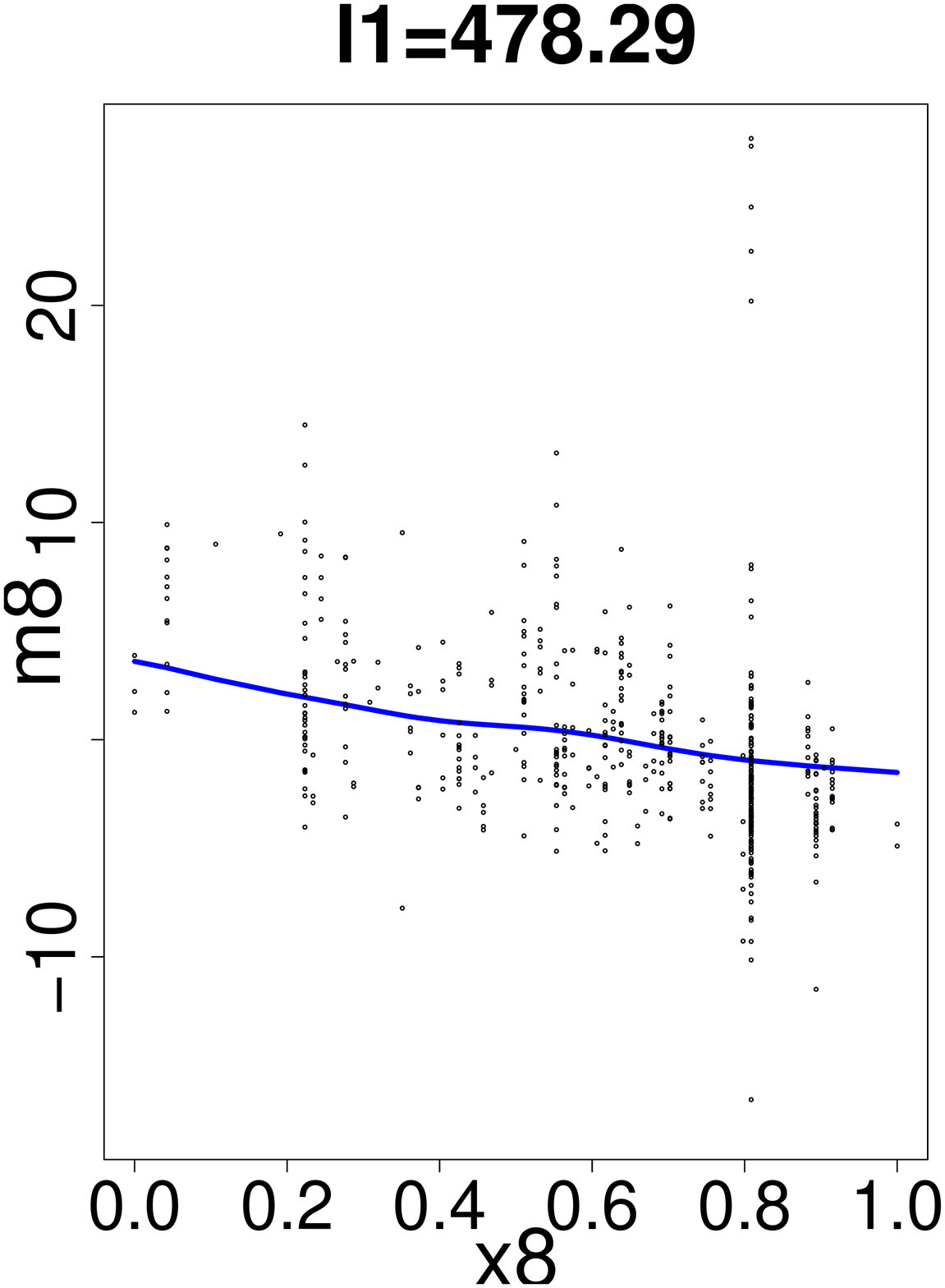} &
\hs\includegraphics[width=.3\textwidth]{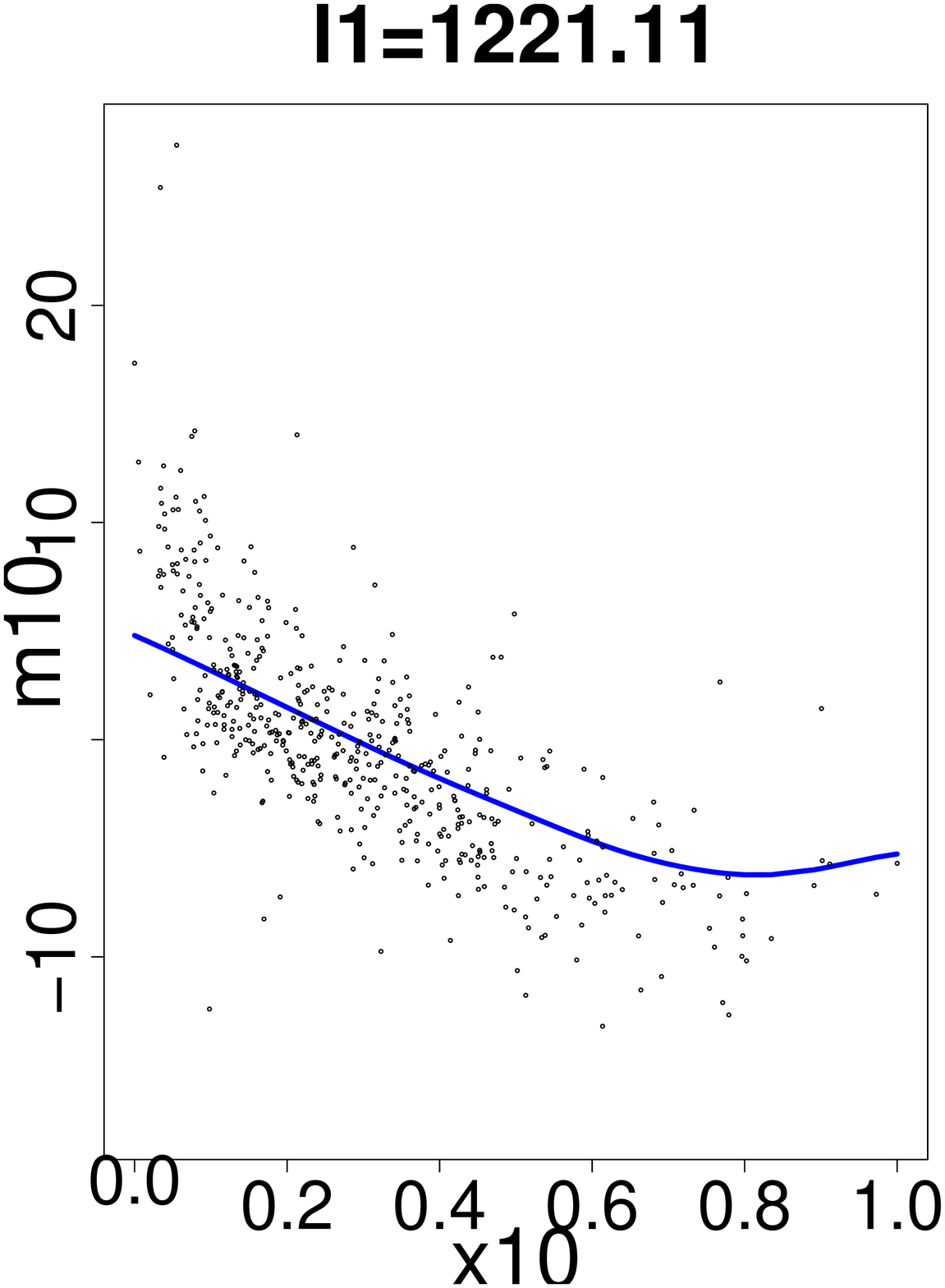}
\end{tabular}
\end{center}
\caption{(Boston housing) Left: The empirical $\ell_2$ norm of
the estimated components versus the regularization parameter 
$\lambda$. Center: The $C_p$ scores against $\lambda$; 
the dashed vertical line corresponds to best
$C_p$ score. Right:  Additive fits for four relevant variables.
}\label{fig.boston}
\end{figure*}

The result of applying SpAM to this 30 dimensional dataset is 
shown in Figure~\ref{fig.boston}.
SpAM identifies 6 nonzero components.  It correctly
zeros out both types of irrelevant variables. From the full
solution path, the important variables are seen to be ${\tt rm}$, {\tt
lstat}, {\tt ptratio}, and {\tt crim}. The importance of
variables {\tt nox} and {\tt b} is borderline. These
results are basically consistent with those obtained by 
other authors \citep{hardle:04}.  However, using
$C_p$ as the selection criterion, the variables ${\tt indux}$, ${\tt
age}$, ${\tt dist}$, and ${\tt tax}$ are estimated to be
irrelevant, a result not seen in other studies.

\def\i{{(i)}}
\def\plotwidth{.08\textwidth}
\def\headwidth{.17\textwidth}

\vspace{.5cm}
{\em SpAM for Spam.} \enspace
Here we consider an email spam classification problem, using the
logistic SpAM backfitting algorithm from Section~\ref{sec:backfitting}.
This dataset has been studied by \cite{Hastie:2001}, using a set of
3,065 emails as a training set, and conducting hypothesis tests to
choose significant variables; there are a total of 4,601 observations
with $p=57$ attributes, all numeric.  The attributes measure the
percentage of specific words or characters in the email, the average
and maximum run lengths of upper case letters, and the total number of
such letters. To demonstrate how SpAM performs with sparse
data, we only sample $n=300$ emails as the training set, with the
remaining $4301$ data points used as the test set.  We also use the
test data as the hold-out set to tune the penalization parameter~$\lambda$.

\begin{figure}[ht]
\begin{center}
\begin{small}
\begin{sc}
\begin{tabular}{c}
\hskip-5pt
\renewcommand{\arraystretch}{1.2}
\begin{tabular}{clccr}
\hline\hline
 $\lambda (\times 10^{-3})$  & Error & \# zeros & selected variables \vspace{.0cm}\\
\hline
5.5    & 0.2009 & 55& \scriptsize\{ 8,54\} \\
5.0    & 0.1725 & 51& \scriptsize \{   8, 9, 27, 53, 54, 57\}\\
4.5    & 0.1354 & 46& \scriptsize \{7, 8, 9, 17, 18, 27, 53, 54, 57, 58\}\\
4.0    & 0.1083 ($\surd$) &20&  \scriptsize \{4, 6--10, 14--22, 26, 27, 38, 53--58\}       \\
3.5    & 0.1117 & 0 & all \\
3.0    & 0.1174 &  0& all \\
2.5    & 0.1251 & 0&  all       \\
2.0    & 0.1259 & 0&  all \\
\hline
\end{tabular}
%\begin{tabular}{clcc}
%\hline\hline
% $\lambda (\times 10^{-3})$  & Error & \# zeros \\
%\hline
%5.5    & 0.2009 & 55 \\
%5.0    & 0.1725 & 51 \\
%4.5    & 0.1354 & 46 \\
%4.0    & 0.1083 ($\surd$) &20 \\
%3.5    & 0.1117 & 0  \\
%3.0    & 0.1174 &  0 \\
%2.5    & 0.1251 & 0 \\
%2.0    & 0.1259 & 0 \\
%\hline
%\end{tabular}
\\
\hskip-5pt
\includegraphics[width=.5\textwidth]{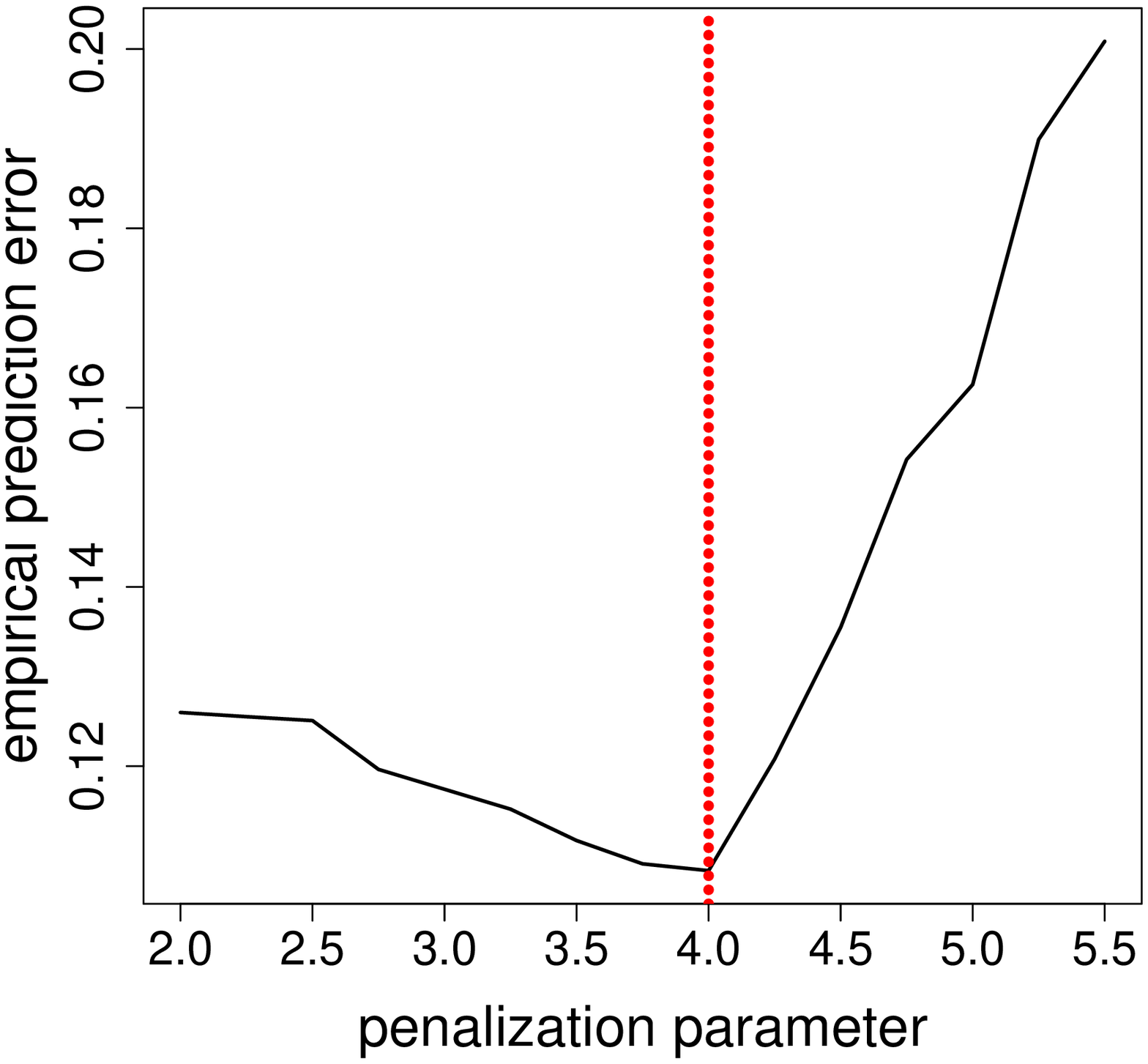}
\end{tabular}
\end{sc}
\end{small}
\end{center}
\caption{(Email spam) Classification accuracies and variable selection for
logistic SpAM.} \label{fig:lrspam} 
\end{figure}

The results of a typical run of logistic SpAM are summarized in
Figure~\ref{fig:lrspam}, using plug-in bandwidths. It is interesting
to note that even with this relatively small sample size, logistic
SpAM recovers a sparsity pattern that is consistent with previous
authors' results.  For example, in the best model chosen by logistic
SpAM, according to error rate, the 33 selected variables cover 80\% of
the significant predictors as determined by \cite{Hastie:2001}.

\vspace{.5cm}
{\em Functional Sparse Coding.}\enspace
\cite{Olshausen:Field:96} propose a method of obtaining
sparse representations of data such as natural images; the motivation
comes from trying to understand principles of neural coding.   In this
example we suggest a nonparametric form of sparse coding.  

Let
$\{y^{\i}\}_{i=1,\ldots, N}$ be the data to
be represented with respect to some learned basis, where each instance
$y^{\i}\in\reals^n$ is an $n$-dimensional vector.  The linear sparse
coding optimization problem is
\begin{eqnarray}
\label{eq:sparsecoding}
 \min_{\beta, X} && \sum_{i=1}^N \left\{\frac{1}{2n}\left\| y^\i - X \beta^\i \right\|^2 + \lambda
  \left\|\beta^\i\right\|_1\right\} \\
\text{such that} && \| X_j\| \leq 1
\end{eqnarray}
Here $X$ is an $n\times p$ matrix  with columns $X_j$, representing 
the ``dictionary'' entries or basis vectors to be learned.  It is
not required that the basis vectors are orthogonal.  The $\ell_1$
penalty on the coefficients $\beta^\i$ encourages sparsity, so that each data
vector $y^\i$ is represented by only a small number of dictionary
elements.  Sparsity allows 
the features to specialize, and to capture salient properties
of the data.

This optimization problem is not jointly convex in $\beta^\i$ and $X$.
However, for fixed $X$, each weight vector $\beta^\i$ is computed
by running the lasso.  
For fixed $\beta^\i$, the optimization is similar to ridge regression, and 
can be solved efficiently.  Thus, an iterative procedure for 
(approximately) solving this optimization problem is easy to derive.

\begin{figure}[t]
\begin{center}
\begin{tabular}{c}
\hskip-5pt
\includegraphics[width=.95\textwidth]{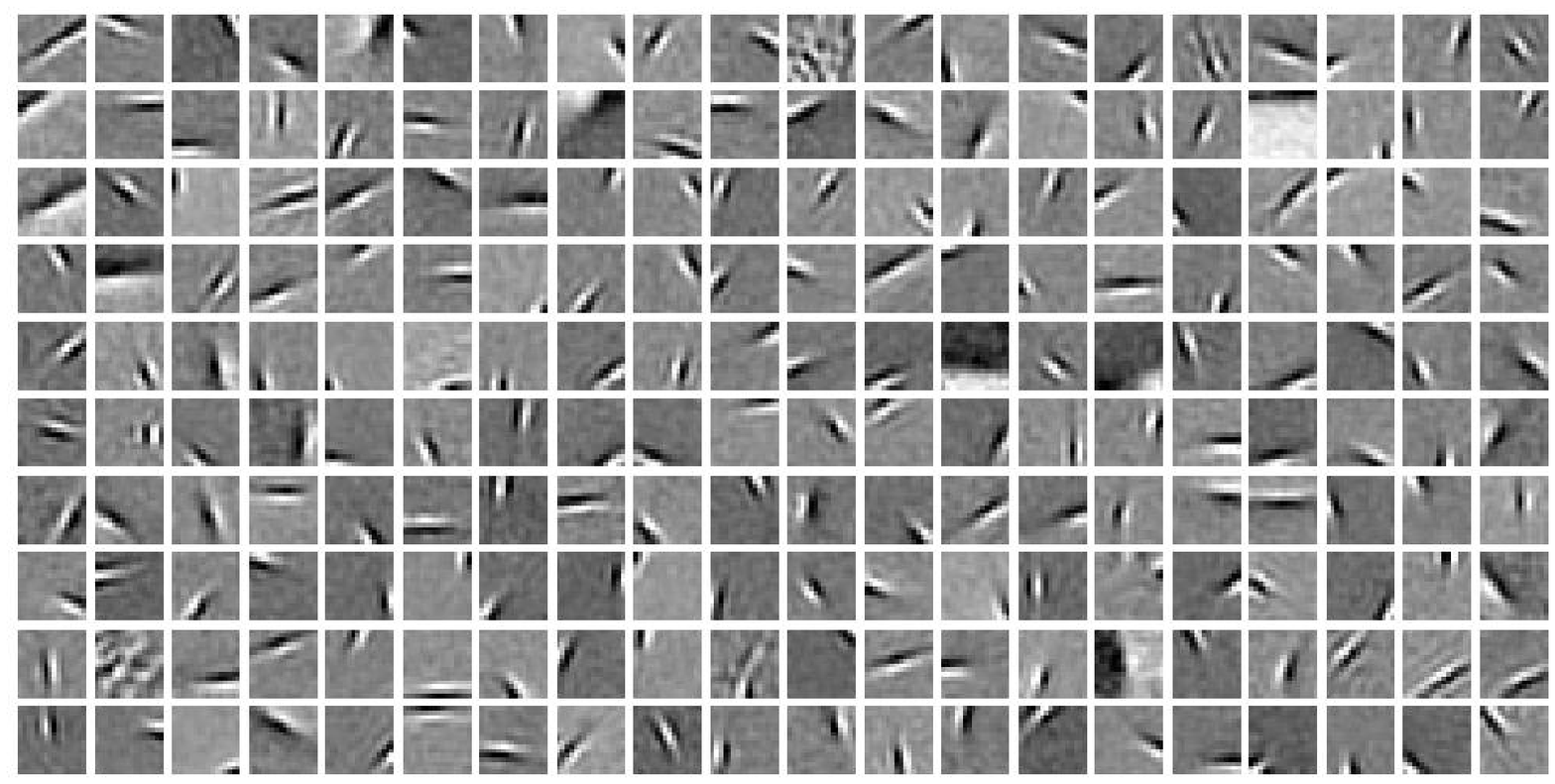} 
\end{tabular}
\end{center}
\caption{A sparse code with 200 codewords, trained on 200,000 patches
  from natural images, each patch $16\times 16$ pixels in size, using
  the lasso and stochastic gradient descent.  The codewords are seen
  to capture edge features at different scales and spatial
  orientations.}
\label{fig:codex}
\end{figure}

In the case of sparse coding of natural images, as in
\cite{Olshausen:Field:96}, the basis vectors $X_j$ encode basic edge
features.  A code with 200 basis vectors, estimated by carrying out
the optimization using the lasso and stochastic gradient descent, is
shown in Figure~\ref{fig:codex}.  The codewords are seen to capture
edge features at different scales and spatial orientations.

In the functional version, 
we no longer assume
a linear, parametric fit between the dictionary $X$ and the data $y$.  Instead,
we model the relationship using an additive model:
\begin{eqnarray}
y^\i = \sum_{j=1}^p \beta_j^\i f_j^\i(X_j) + \epsilon^\i
\end{eqnarray}
where $X_j\in \reals^p$ is a dictionary vector and $\epsilon^\i\in\reals^n$
is a noise vector.  This leads to the following optimization problem
for functional sparse coding:
\begin{eqnarray}
 \min_{f, X} && \sum_{i=1}^N \left\{\frac{1}{2n}\left\| y^\i - \textstyle \sum_{j=1}^p
     f_j^\i (X_j) \right\|^2 + \lambda \sum_{j=1}^p
  \left\|f_j^\i\right\|\right\} \\
\text{such that} && \| X_j\| \leq 1,\;j=1,\ldots, p .
\end{eqnarray}

\begin{figure*}[ht]
\begin{center}
\begin{tabular}{c|c}
\begin{tabular}{c}
\begin{tabular}{cc}
\multicolumn{2}{c}{Lasso} \\
\small Original patch & \small $\text{RSS}=0.0561$ \\
\hskip-5pt
\includegraphics[width=\headwidth]{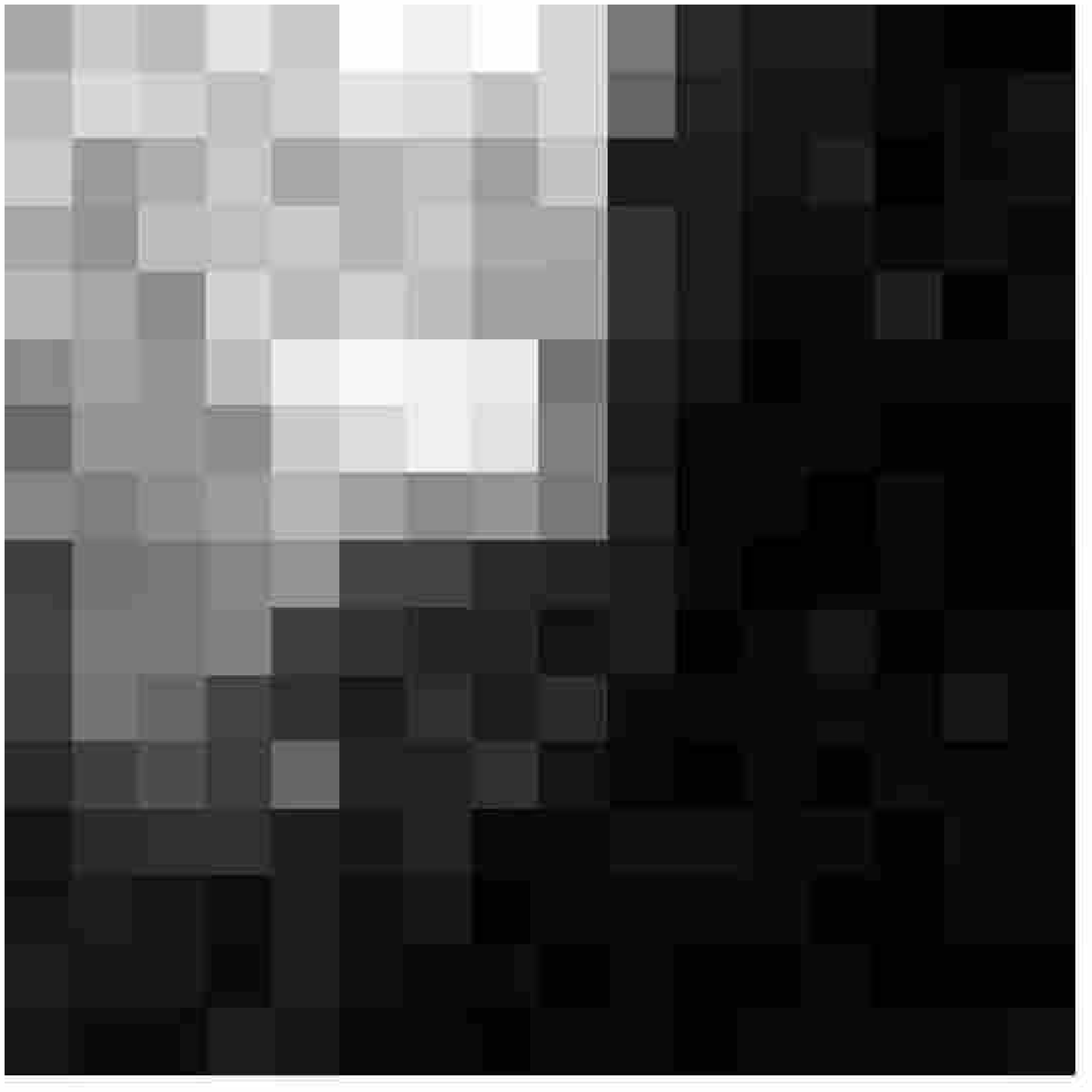} &
\includegraphics[width=\headwidth]{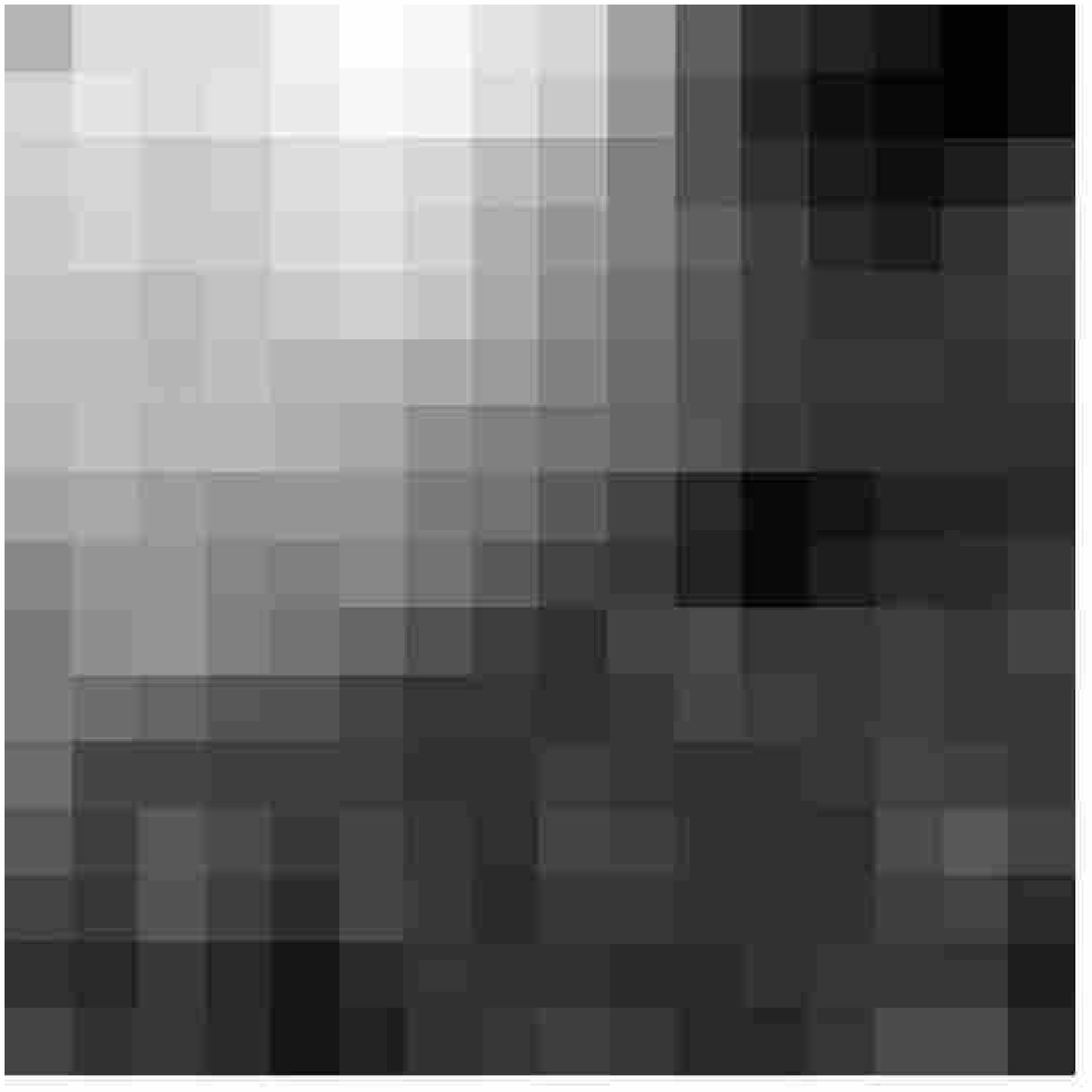}
\\[10pt]
\end{tabular}
\\
\hskip-14pt
\begin{tabular}{cccc}
\includegraphics[width=\plotwidth]{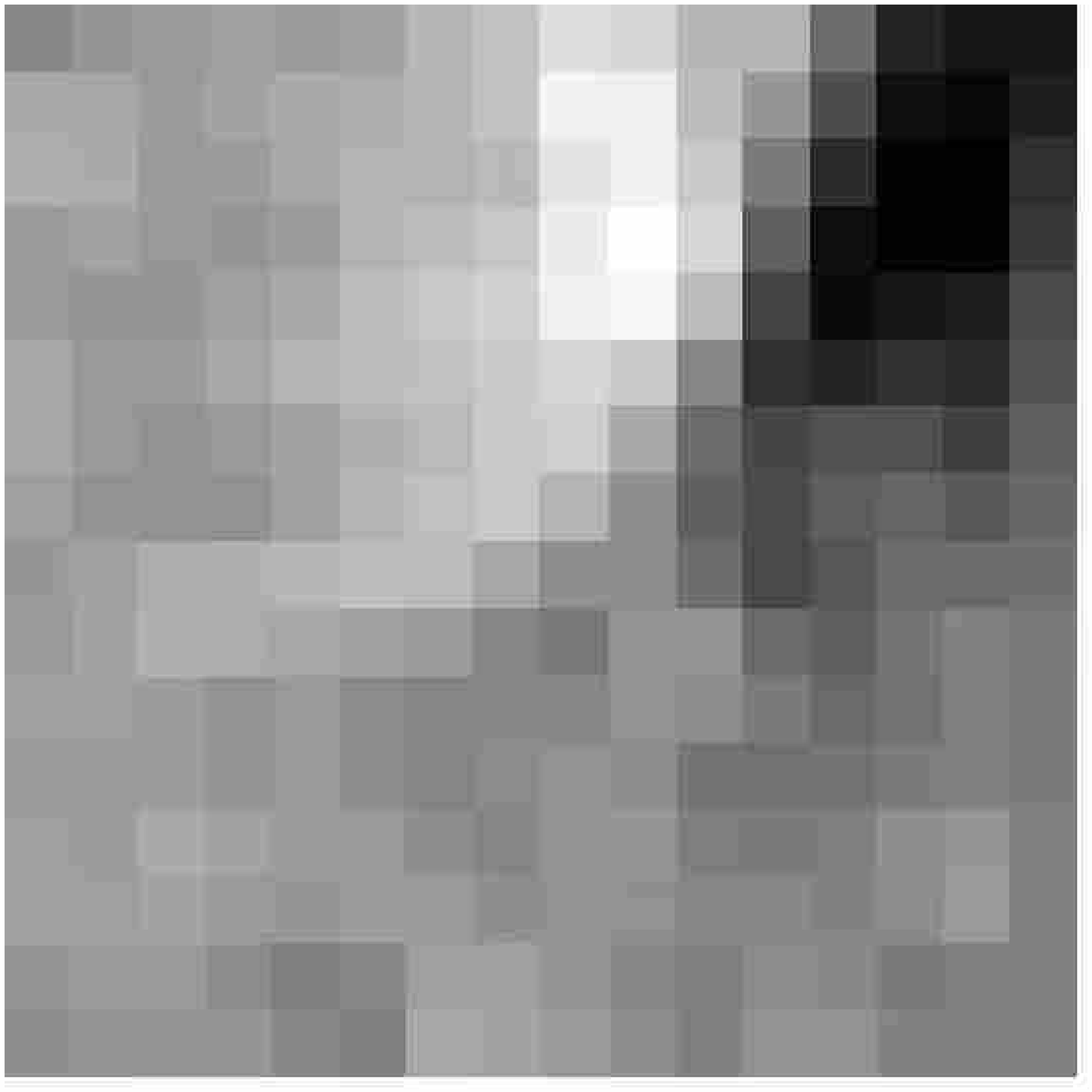}&
\includegraphics[width=\plotwidth]{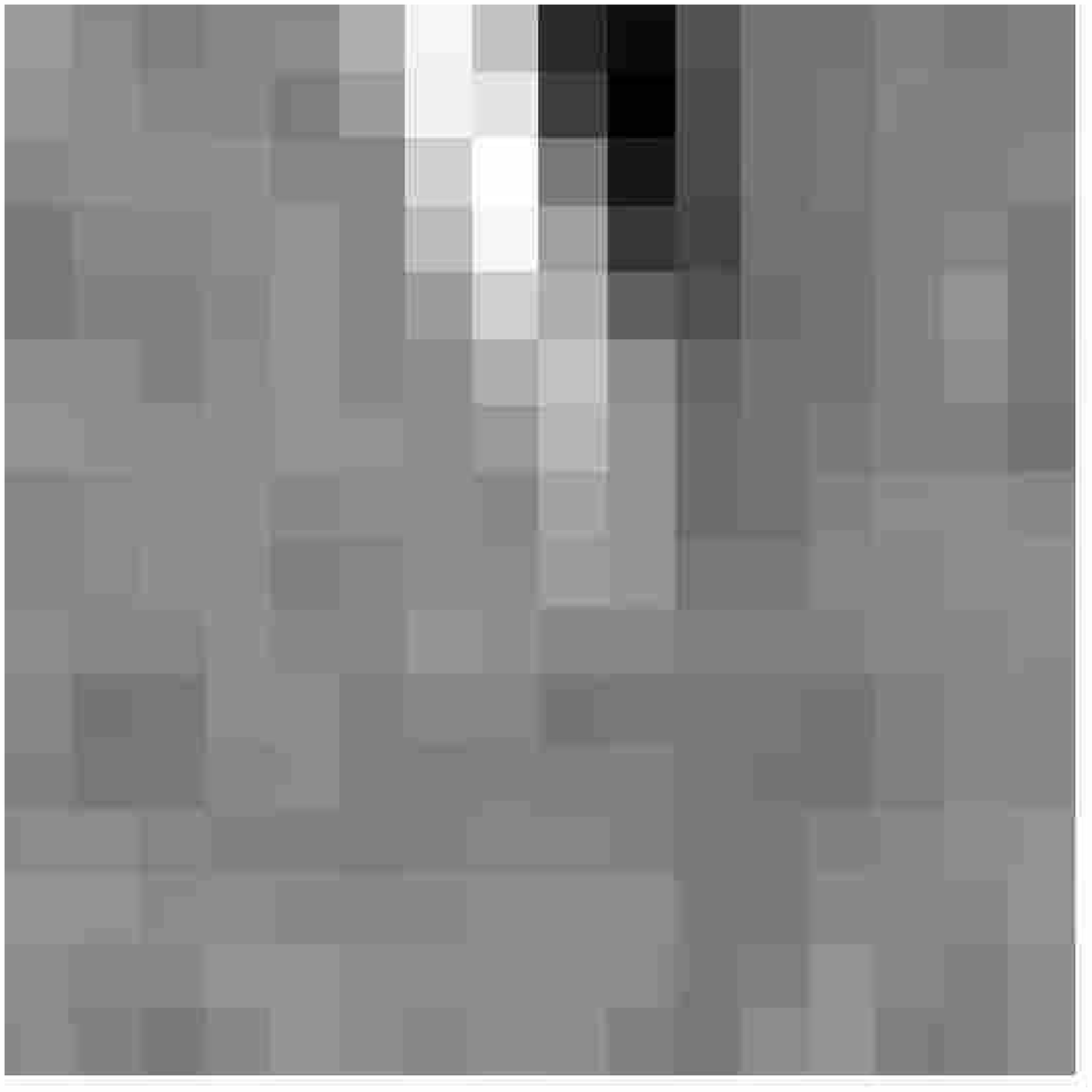}&
\includegraphics[width=\plotwidth]{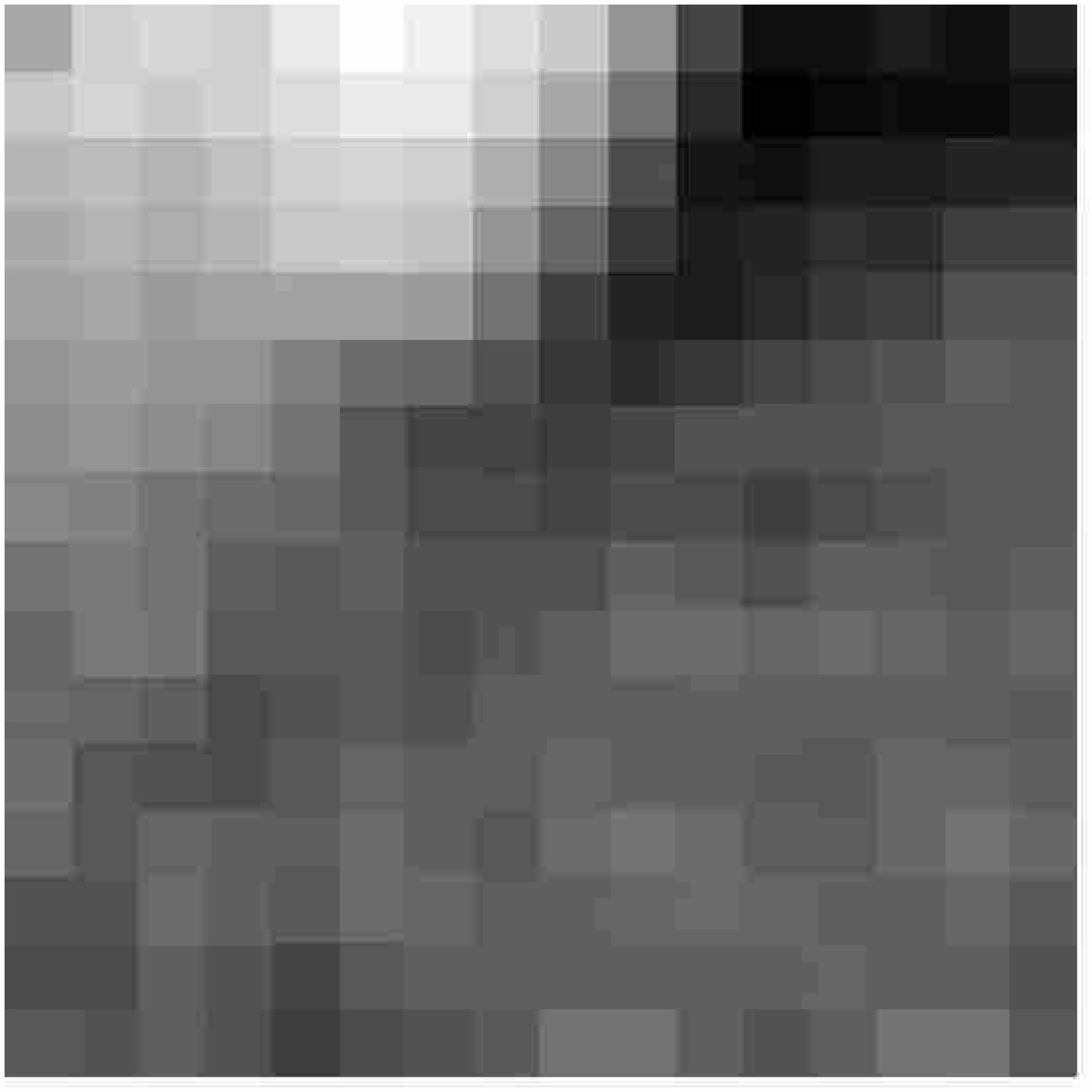}&
\includegraphics[width=\plotwidth]{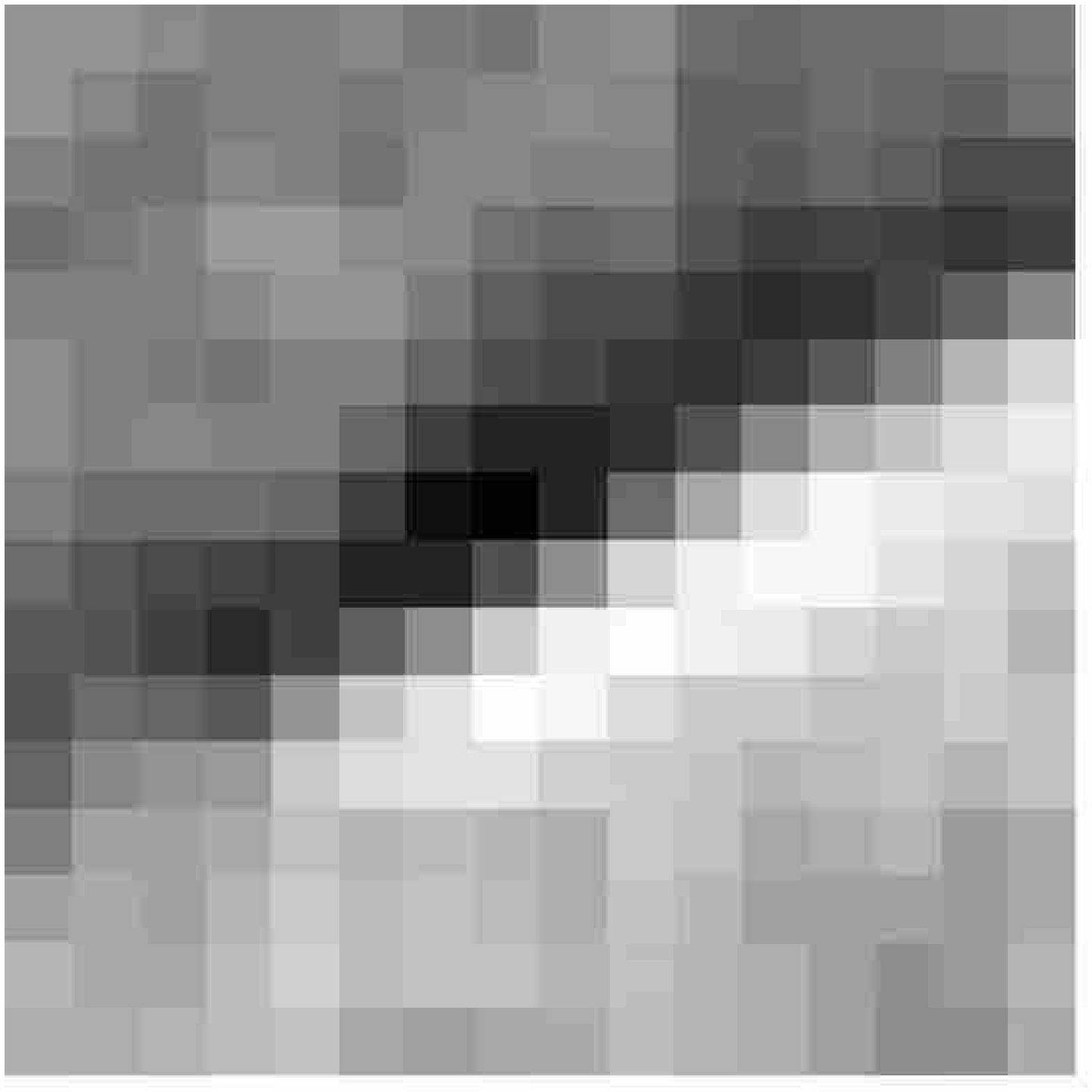} \\
\includegraphics[width=\plotwidth]{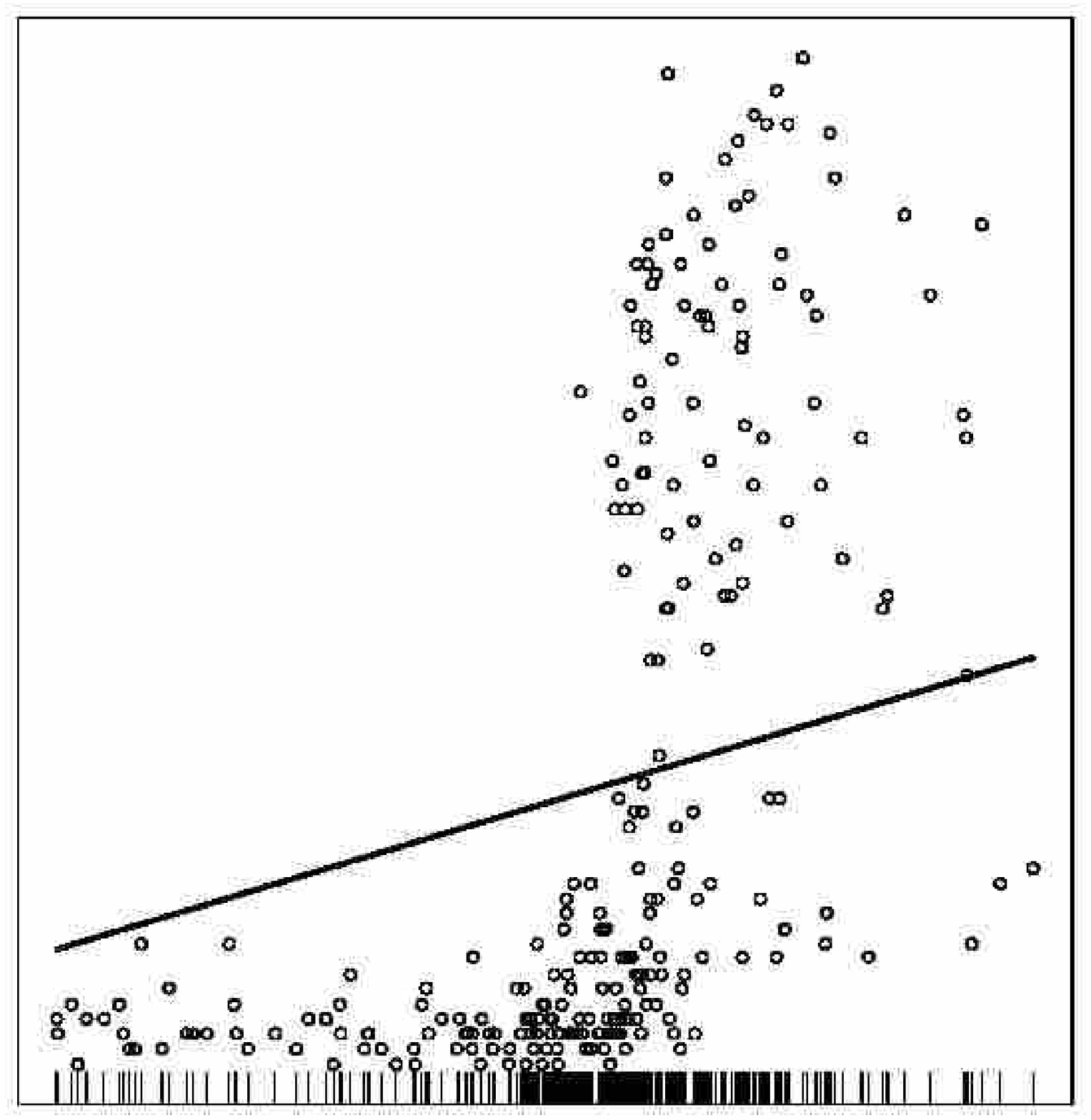}&
\includegraphics[width=\plotwidth]{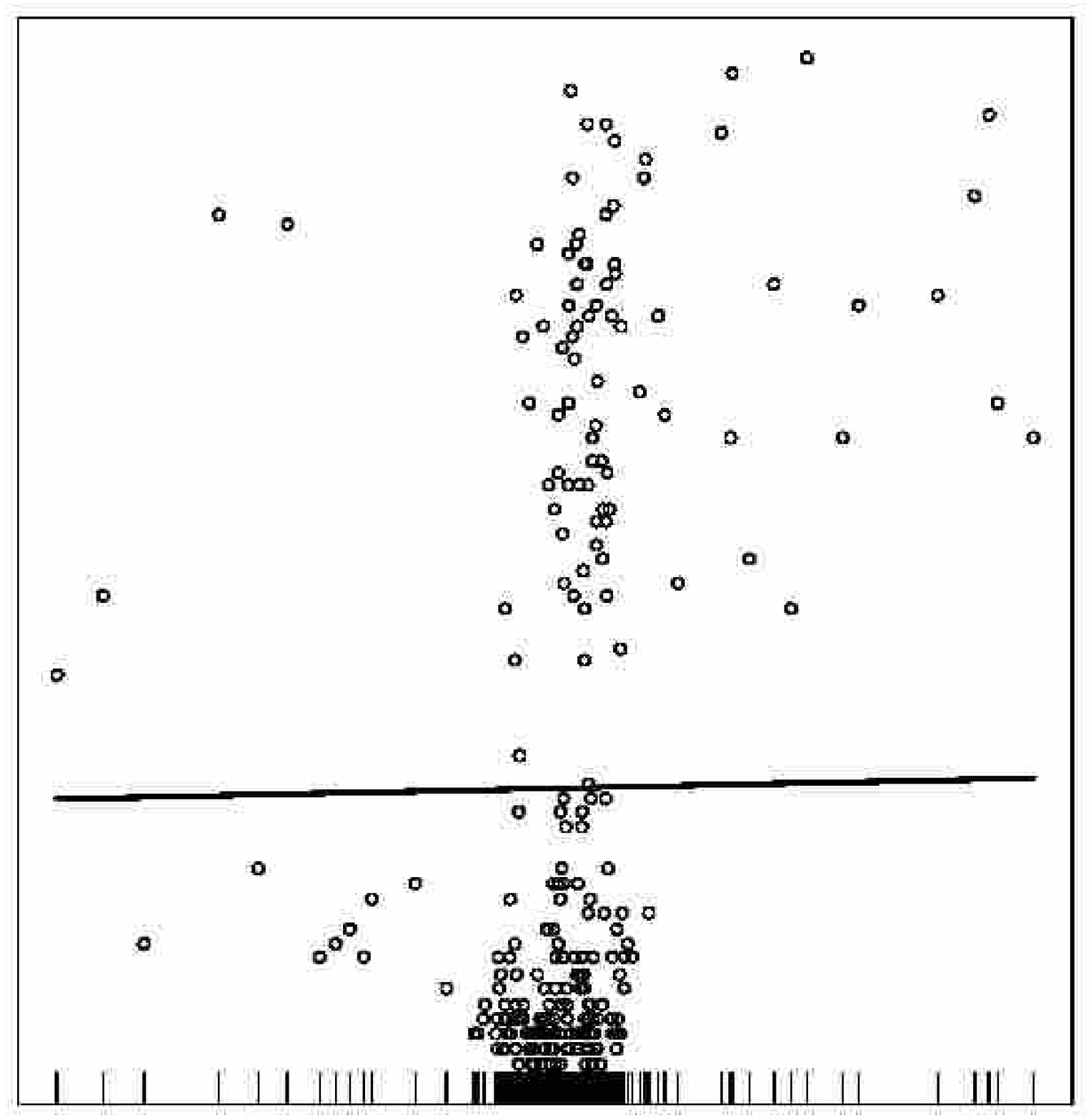}&
\includegraphics[width=\plotwidth]{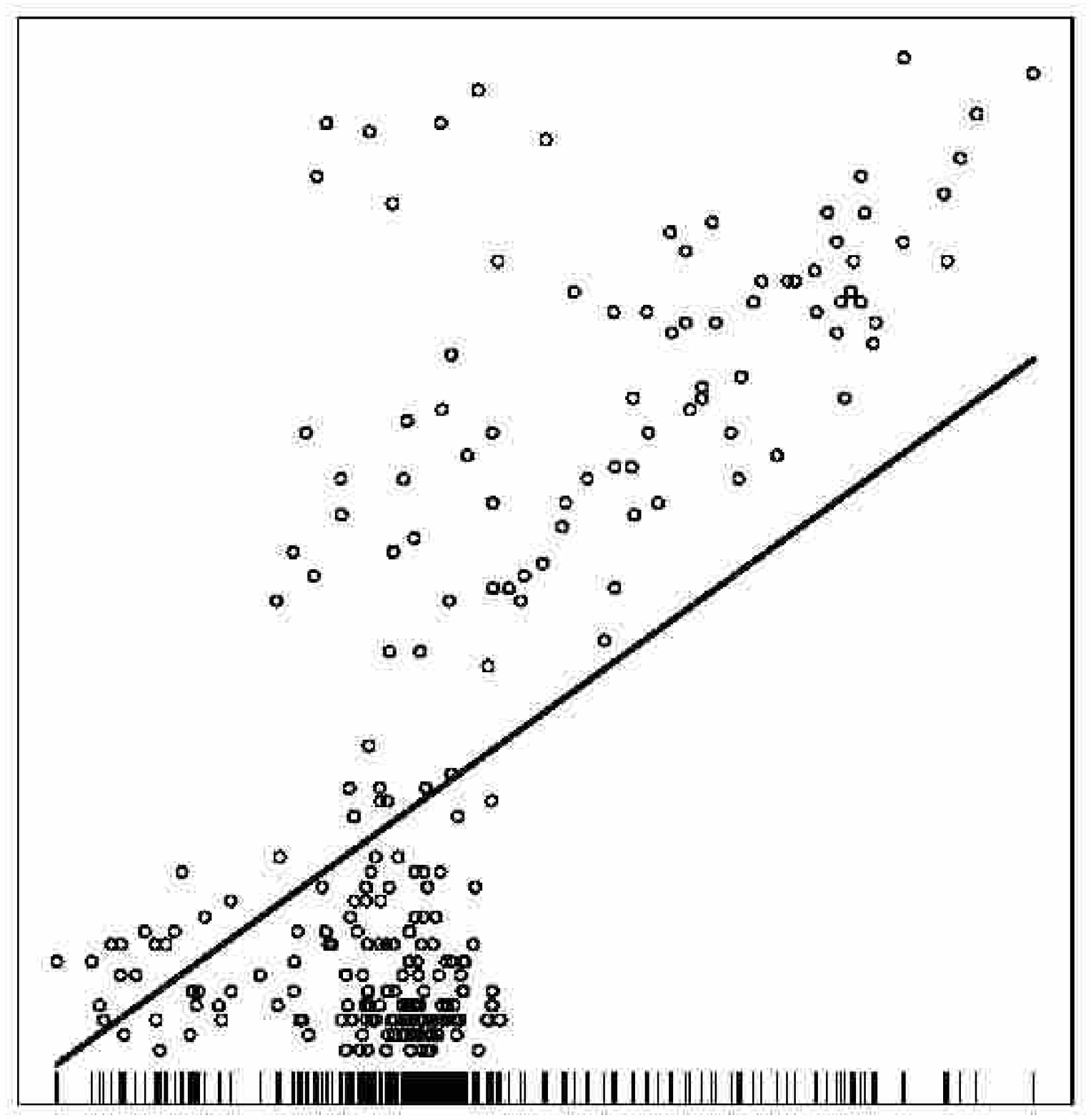}&
\includegraphics[width=\plotwidth]{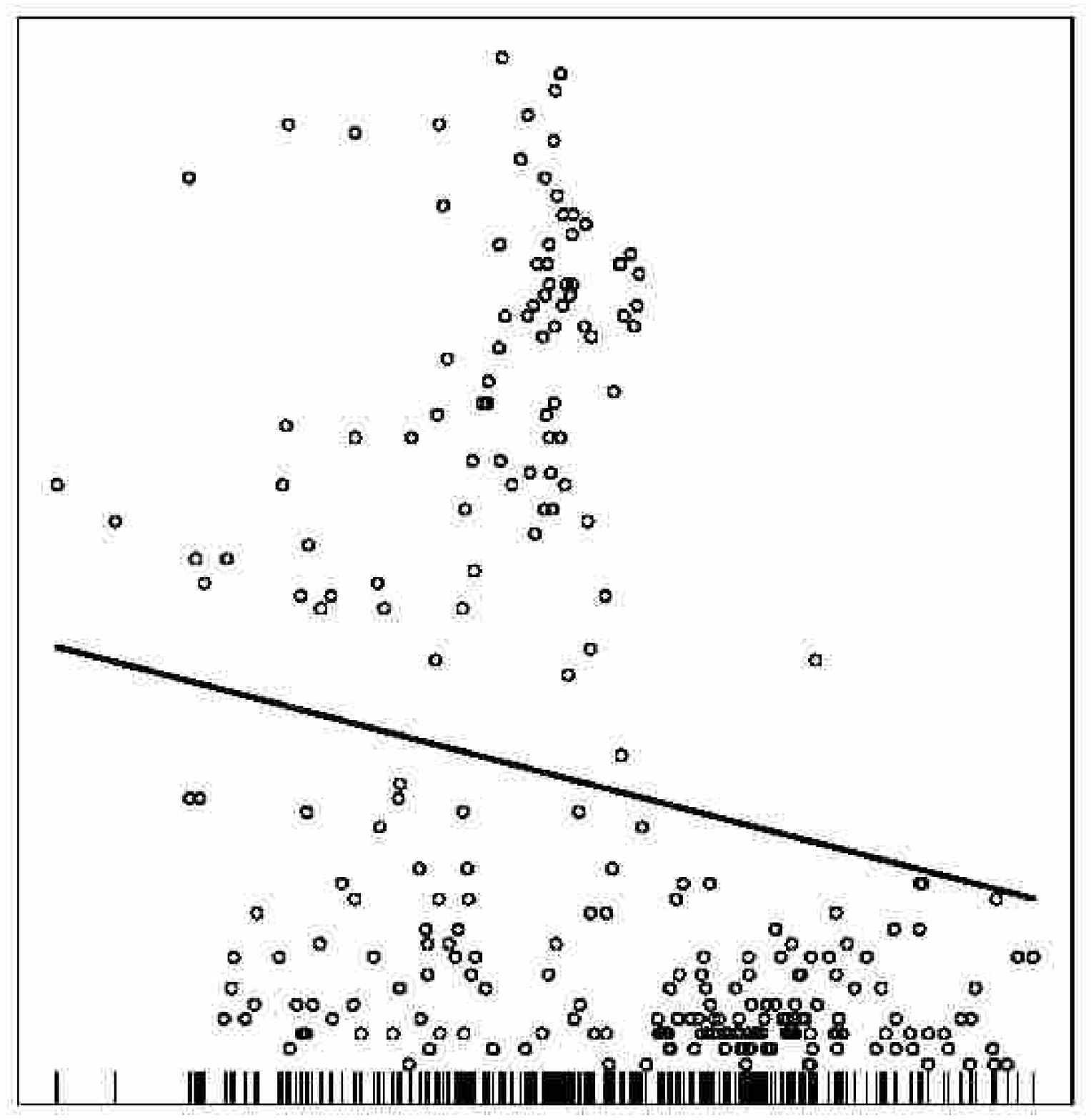} \\
\includegraphics[width=\plotwidth]{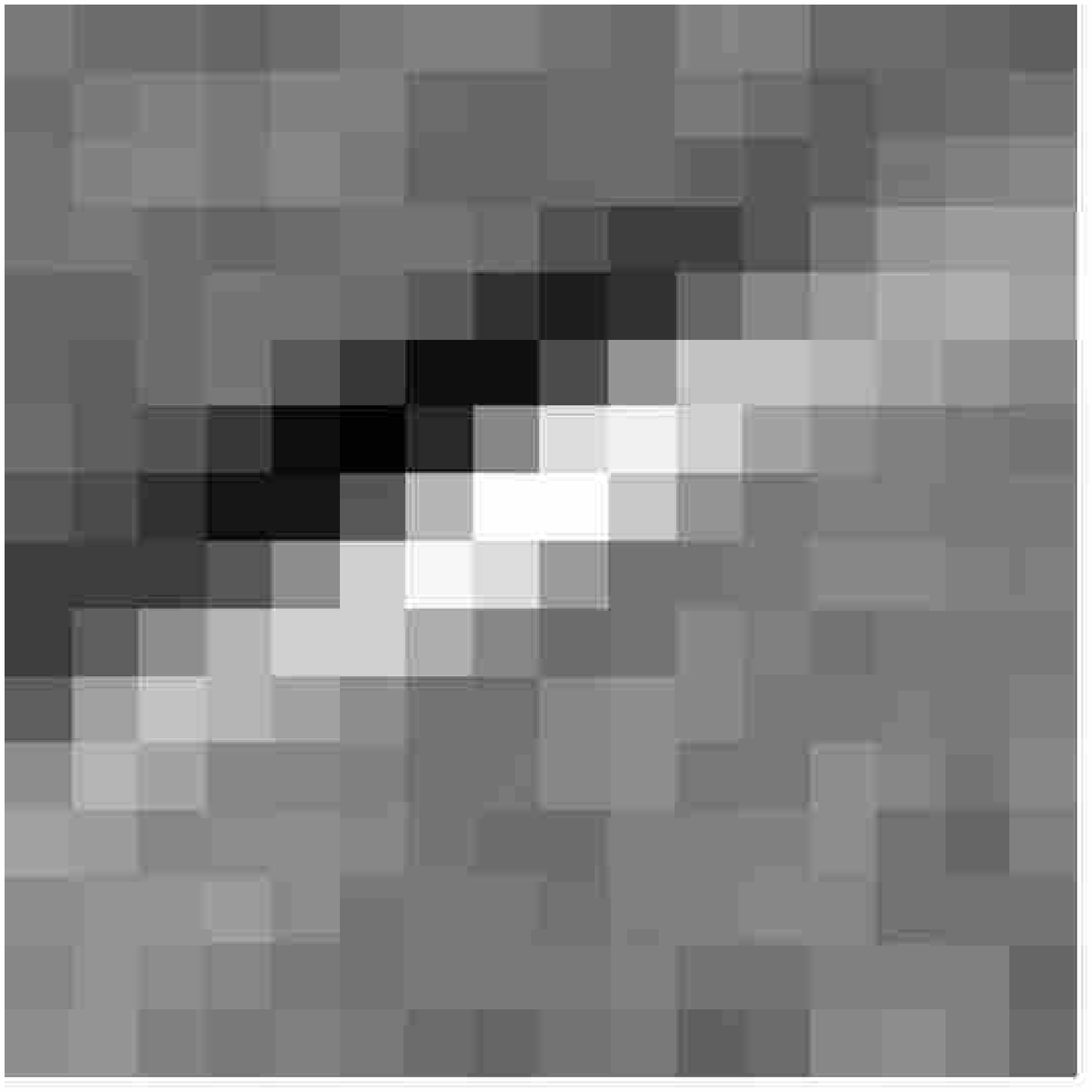}&
\includegraphics[width=\plotwidth]{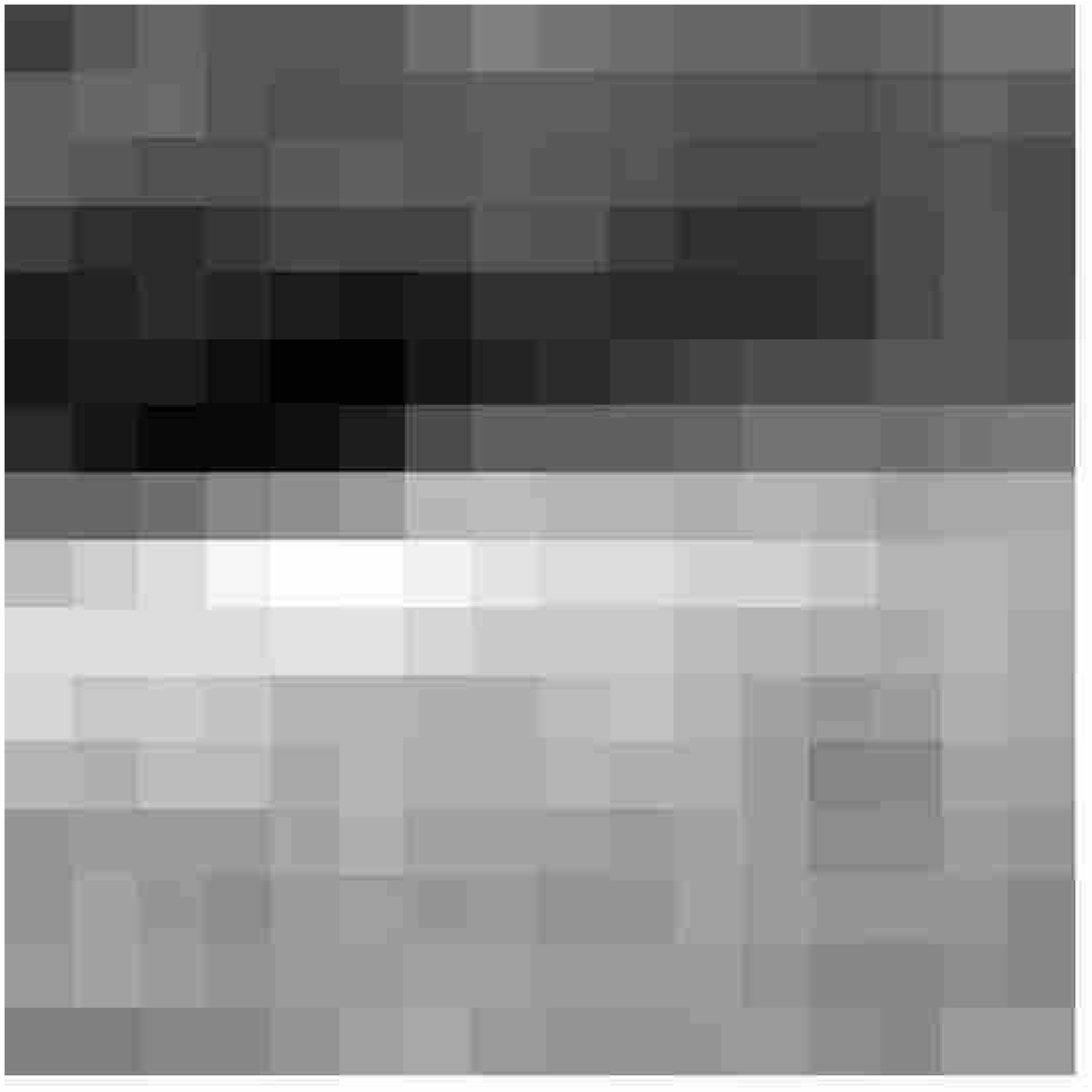}&
\includegraphics[width=\plotwidth]{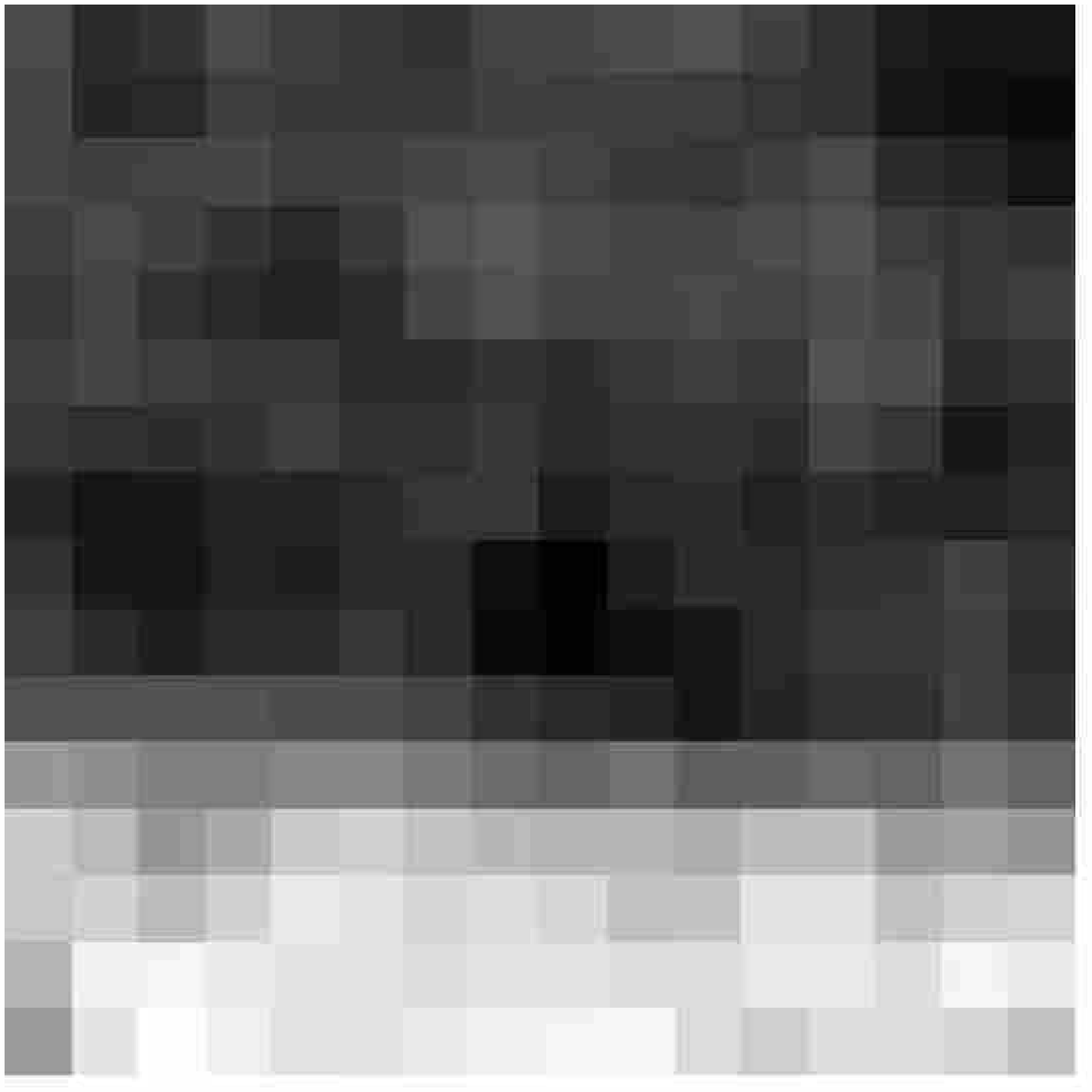}&
\includegraphics[width=\plotwidth]{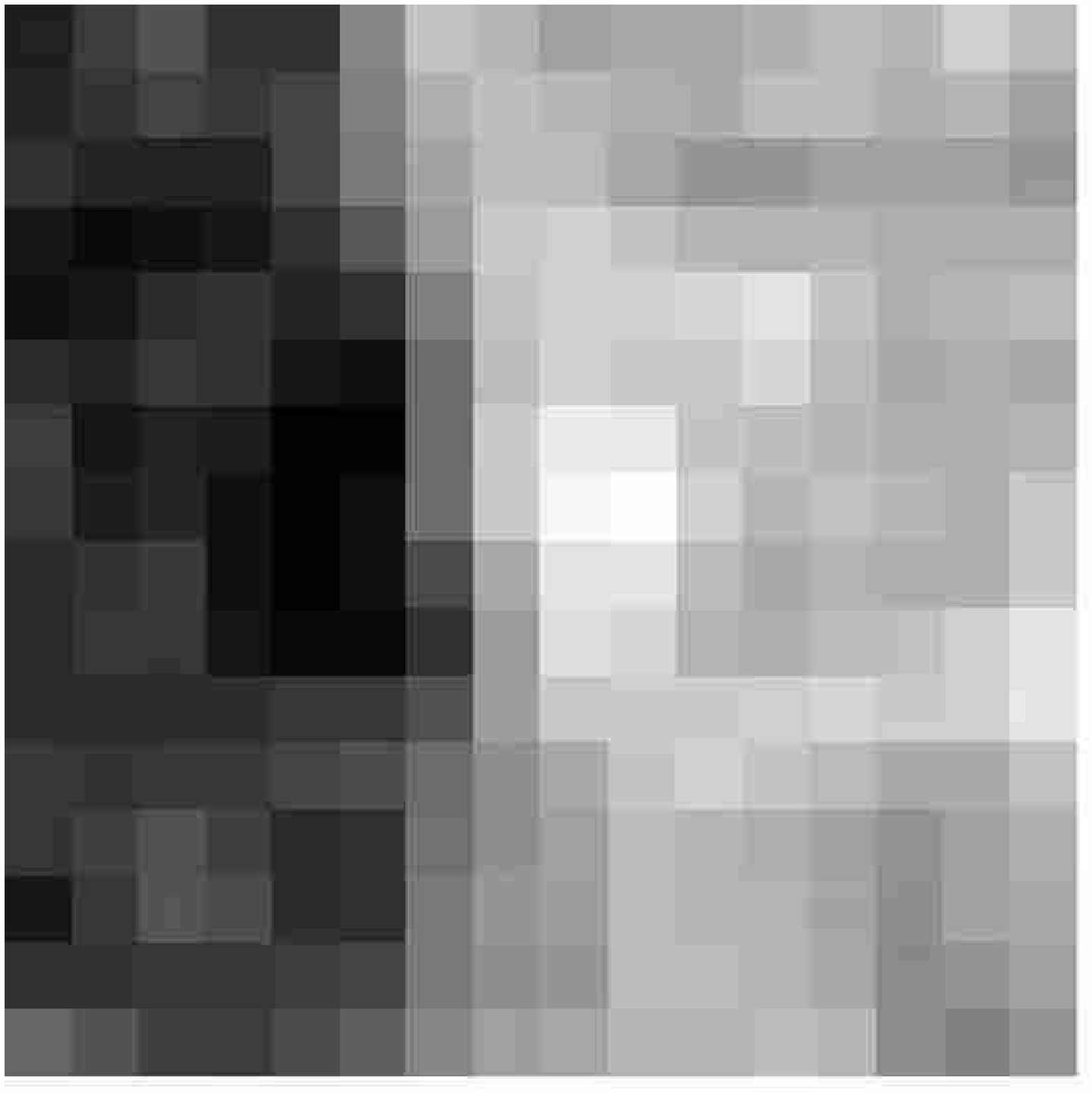}\\
\includegraphics[width=\plotwidth]{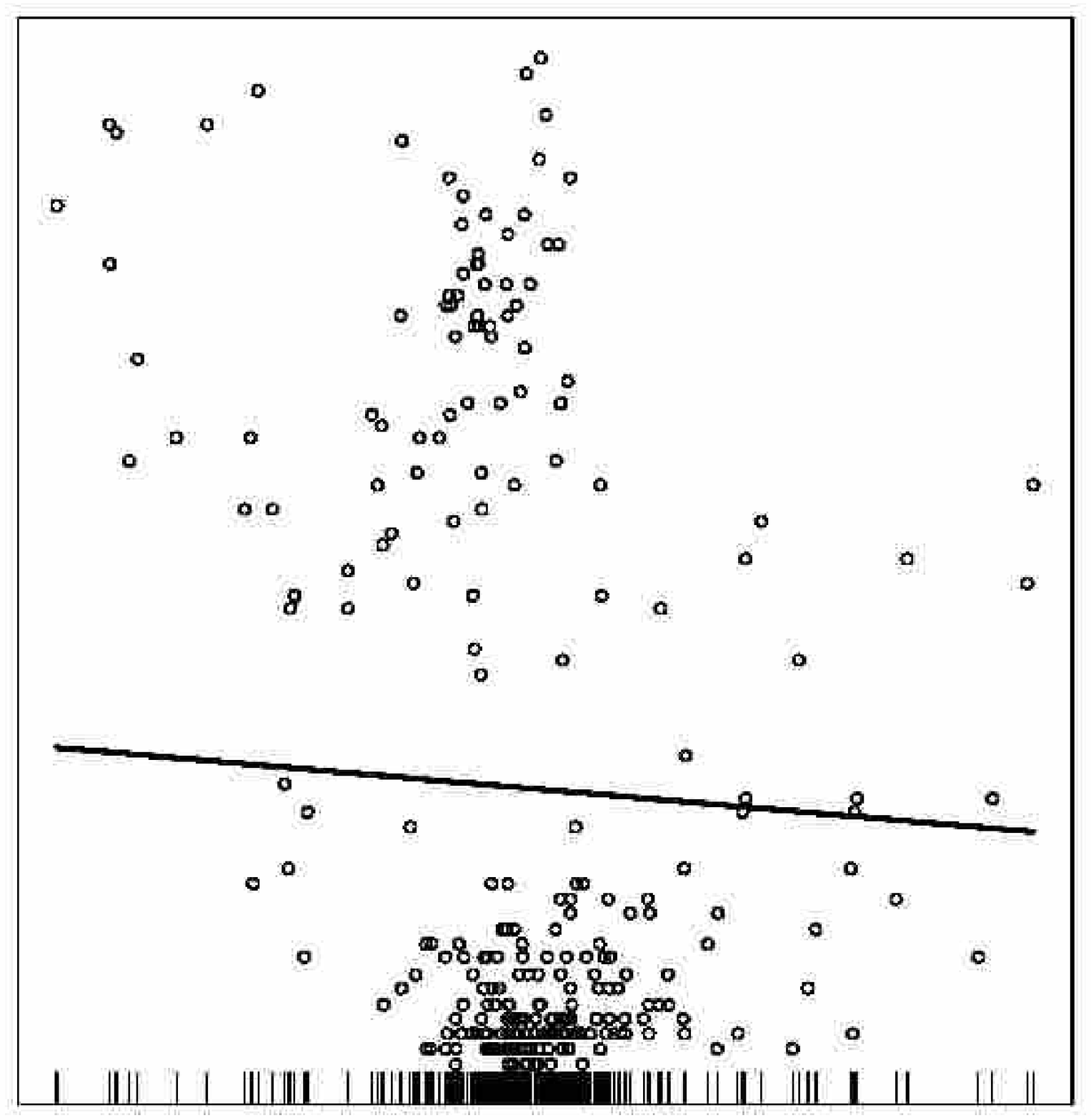}&
\includegraphics[width=\plotwidth]{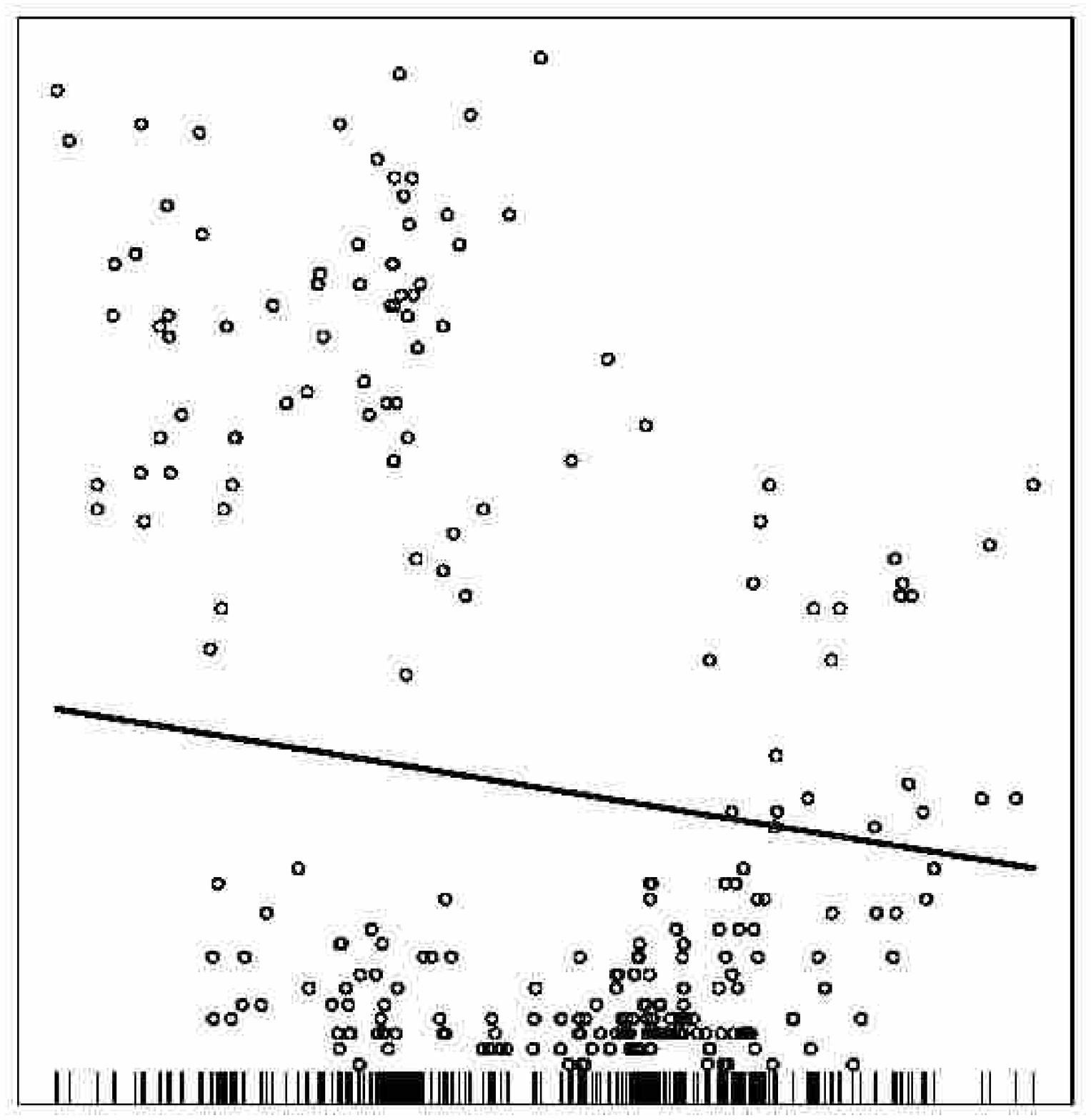}&
\includegraphics[width=\plotwidth]{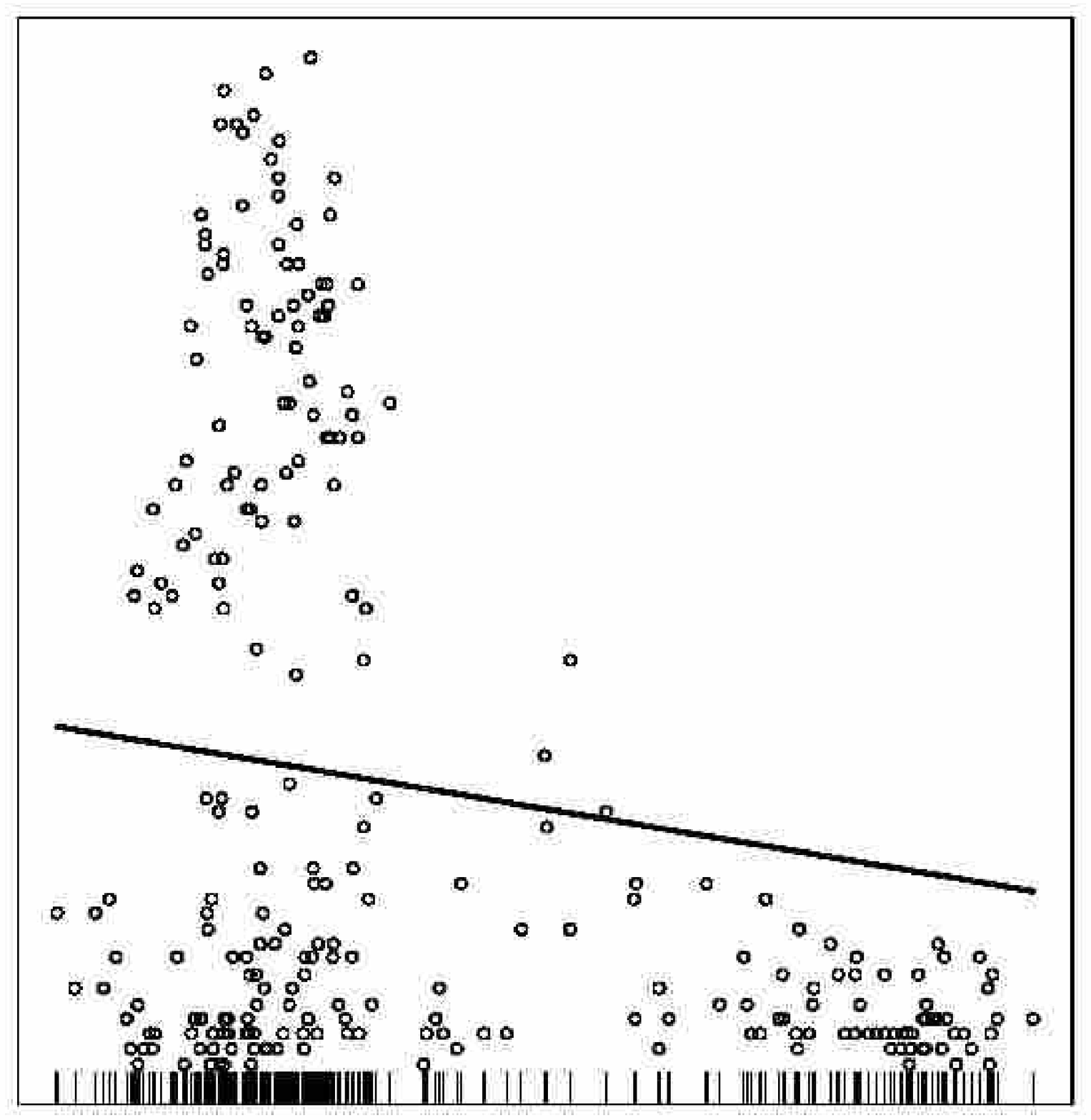} &
\includegraphics[width=\plotwidth]{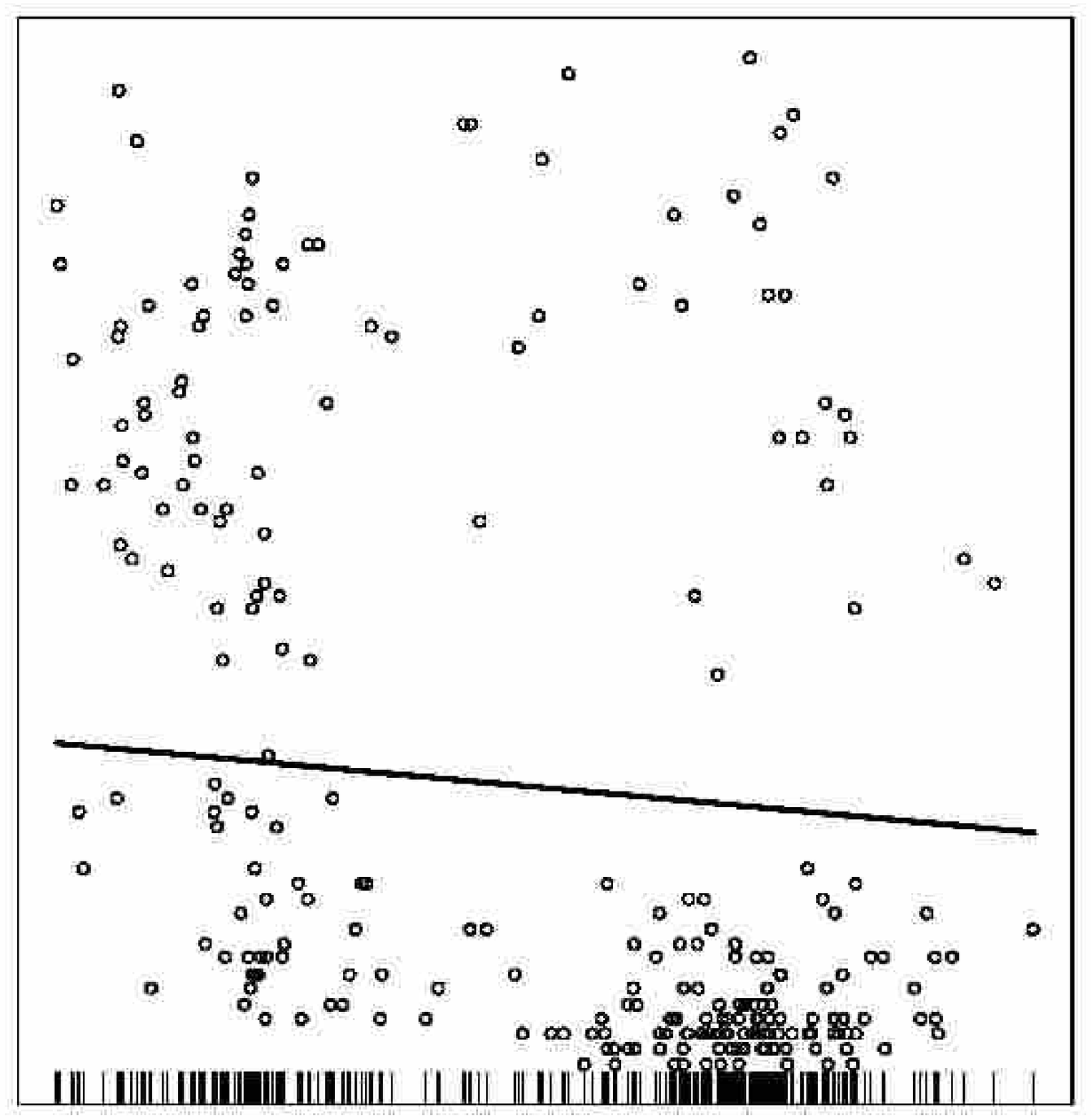}\\[10pt]
\end{tabular}
\end{tabular}
&
\begin{tabular}{c}
\begin{tabular}{cc}
\multicolumn{2}{c}{SpAM} \\
\small Original patch & \small $\text{RSS}=0.0206$ \\
\hskip-5pt
\includegraphics[width=\headwidth]{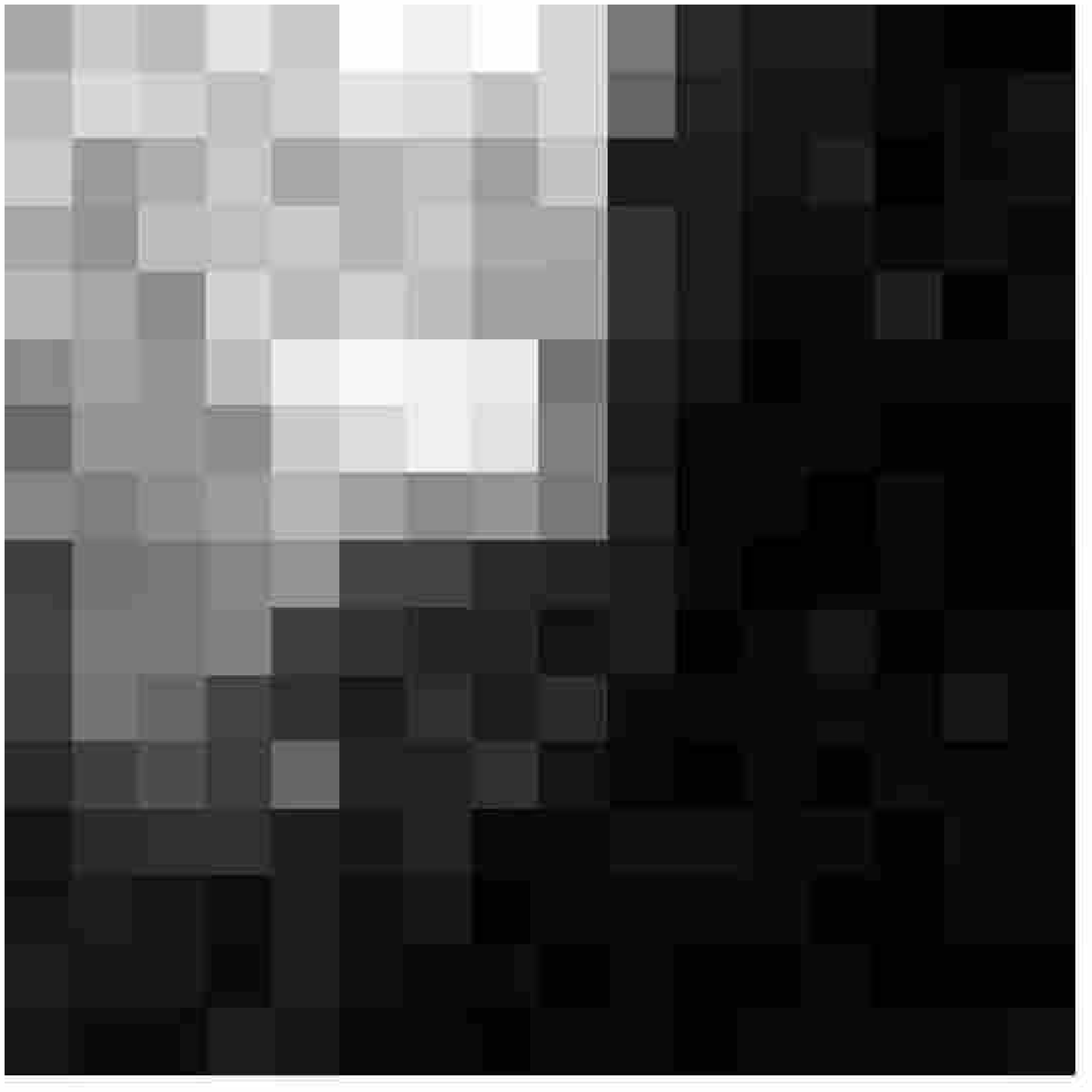} &
\includegraphics[width=\headwidth]{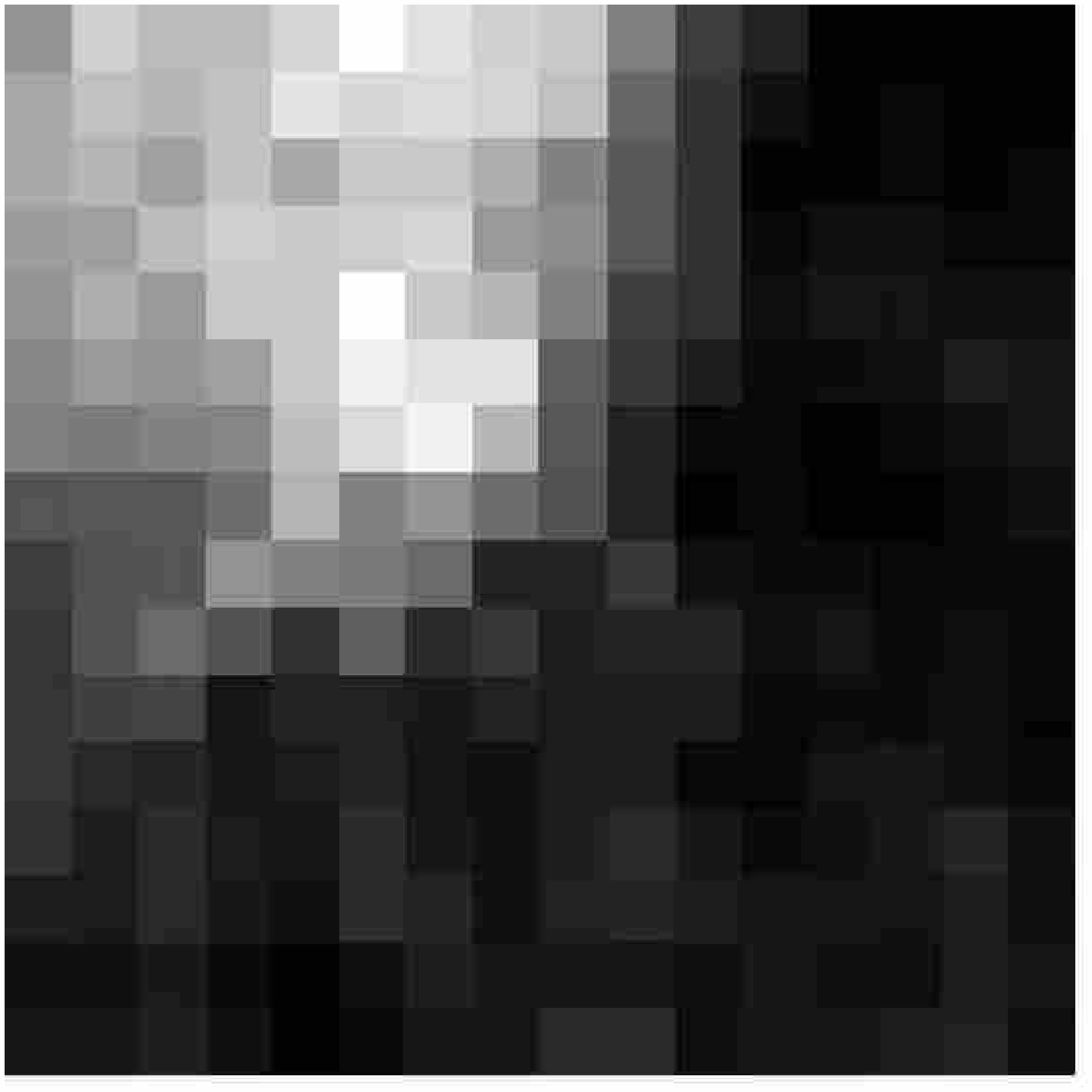}
\\[10pt]
\end{tabular}
\\
\hskip-14pt
\begin{tabular}{cccc}
\includegraphics[width=\plotwidth]{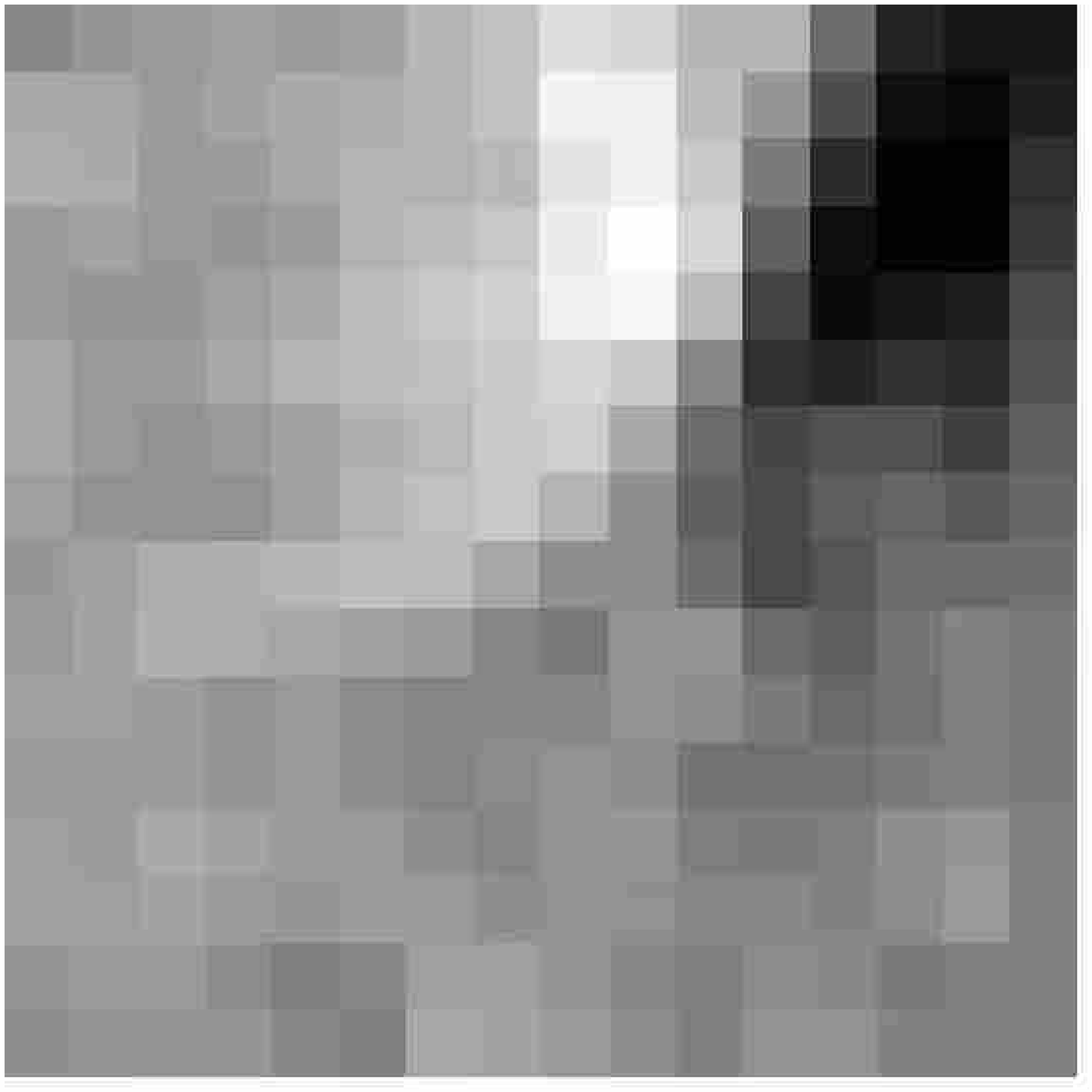}&
\includegraphics[width=\plotwidth]{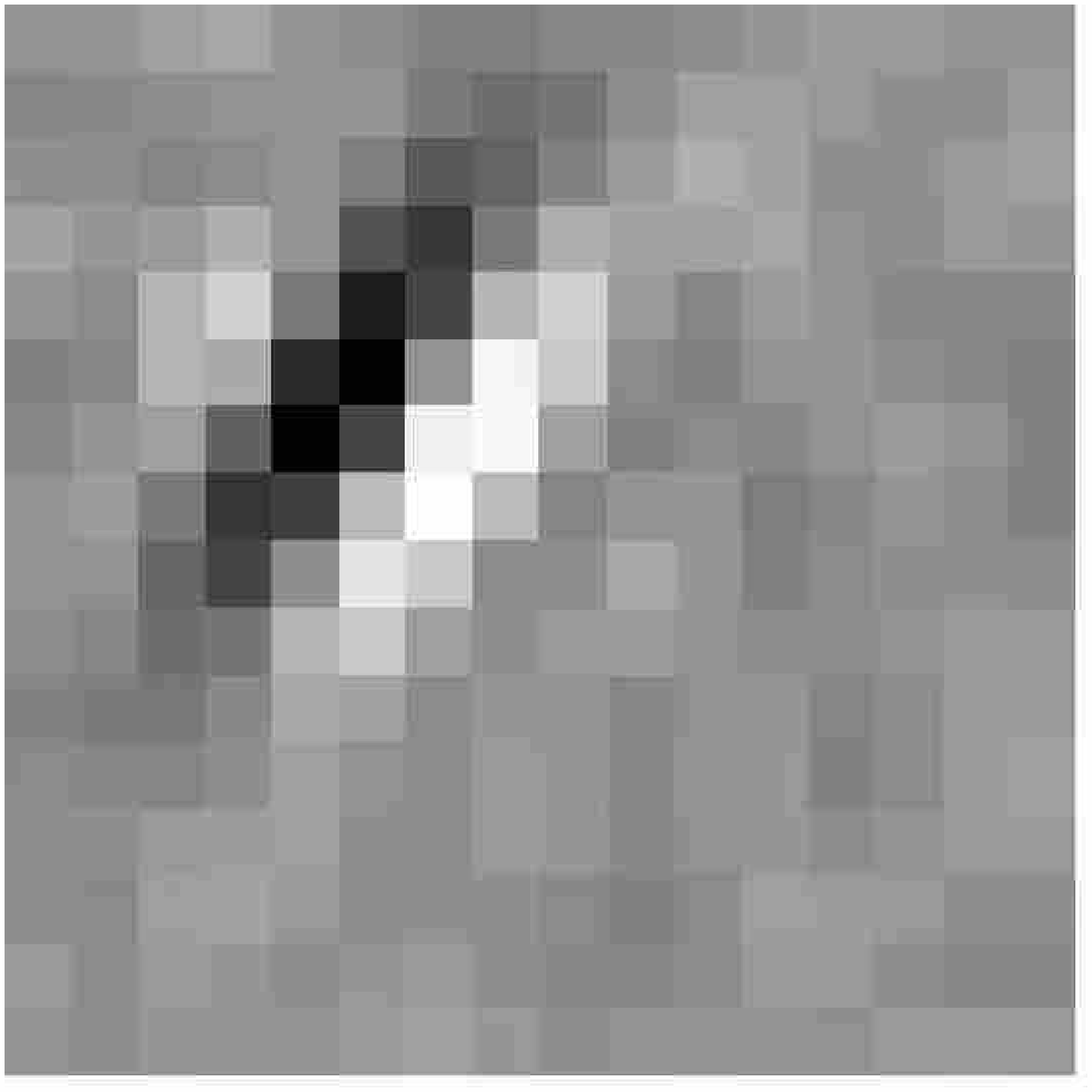}&
\includegraphics[width=\plotwidth]{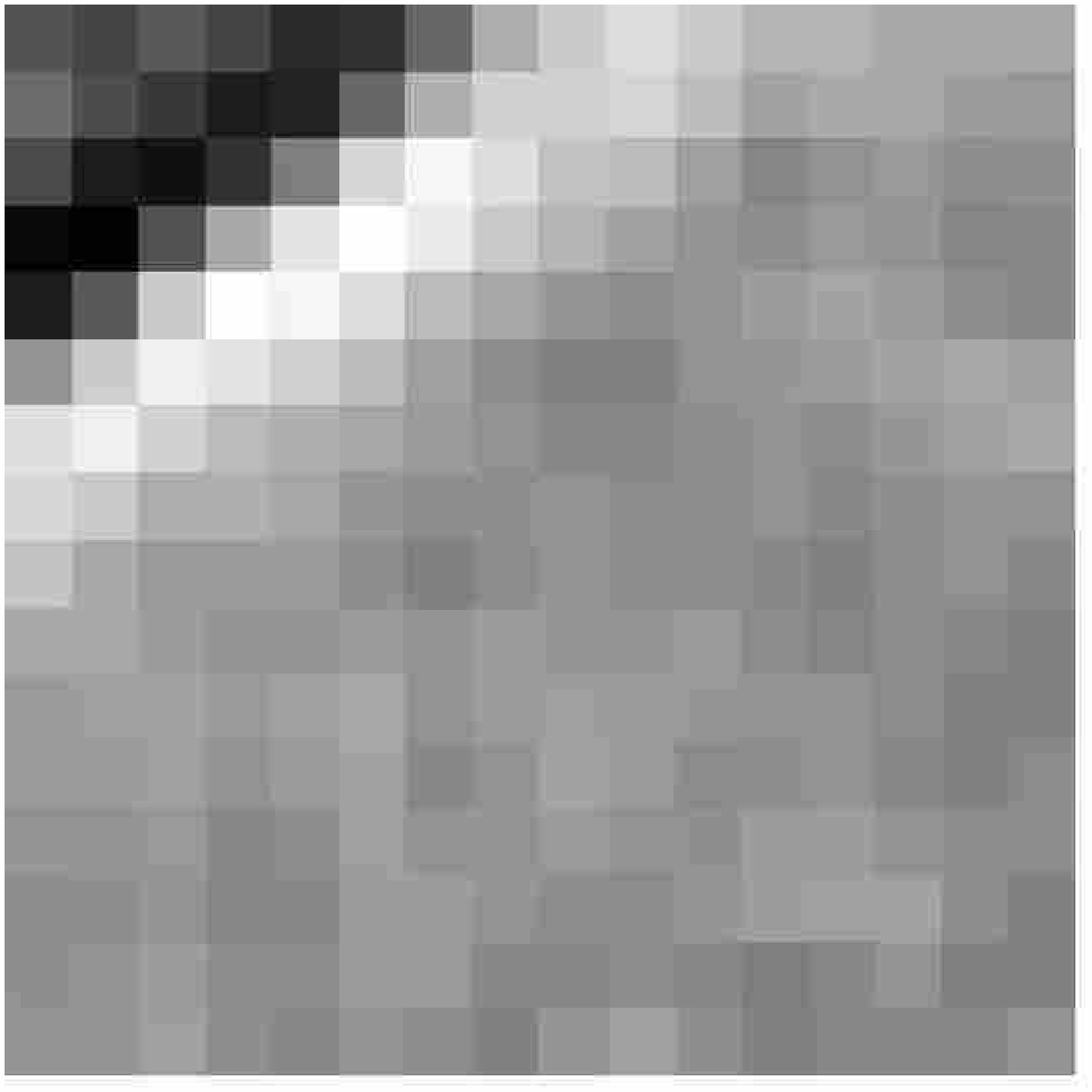}&
\includegraphics[width=\plotwidth]{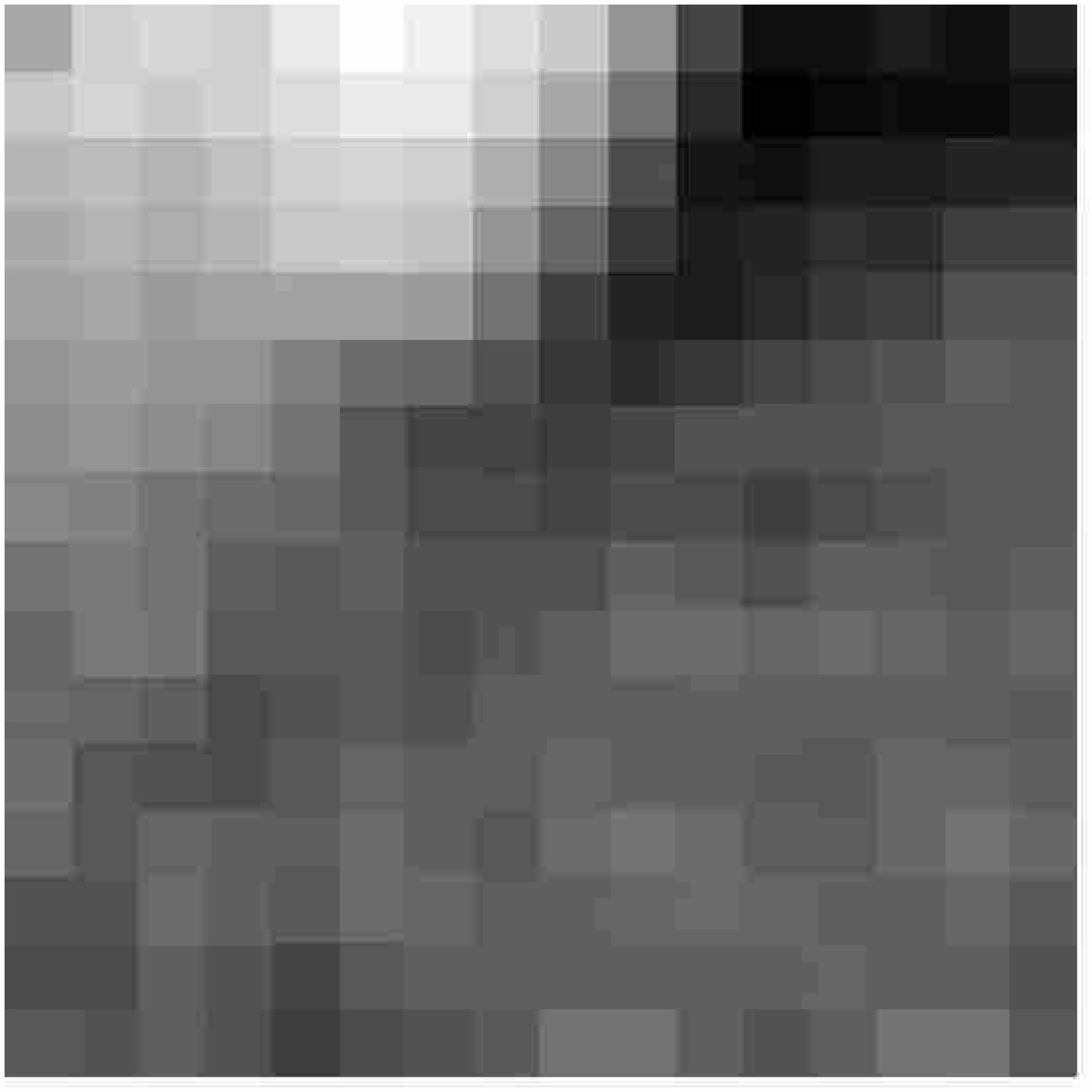}\\
\includegraphics[width=\plotwidth]{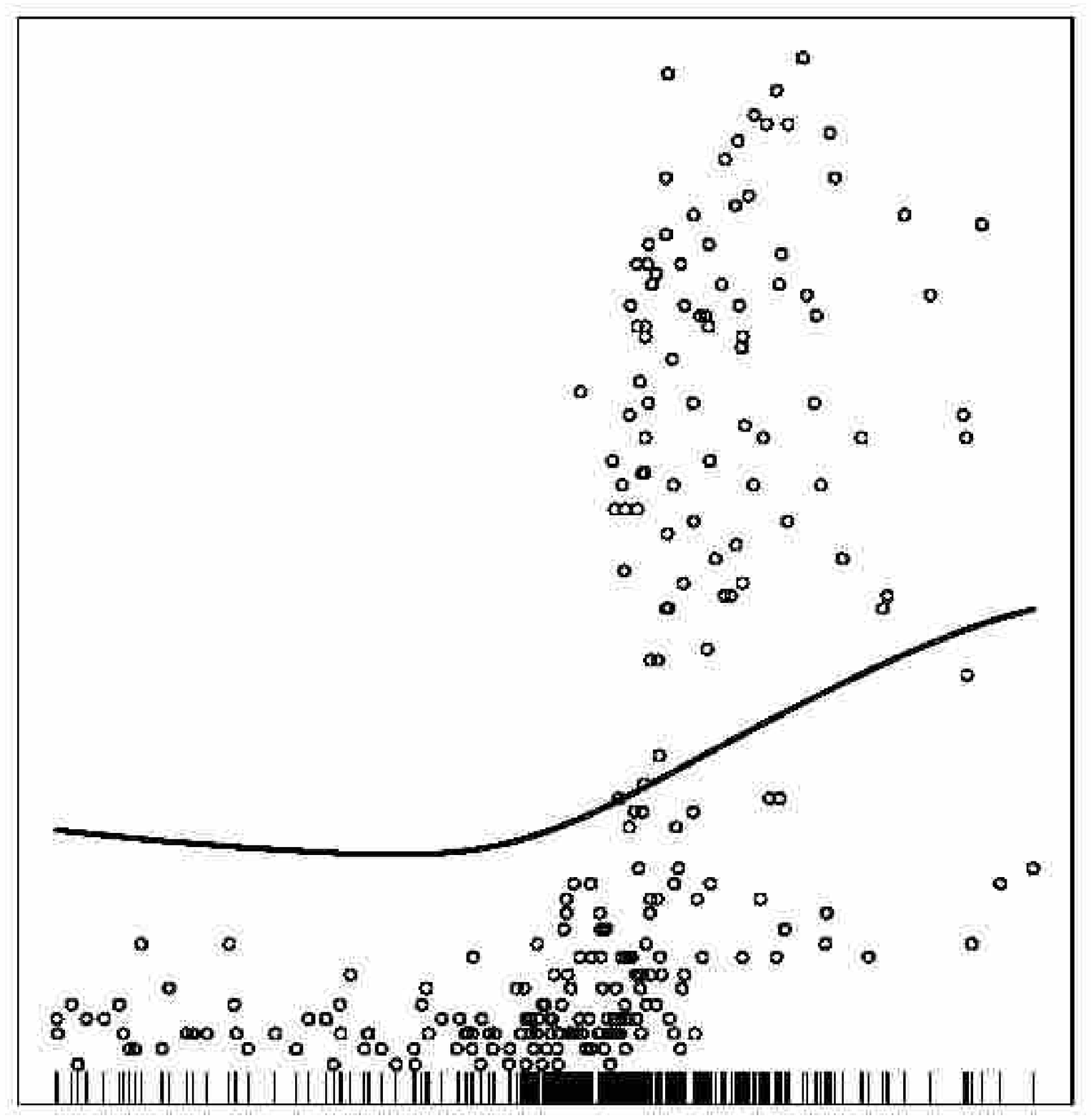}&
\includegraphics[width=\plotwidth]{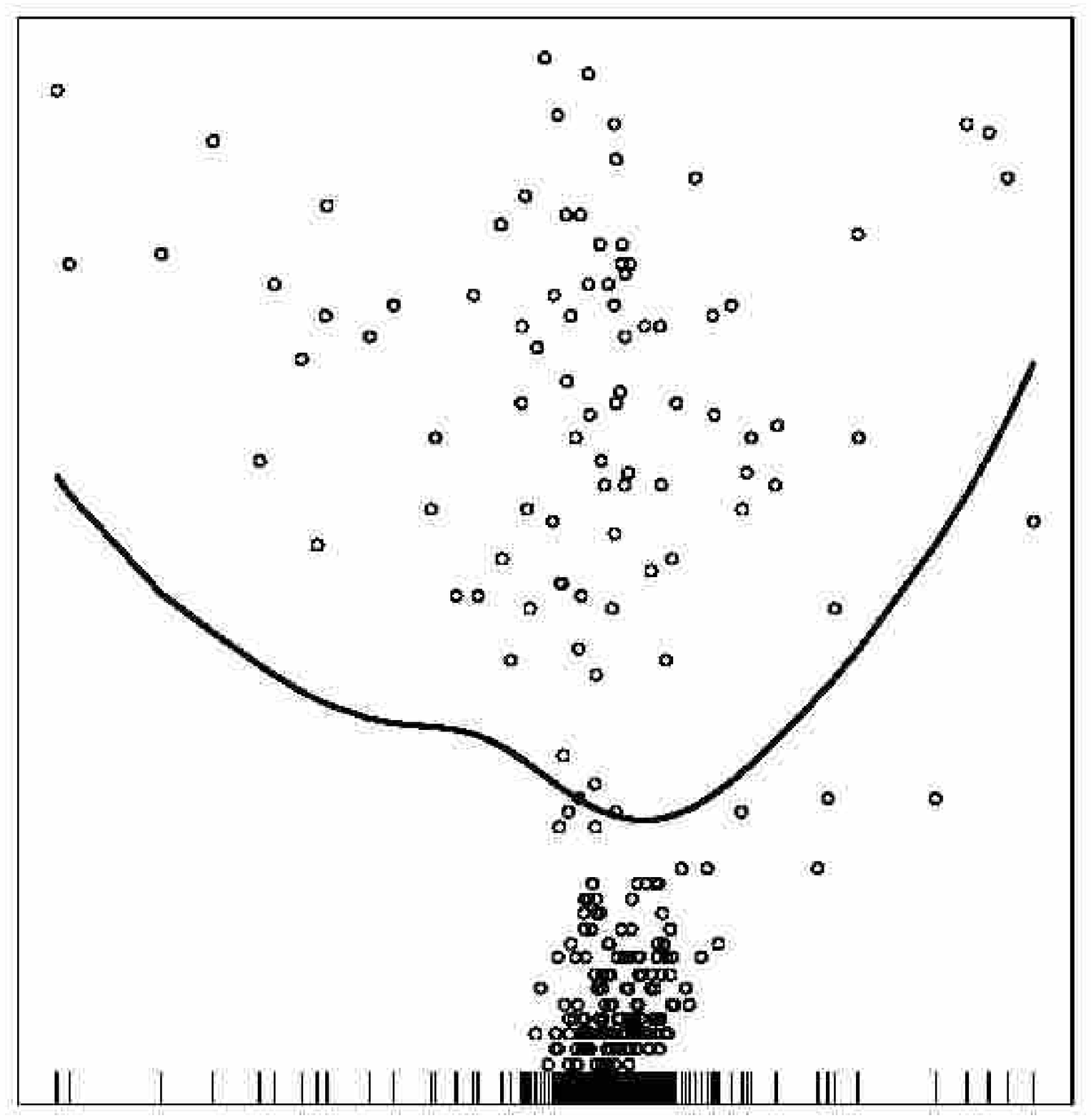}&
\includegraphics[width=\plotwidth]{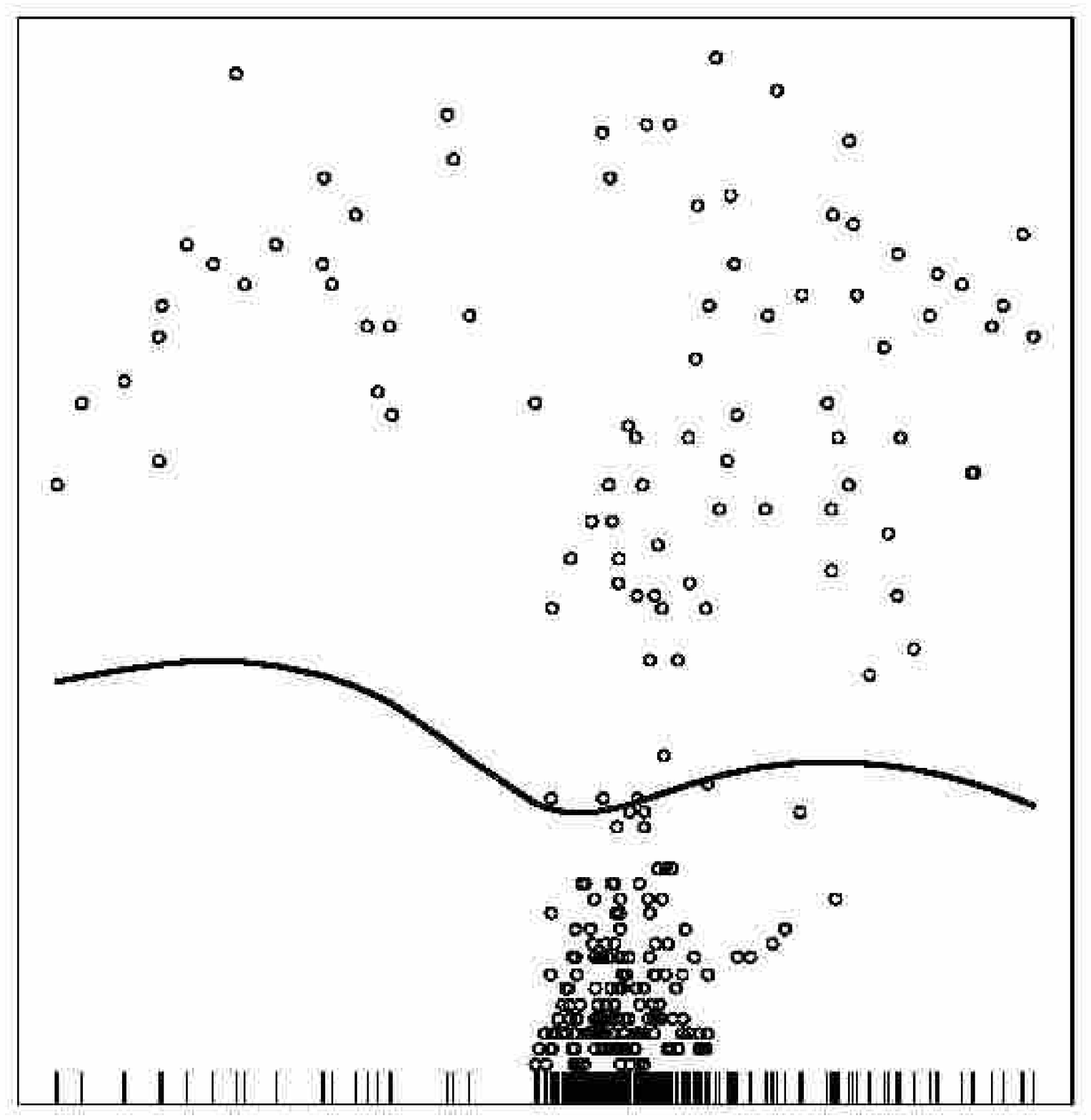}&
\includegraphics[width=\plotwidth]{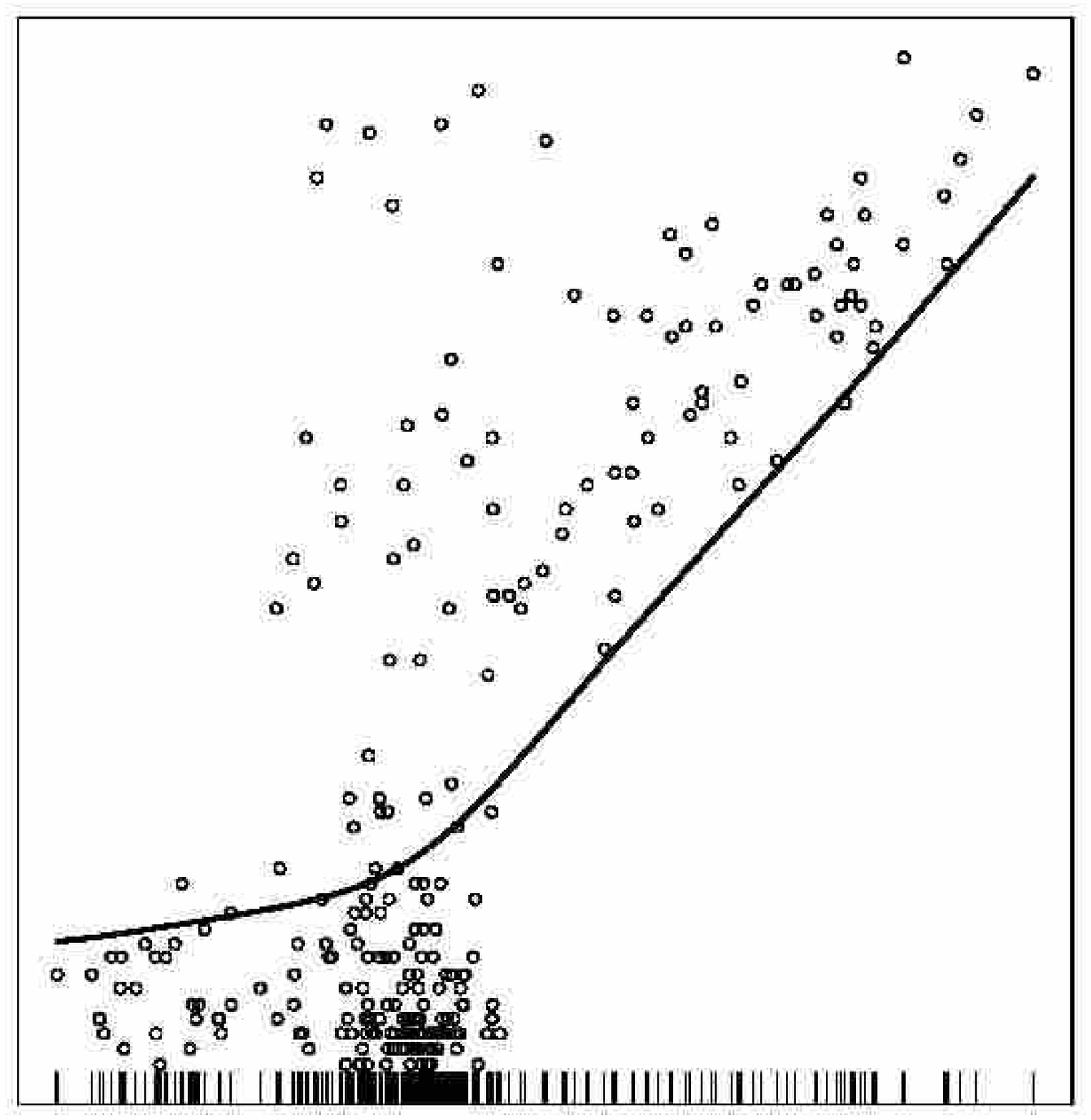}\\
\includegraphics[width=\plotwidth]{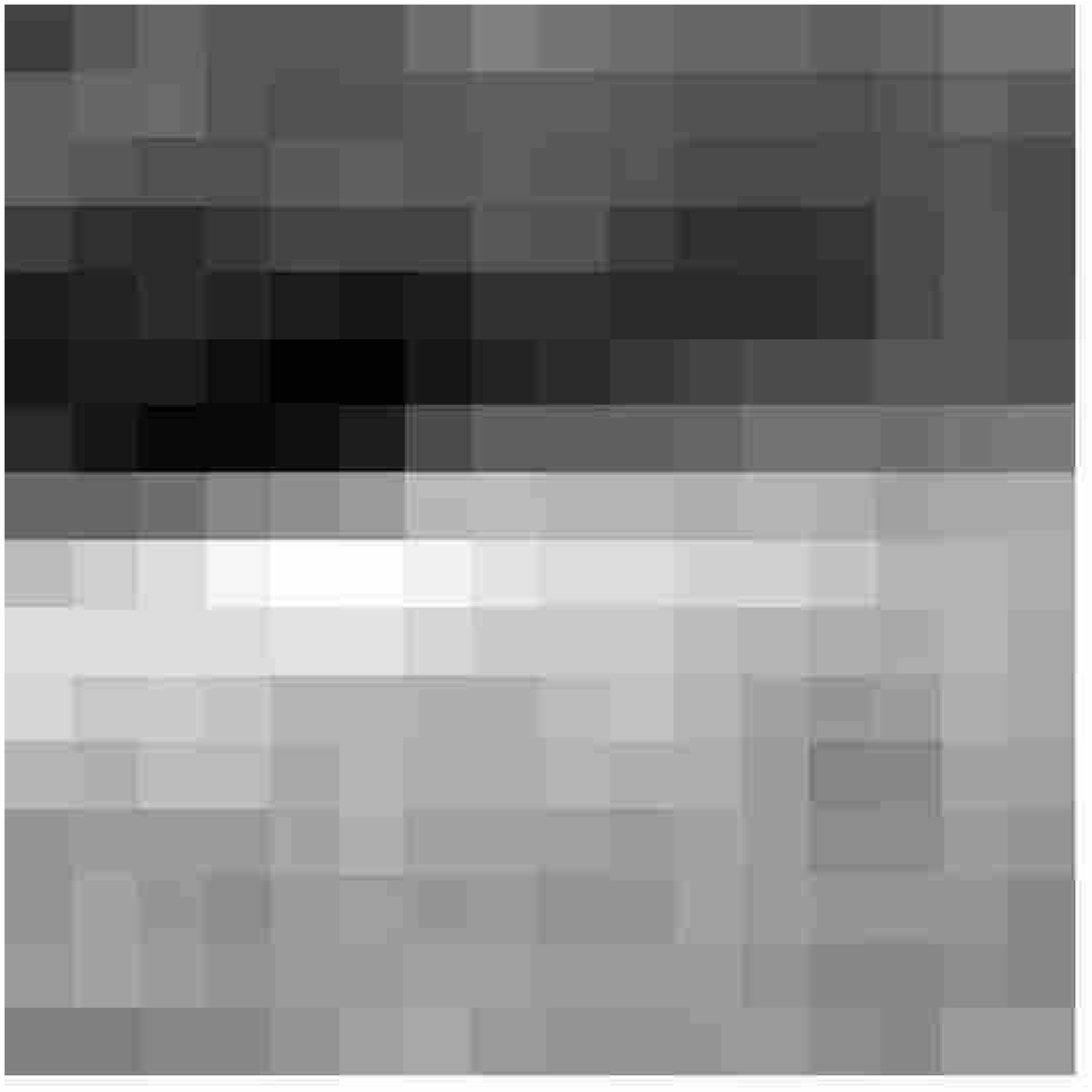}&
\includegraphics[width=\plotwidth]{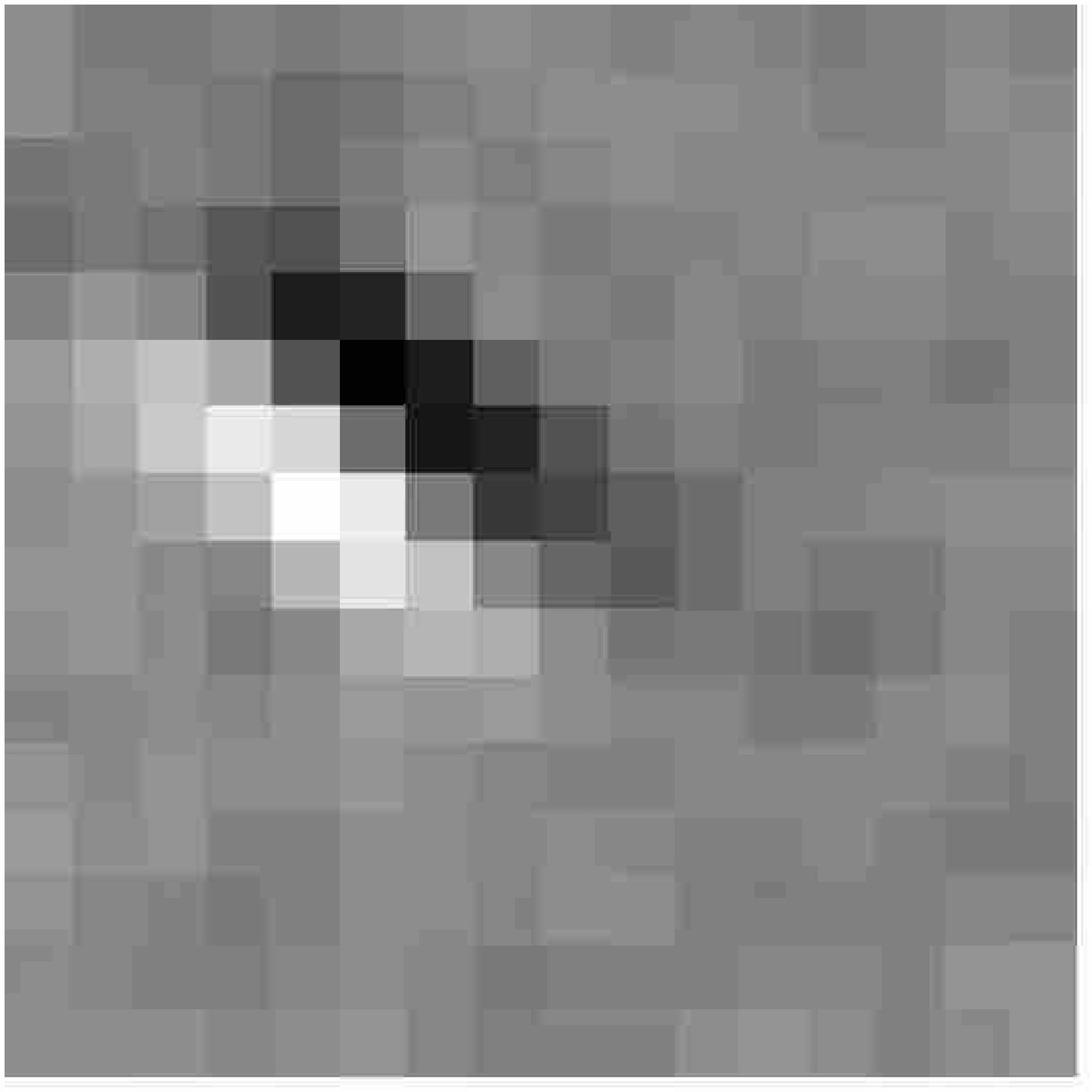}&
\includegraphics[width=\plotwidth]{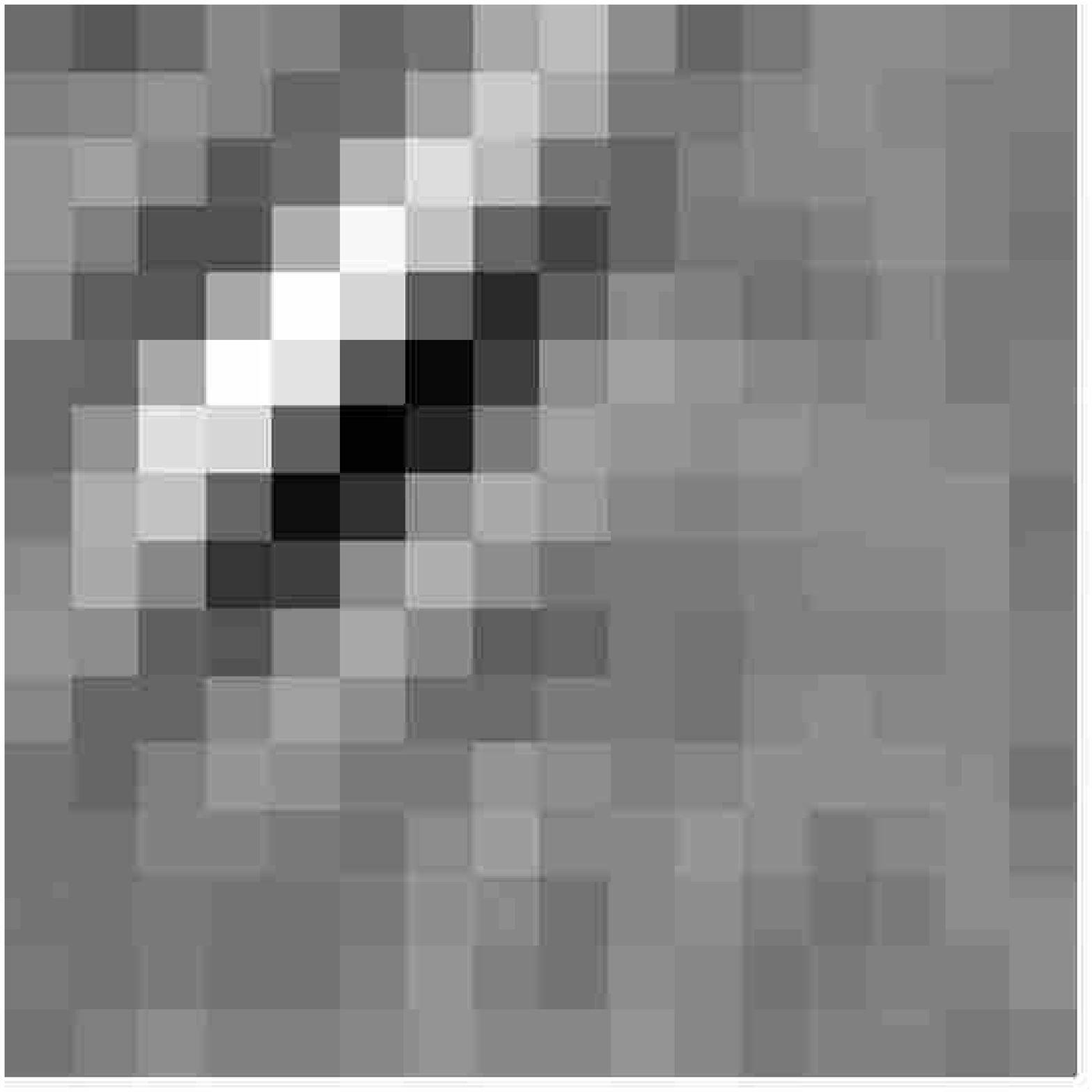}&\\
\includegraphics[width=\plotwidth]{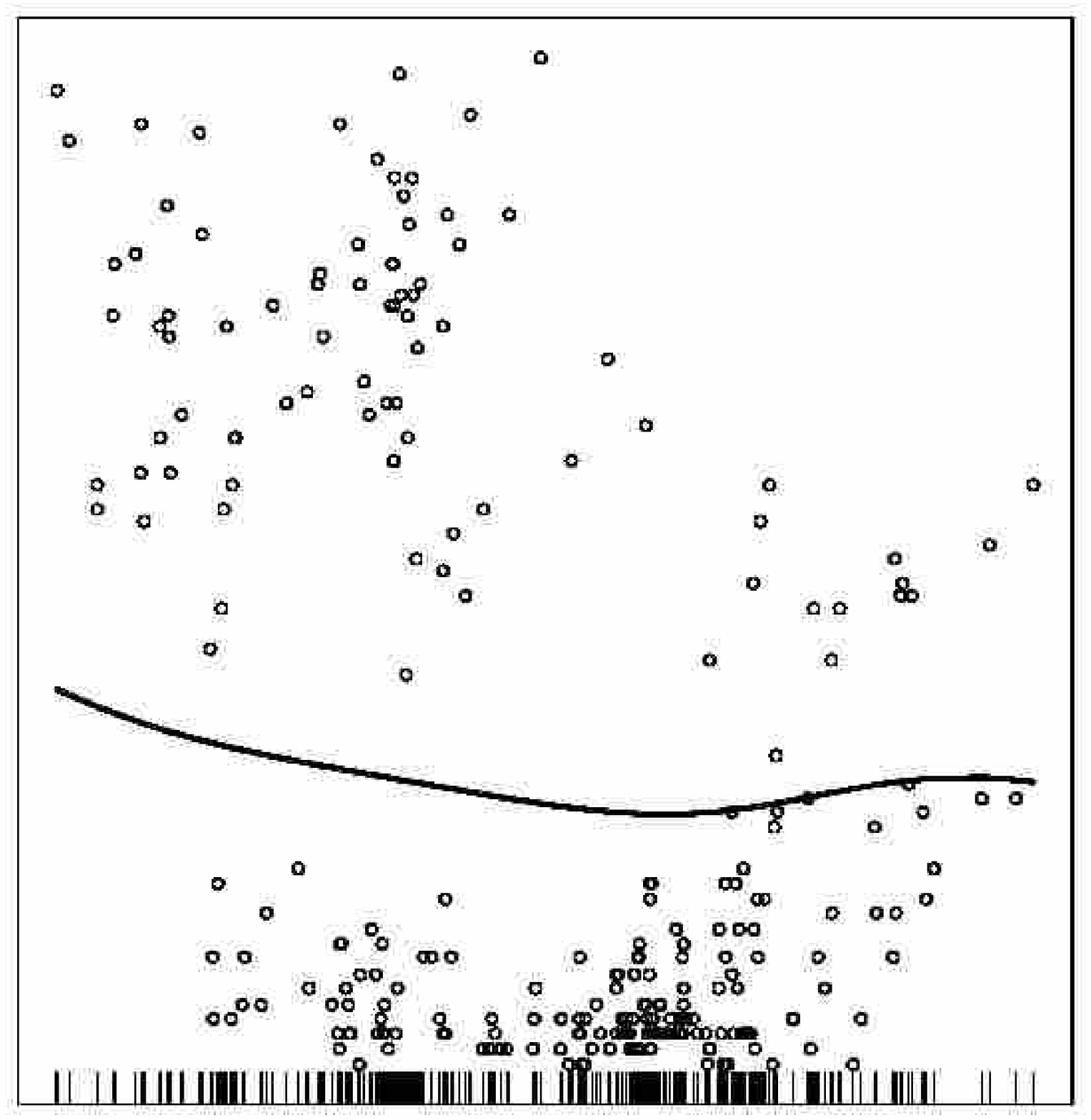}&
\includegraphics[width=\plotwidth]{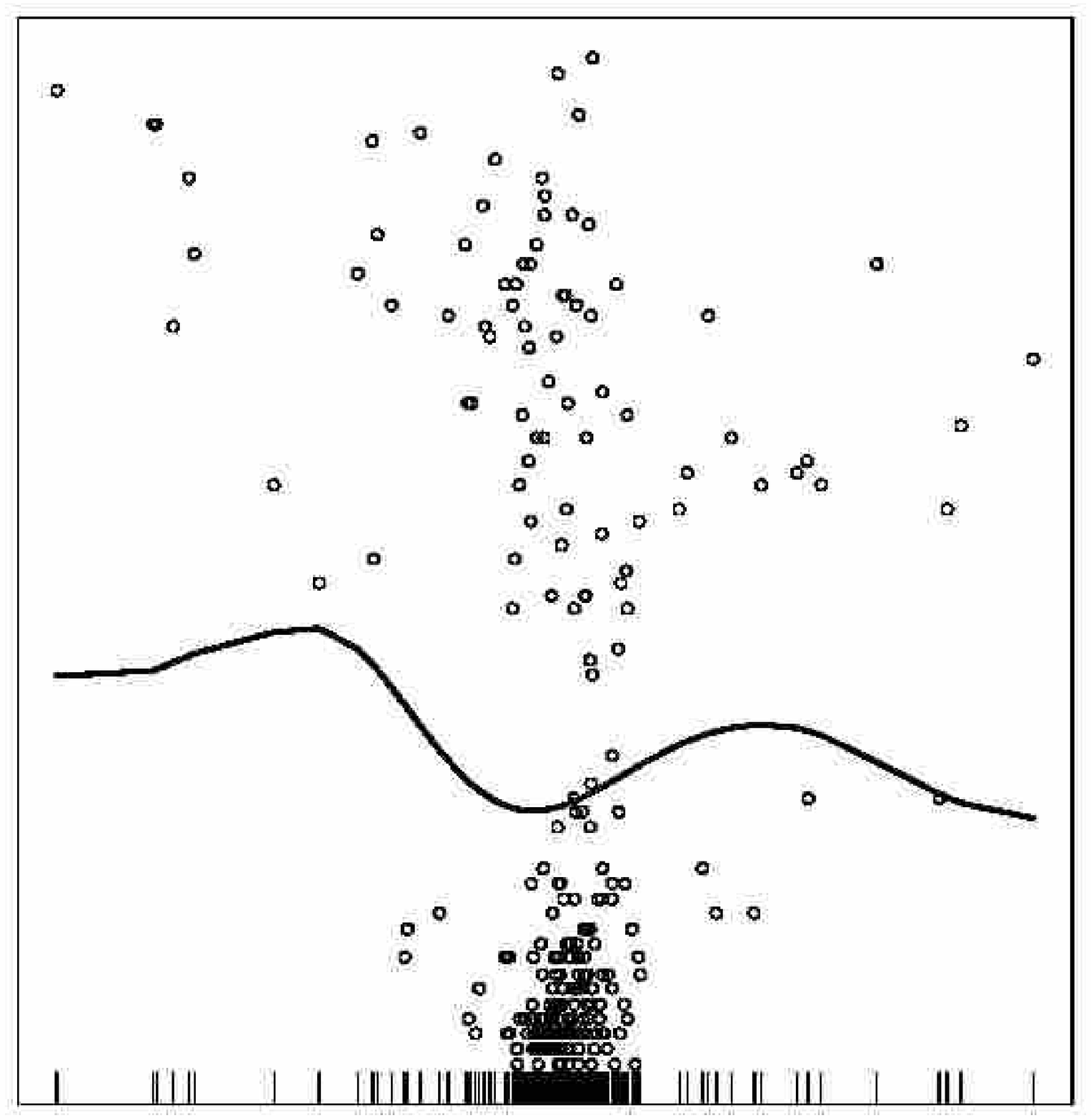}&
\includegraphics[width=\plotwidth]{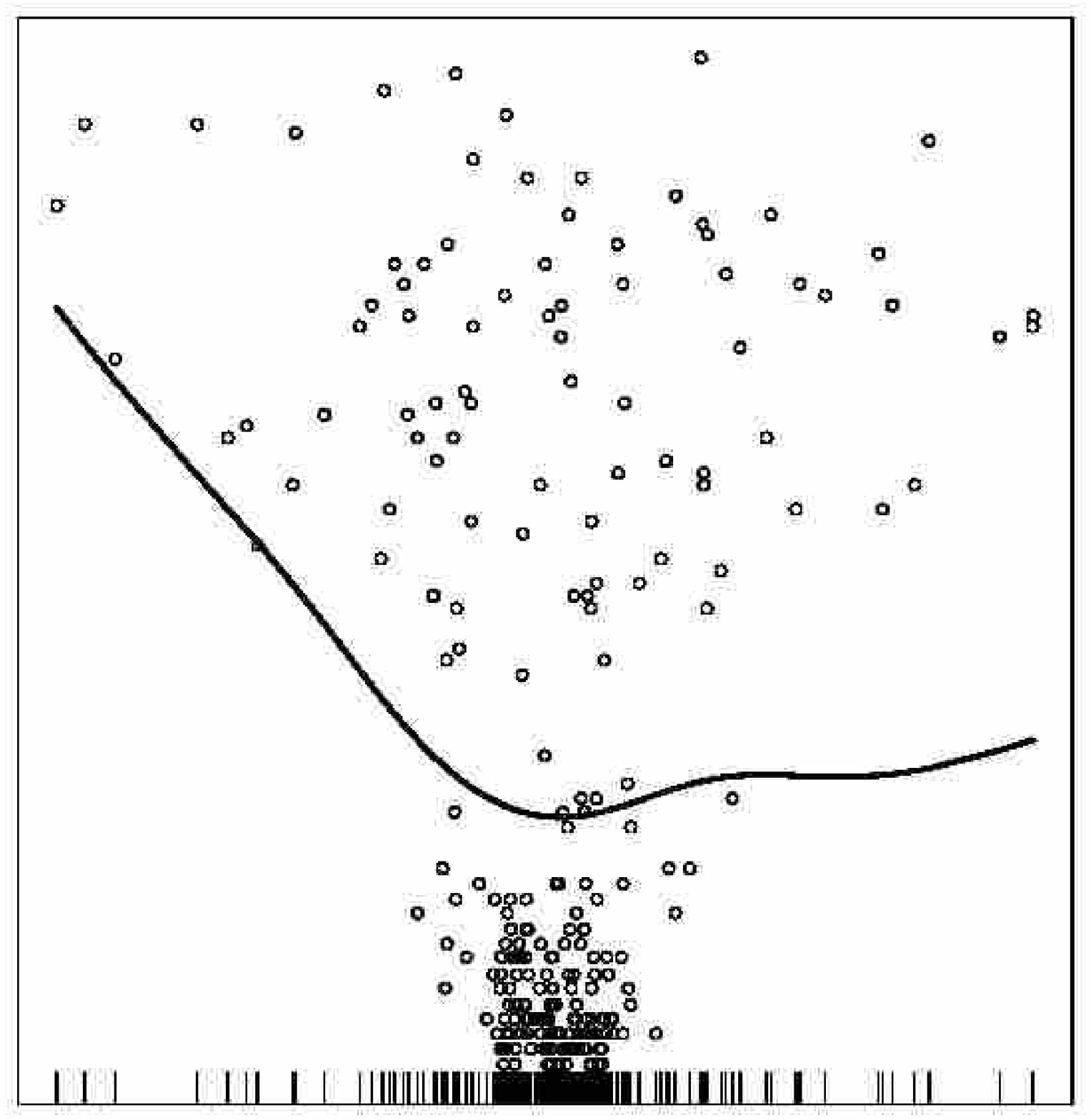}&\\[10pt]
\end{tabular}
\end{tabular}
\end{tabular}
\end{center}
\caption{Comparison of sparse reconstruction using the 
lasso (left) and SpAM (right).}
\label{fig:compare2}
\end{figure*}

\begin{figure*}
\begin{center}
\begin{tabular}{c|c}
\begin{tabular}{c}
\begin{tabular}{cc}
\multicolumn{2}{c}{Lasso} \\
\small Original patch & \small $\text{RSS}=0.2215$ \\
\hskip-5pt
\includegraphics[width=\headwidth]{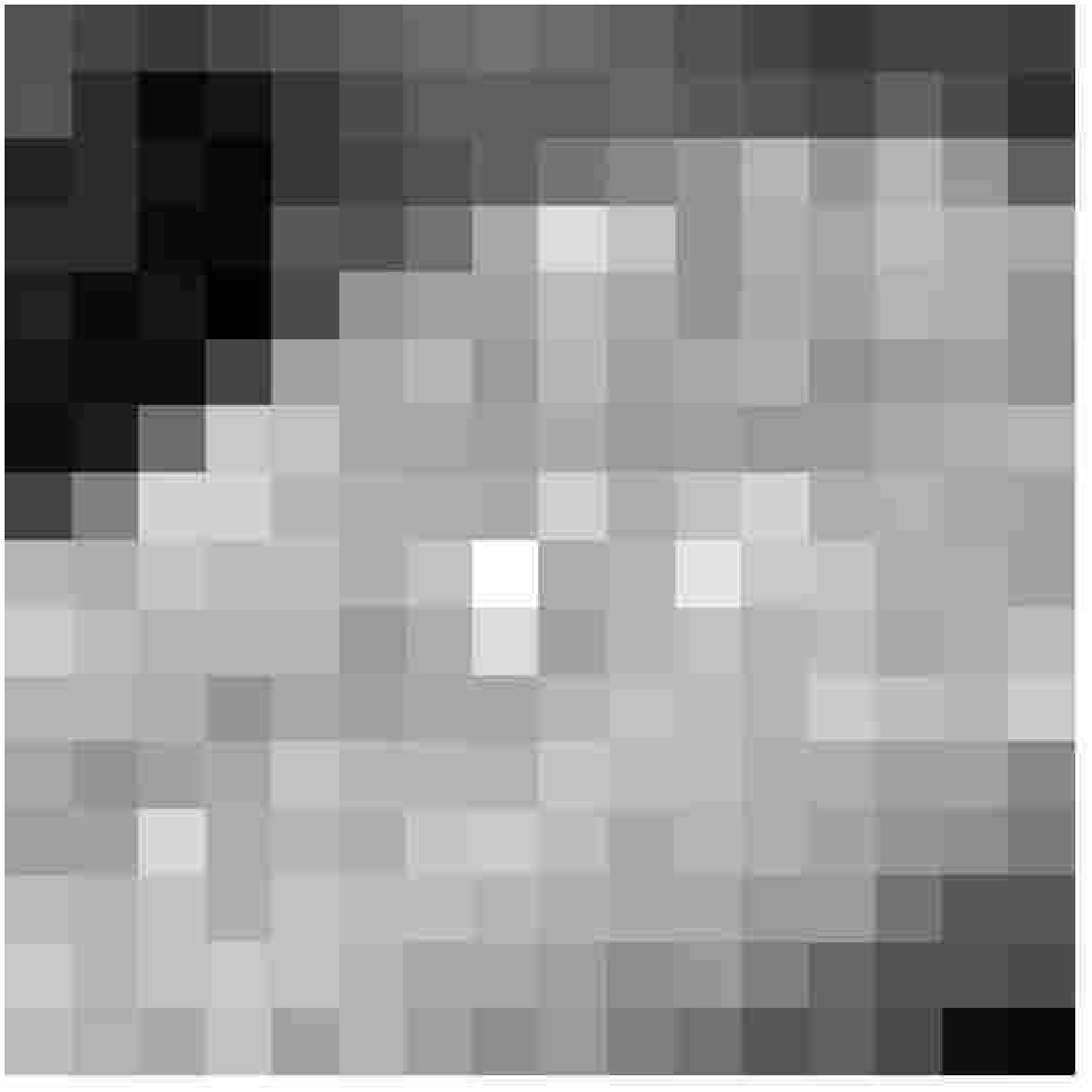} &
\includegraphics[width=\headwidth]{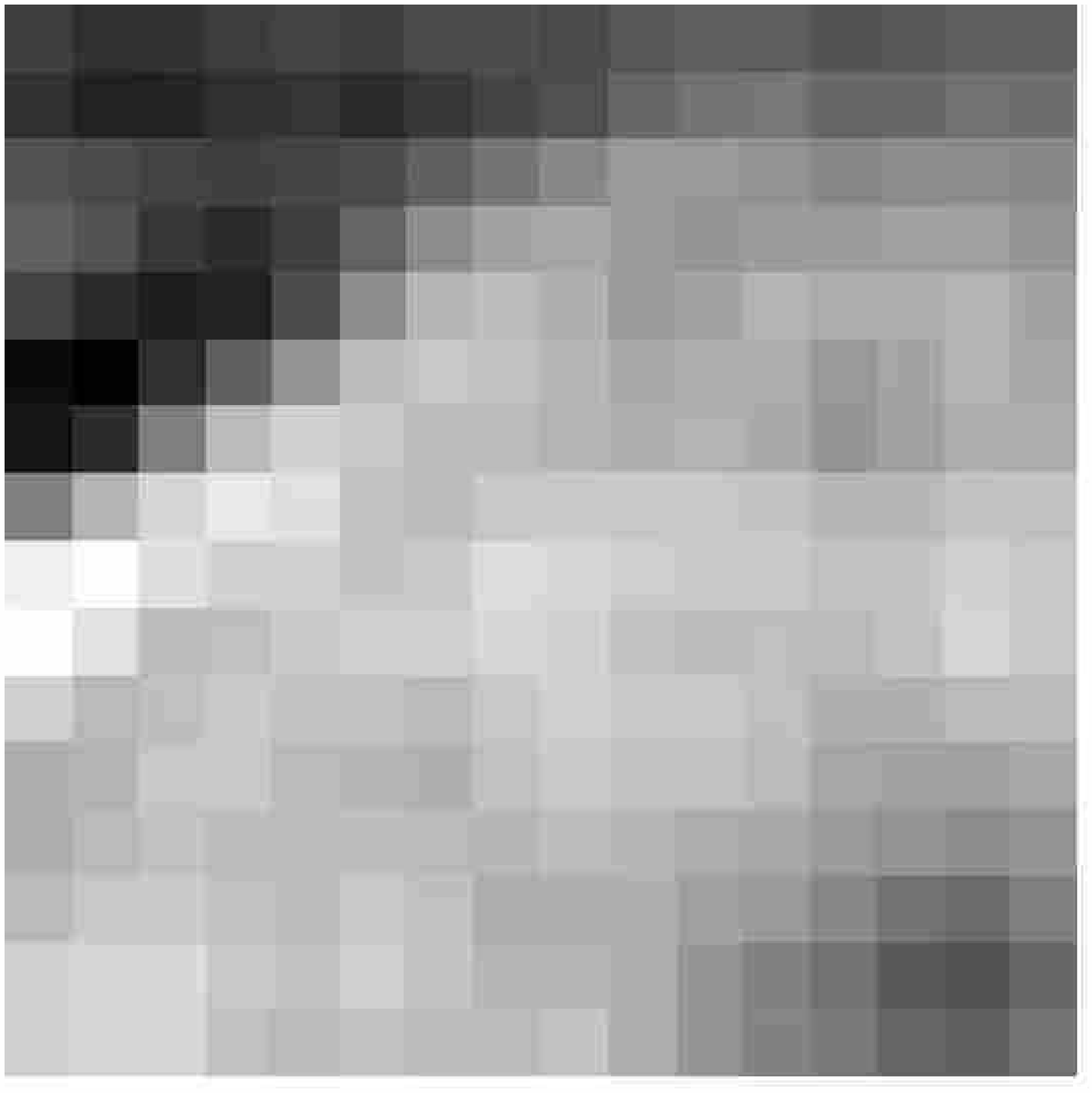}
\\[10pt]
\end{tabular}
\\
\hskip-14pt
\begin{tabular}{cccc}
\includegraphics[width=\plotwidth]{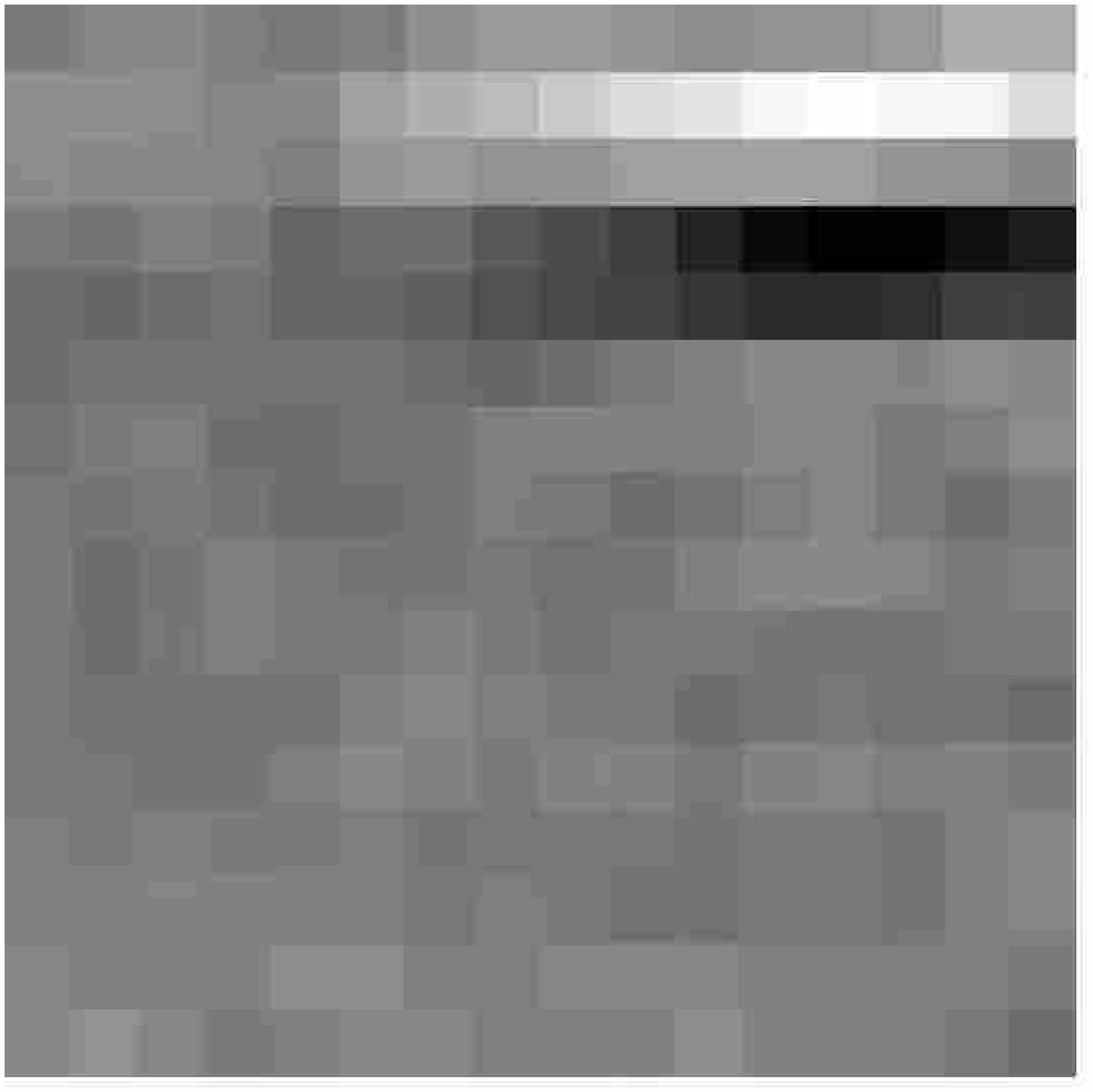}&
\includegraphics[width=\plotwidth]{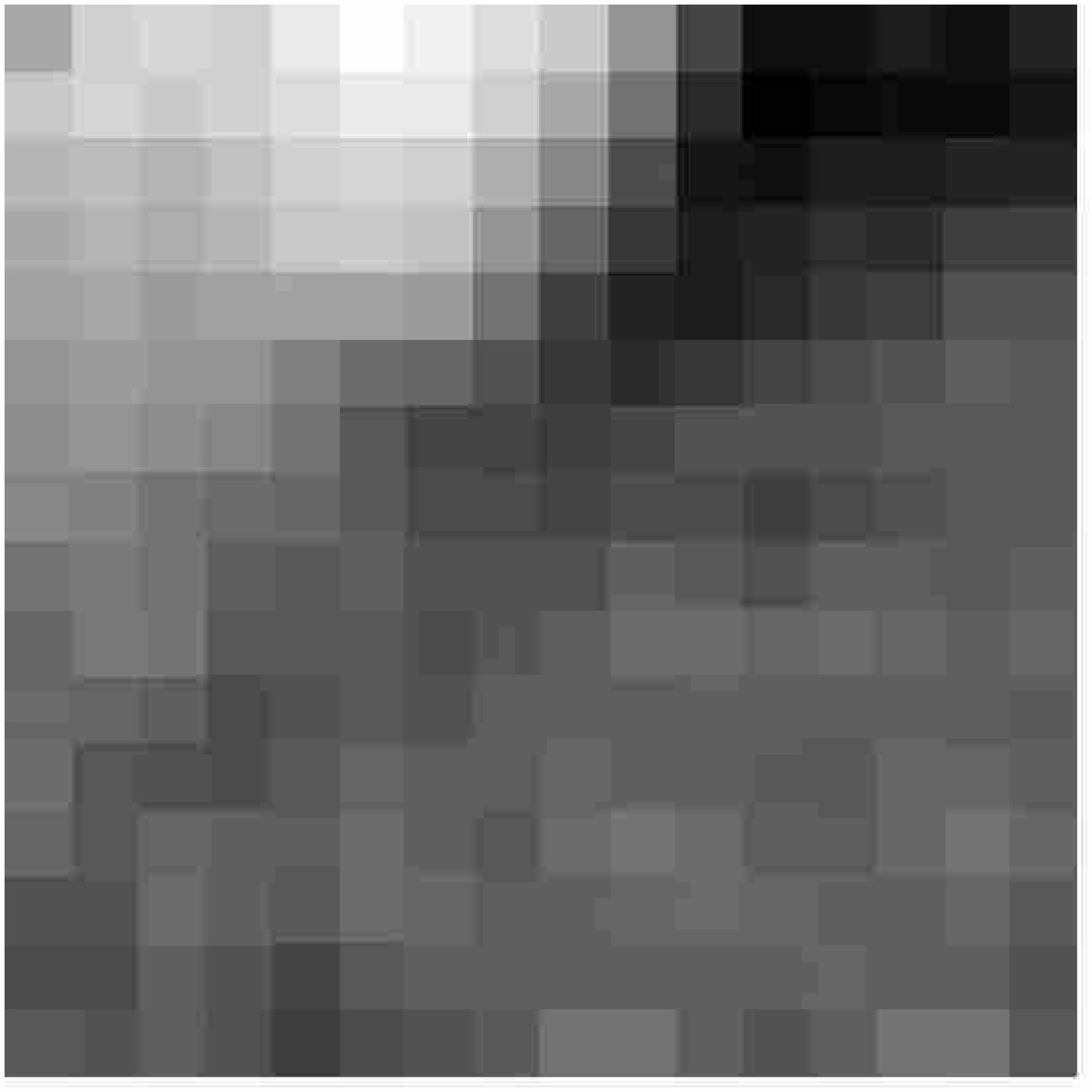}&
\includegraphics[width=\plotwidth]{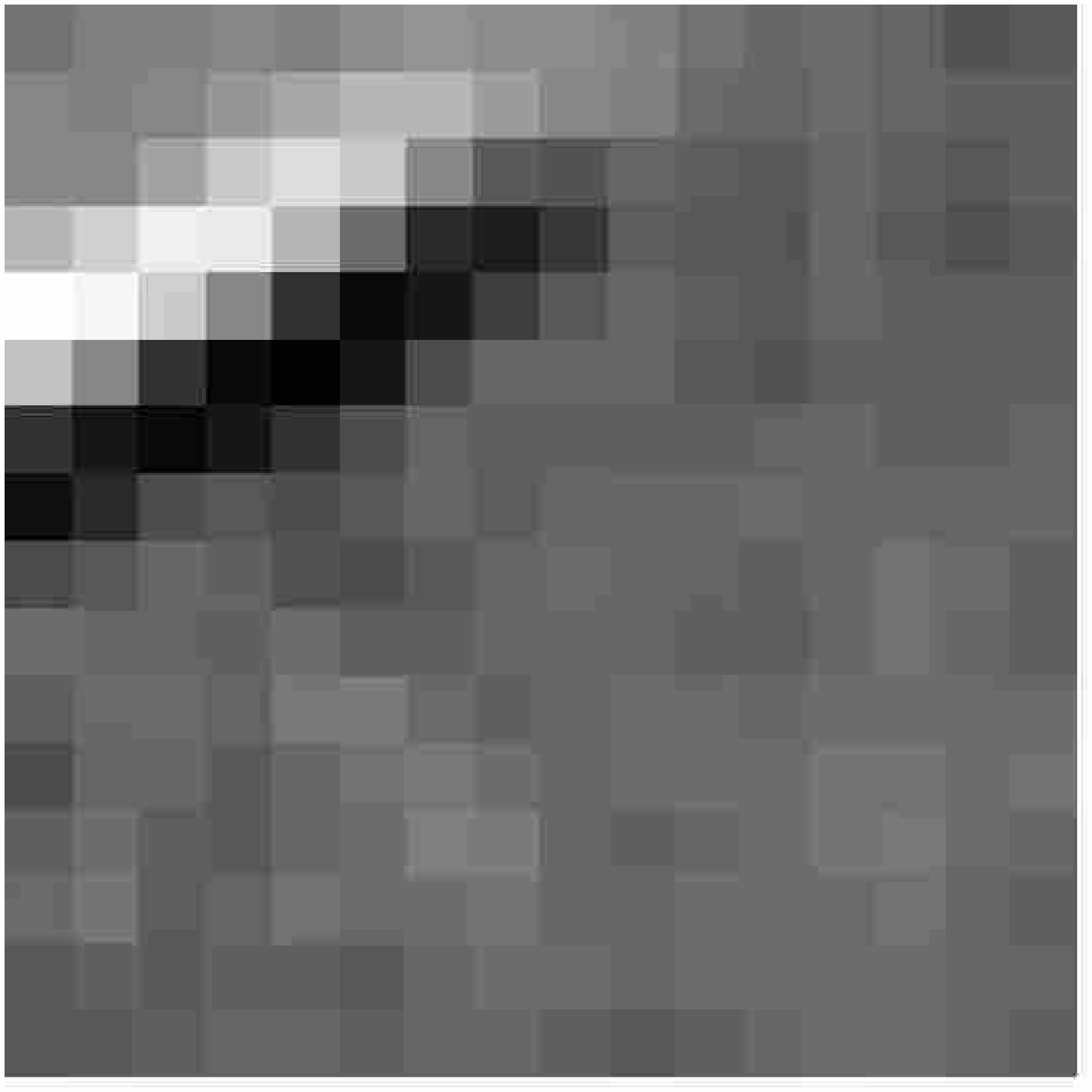}&
\includegraphics[width=\plotwidth]{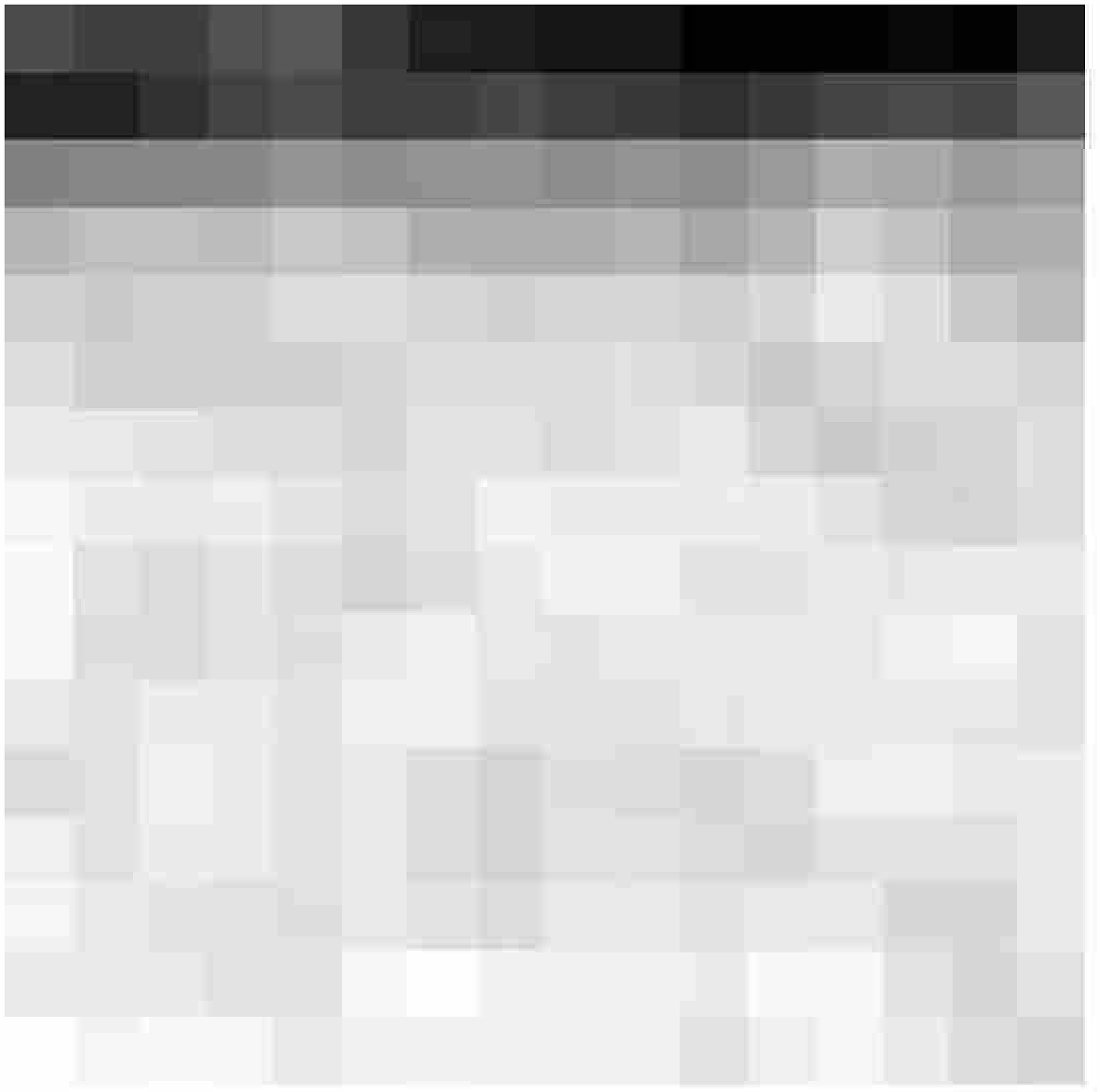} \\
\includegraphics[width=\plotwidth]{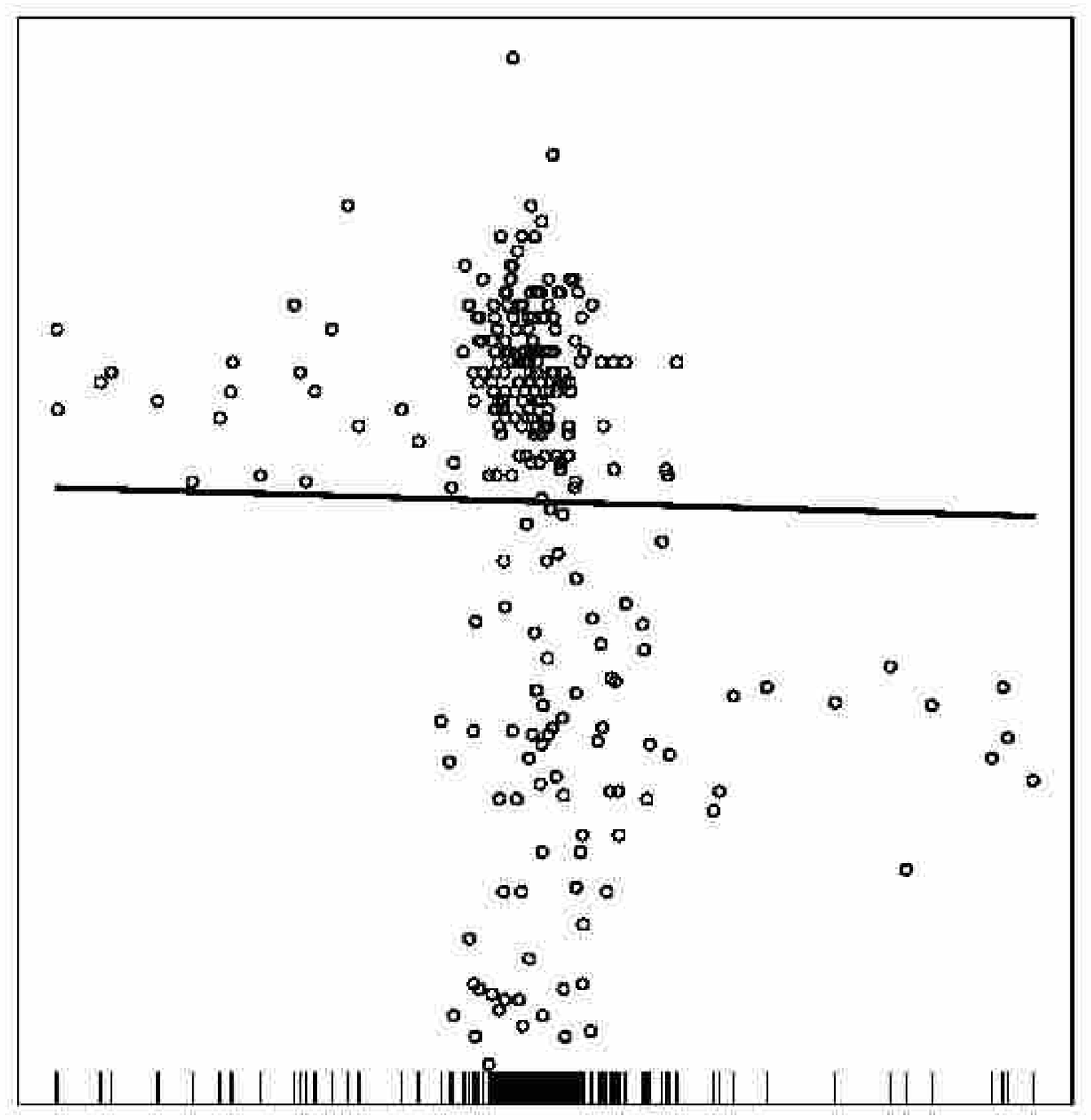}&
\includegraphics[width=\plotwidth]{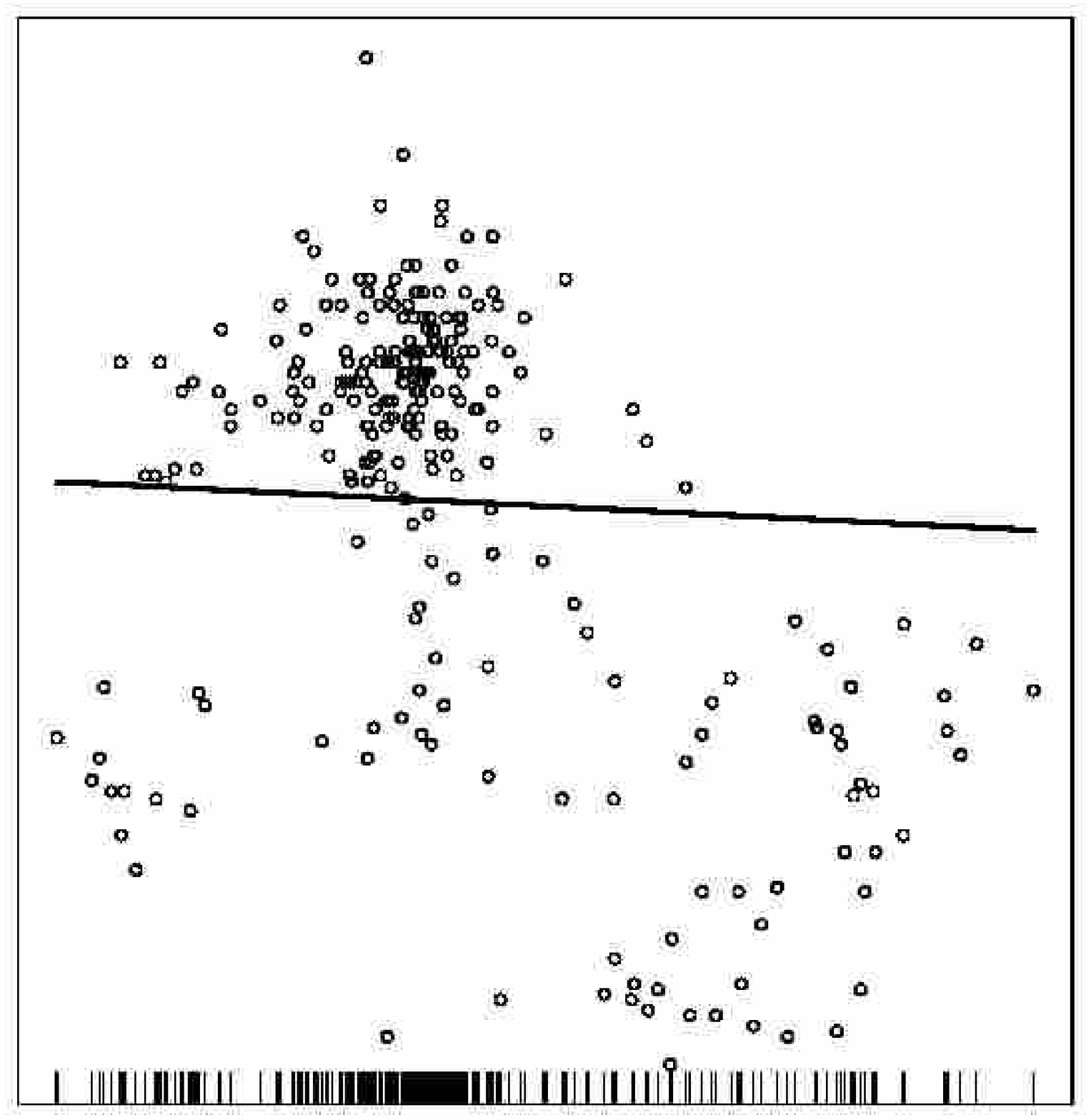}&
\includegraphics[width=\plotwidth]{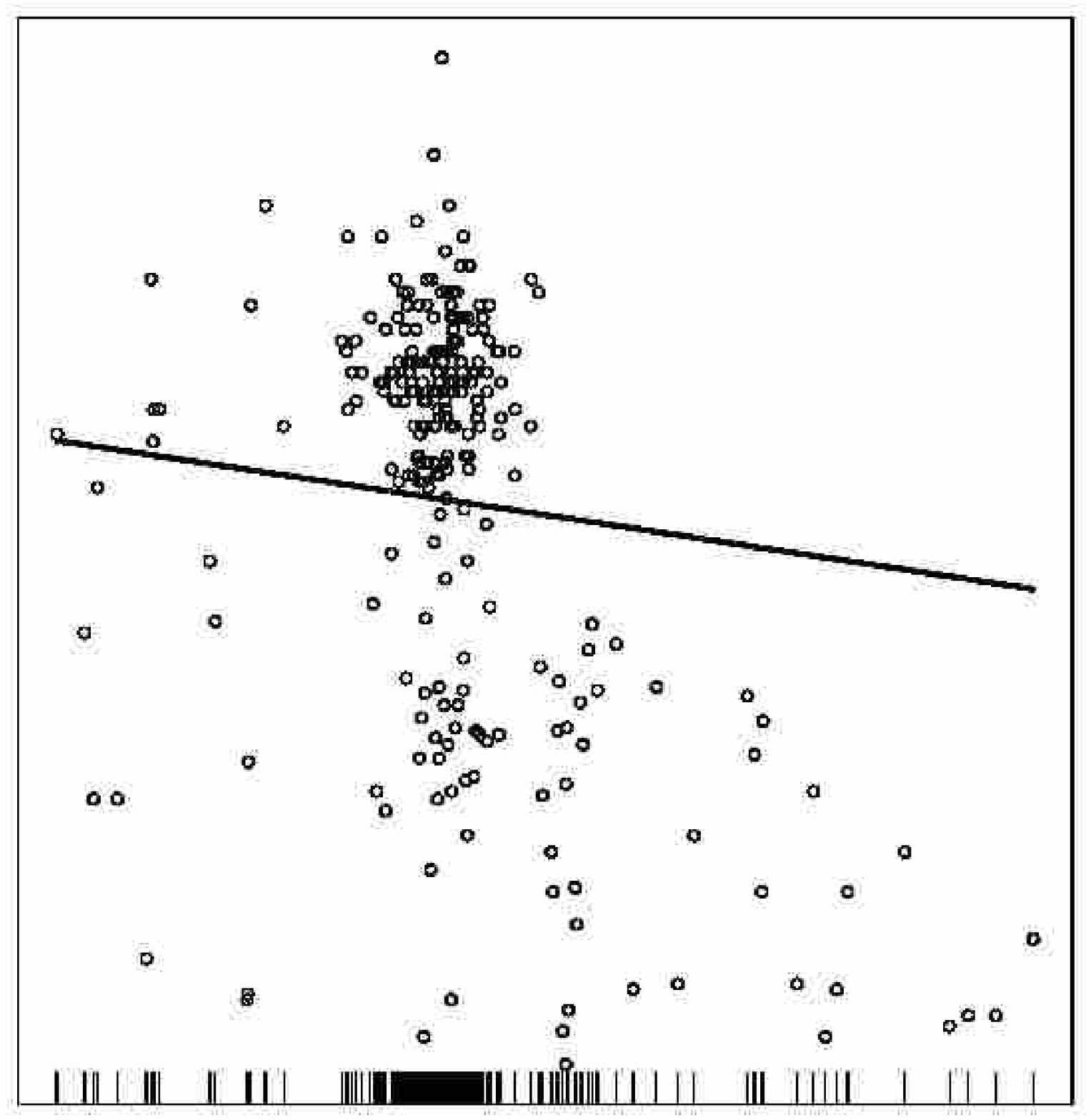}&
\includegraphics[width=\plotwidth]{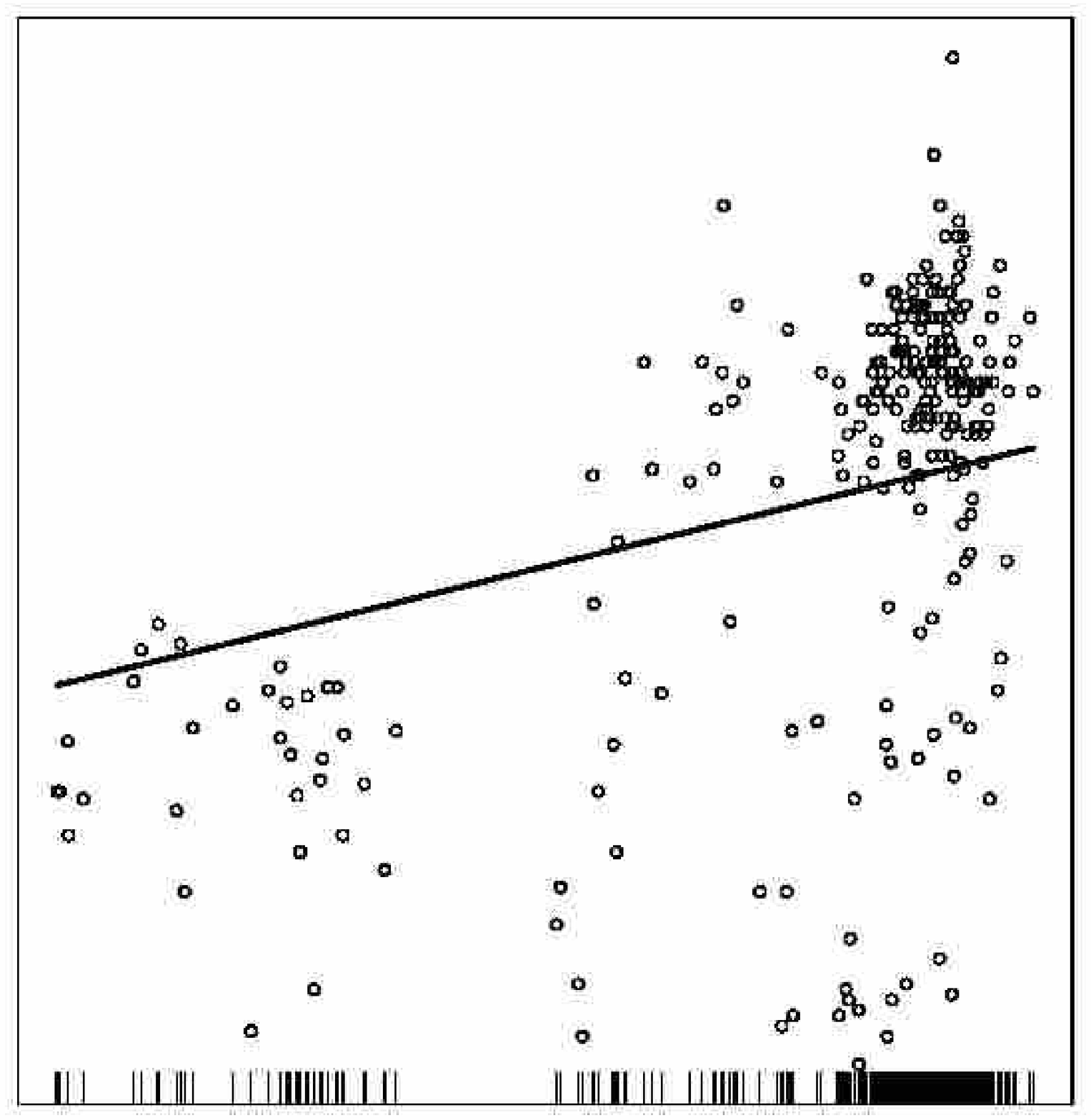} \\
\includegraphics[width=\plotwidth]{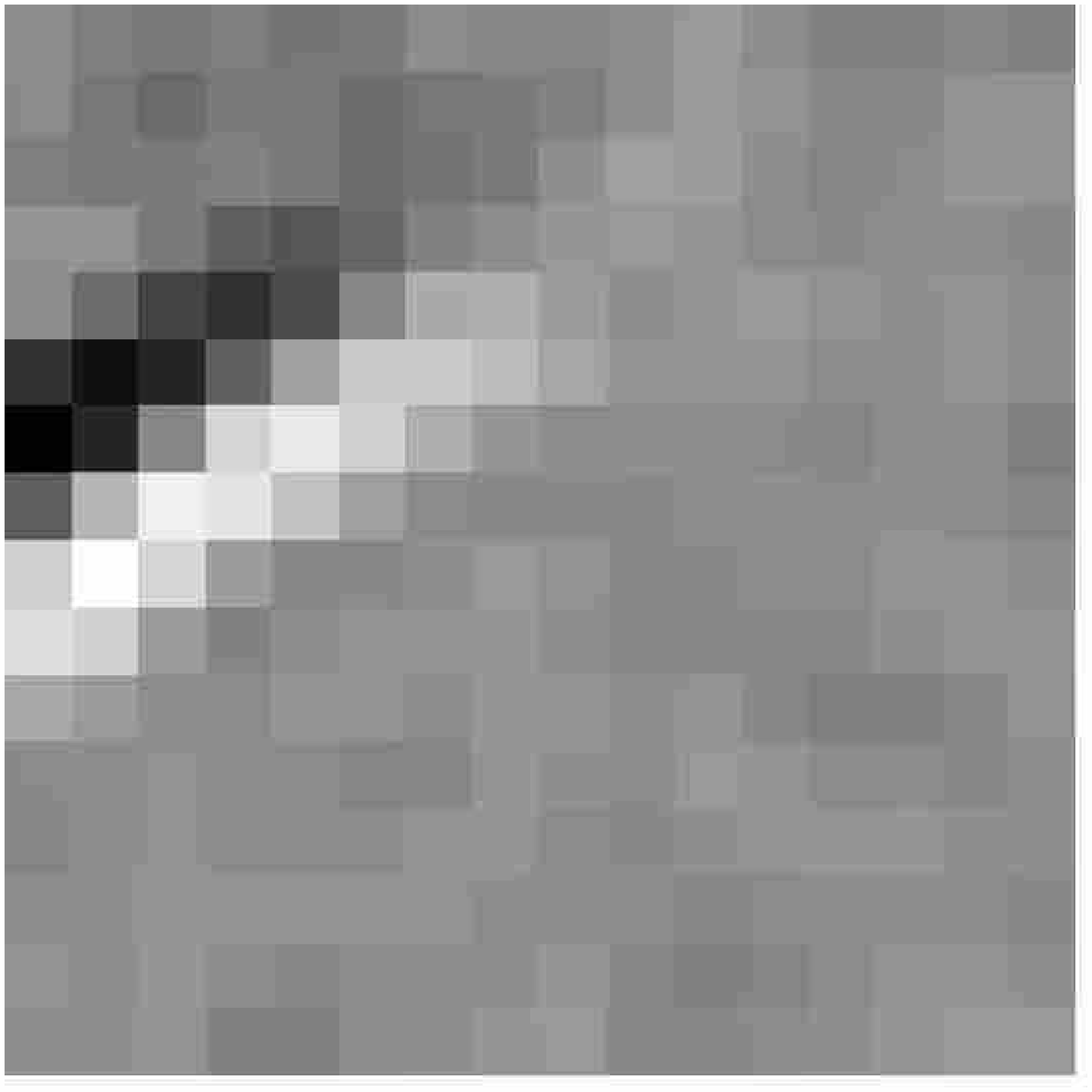}&
\includegraphics[width=\plotwidth]{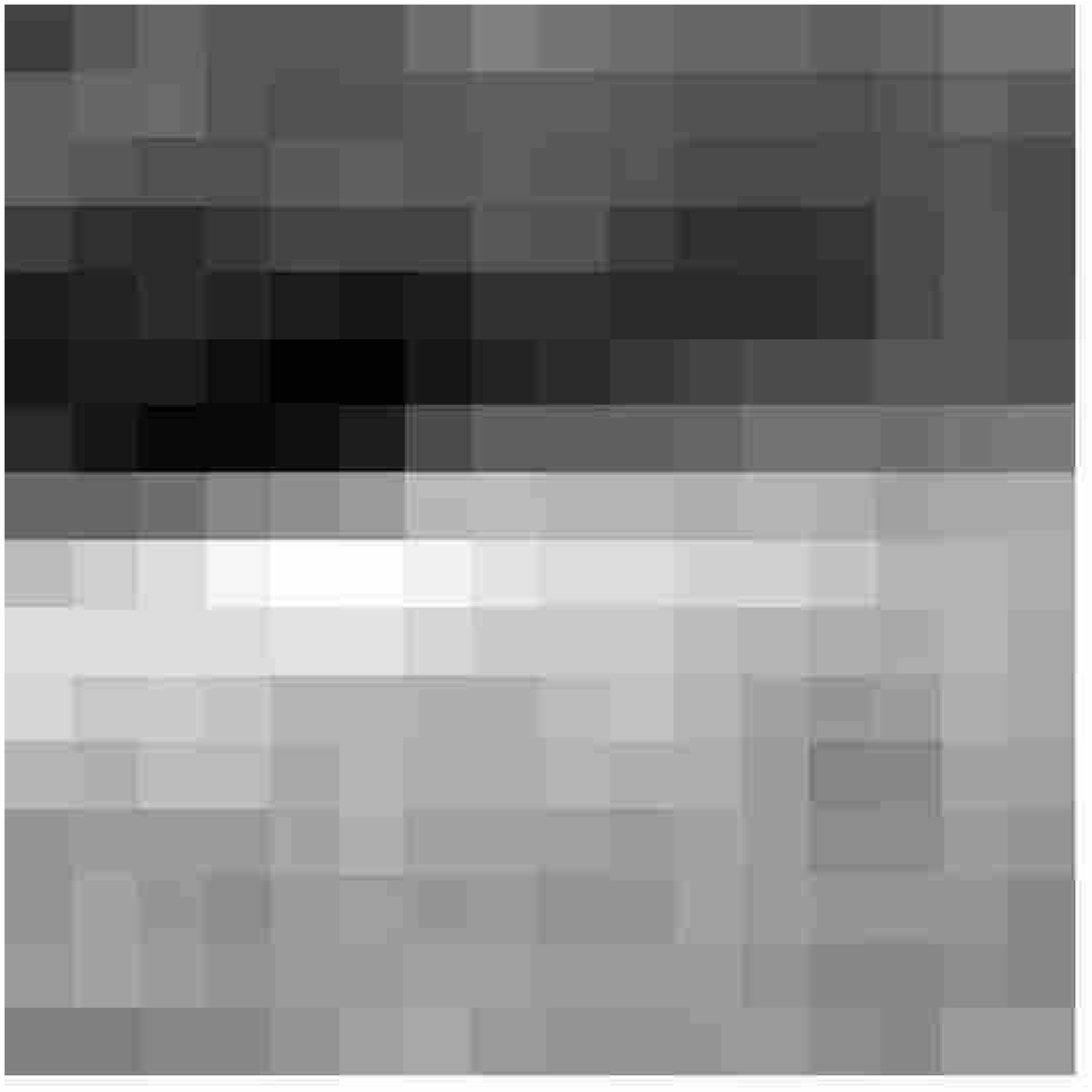}&
\includegraphics[width=\plotwidth]{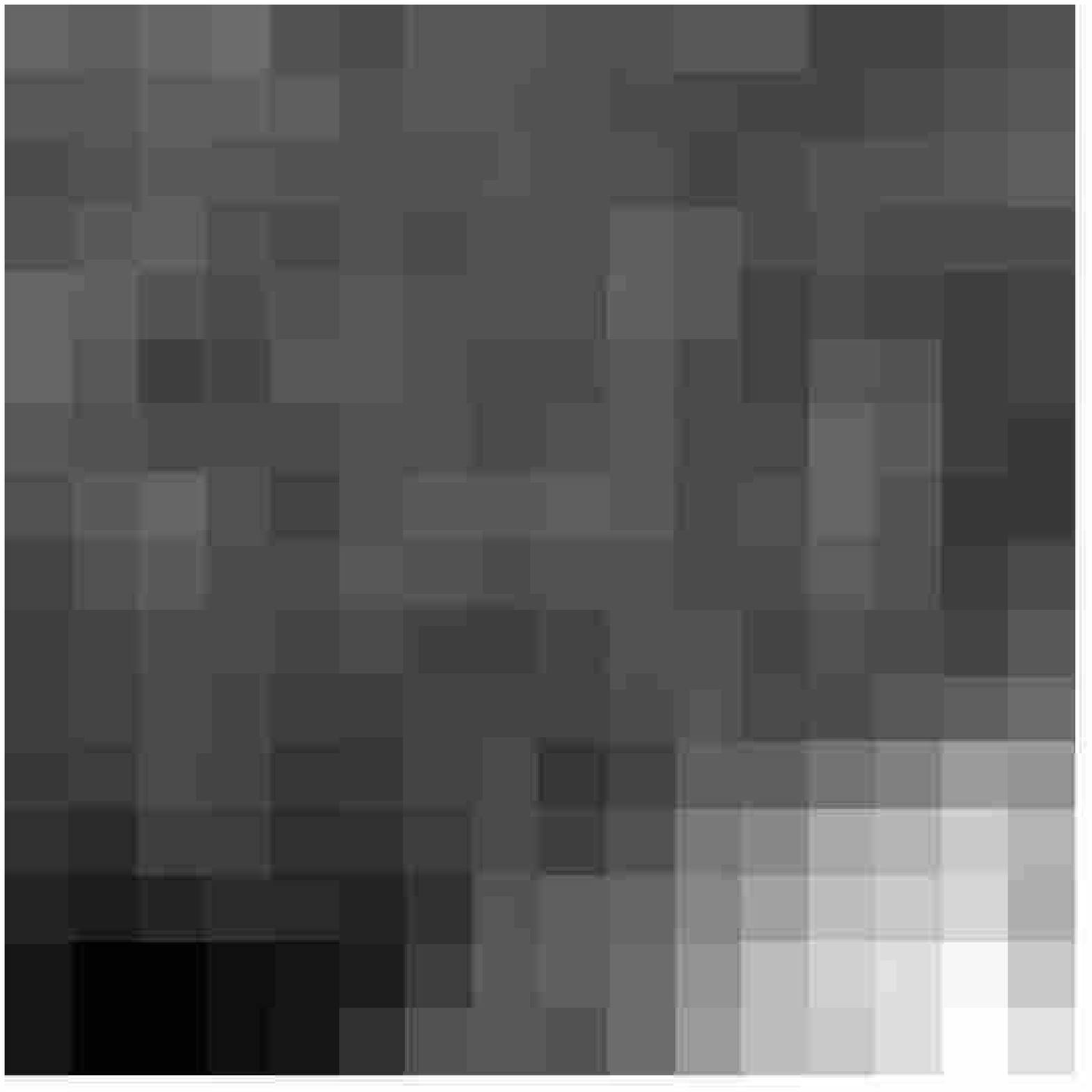}&
\includegraphics[width=\plotwidth]{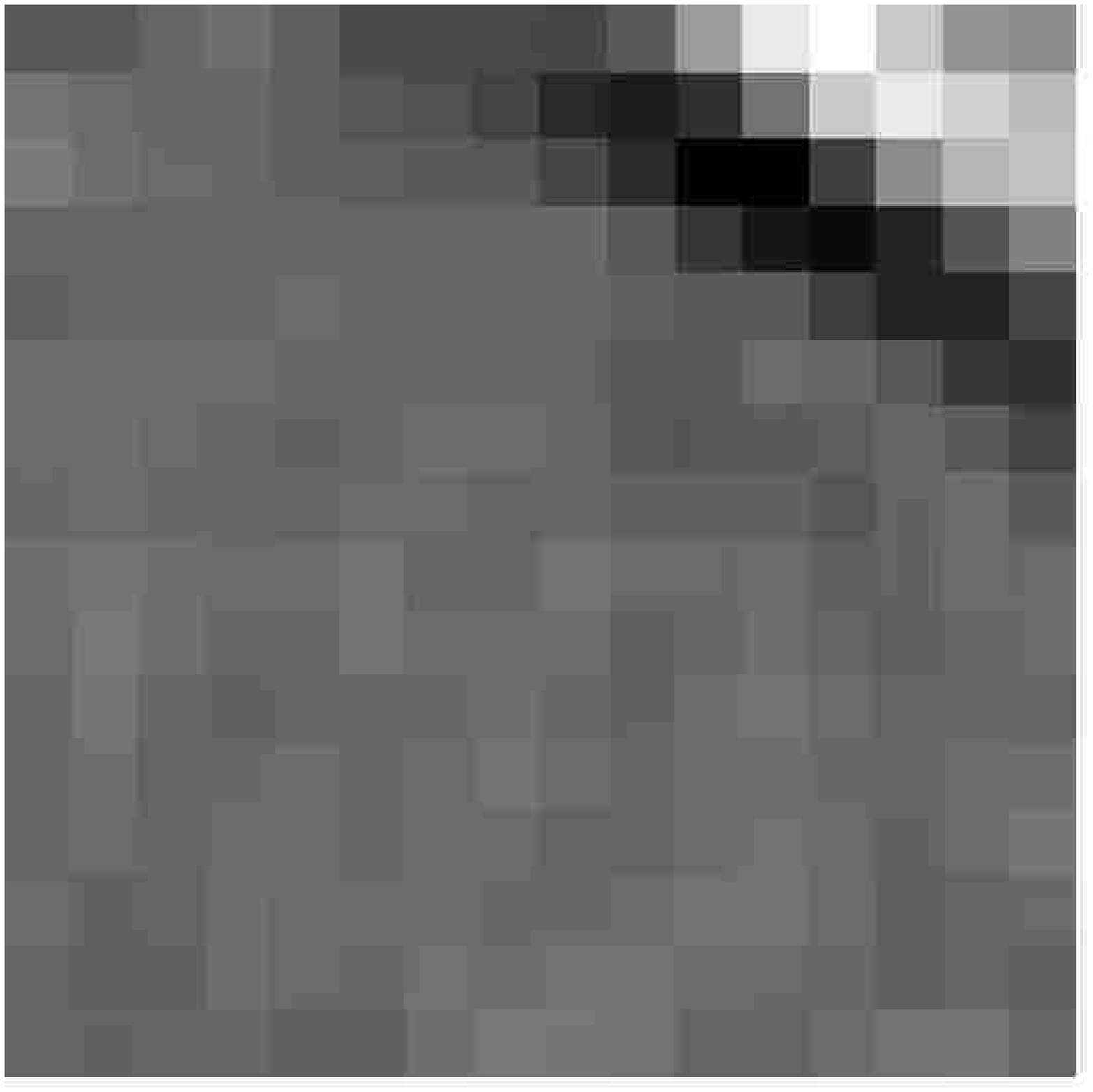}\\
\includegraphics[width=\plotwidth]{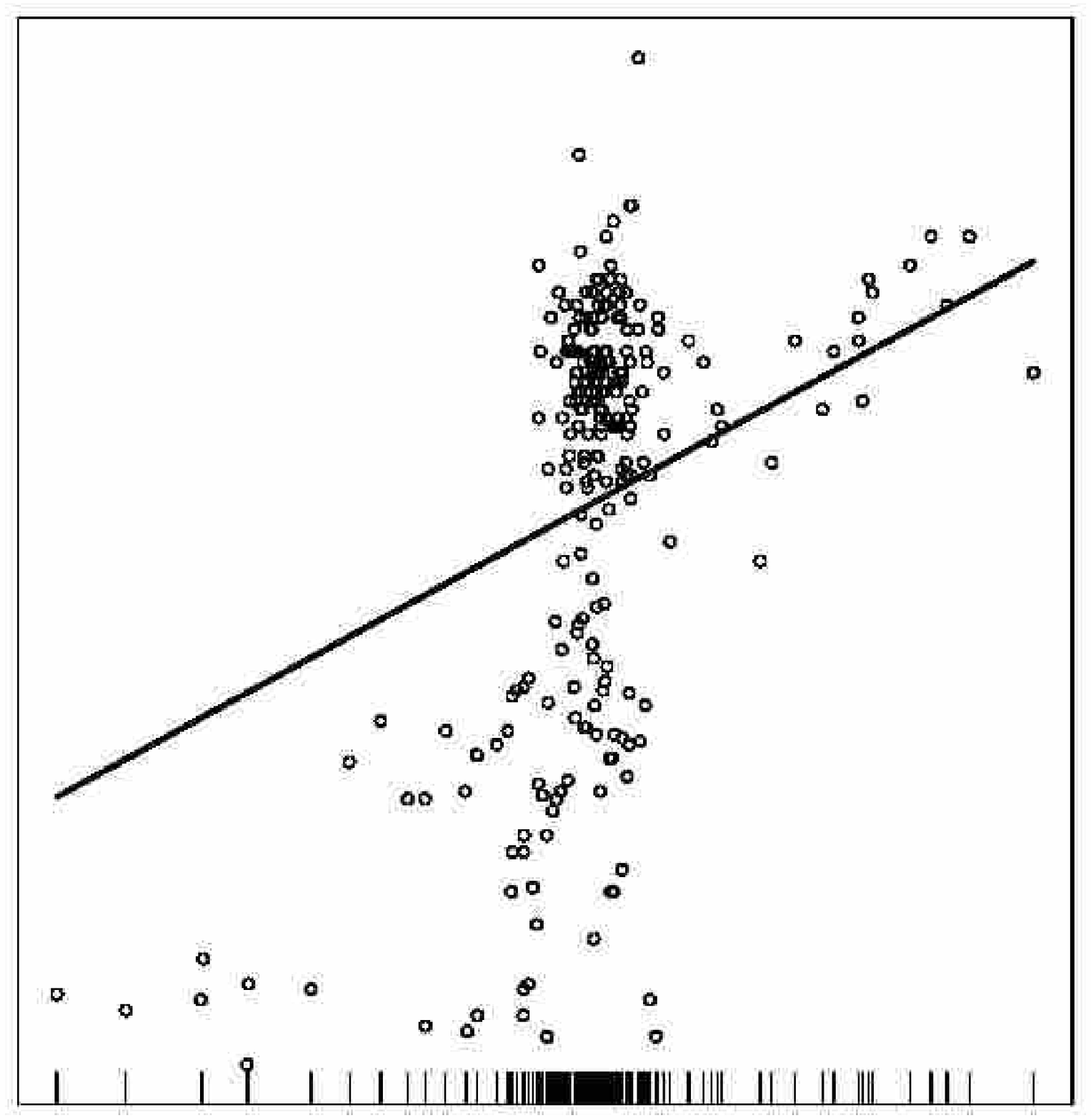}&
\includegraphics[width=\plotwidth]{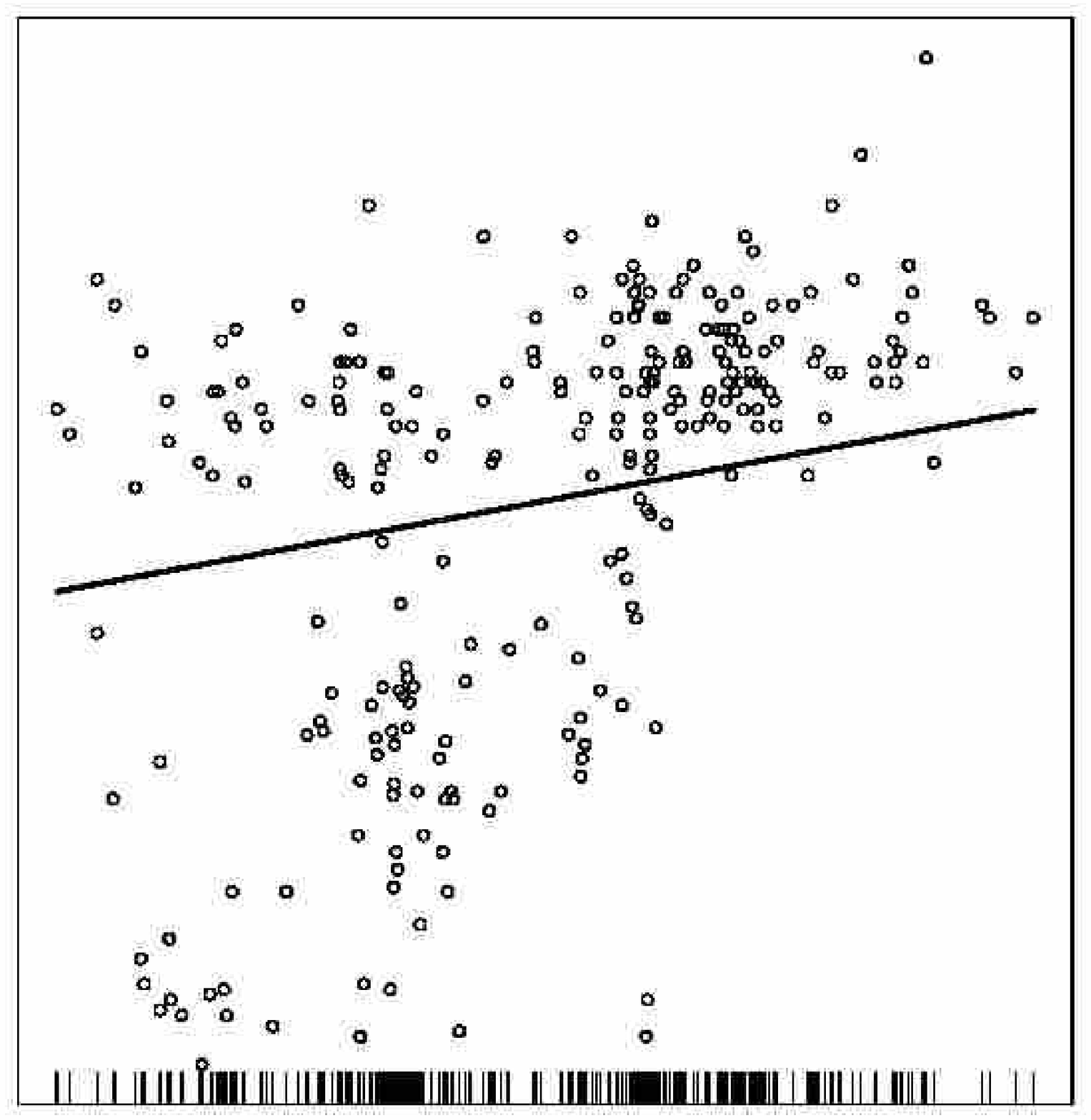}&
\includegraphics[width=\plotwidth]{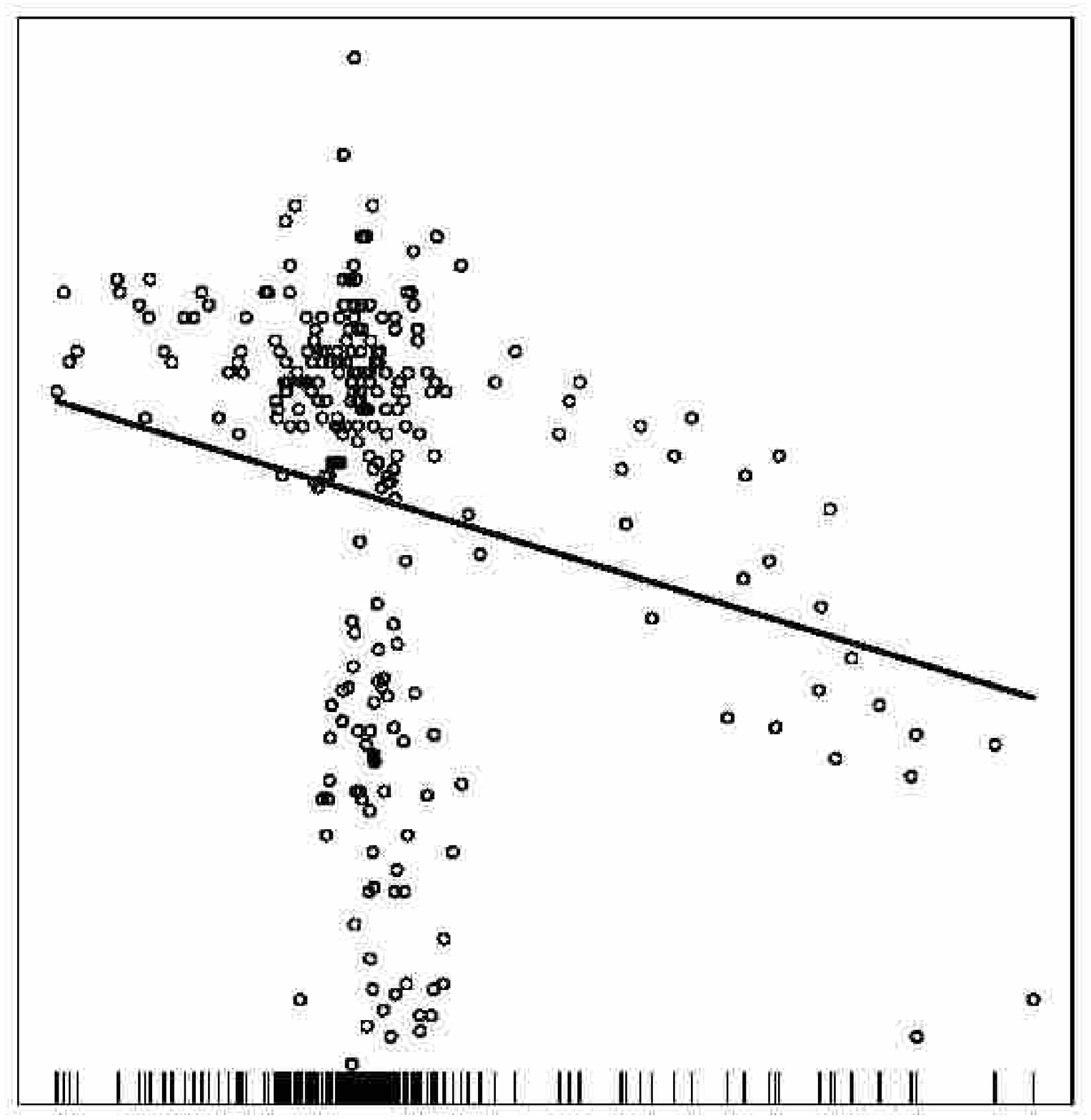} &
\includegraphics[width=\plotwidth]{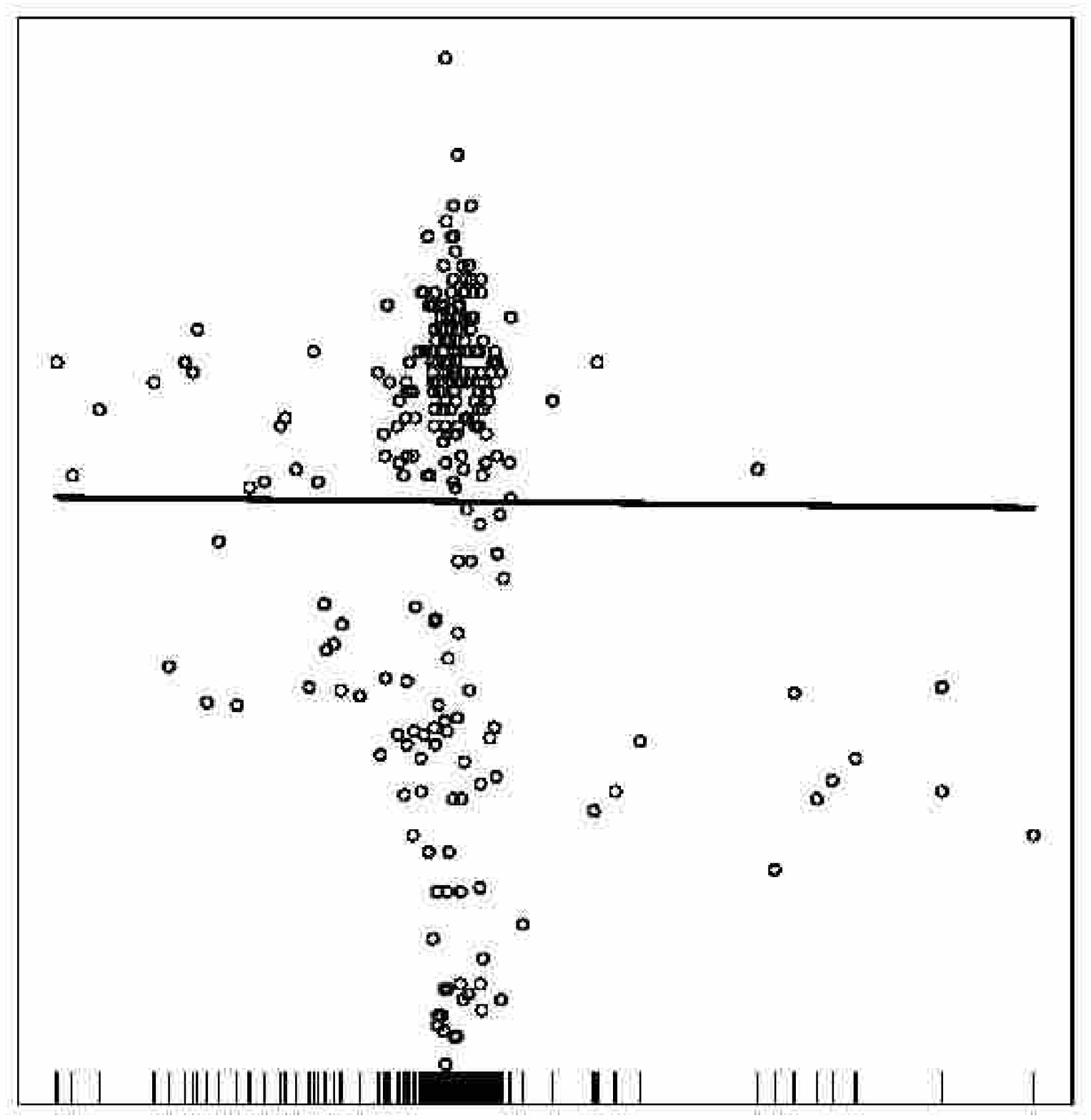}\\[10pt]
\includegraphics[width=\plotwidth]{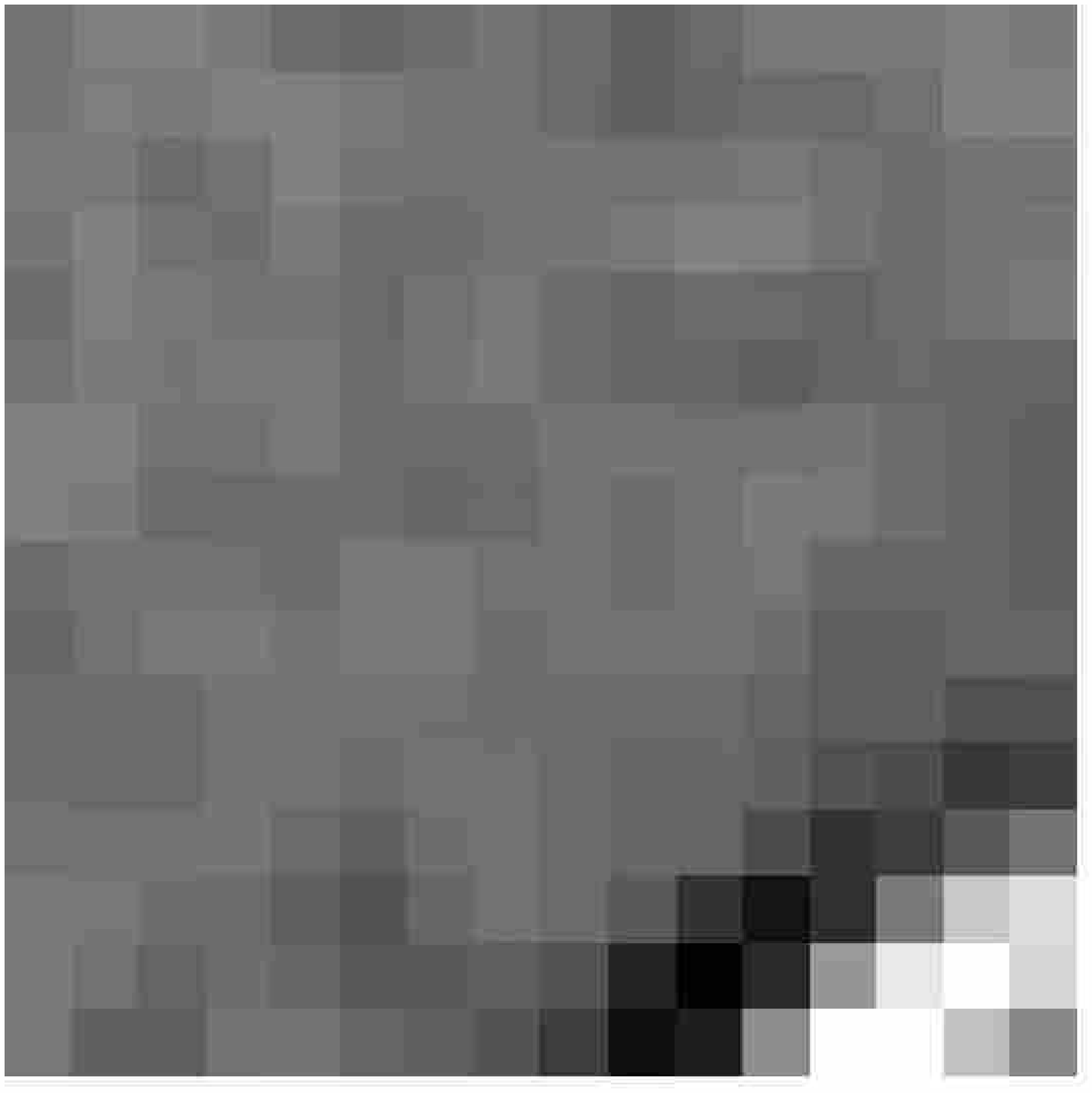}&
\includegraphics[width=\plotwidth]{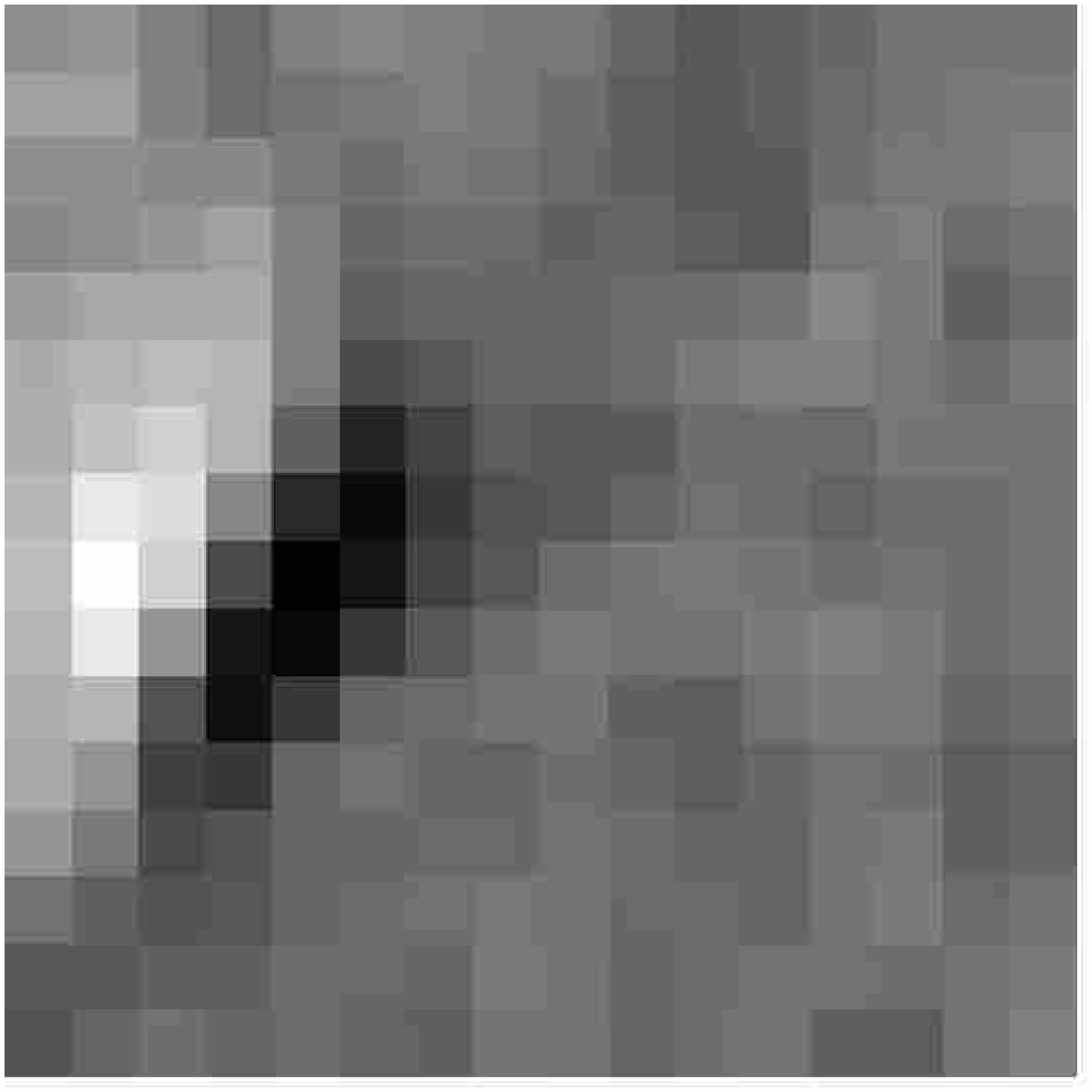}&
\includegraphics[width=\plotwidth]{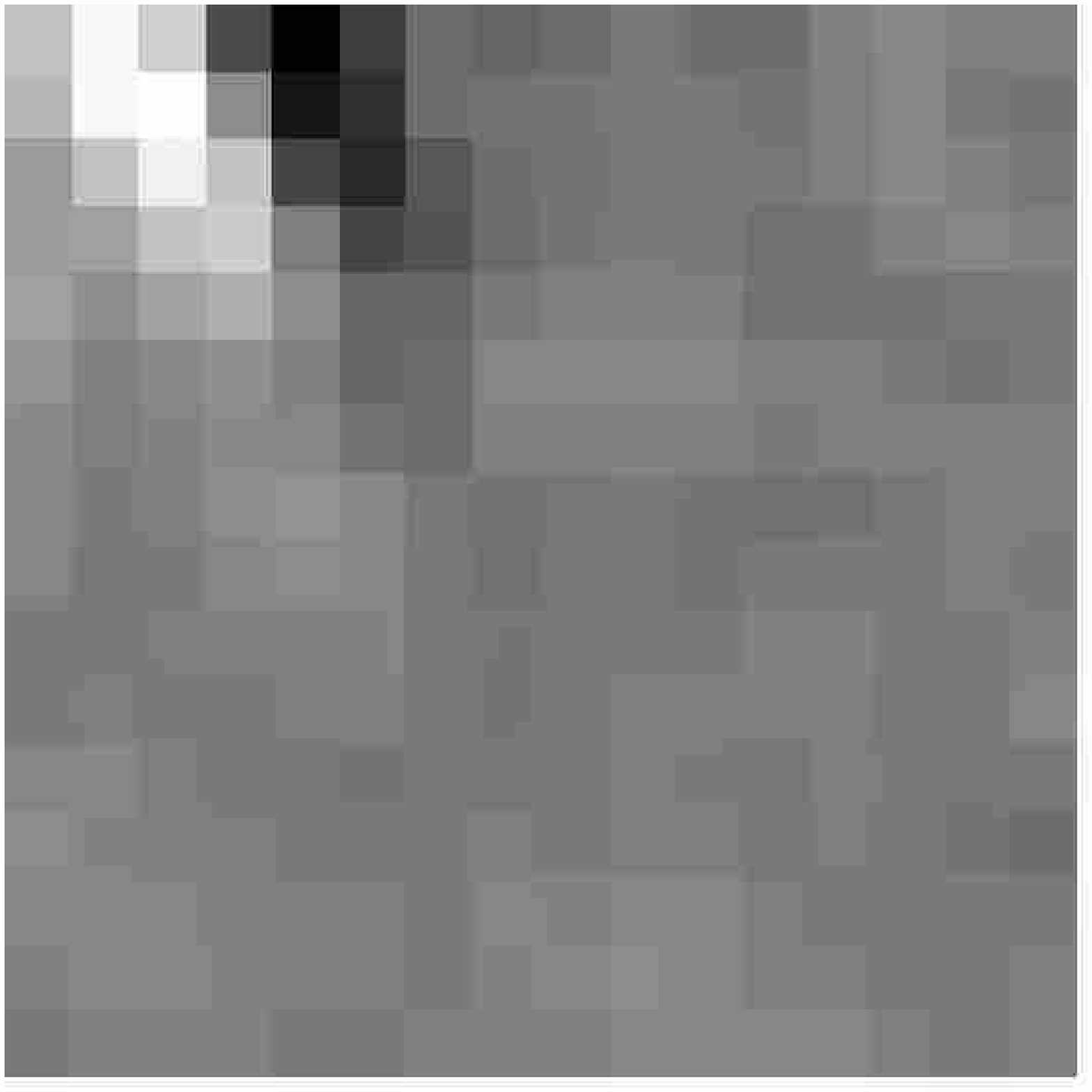}&
\includegraphics[width=\plotwidth]{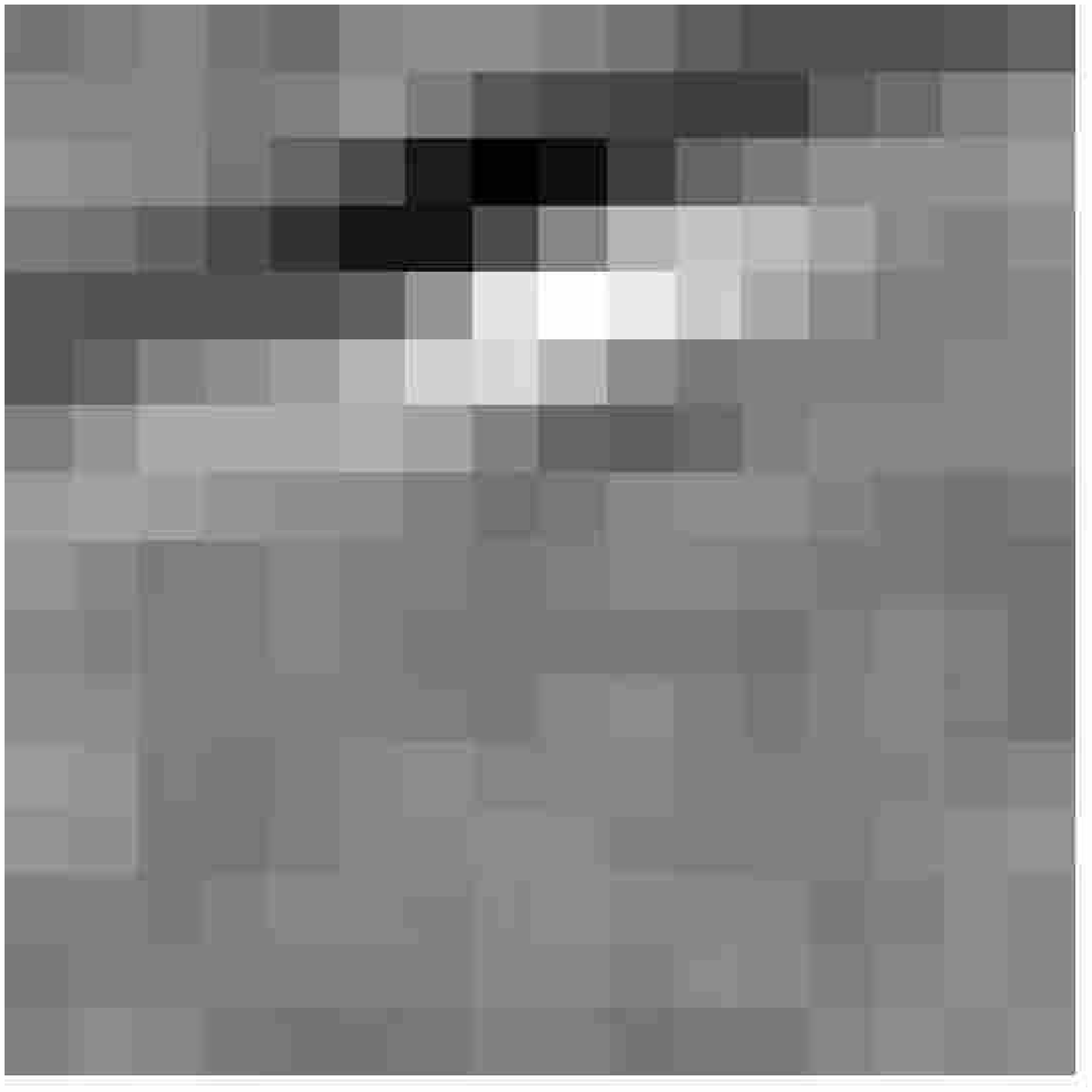}\\
\includegraphics[width=\plotwidth]{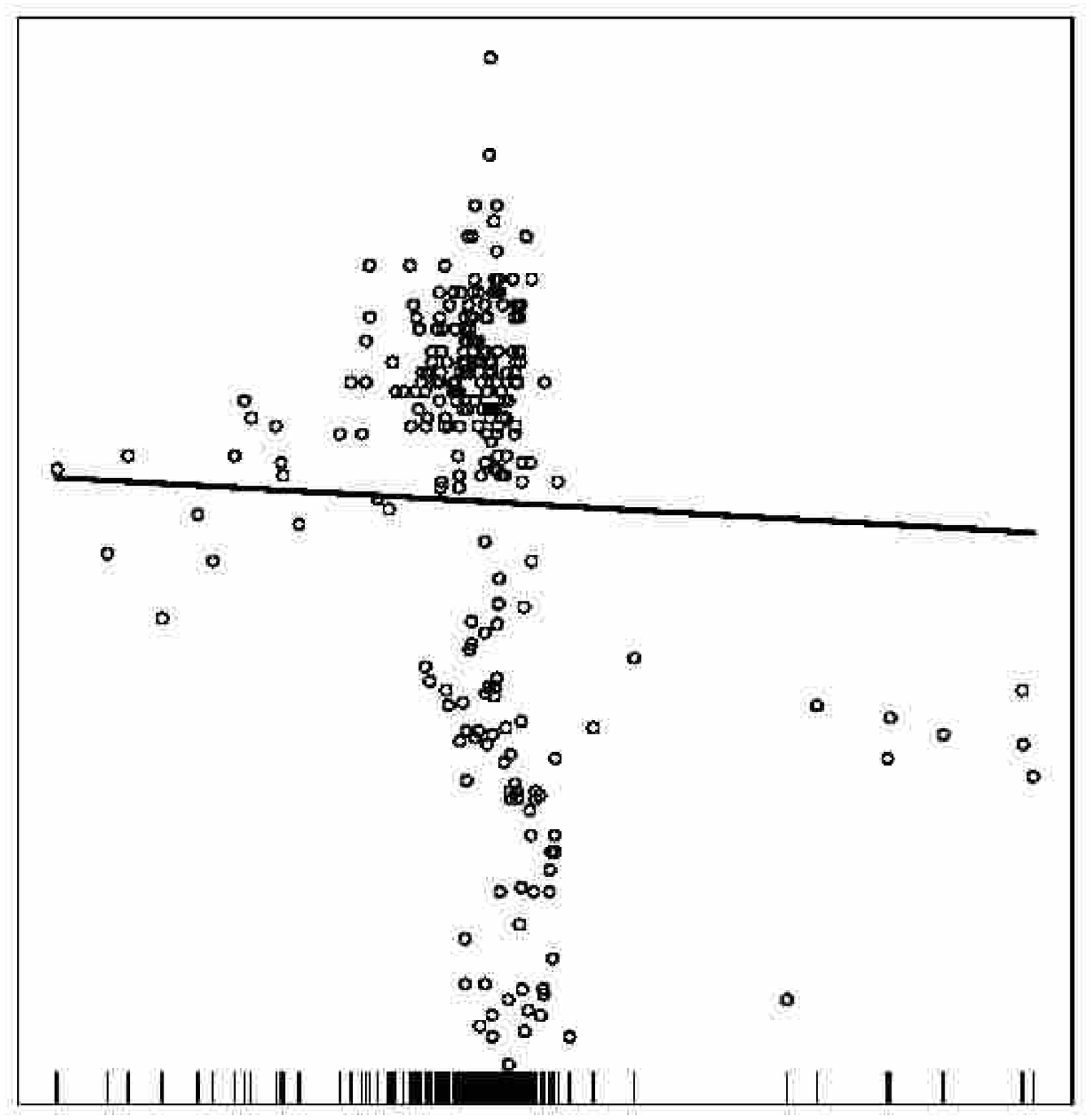}&
\includegraphics[width=\plotwidth]{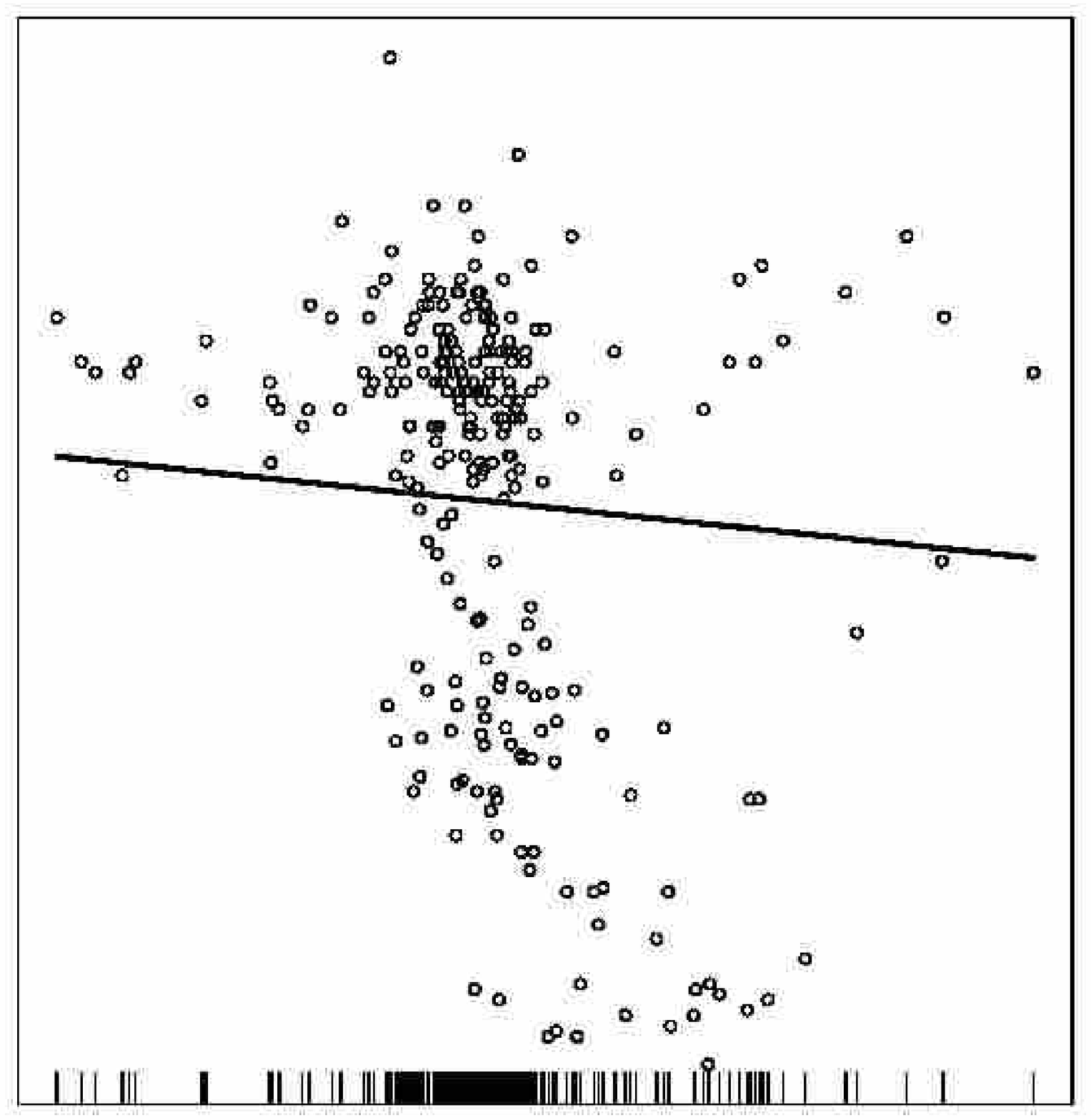}&
\includegraphics[width=\plotwidth]{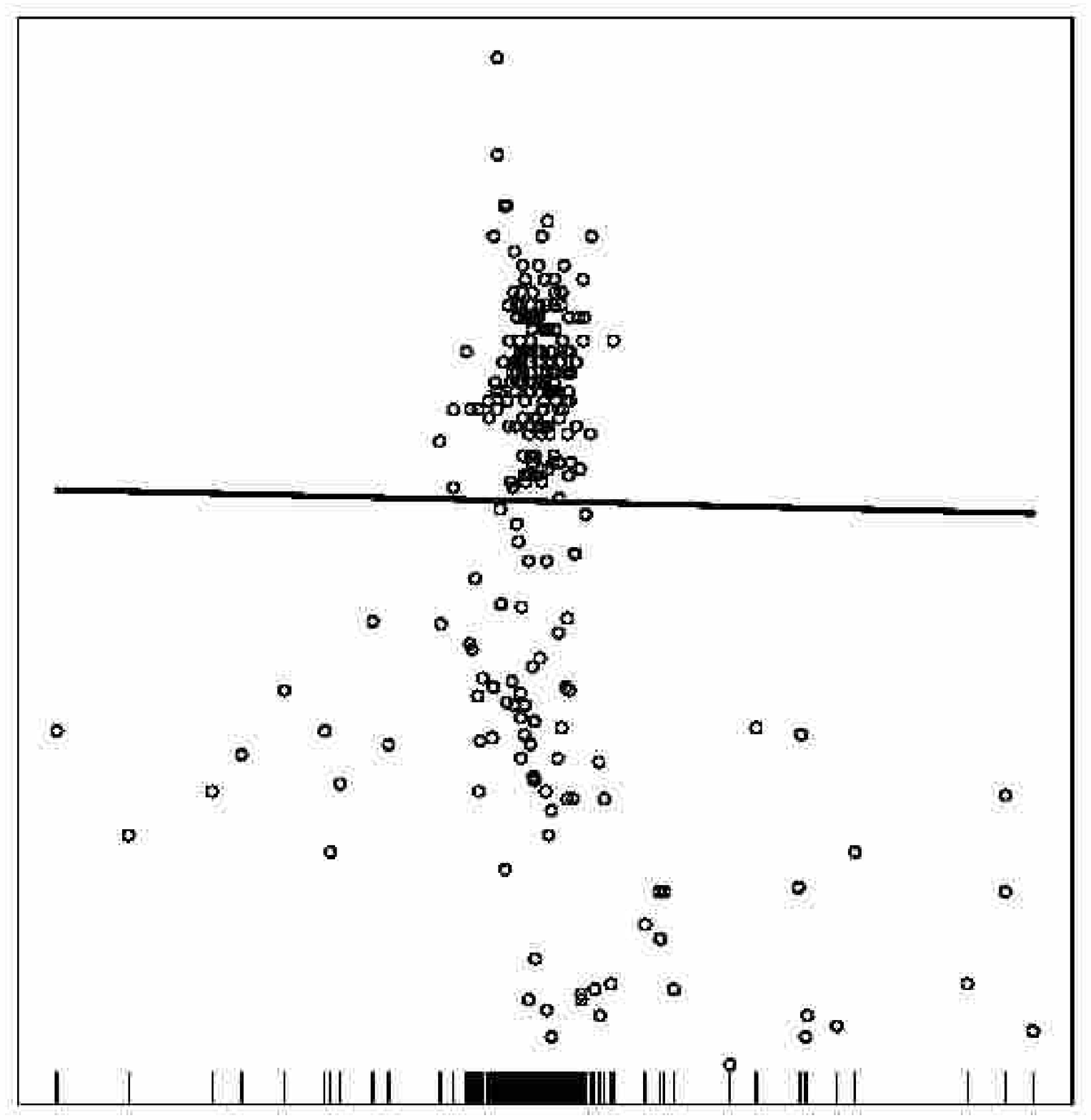}&
\includegraphics[width=\plotwidth]{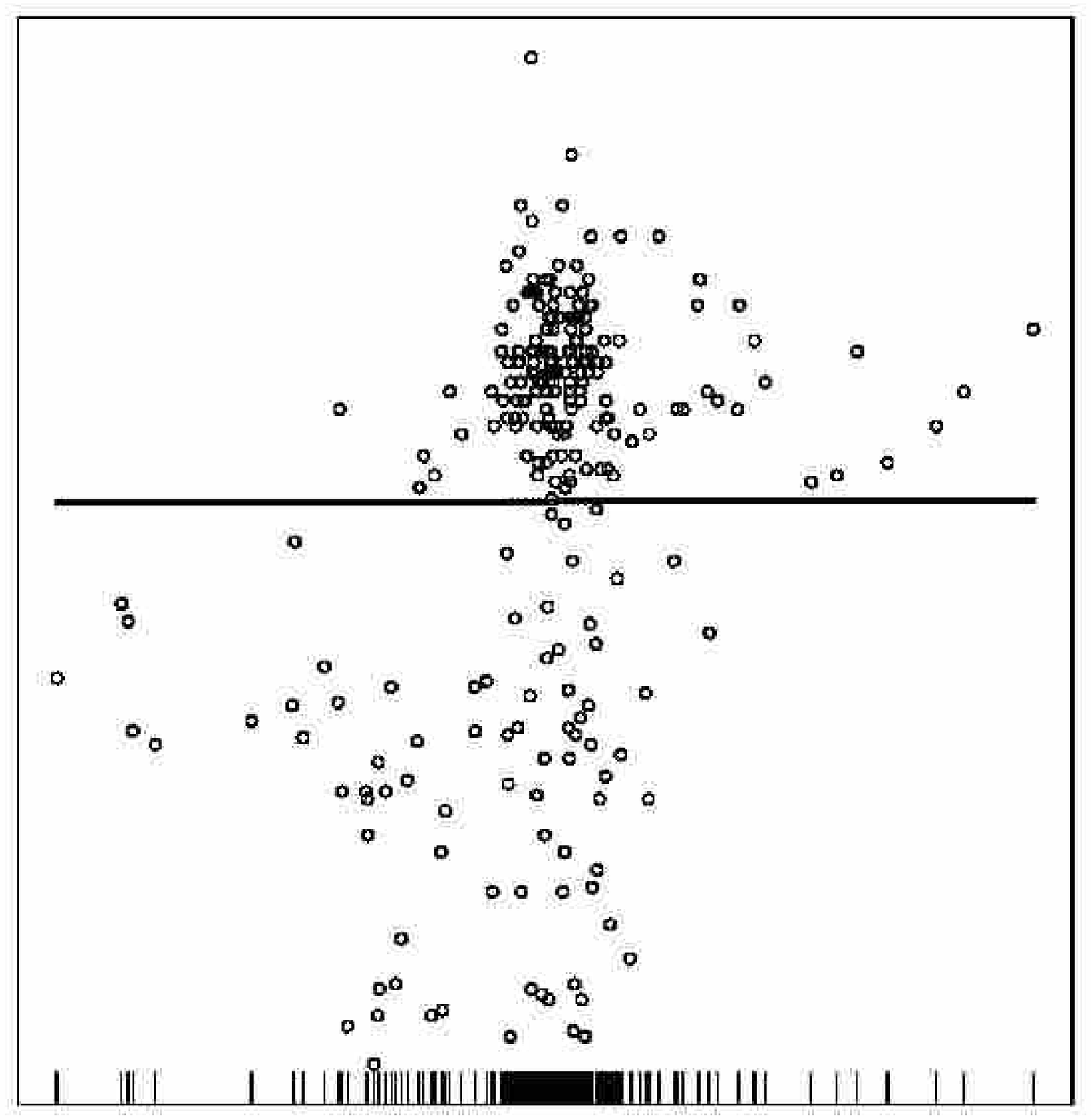}\\
\includegraphics[width=\plotwidth]{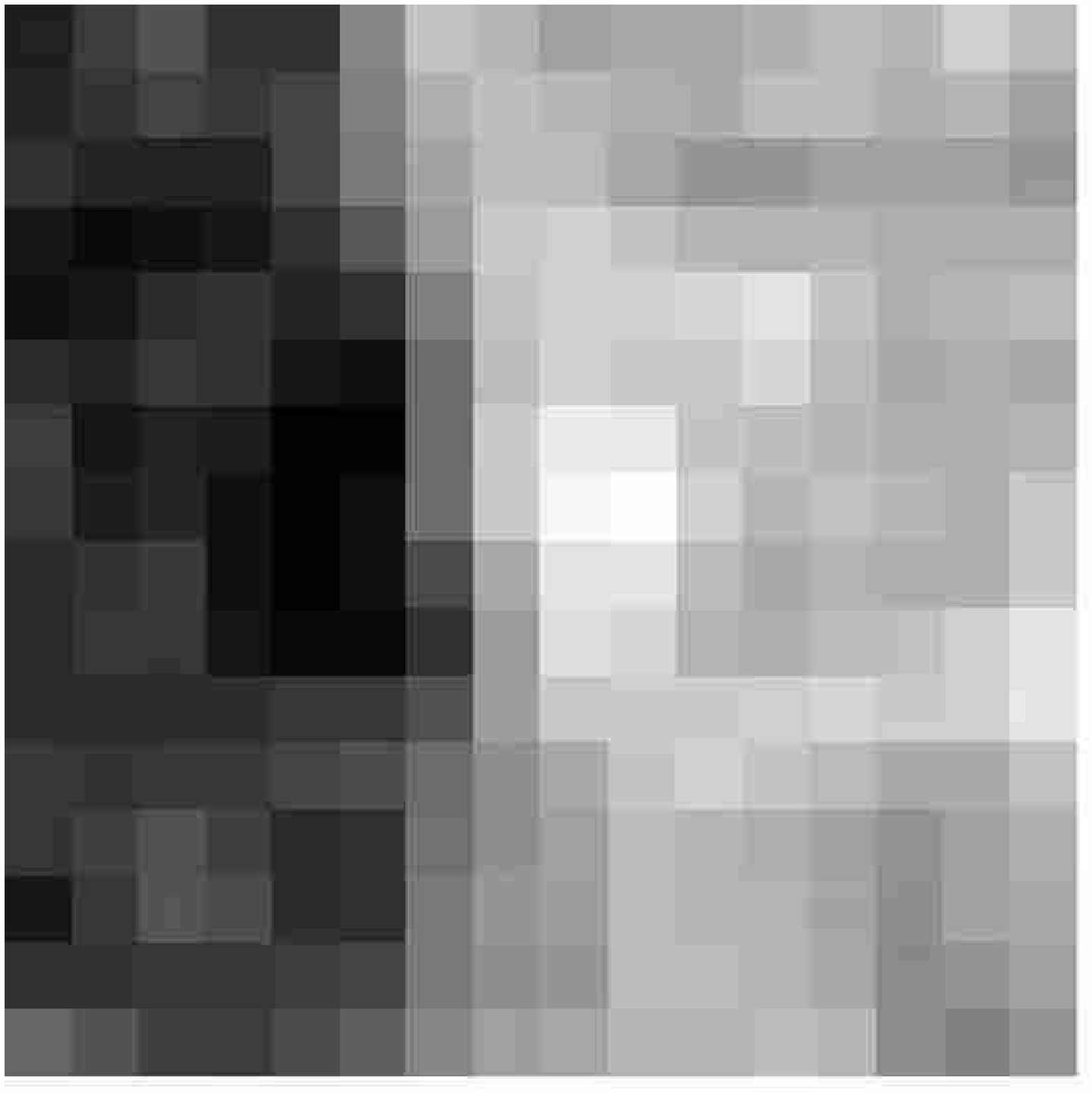}&
\includegraphics[width=\plotwidth]{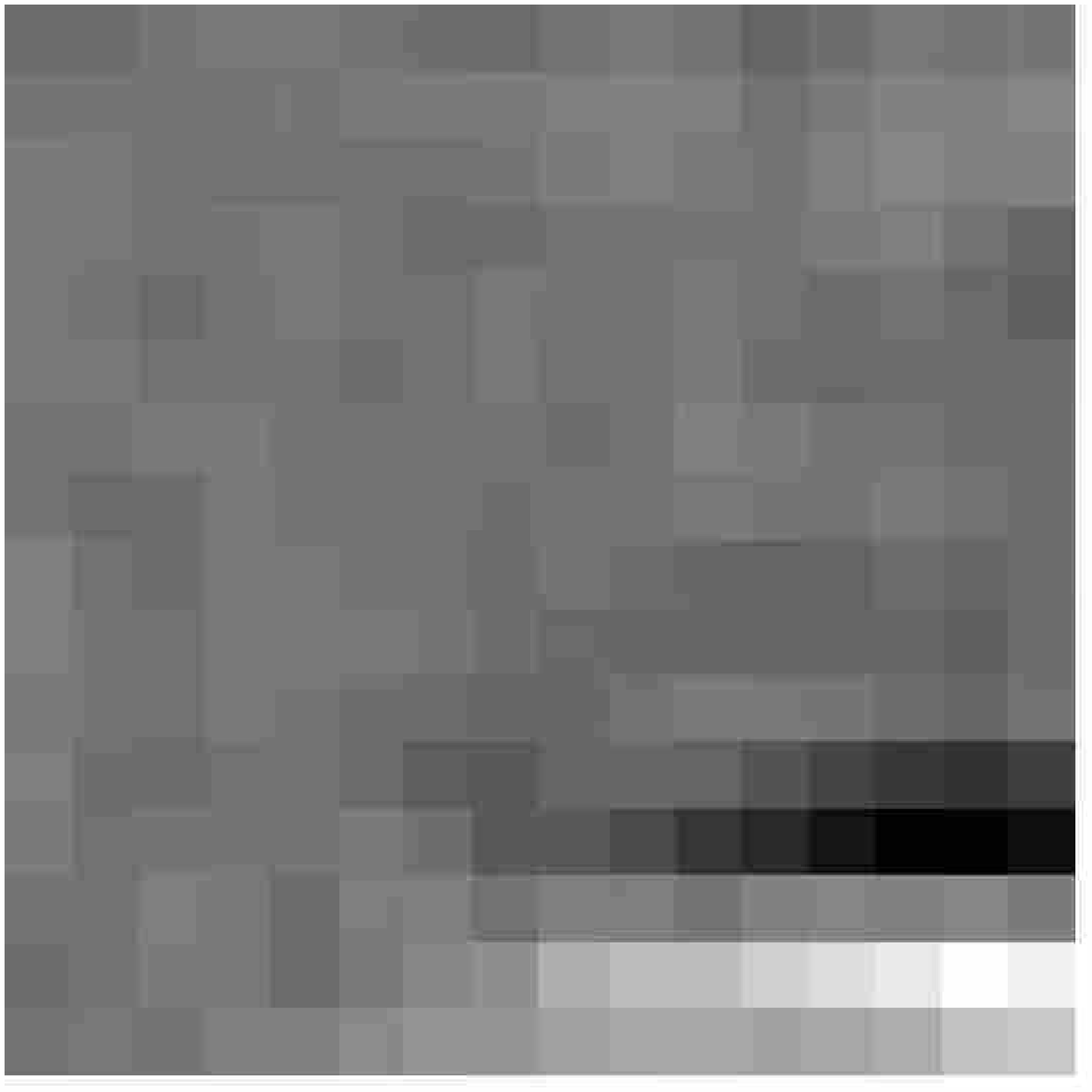}&&\\
\includegraphics[width=\plotwidth]{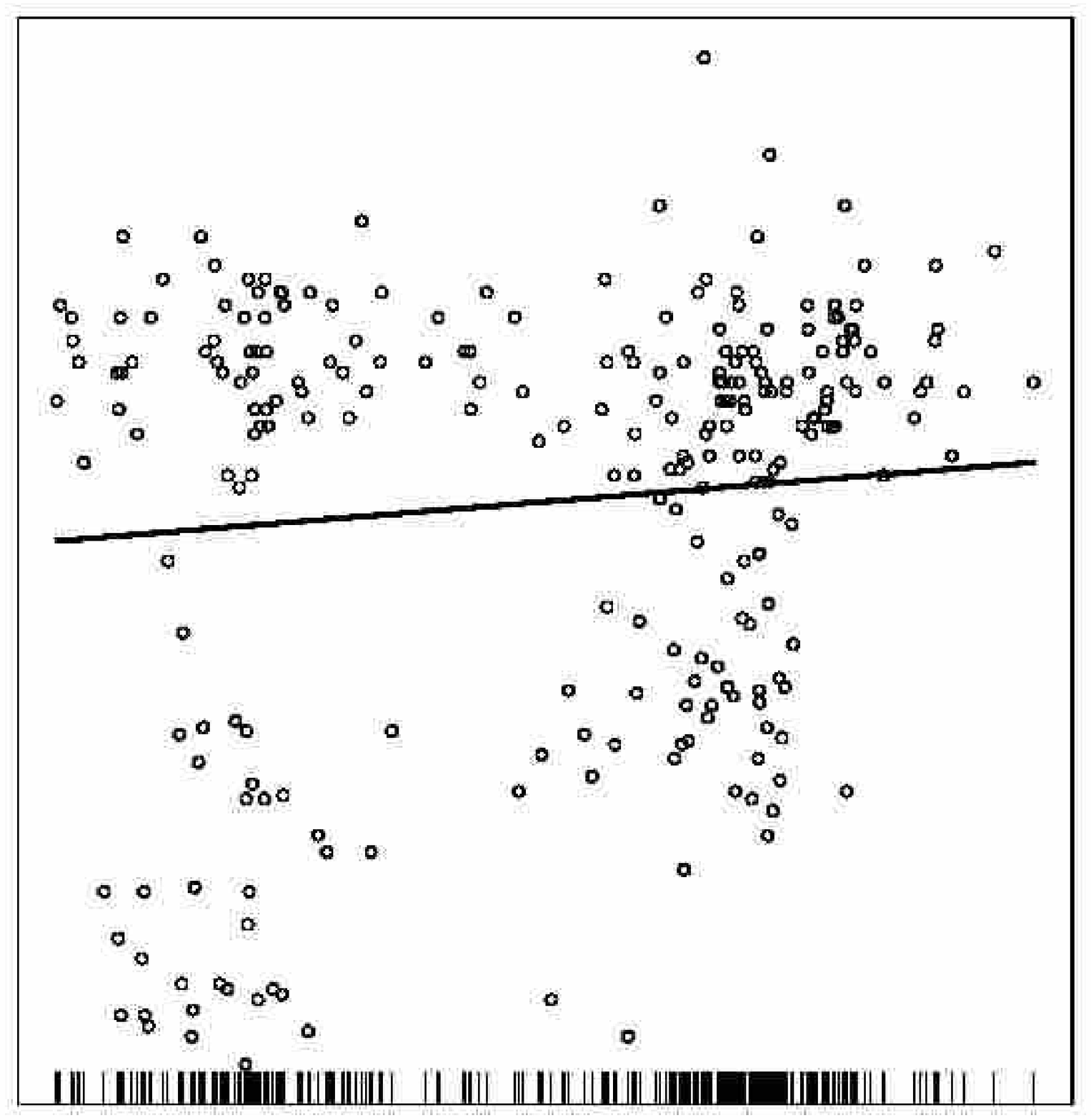}&
\includegraphics[width=\plotwidth]{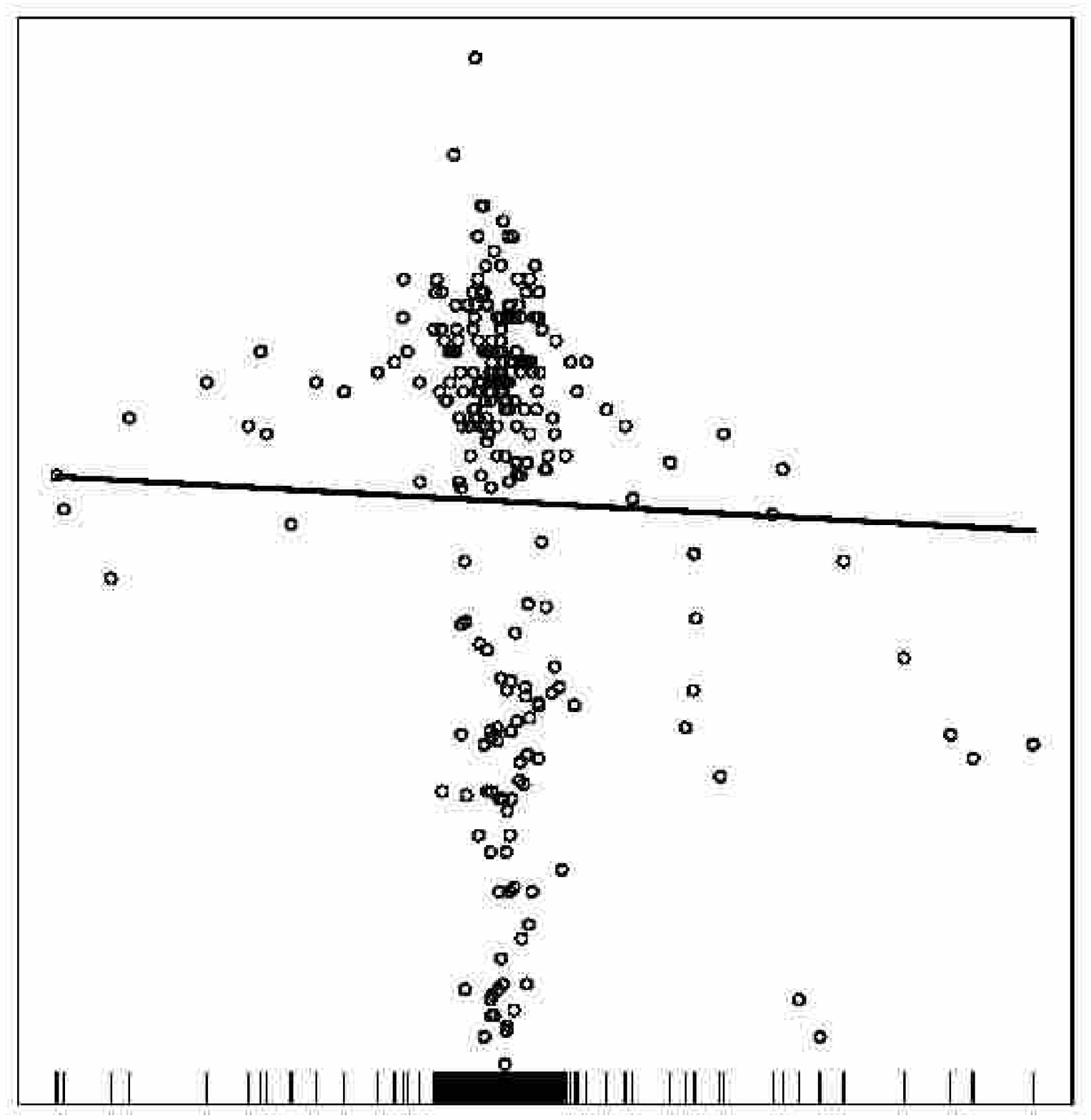}&&\\
\end{tabular}
\end{tabular}
&
\begin{tabular}{c}
\begin{tabular}{cc}
\multicolumn{2}{c}{SpAM} \\
\small Original patch & \small $\text{RSS}=0.0887$ \\
\hskip-5pt
\includegraphics[width=\headwidth]{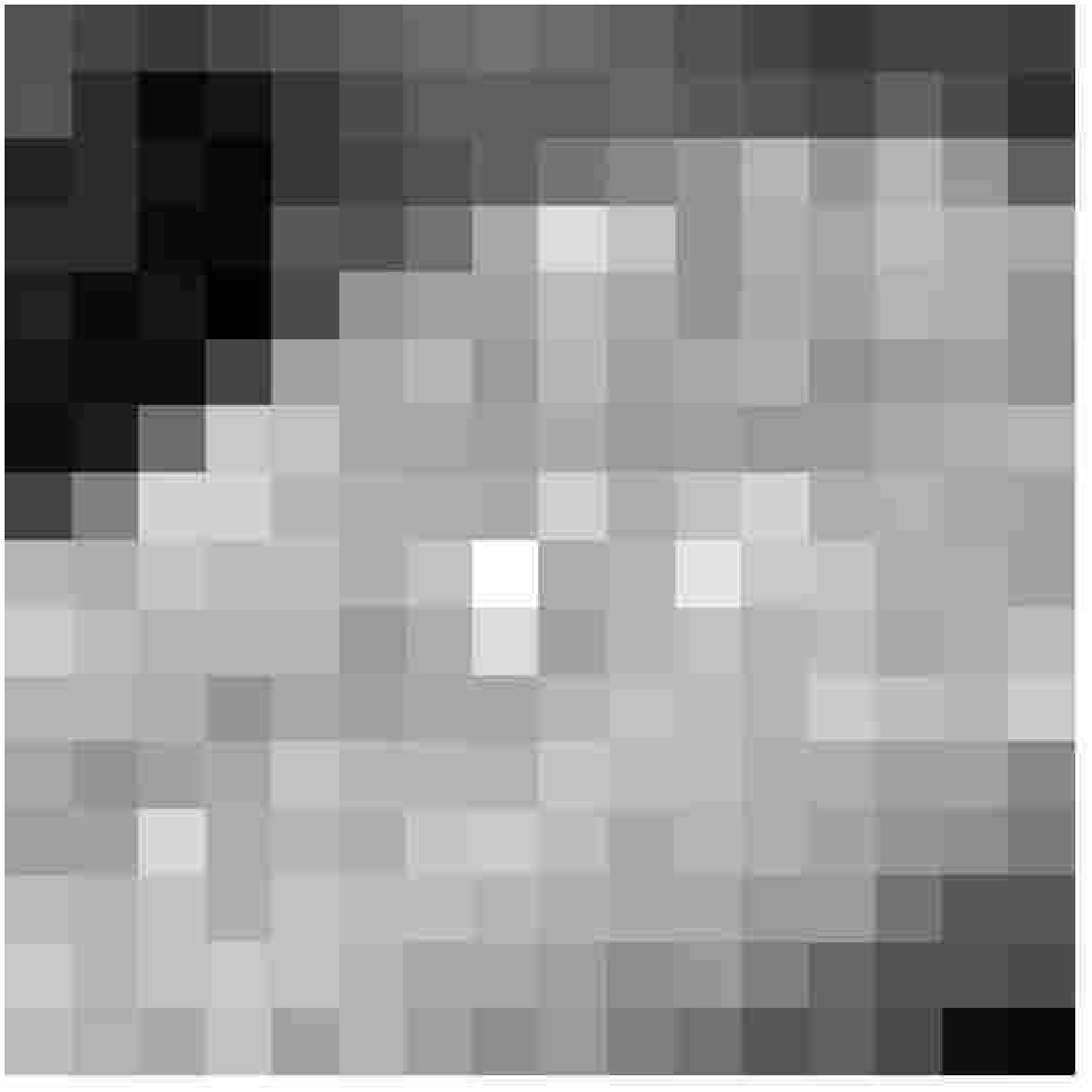} &
\includegraphics[width=\headwidth]{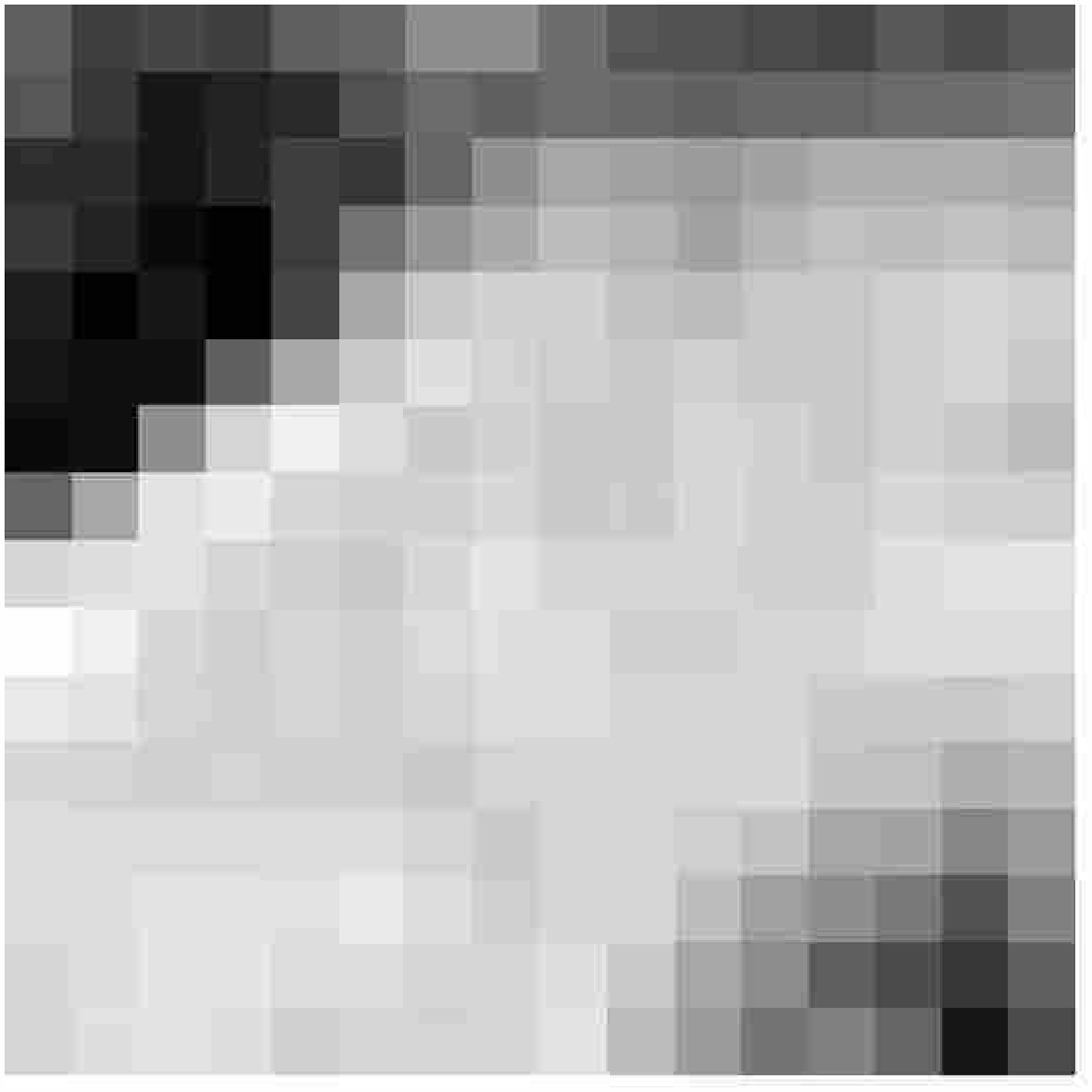}
\\[10pt]
\end{tabular}
\\
\hskip-10pt
\begin{tabular}{cccc}
\includegraphics[width=\plotwidth]{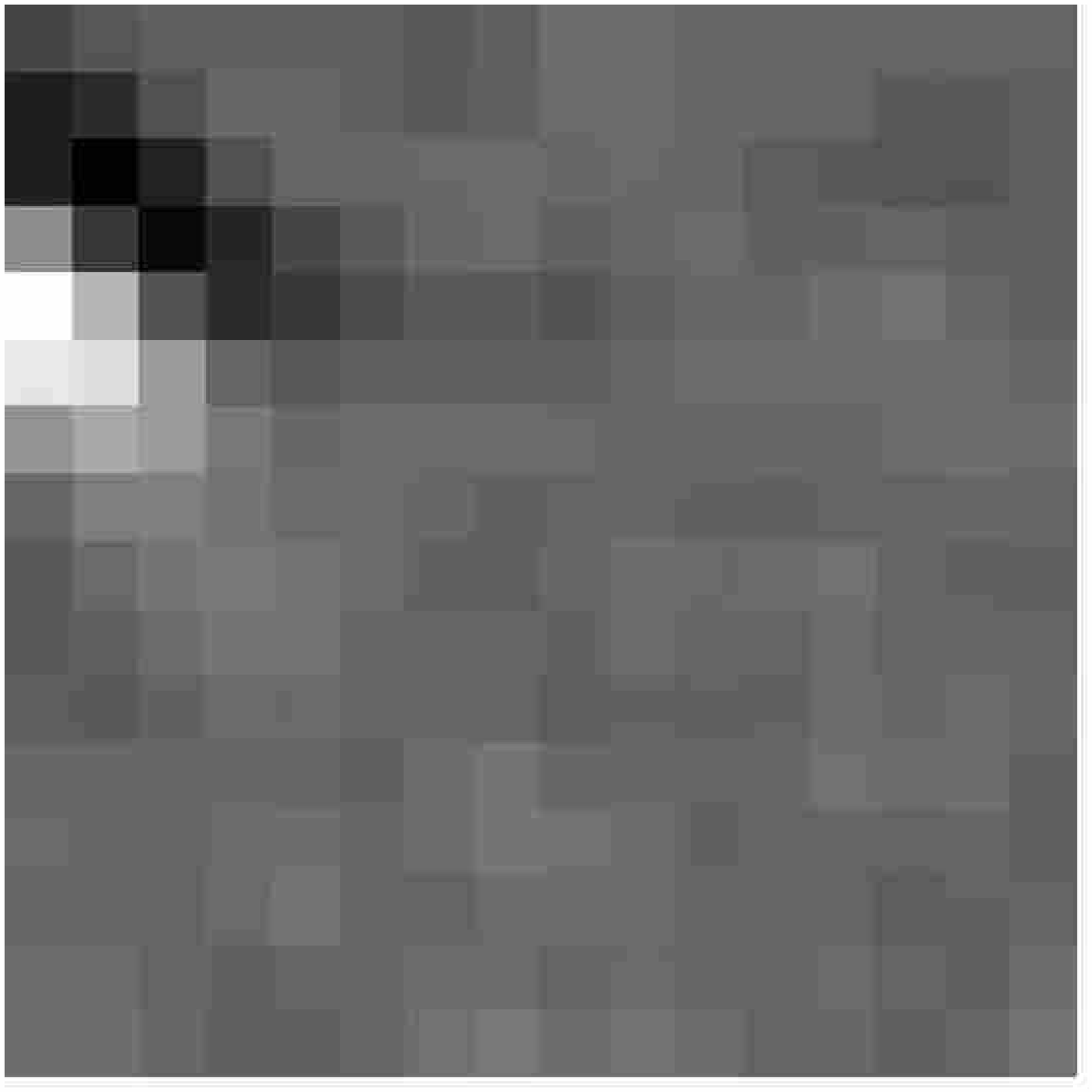}&
\includegraphics[width=\plotwidth]{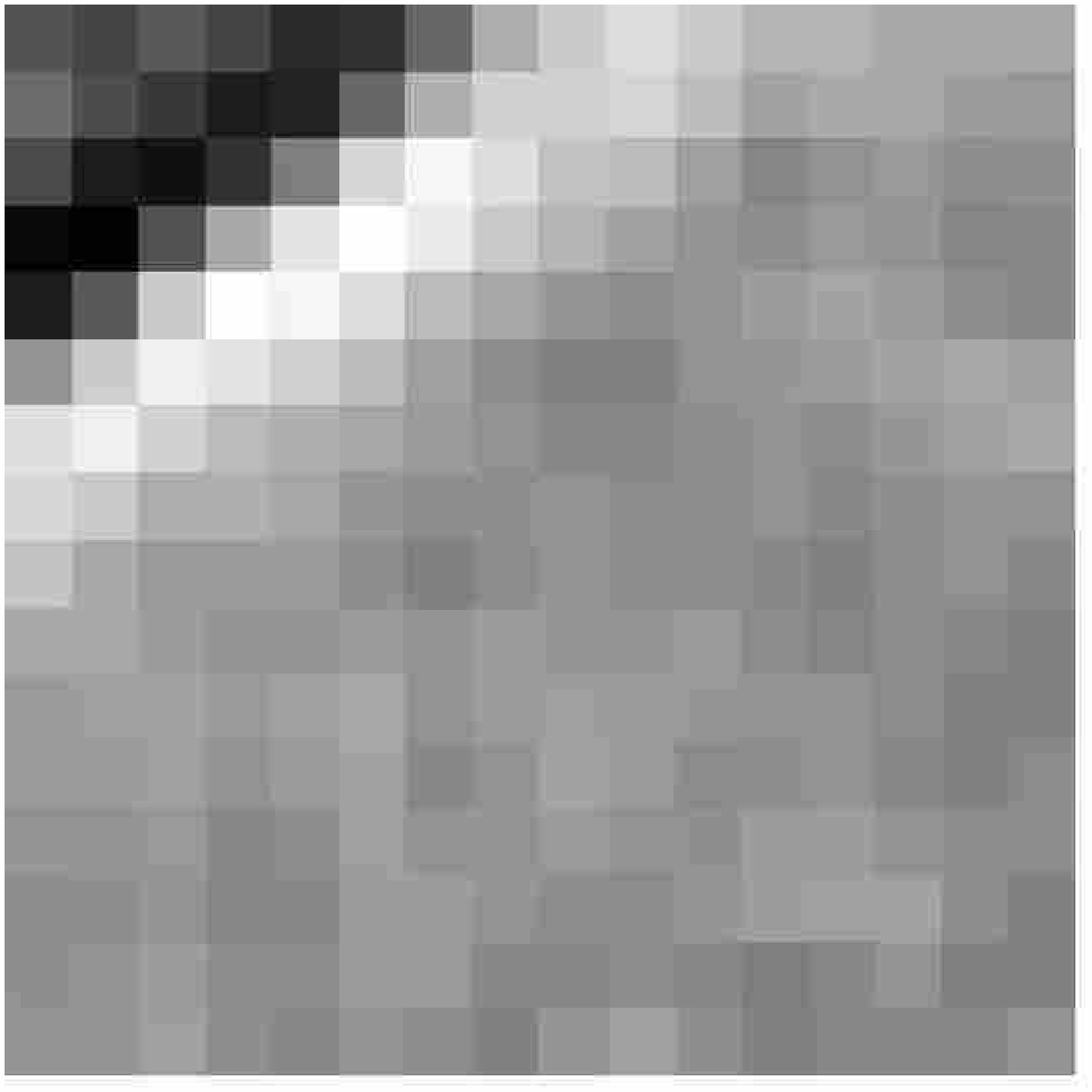}&
\includegraphics[width=\plotwidth]{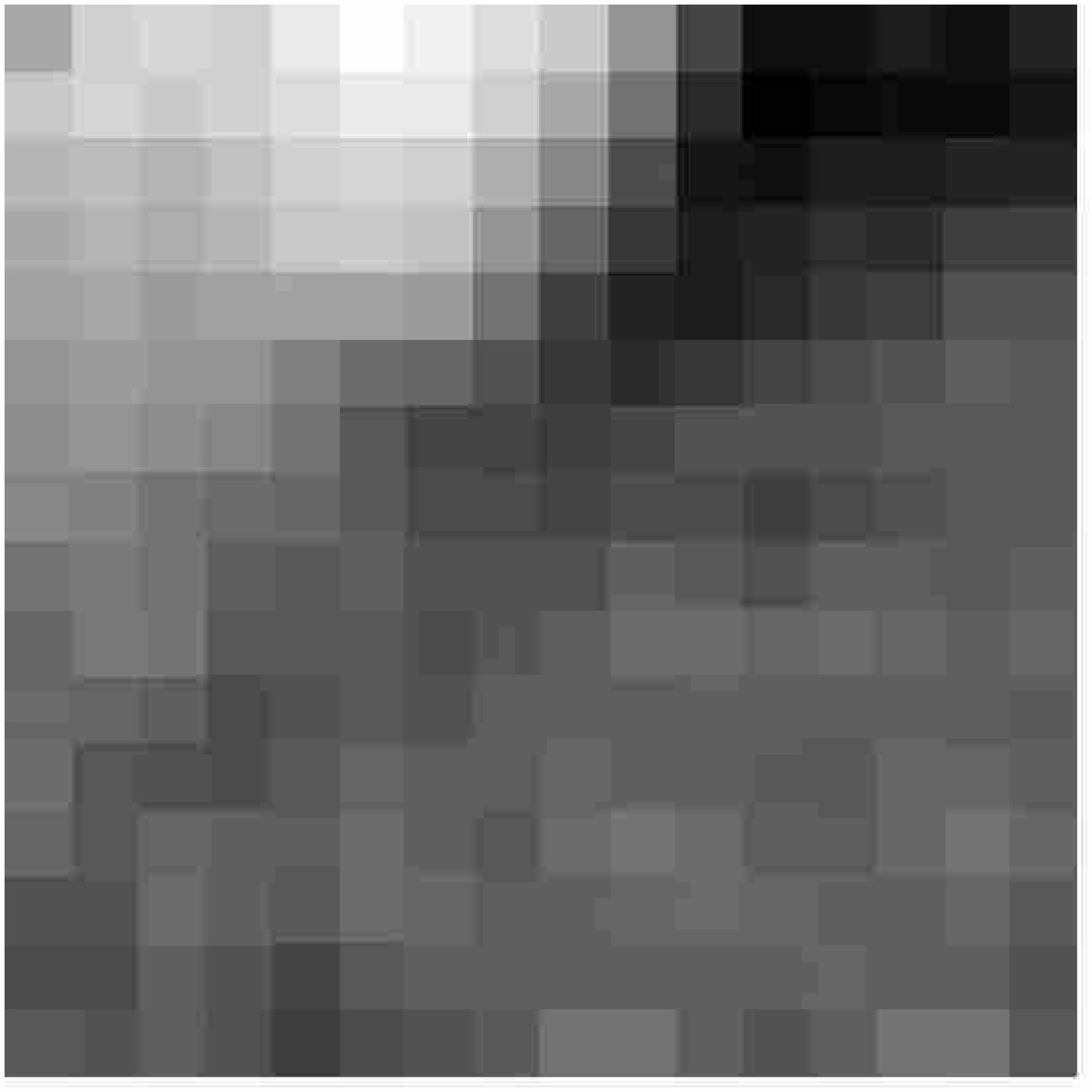}&
\includegraphics[width=\plotwidth]{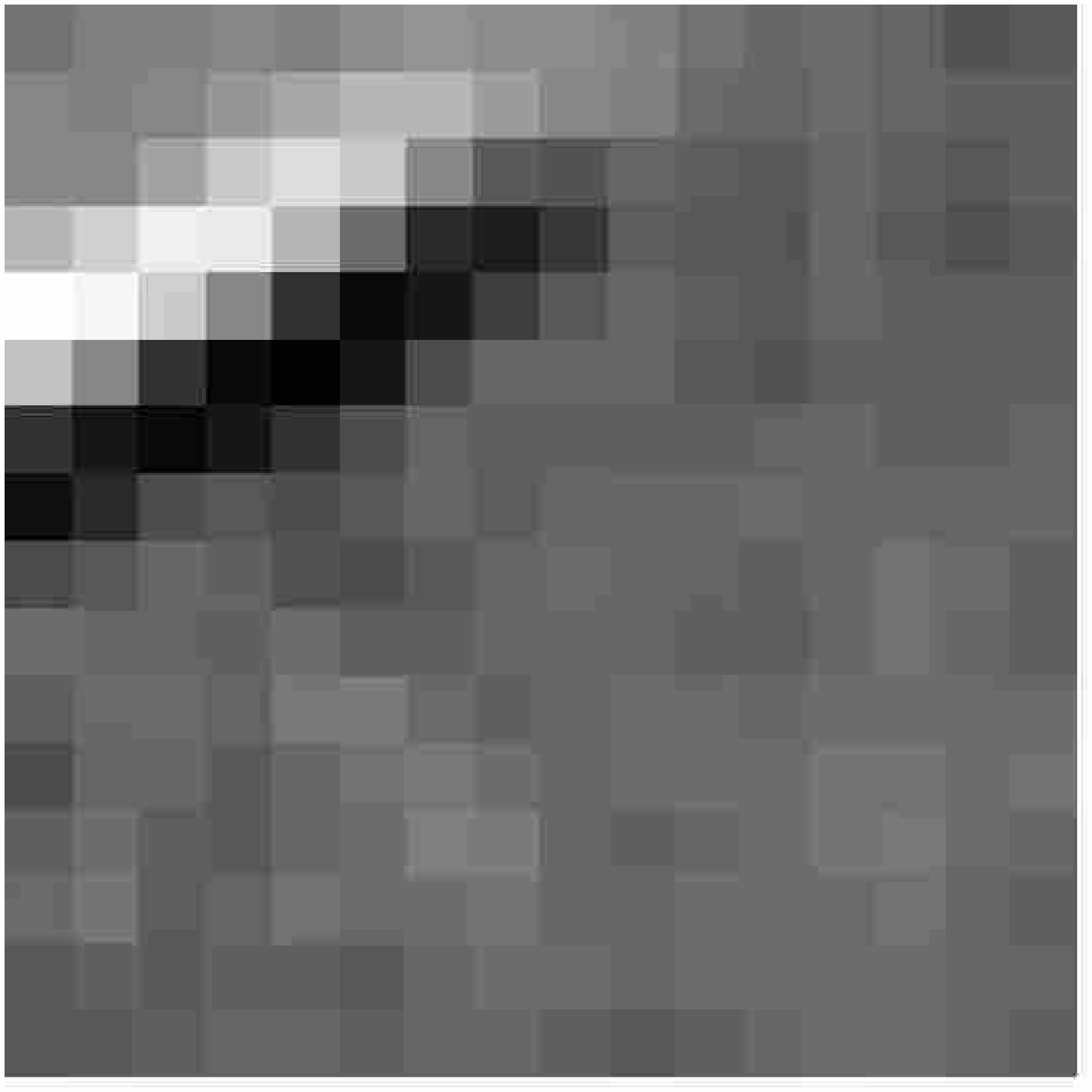}\\
\includegraphics[width=\plotwidth]{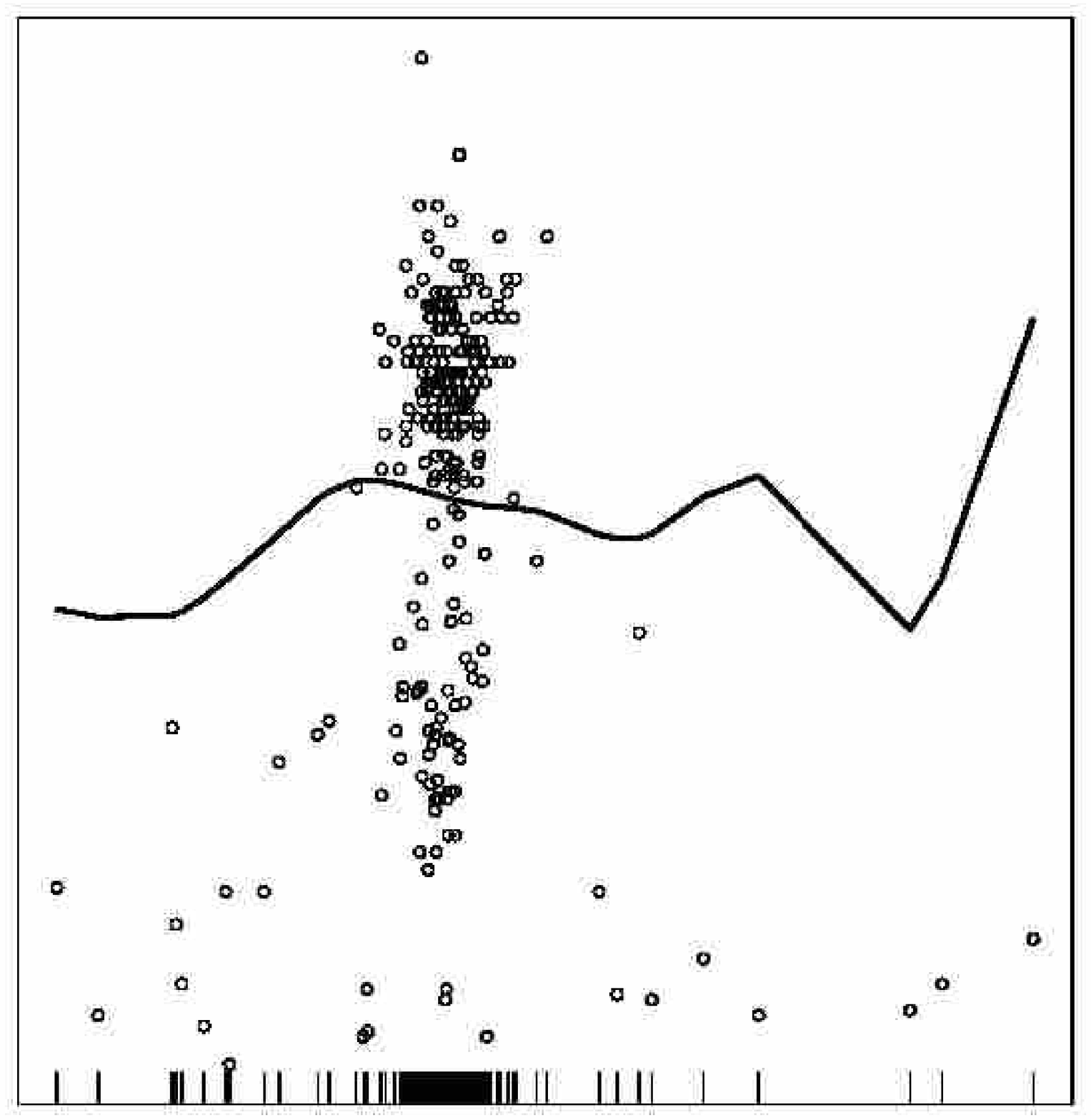}&
\includegraphics[width=\plotwidth]{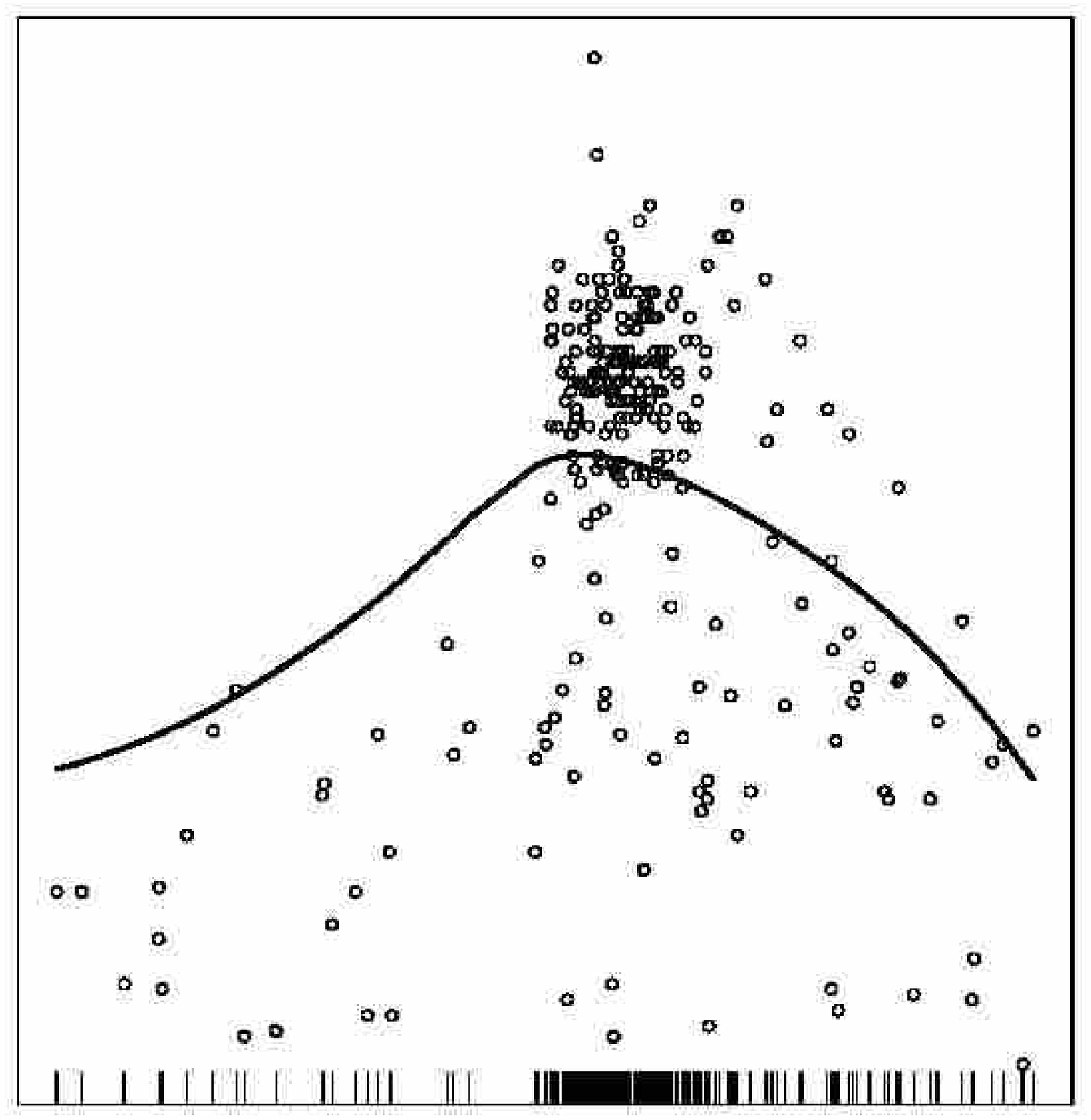}&
\includegraphics[width=\plotwidth]{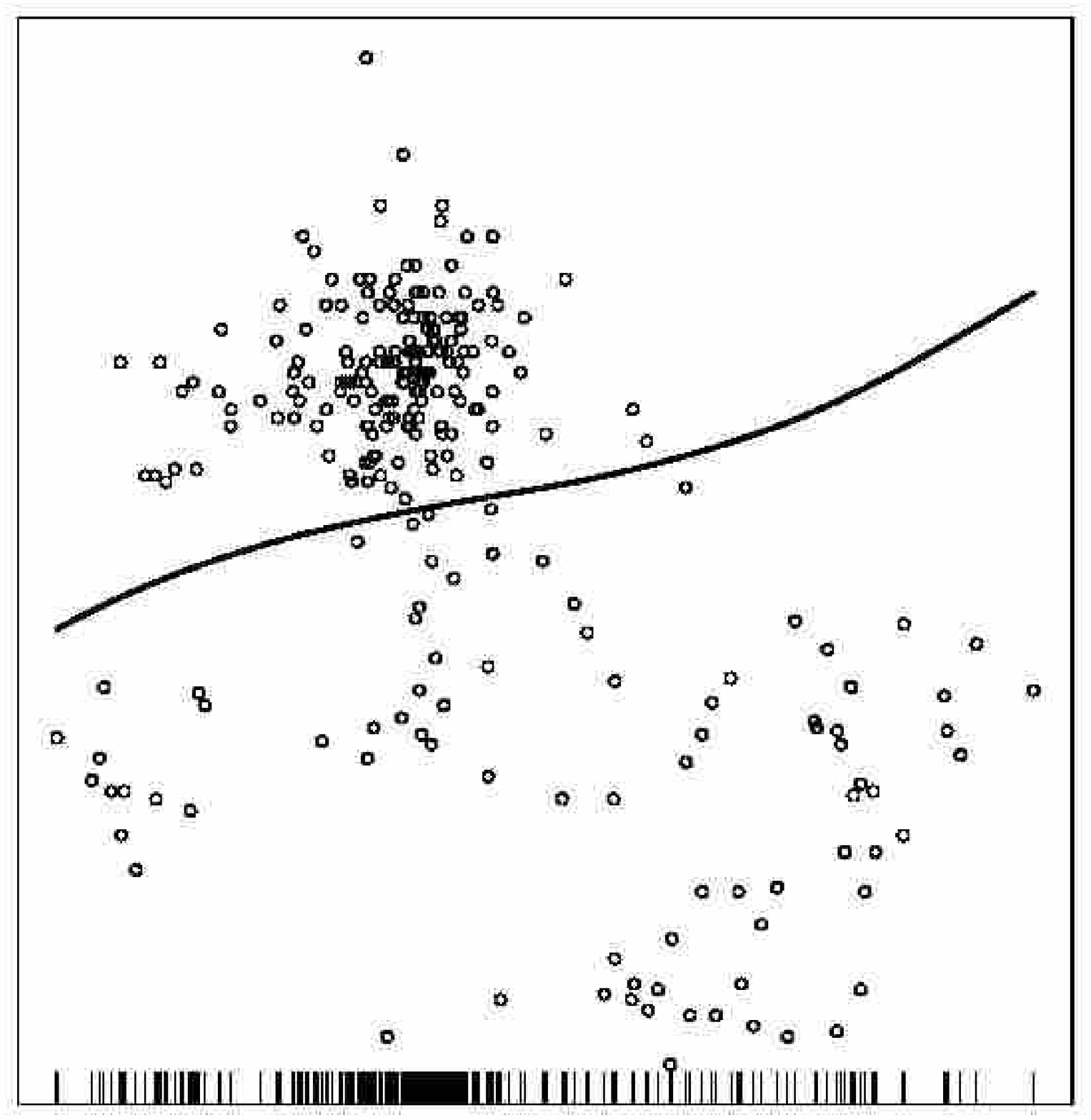}&
\includegraphics[width=\plotwidth]{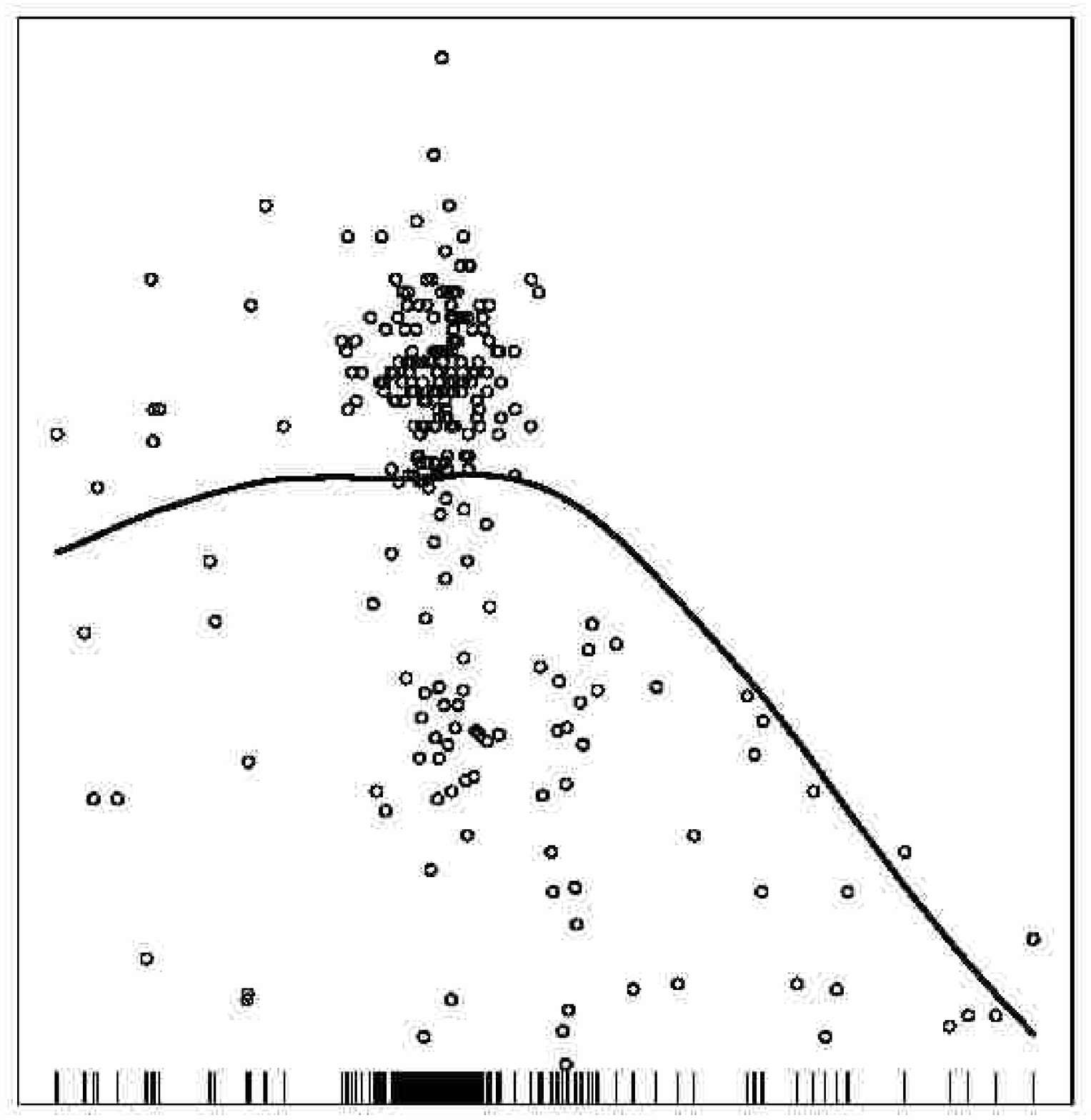}\\
\includegraphics[width=\plotwidth]{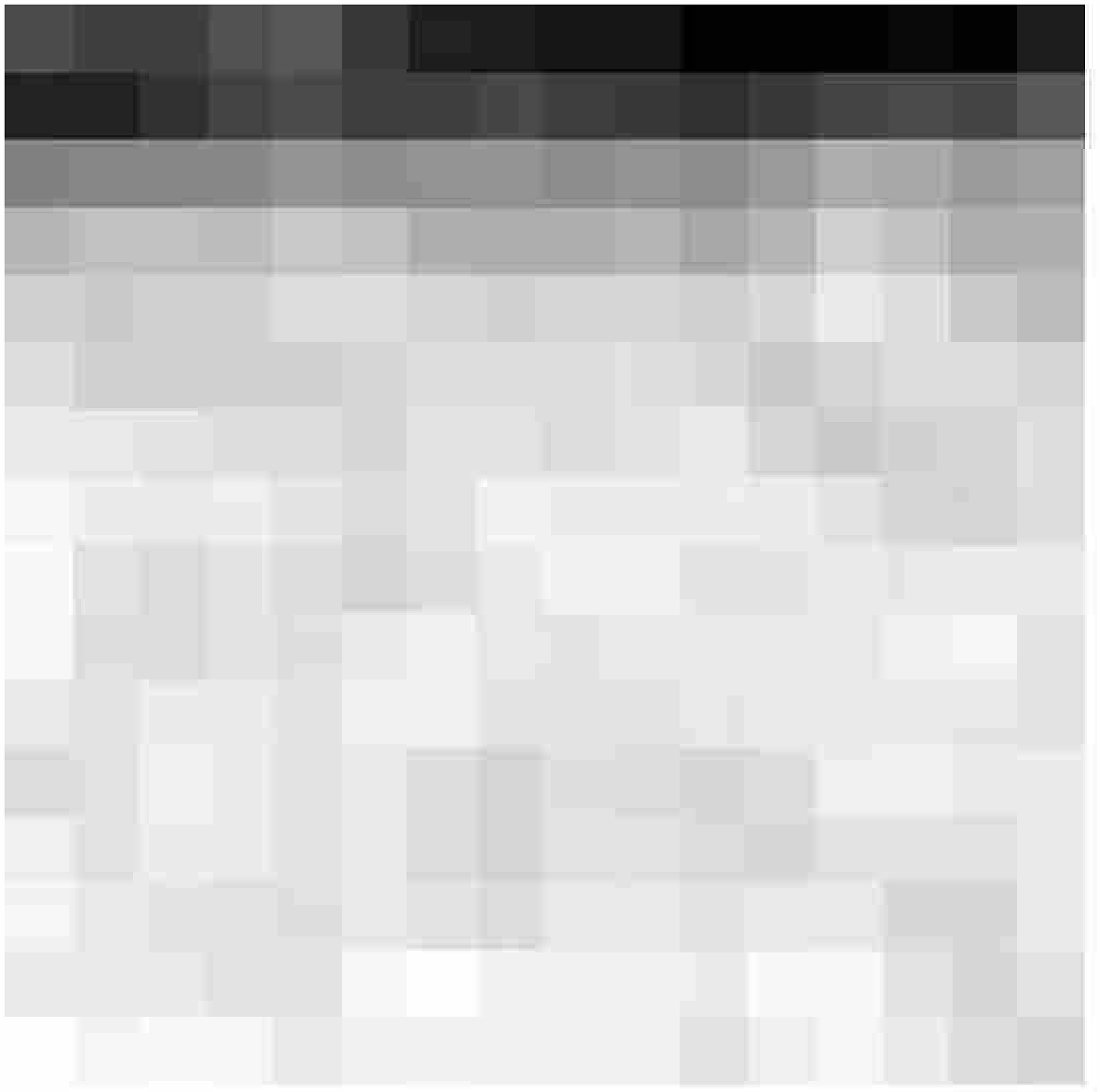}&
\includegraphics[width=\plotwidth]{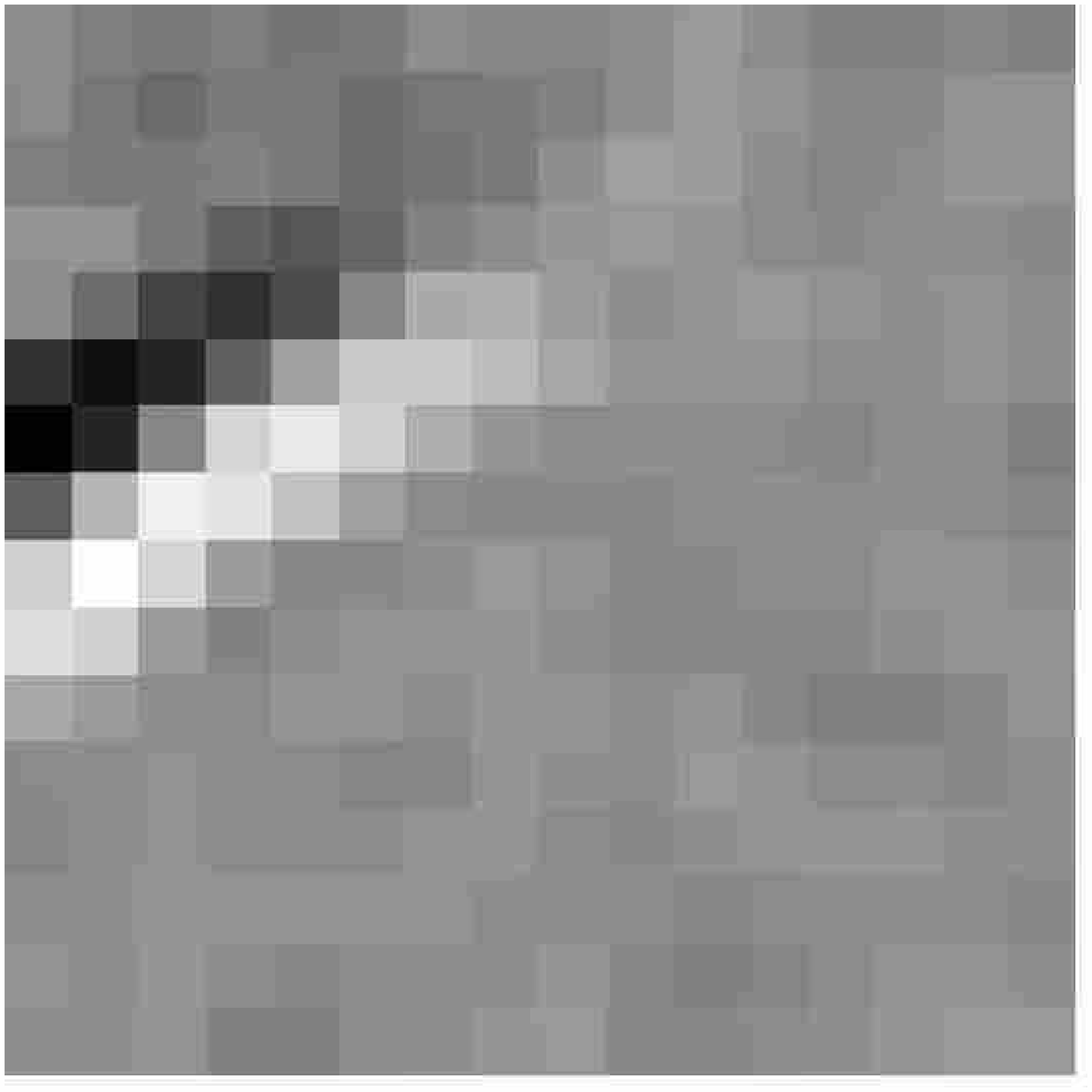}&
\includegraphics[width=\plotwidth]{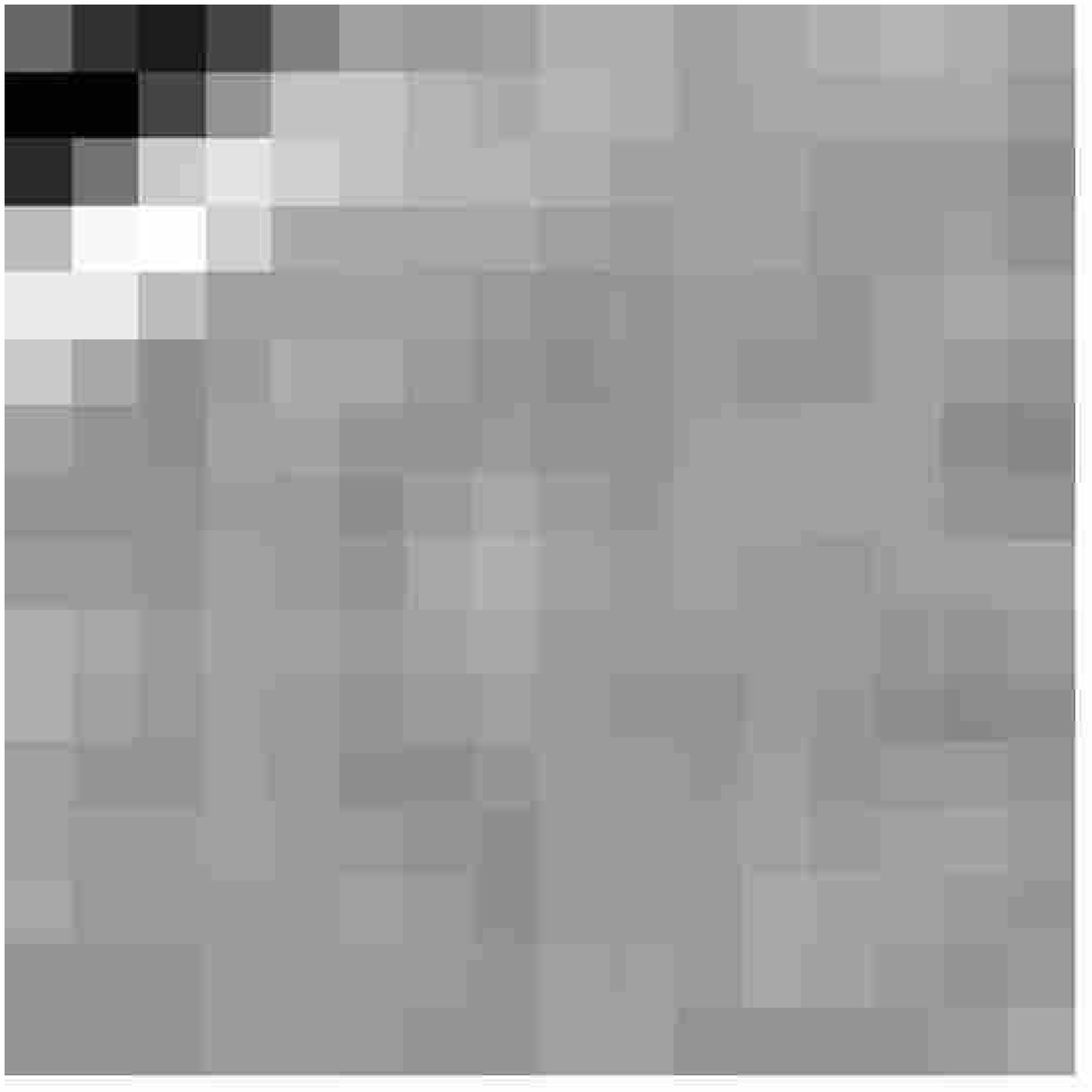}&
\includegraphics[width=\plotwidth]{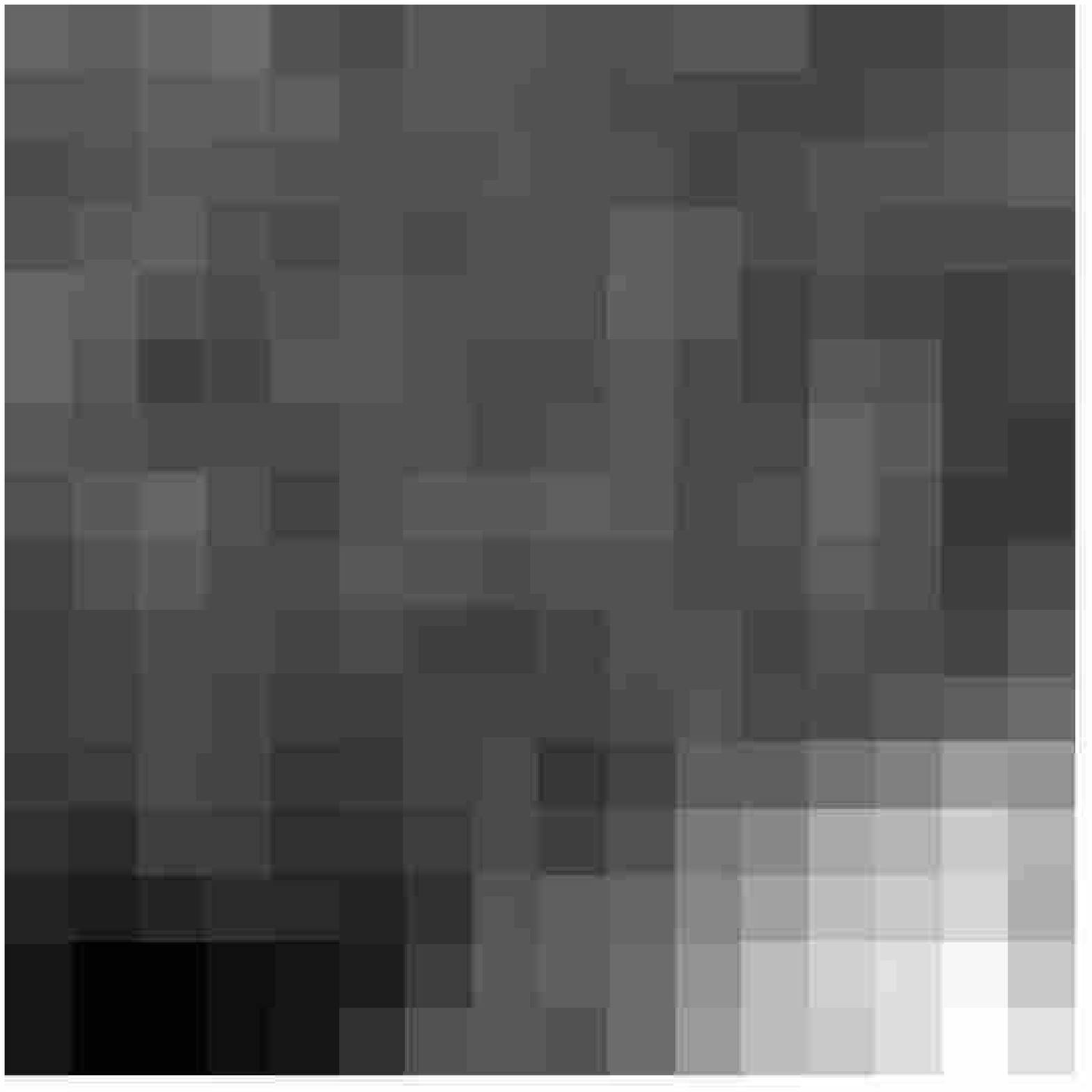}\\
\includegraphics[width=\plotwidth]{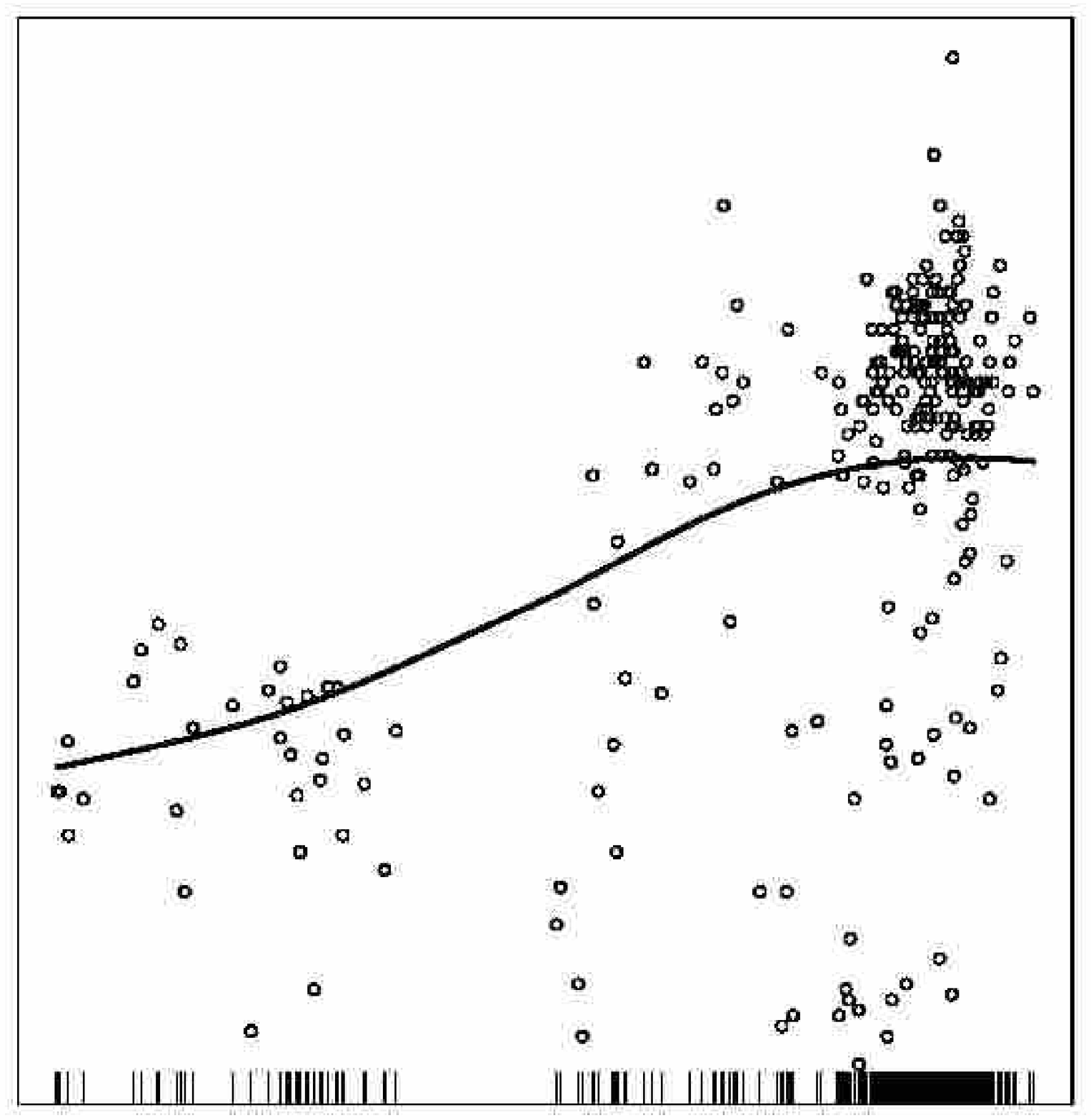}&
\includegraphics[width=\plotwidth]{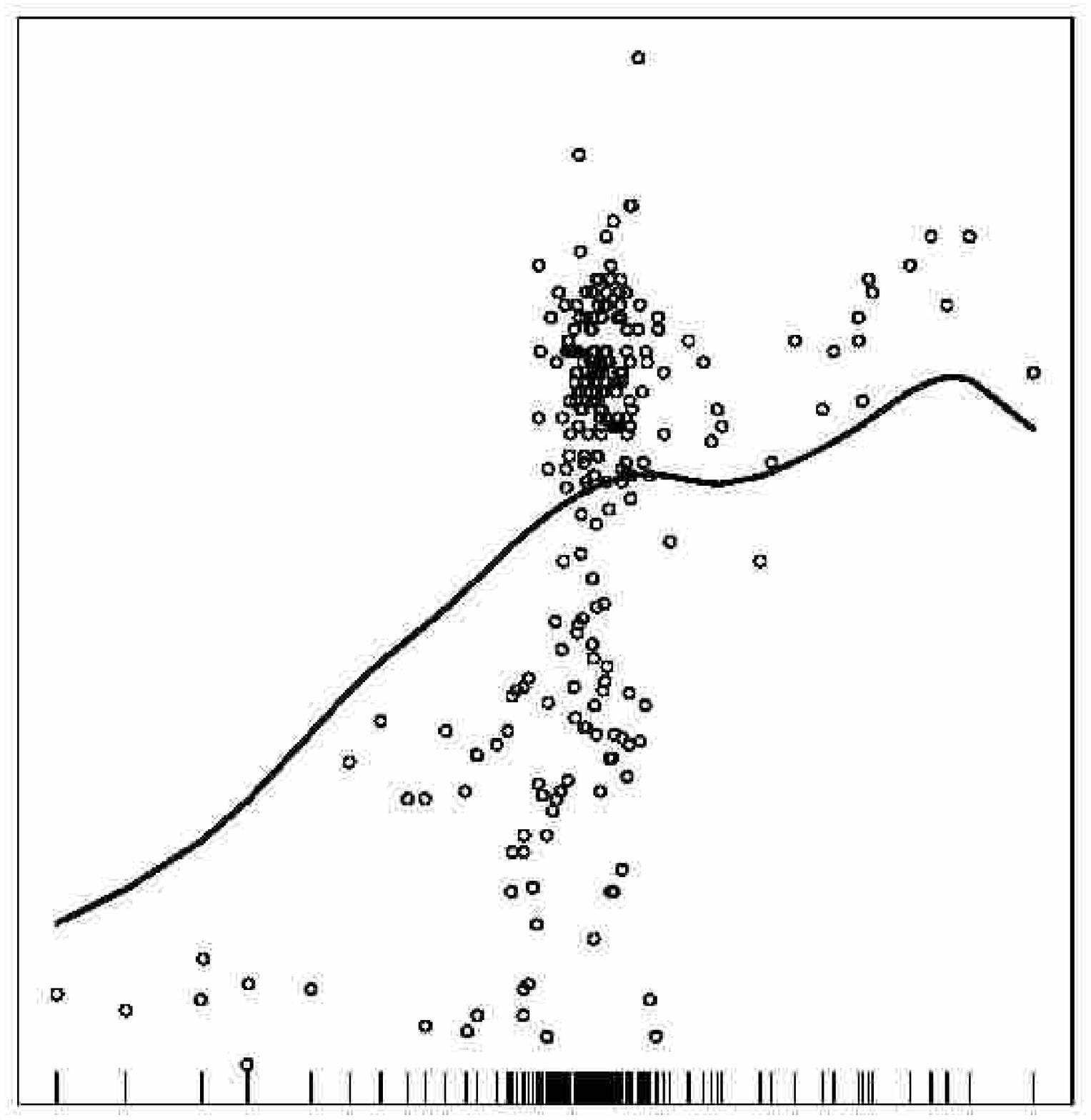}&
\includegraphics[width=\plotwidth]{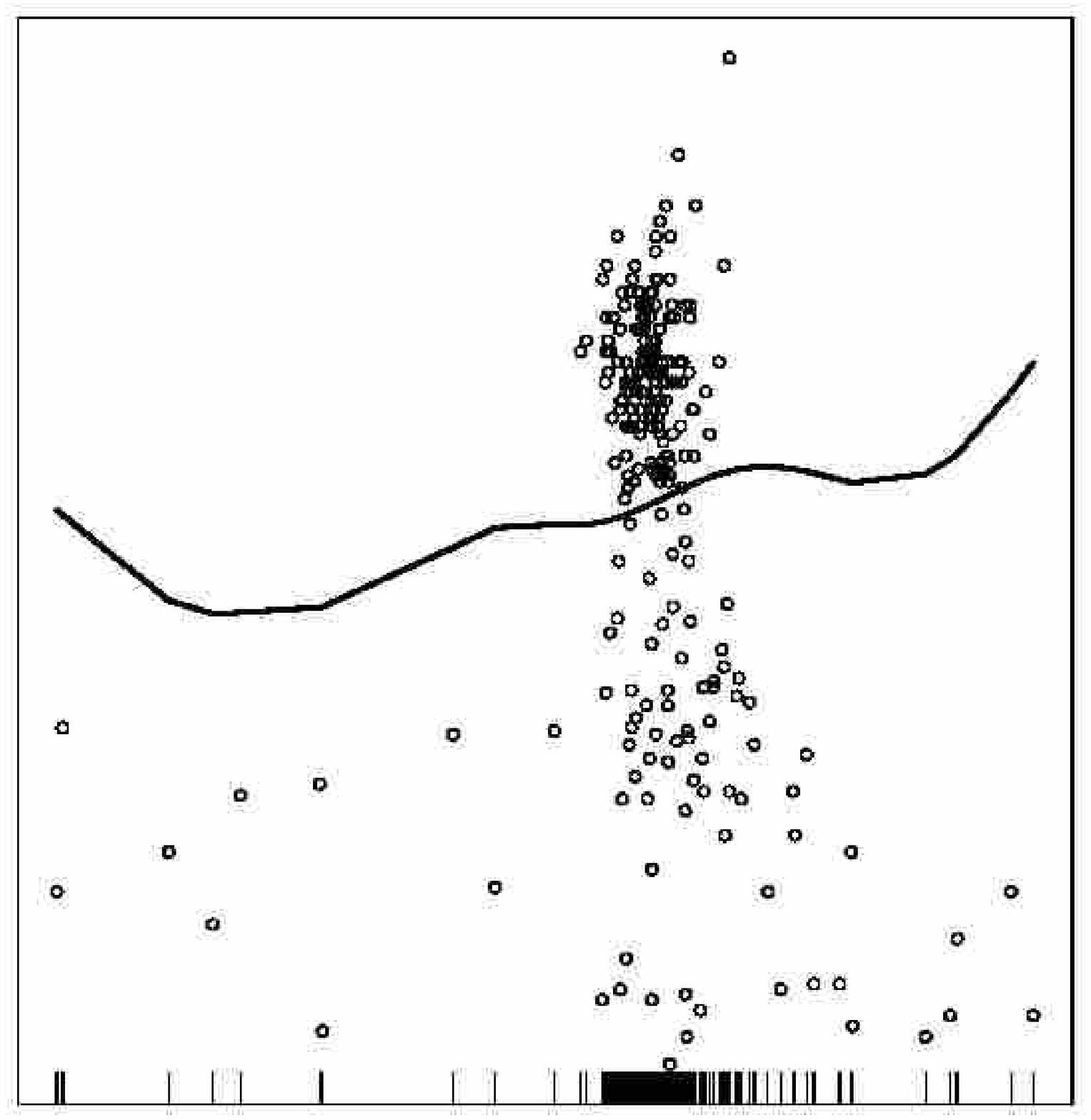}&
\includegraphics[width=\plotwidth]{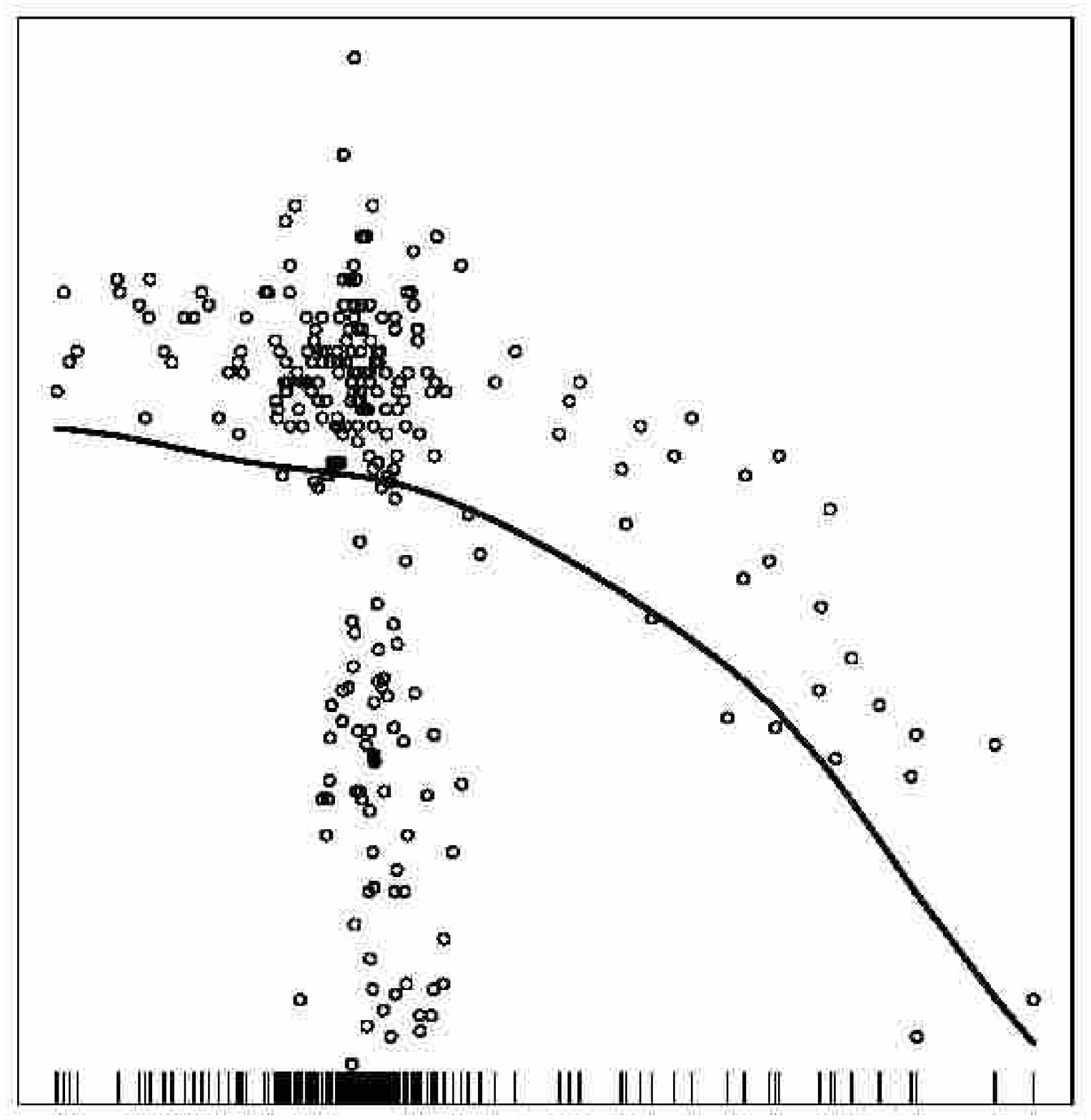}\\[10pt]
\includegraphics[width=\plotwidth]{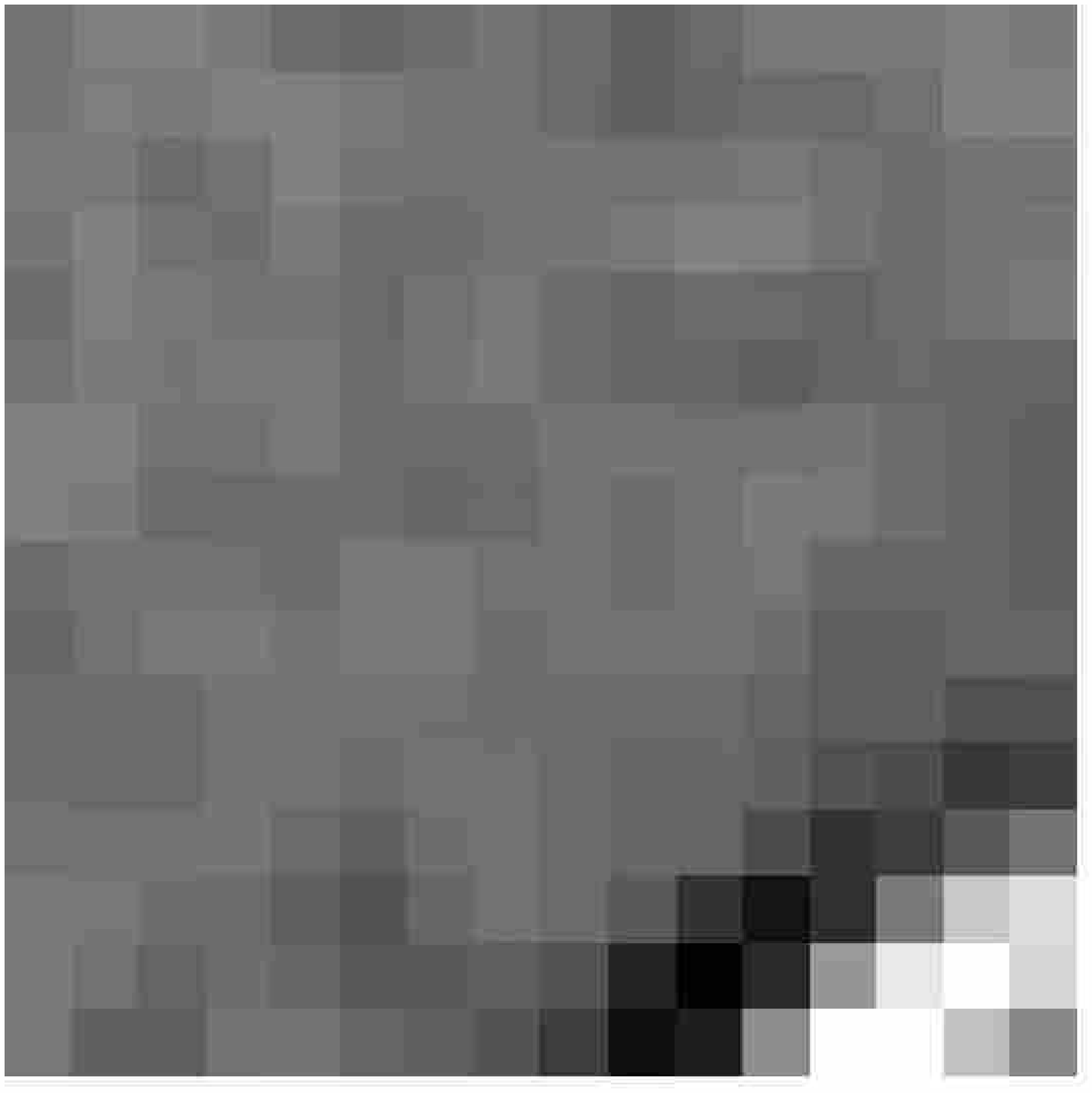}&&&\\
\includegraphics[width=\plotwidth]{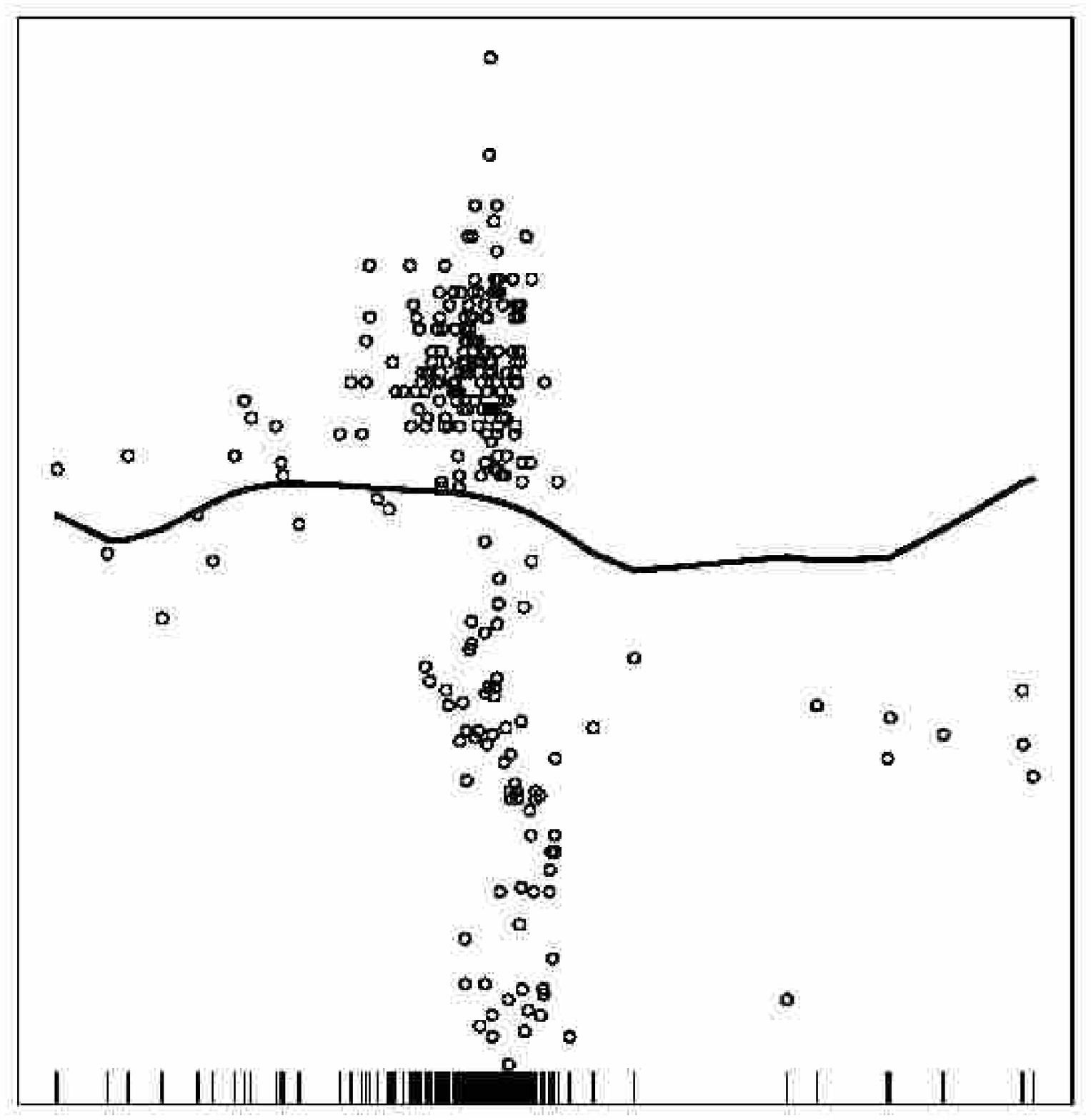}&&&\\[70pt]
\end{tabular}
\end{tabular}
\end{tabular}
\end{center}
\caption{Comparison of sparse reconstruction using the 
lasso (left) and SpAM (right).}
\label{fig:compare5}
\end{figure*}

\begin{figure*}
\begin{center}
\begin{tabular}{cc}
%\multicolumn{2}{c}{\small Original image} \\
%\multicolumn{2}{c}{\includegraphics[width=.40\textwidth]{figs/image1-lambda0p17-h0p05-original.ps}} 
%\\[20pt]
\small Lasso reconstruction & \small SpAM reconstruction \\
\hskip-5pt
\includegraphics[width=.40\textwidth]{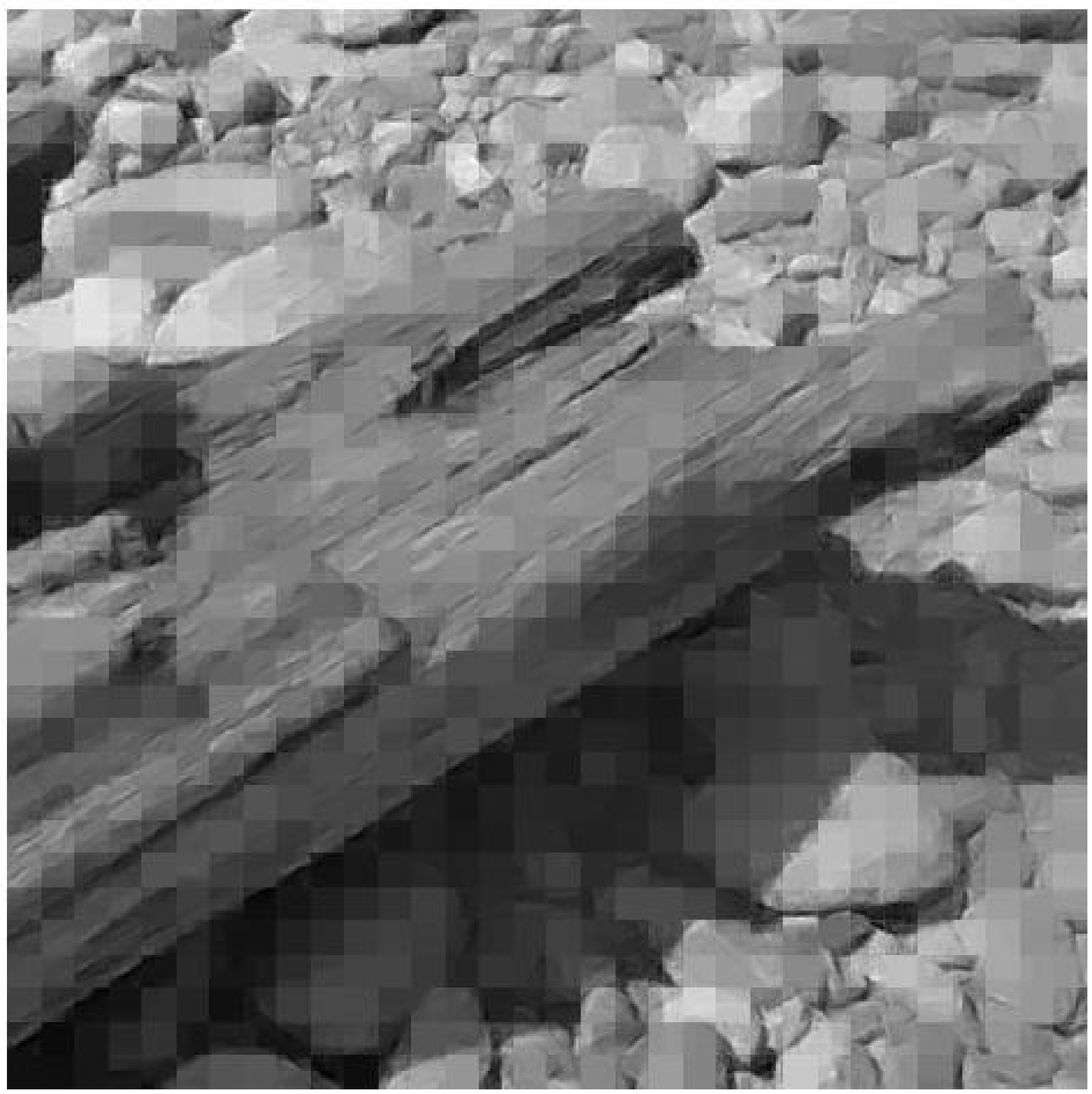}&
\includegraphics[width=.40\textwidth]{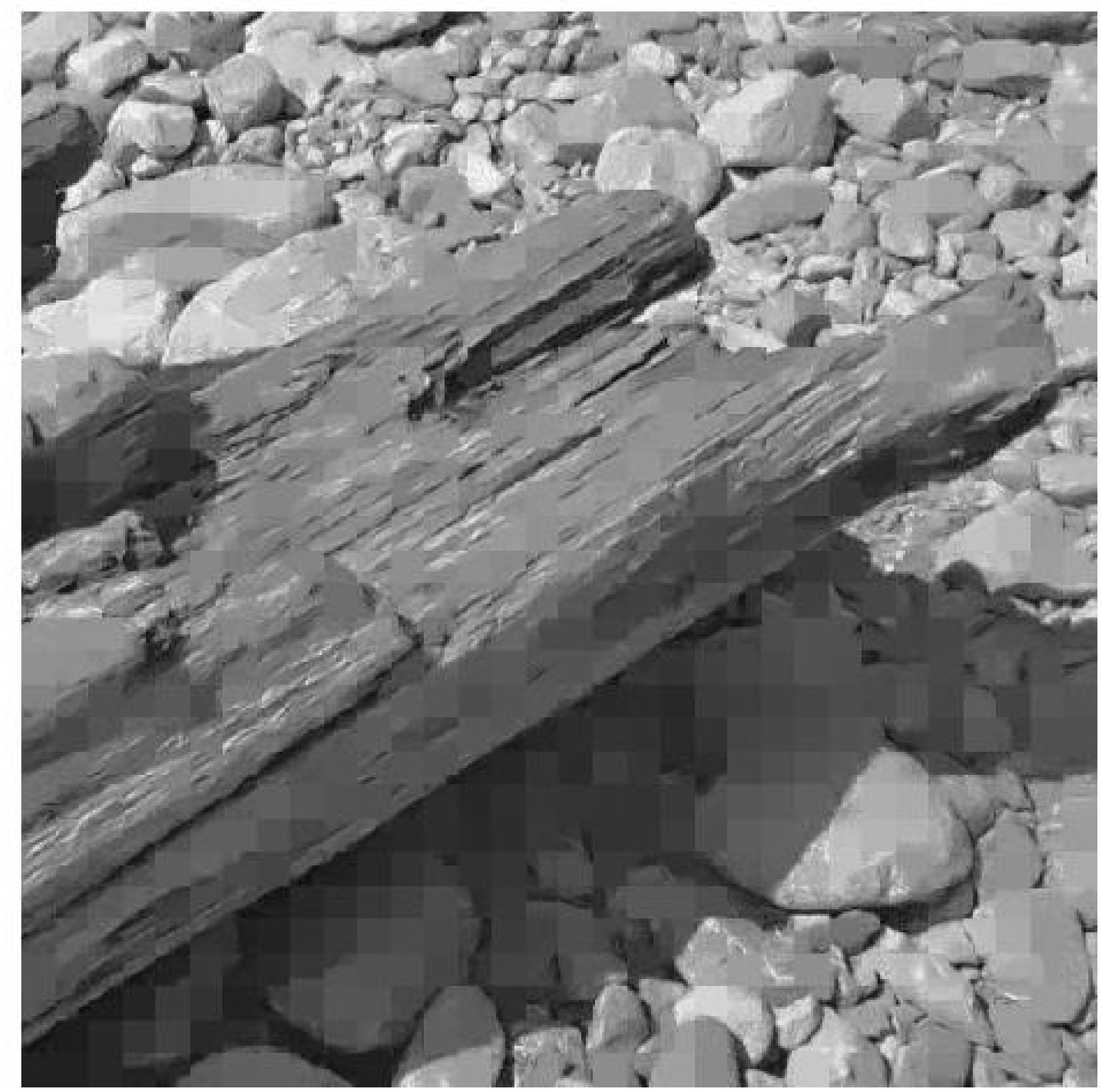}
\\[10pt]
{\small Codewords/patch 8.14, RSS 0.1894} & 
%%{\small Codewords/patch 8.14, RSS 0.1170} & 
{\small Codewords/patch 8.08, RSS 0.0913} \\
\end{tabular}
\end{center}
\caption{Image reconstruction using the lasso (left) and SpAM
  (right).  The regularization parameters were set so that the number
  of codewords used in each reconstruction was approximately equal.
  To achieve the same residual sum of squares (RSS), the
  linear model requires an average of more than 26 codewords per patch.}
\end{figure*}

Figures \ref{fig:compare2} and \ref{fig:compare5} illustrate the
reconstruction of different image patches using the sparse linear
model compared with the sparse additive model.  The codewords $X_j$
are those obtained using the Olshausen-Field procedure; these 
become the design points in the regression estimators.  Thus,
a codeword for a $16\times 16$ patch corresponds to a vector 
$X_j$ of dimension $256$, with each $X_{ij}$ the gray level for a
particular pixel.  

It can be seen how the functional fit achieves a more accurate
approximation to the image patch, with fewer codewords.  Each set of
plots shows the original image patch, the reconstruction, and the
codewords that were used in the reconstruction, together with their
marginal fits to the data.  For instance, in Figure~\ref{fig:compare2}
it is seen that the sparse linear model uses eight codewords, with a
residual sum of squares (RSS) of 0.0561, while the sparse additive
model uses seven codewords to achieve a residual sum of squares of
0.0206.  It can also be seen that the linear and nonlinear fits use
different sets of codewords.  Local linear smoothing was used with a
Gaussian kernel having fixed bandwidth $h=0.05$ for all patches and
all codewords.  

These results are obtained using the set of codewords obtained under
the sparse linear model.  The codewords can also be learned using the
sparse additive model; this will be reported in a separate paper.

\section{Theoretical Properties}

\subsection{Sparsistency}
\label{sec:sparsistency}

In the case of linear regression, with $f_j(X_j) = \beta_j^{*T} X_j$, several authors have shown that,
under certain conditions on $n$, $p$, the number of relevant
variables $s = |\supp(\beta^*)|$, and the
design matrix $X$, the lasso recovers the sparsity pattern
asymptotically; that is, the lasso estimator $\hat\beta_n$ is \textit{sparsistent}:
\begin{equation}
\P\left(\supp(\beta^*) = \supp(\hat\beta_n)\right) \rightarrow 1.
\end{equation}
Here, $\supp(\beta) = \left\{j:\ \beta_j \neq \bm{0} \right\}$.
References include \cite{Wai06},
\cite{MB}, \cite{Zou} and \cite{ZY07}.
We show a similar result for sparse additive models under orthogonal function regression.

In terms of an orthogonal basis $\psi$, we can write
\begin{equation}
\label{eq:spmmdl}
Y_i = \sum_{j=1}^p \sum_{k=1}^\infty \beta_{jk}^*\psi_{jk}(X_{ij}) + \epsilon_i.
\end{equation}

To simplify notation, let $\beta_{j}$ be the $d_n$ dimensional vector
$\{\beta_{jk},\,k = 1,\hdots,d_n\}$ and let $\Psi_{j}$ be the $n \times
d_{n}$ matrix $\Psi_{j}[i,k] = \psi_{jk}(X_{ij})$.
If $A\subset \{1,\ldots,p\}$, we denote by $\Psi_A$ the
$ n\times d|A|$ matrix where for each $j\in A$, $\Psi_j$ appears
as a submatrix in the natural way.

We now analyze the sparse backfitting Algorithm (\ref{fig:backfitting:algo})
assuming an orthogonal series smoother is used to estimate the conditional
expectation in its Step (2). As noted
earlier, an orthogonal series smoother for a predictor $X_j$ is the least squares
projection onto a truncated set of basis functions
$\{\psi_{j1},\ldots, \psi_{jd}\}$. Combined with the soft-thresholding
step, the update for $f_j$ in Algorithm (\ref{fig:backfitting:algo}) can
thus be seen to solve the following problem,
\begin{align*}
\min_{\beta} \frac{1}{2n}\|R_{j} - \Psi_{j}\beta_{j}\|_{2}^2
+ \lambda_n \sqrt{\frac{1}{n} \beta_j^T \Psi_j^T \Psi_j \beta_j} \\
\end{align*}
where $\|v\|_{2}^{2}$ denotes $\sum_{i=1}^{n}v_{i}^{2}$ and $R_j = Y - \sum_{l \neq j} \Psi_{l} \beta_{l}$ is
the residual for $f_j$. The sparse backfitting algorithm can then be seen to solve,
\begin{eqnarray}
\label{eq::objective}
\nonumber
\lefteqn{\min_{\beta} \left\{R_n(\beta) + \lambda_n \Omega(\beta) \right\}}
&& \\
&=& \min_{\beta} \frac{1}{2n}\left\|Y - \sum_{j=1}^{p}\Psi_{j}\beta_{j}\right\|_{2}^{2} +
\lambda \sum_{j=1}^p \sqrt{\frac{1}{n} \beta_j^T \Psi_j^T \Psi_j \beta_j} \\
&=&\min_{\beta} \frac{1}{2n}\left\|Y - \sum_{j=1}^{p}\Psi_{j}\beta_{j}\right\|_{2}^2 +
\lambda \sum_{j=1}^{p} \left\|\frac{1}{\sqrt{n}}\Psi_j \beta_j\right\|_2.
\label{eq:orthospam}
\end{eqnarray}
where $R_n$ denotes the squared error term and $\Omega$ denotes
the regularization term, and each $\beta_j$ is a $d_n$-dimensional vector. 
Let $S$ denote the true set of variables
$\{j: f_j\neq 0\}$, with $s = |S|$, and let $S^c$ denote its complement.  
Let $\hat S_{n} = \{j: \hat \beta_j \neq 0\}$ denote
the estimated set of variables from the minimizer $\hat\beta_n$ of 
\eqref{eq:orthospam}.
For the results in this section, we will treat the covariates as fixed.

\begin{theorem} 
\label{thm:sparsistence}
Suppose that the following conditions hold on the design matrix $X$ in the
orthogonal basis $\psi$:
\begin{gather}
\label{eq:design1}
\Lambda_{\max}\left(\frac{1}{n} \Psi_S^T \Psi_S\right) \;\leq\;
  C_{\max}  \;<\; \infty   \\
\label{eq:design2}
\Lambda_{\min}\left(\frac{1}{n} \Psi_S^T \Psi_S\right) \;\geq\;
  C_{\min}  \;>\; 0 \\
\label{eq:design3}
\max_{j\in S^c} \; \left\|  \left( \textstyle \frac{1}{n} \Psi_{j}^T \Psi_S\right)
  \left(\textstyle \frac{1}{n}
    \Psi_{S}^T \Psi_S\right)^{-1} \right\| \;\leq\;
\sqrt{\frac{C_{\min}}{C_{\max}}} \,\frac{1-\delta}{\sqrt{s}}, \;\;\text{for some $0 < \delta\leq 1$}.
\end{gather}
Assume that the truncation dimension $d_n$
satisfies $d_n \rightarrow \infty$ and $d_n = o(n)$.  Furthermore, suppose the following conditions,
which relate the regularization parameter $\lambda_n$ to the
design parameters $n, p$, the number of relevant variables $s$, and the truncation
size $d_n$:
\begin{gather}
\label{eq:conda}
\frac{s}{d_n \lambda_n} \;\longrightarrow\; 0 \\[5pt]
\label{eq:condb}
\frac{d_n \log \left(d_n (p-s)\right)}{n\lambda_n^2} \;\longrightarrow\; 0 \\[5pt]
\label{eq:condc}
\frac{1}{\rho_n^*}\left(
\sqrt{\frac{\log(sd_n)}{n}}
+ \frac{s^{3/2}}{d_n} + \lambda_n \sqrt{s d_n} 
\right) \;\longrightarrow\; 0
\end{gather}
where $\rho_n^* = \min_{j\in S} \|\beta_j^*\|_\infty$.
Then SpAM is sparsistent:  
$\P\left(\hat S_n = S\right) \rightarrow 1$.
\end{theorem}

The proof is given in an appendix.
Note that condition \eqref{eq:design3} implies that
\begin{eqnarray}
\left\|\Psi_{S^c}^T \Psi_S \left(\Psi_S^T \Psi_S\right)^{-1}\right\|_\infty &=&
\max_{j\in S^c} \left\|\Psi_{j}^T \Psi_S \left(\Psi_S^T \Psi_S\right)^{-1}\right\|_\infty \\
&\leq& \sqrt{\frac{C_{\min} \, d_n}{C_{\max}}} \, (1-\delta)
\end{eqnarray}
since $\frac{1}{\sqrt{n}} \|A\|_\infty \leq \|A\| \leq \sqrt{m} \, \|A\|_\infty$ for an
$m\times n$ matrix $A$.  This relates the condition to previous $\infty$-norm
incoherence conditions that have been used for sparsistency in the linear case
\citep{Wai06}.

For $\nu=2$ we take $d_{n} = n^{1/5}$,
which achieves the minimax error rate in the one-dimensional
case. The theorem under this design setting, with the simplifying
assumption that $s = O(1)$, gives the following

\begin{corollary}
Suppose that $s = O(1)$, and $d_n = O(n^{1/5})$.  Assume the design
conditions \eqref{eq:design1}, \eqref{eq:design2} and \eqref{eq:design3}.
Suppose the penalty $\lambda_n$ is chosen to satisfy
\begin{equation}
n^{1/5} \lambda_n \;\rightarrow\; \infty,\ \ \ 
\frac{\log(np)}{n^{4/5} \lambda_n^2} \;\rightarrow\; 0,\ \ \ 
\frac{\lambda_n \, n^{1/10}}{\rho^*_n} \;\rightarrow\; 0.
\end{equation}
Then $\mathbb{P}\left(\hat{S}_{n} = S\right) \rightarrow 1$.
\end{corollary}

For example, under these conditions and $\rho_n^* \asymp 1$, the dimension $p_n$
can be taken as large as $e^{n^{3/5}}$; a suitable choice for the regularization parameter 
in this case would be $\lambda_n = {C n^{-\frac{1}{10} + \delta} \log n}$ for some
$\delta > 0$.

\subsection{Persistence}
\label{sec:persistence}

The previous assumptions are very strong.
They can be weakened at the expense of getting
weaker results.
In particular, in the section we do not assume that
the true regression function is additive.
We use arguments like those
in \cite{JN:00} and \cite{GR:04}
in the context of linear models.
In this section we treat $X$ as random
and we use triangular array asymptotics, that is,
the joint distribution for the data can change with $n$.
Let $(X,Y)$ denote a new pair
(independent of the observed data) and define the predictive
risk when predicting $Y$ with $v(X)$ by
\begin{equation}
R(v) = \mathbb{E}(Y - v(X))^2.
\end{equation}
When $v(x)=\sum_j \beta_j g_j(x_j)$ we also write the risk as
$R(\beta,g)$
where $\beta = (\beta_1,\ldots, \beta_p)$ and
$g = (g_1,\ldots, g_p)$.
Following \cite{GR:04} we say
that an estimator $\hat{m}_n$ is persistent (risk consistent)
relative to
a class of functions ${\cal M}_n$,
if
\begin{equation}
R(\hat{m}_n) - R(m_n^*) \stackrel{P}{\to} 0
\end{equation}
where
\begin{equation}
m_n^* = \argmin_{v\in {\cal M}_n} R(v)
\end{equation}
is the predictive oracle.
\cite{GR:04} show that the lasso 
is persistent
for
${\cal M}_n = \{ \ell(x) = x^T \beta :\ \left\|\beta\right\|_1 \leq L_n\}$
and $L_n = o( (n/\log n)^{1/4})$.
Note that $m_n^*$ is the best linear approximation (in prediction risk)
in ${\cal M}_n$ but the true regression function is not assumed to be linear.
Here we show a similar result for SpAM.

In this section, we assume that the SpAM estimator $\hat{m}_n$
is chosen to minimize
\begin{equation}
\frac{1}{n}\sum_{i=1}^n(Y_i - \sum_j \beta_j g_j(X_{ij}))^2
\end{equation}
subject to
$\left\|\beta\right\|_1 \leq L_n$ and
$g_j\in {\cal T}_j$.
We make no assumptions about the design matrix.
Let ${\cal M}_n \equiv {\cal M}_n(L_n)$ be defined by
\begin{equation}
{\cal M}_n = \Biggl\{ m:\ 
m(x) = \sum_{j=1}^{p_n}\beta_j g_j(x_j):\ 
\E(g_j) = 0,\ \E(g_j^2) = 1,\ \sum_j | \beta_j| \leq L_n\Biggr\}
\end{equation}
and let
$m_n^* = \argmin_{v\in {\cal M}_n} R(v)$.

\begin{theorem}\label{thm::persist}
Suppose that $p_n \leq e^{n^\xi}$ for some
$\xi < 1$.
Then,
\begin{equation}
R(\hat{m}_n) - R(m_n^*) = O_P\left(\frac{L_n^2}{n^{(1-\xi)/2}}\right)
\end{equation}
and hence , if $L_n = o(n^{(1-\xi)/4})$ then 
SpAM is persistent.
\end{theorem}

\section{Discussion}
\label{sec:discussion}

The results presented here show how many of the recently established
theoretical properties of $\ell_1$ regularization for linear models 
extend to sparse additive models.  The sparse backfitting algorithm we
have derived is attractive because it decouples smoothing and
sparsity, and can be used with any nonparametric smoother.  It thus
inherits the nice properties of the original backfitting procedure.
However, our theoretical analyses have made use of a particular
form of smoothing, using a truncated orthogonal basis.  An 
important problem is thus to extend the theory to cover
more general classes of smoothing operators.

An additional direction for future work is to develop procedures
for automatic bandwidth selection in each dimension.  We have
used plug-in bandwidths and truncation dimensions $d_n$
in our experiments and theory. It is of particular interest 
to develop procedures that are adaptive to different levels
of smoothness in different dimensions.

Finally, we note that while we have considered basic additive models
that allow functions of individual variables, it is natural to consider
interactions, as in the functional ANOVA model. One challenge
is to formulate suitable incoherence conditions on the functions that enable
regularization based procedures or greedy algorithms to recover the
correct interaction graph.  In the parametric setting, one
result in this direction is \cite{Wainwright:07}.

\section{Proofs}
\label{sec:proofs}

\begin{proof}[{\it Proof of} Theorem \ref{thm::backfit}]
Consider the minimization of the Lagrangian
\begin{equation}
\min_{\{f_j \in \H_j\}}\L(f, \lambda) \equiv \frac{1}{2} \E\left(Y-\textstyle \sum_{j=1}^p
  f_j(X_j)\right)^2  + \lambda \sum_{j=1}^p \sqrt{\E(f_j^2(X_j))}
\end{equation}
with respect to $f_j \in \H_j$, holding
the other components $\{f_{k},\,k \neq j\}$ fixed. The stationary condition is 
obtained by setting to zero the Fr\'echet directional derivative with respect to $f_j$,
denoted $\partial \L(f,\lambda;\eta_j)$,  for all feasible directions 
$\eta_j(X_j) \in \cH_j\, (\E(\eta_j) = 0,\,\E(\eta_{j}^2) < \infty)$.
This leads to the condition
\begin{align}
\partial \L(f,\lambda;\eta) = \frac{1}{2} \E\left[(f_j - R_j + \lambda v_{j})~\eta_j\right] = 0
\end{align}
where $R_j = Y - \sum_{k \neq j} f_{k}$ is the residual for $f_j$,
and $v_j$ is an element of the subgradient $\partial \sqrt{\E(f_{j}^2)}$,
satisfying $v_j \in \H_j$; $v_{j} = f_{j}/\sqrt{\E(f_j)^2}$ if $\E(f_{j}^2) \neq 0$ and belonging
to the set $\{u_j \in \H_j |~\E(u_{j}^{2}) \le 1\}$ otherwise.

Using iterated expectations, the above condition can be rewritten as
\begin{align}
\E\left[(f_{j} + \lambda v_{j} - \E(R_{j}|X_j))~\eta_j \right] = 0.
\end{align}
%Let $d_j(X_j) = \delta(X_j - x_j) - p_j(x_j)$ be a translated delta function,
%where $p_j(x_j)$ is the marginal sample density of $X_j$ at $x_j$.
%Clearly $d_j \in \H_j$ since
%\begin{align}
%\E(\eta_j(X_j)) &= \E(\delta(X_j - x_j) - p_j(x_j)\\
%\nn            &= p_j(x_j) - p_j(x_j) = 0\\[10pt]
%\E(\eta_{j}^{2}(X_j)) &= \E(\delta(X_j - x_j) - p_j(x_j))^2\\
%\nn            &= p_j(x_j)(1 - p_j(x_j)) < \infty
%\end{align}
%For the direction $\eta_j = d_j$, the condition 
%$\partial \L(f,\lambda;\eta) =0$ becomes
%\begin{eqnarray}
%\lefteqn{\left(f_j(x_j) - \E(R_j|X_j = x_j) + \lambda
%    v_j(x_j)\right) p_j(x_j) \;=\;}&& \\
%&& \E\left(f_j(x_j) - \E(R_j|X_j = x_j) + \lambda v_j(x_j)\right)p_j(x_j) .
%\end{eqnarray}
But since $f_j - \E(R_j|X_j) + \lambda v_j \in \H_j$, we can compute the
derivative in the direction $\eta_j = f_j - \E(R_j|X_j) + \lambda v_j \in \H_j$, implying that
\begin{align}
\E\left[\left(f_j(x_j) - \E(R_j|X_j = x_j) + \lambda v_j(x_j)\right)^2\right] = 0;
\end{align}
that is,
\begin{align}\label{eq:popstatcond}
f_{j} + \lambda v_{j} = \E(R_{j}|X_j)\quad \text{a.e.}
\end{align}
Denote the conditional expectation $\E(R_{j}|X_{j})$---also the 
projection of the residual $R_j$ onto $\H_j$---by $P_{j}$. 
Now if $\E(f_{j}^2) \neq 0$, then $v_{j} = \frac{f_{j}}{\sqrt{\E(f_j^2)}}$, 
which from condition~\eqref{eq:popstatcond} implies
\begin{eqnarray}
\sqrt{\E(P_j^2)} &=& \sqrt{\E\left[\left(f_{j} +\lambda {f_{j}}/{\sqrt{\E(f_j^2)}}\right)^2\right]}\\
                &=& \left(1 + \frac{\lambda}{\sqrt{\E(f_j^2)}}\right)\sqrt{\E(f_j^2)}\\
\label{eq:normproj}             &=& \sqrt{\E(f_j^2)} + \lambda\\
                &\ge& \lambda.
\end{eqnarray}
If $\E(f_{j}^2) = 0$, then $f_j = 0$ a.e., and $\sqrt{\E(v_{j}^2)} \le 1$.
\eqref{eq:popstatcond} implies that
\begin{align}
\sqrt{\E(P_j^2)} = \lambda \sqrt{\E(v_{j}^2)} \le \lambda .
\end{align}
We thus obtain the equivalence
\begin{align}\label{eq:subgradequiv}
\sqrt{\E(P_j^2)} \le \lambda \;\iff\; f_{j} = 0 \quad \text{a.s.}
\end{align}
Rewriting equation~\eqref{eq:popstatcond} in light of \eqref{eq:subgradequiv}, we obtain
\begin{eqnarray*}
\left(1 + \frac{\lambda}{\sqrt{\E(f_{j}^2)}}\right)f_{j} =
P_{j} &&
\text{if $\sqrt{\E(P_{j}^2)} > \lambda$} \\
f_j = 0 && \text{otherwise}.
\end{eqnarray*}
Using \eqref{eq:normproj}, we thus arrive at the 
soft thresholding update for $f_j$:
\begin{equation}
\label{eq:softy}
f_j = \left[1 - \frac{\lambda}{\sqrt{\E(\Proj_j^2)}}\right]_+ \Proj_j
\end{equation}
where $[\cdot]_+$ denotes the positive part and
$\Proj_j = \E[\Res_j\given X_j]$.
\end{proof}

\begin{proof}[{\it Proof of} Theorem \ref{thm:sparsistence}]

A vector $\hat{\beta} \in \mathbb{R}^{d_{n}p}$ is an optimum of the objective function in
\eqref{eq::objective} if and only if there exists a subgradient $\hat{g} \in \partial
\Omega(\hat{\beta})$, such that
\begin{equation}
\label{eq:dbstat}
\frac{1}{n}\Psi^{\mytop}\left(\sum_{j}\Psi_{j}\hat\beta_{j} - Y\right) + \lambda_{n}\hat{g} = \bm{0}.
\end{equation}
The subdifferential $\partial \Omega(\beta)$ is the set of vectors $g\in
\reals^{pd_n}$ satisfying
\begin{eqnarray*}
g_{j} \;=\;
\frac{\frac{1}{n} \Psi_j^T \Psi_j \beta_j}{\sqrt{\frac{1}{n} \beta_j^T
      \Psi_j^T \Psi_j \beta_j}} &&
\text{if $\beta_{j} \neq \bm{0}$} \\[10pt]
g_j^T \left(\frac{1}{n} \Psi_j^T \Psi_j\right)^{-1} g_j \;\le\; 1 &&\text{if $\beta_{j} = \bm{0}$}.
\end{eqnarray*}

Our argument closely follows the approach of \cite{Wai06} in the linear case.  In particular, we
proceed by a ``witness'' proof technique, to show the existence of a
coefficient-subgradient pair $(\beta, g)$ for which $\supp(\beta) = \supp(\beta^*)$.  To
do so, we first set $\hat{\beta}_{\Sc} = \bm{0}$ and $\hat{g}_{S} = \partial
\Omega(\beta^*)_S$, and we then obtain $\hat{\beta}_{S}$ and $\hat{g}_{\Sc}$ from the
stationary conditions in~\eqref{eq:dbstat}.  By showing that, with high probability,
\begin{eqnarray*}
\hat{\beta}_{j} &\neq& \bm{0} \;\; \text{for $j \in S$}\\[5pt]
g_j^T \left(\frac{1}{n} \Psi_j^T \Psi_j\right)^{-1} g_j &\le& 1 \;\;\text{for $j \in \Sc$},
\end{eqnarray*}
this demonstrates that with high probability there exists an optimal
solution to the optimization problem in~\eqref{eq::objective} that has
the same sparsity pattern as the true model.

Setting $\hat{\beta}_{\Sc} = \bm{0}$ and 
\begin{eqnarray}
\hat{g}_{j} &=& \frac{\frac{1}{n} \Psi_j^T \Psi_j \beta_j^*}{\sqrt{\frac{1}{n} \beta_j^{*T}
      \Psi_j^T \Psi_j \beta_j^*}} 
\end{eqnarray}
for $j\in S$, the stationary condition for
$\beta_{S}$ is
\begin{equation}
\frac{1}{n}\Psi_{S}^{\mytop}\left(\Psi_{S}\beta_{S} - Y\right) +
\lambda_n \hat g_{S} = \bm{0}.
\end{equation}
Let $V = Y - \Psi_{S}\beta^{*}_{S} - W$ denote the error due to
finite truncation of the orthogonal basis, where $W =
(\epsilon_1,\ldots, \epsilon_n)^T$.  Then the stationary 
condition can be written as
\begin{eqnarray}
\frac{1}{n} \Psi_S^T \Psi_S \left(\beta_S - \beta_S^*\right)  - 
\frac{1}{n} \Psi_S^T W - 
\frac{1}{n} \Psi_S^T V +
\lambda_n \hat g_S = \bm{0}
\end{eqnarray}
or
\begin{eqnarray}
\label{eq:bseq}
\beta_S - \beta_S^* = \left(\frac{1}{n} \Psi_S^T \Psi_S\right)^{-1}
\left(\frac{1}{n} \Psi_S^T W + \frac{1}{n} \Psi_S^T V - \lambda_n \hat g_S\right)
\end{eqnarray}
assuming that $\frac{1}{n} \Psi_S^T \Psi_S$ is nonsingular.
Recalling our definition
\begin{equation}
\rho_n^* = \min_{j \in S}\|\beta_j^*\|_\infty > 0.
\end{equation}
it suffices to show that
\begin{equation}
\|\beta_{S} - \beta^{*}_{S}\|_{\infty} < \frac{\rho_n^*}{2}
\end{equation}
in order to ensure that $\supp(\beta_{S}^*) = \supp(\beta_S) =
\left\{j\,:\, \|\beta_j\|_\infty \neq 0\right\}$.

Using $\Sigma_{SS} =
\frac{1}{n}\left(\Psi_{S}^{\mytop}\Psi_{S}\right)$ to simplify
notation, we have the $\ell_\infty$ bound
\begin{eqnarray}
\label{eq:linftybnd} 
\lefteqn{\|\beta_{S} - \beta^{*}_{S}\|_{\infty} \;\leq\;} && \\
&& \left\|\Sigma_{SS}^{-1}\left(\onen\Psi_{S}^{\mytop}W\right)\right\|_{\infty} +
\left\|\Sigma_{SS}^{-1}\left(\onen\Psi_{S}^{\mytop}V\right)\right\|_{\infty} +
\lambda_n \left\|\Sigma_{SS}^{-1}\hat g_{S}\right\|_{\infty}.
\end{eqnarray}
We now proceed to bound the quantities above.  First note that for
$j\in S$, 
\begin{eqnarray}
1 &=& g_j^T \left(\onen \Psi^T_j \Psi_j \right)^{-1} g_j \\
  &\geq& \frac{1}{C_{\max}} \|g_j\|^2
\end{eqnarray}
and thus $\|g_j\| \leq \sqrt{C_{\max}}$.  Therefore, 
since 
\begin{equation}
\|g_{S}\|_{\infty} = \max_{j \in S}\|g_{j}\|_{\infty}
\le \max_{j \in S}\|g_{j}\|_{2} \leq \sqrt{C_{\max}} 
\end{equation}
we have that
\begin{equation}
\left\|\Sigma_{SS}^{-1} \hat g_S\right\|_{\infty} \le
\sqrt{C_{\max}} \left\|\Sigma_{SS}^{-1}\right\|_{\infty} .
\end{equation}

Now, to bound $\left\|\onen \Psi_S^T V\right\|_\infty$, first note
that, as we are working over the Sobolev spaces $\mathcal S_j$ of order two,
\begin{eqnarray}
|V_{i}| &=& \left|\sum_{j\in S}\sum_{k=d_n+1}^\infty
  \beta^*_{jk}\Psi_{jk}(X_{ij}) \right| 
\leq  B \sum_{j\in S}\sum_{k=d_n+1}^\infty \left|\beta_{jk}^*\right| \\
&=& B \sum_{j\in S}\sum_{k=d_n+1}^\infty \frac{\left|\beta_{jk}^*\right|k^2 }{k^2}
\leq  B \sum_{j\in S} \sqrt{\sum_{k=d_n+1}^\infty \beta_{jk}^{*2} k^4}
\sqrt{\sum_{k=d_n+1}^\infty \frac{1}{k^4}}\\
& \leq & s B C \sqrt{\sum_{k=d_n+1}^\infty \frac{1}{k^4}} \leq  \frac{s B'}{d_n^{3/2}}
\end{eqnarray}
for some constant $B' > 0$.
Therefore,
\begin{eqnarray}
\label{eq:Vinftybnd}\|V\|_{\infty} &\le& \frac{B' s}{d_n^{3/2}},
\end{eqnarray}
and also
\begin{equation}
\left|\frac{1}{n}\Psi_{jk}^{\mytop}V\right| \le
\left|\frac{1}{n}\sum_{i}\Psi_{jk}(X_{ij})\right|\|V\|_{\infty}
\le \frac{D s}{d_n^{3/2}}
\end{equation}
where $D$ denotes a generic constant.  Together then, we have that
\begin{eqnarray}
\nonumber
\lefteqn{\|\beta_{S} - \beta^{*}_{S}\|_{\infty} \;\leq\;} && \\
&& \left\|\Sigma_{SS}^{-1}\left(\onen\Psi_{S}^{\mytop}W\right)\right\|_{\infty} +
\left\|\Sigma_{SS}^{-1}\right\|_\infty \left( \frac{Ds}{d_n^{3/2}} +
  \lambda_n \sqrt{C_{\max}}\right)
\label{eq:linftybnd2} 
\end{eqnarray}

Finally, consider $Z := \Sigma_{SS}^{-1}\left(\frac{1}{n}\Psi_{S}^{\mytop}W\right)$.
Note that $W \sim N(0,\sigma^{2}I)$, so that $Z$ is Gaussian as well,
with mean zero. Consider its $l$-th component, $Z_{l} = e_{l}^{\mytop}Z$.
Then $\E[Z_{l}] = 0$, and 
\begin{equation}
\Var(Z_{l}) =
\frac{\sigma^{2}}{n}e_{l}^{\mytop}\Sigma_{SS}^{-1}e_{l} \leq
\frac{\sigma^{2}}{C_{\min}n}.
\end{equation}
By Gaussian comparison results
\citep{Ledoux:Talagrand:91}, we have then that
\begin{equation}
\E\left[\|Z\|_{\infty}\right] \;\le\; 3 \sqrt{\log(sd_n) \left\|\Var(Z)\right\|_\infty }
\;\le\;  3\sigma \sqrt{\frac{\log(sd_n)}{n C_{\min}}}.
\end{equation}
An application of Markov's inequality then gives that
\begin{eqnarray}
\nonumber
\lefteqn{\P\left(\|\beta_{S} - \beta^{*}_{S}\|_{\infty} >
    \frac{\rho_n^*}{2}\right) \;\le\;} \\
&&
\P\left(\|Z\|_{\infty} + 
\left\|\Sigma_{SS}^{-1}\right\|_\infty \left( {Ds}{d_n^{-3/2}} +
  \lambda_n \sqrt{C_{\max}} \right) > \frac{\rho_n^*}{2} \right) \\
&\le& \frac{2}{\rho_n^*}\left\{\E\left[\|Z\|_{\infty}\right]
+ \left\|\Sigma_{SS}^{-1}\right\|_\infty \left( {Ds}{d_n^{-3/2}} +
  \lambda_n \sqrt{C_{\max}} \right) \right\} \\
&\le& \frac{2}{\rho_n^*}\left\{
3 \sigma \sqrt{\frac{\log(sd_n)}{n C_{\min}}} 
+ \left\|\Sigma_{SS}^{-1}\right\|_\infty \left( \frac{Ds}{d_n^{3/2}} +
  \lambda_n \sqrt{C_{\max}} \right) \right\} 
\end{eqnarray}
which converges to zero under the condition that
\begin{equation}
\label{eq:thiscond}
\frac{1}{\rho_n^*}\left\{
\sqrt{\frac{\log(sd_n)}{n}}
+ \left\|\left(\onen \Psi_S^T \Psi_ S\right)^{-1}\right\|_\infty \left( \frac{s}{d_n^{3/2}} +
  \lambda_n\right) \right\} \longrightarrow 0.
\end{equation}
Now, since 
\begin{eqnarray}
\left\|\left(\onen \Psi_S^T \Psi_ S\right)^{-1}\right\|_\infty &\leq& 
\frac{\sqrt{s d_n}}{C_{\min}}
\end{eqnarray}
condition \eqref{eq:thiscond} holds in case
\begin{equation}
\frac{1}{\rho_n^*}\left(
\sqrt{\frac{\log(sd_n)}{n}}
+ \frac{s^{3/2}}{d_n} + \lambda_n \sqrt{s d_n} 
\right) \;\longrightarrow\; 0,
\end{equation}
which is condition \eqref{eq:condc} in the statement of the theorem.

We now analyze $\hat{g}_{\Sc}$. Recall that we have set 
$\hat \beta_{S^c} = \beta^*_{S^c} = \bm{0}$.  The stationary condition
for $j \in \Sc$ is thus given by
\begin{equation}
\frac{1}{n}\Psi_{j}^{\mytop}\left(\Psi_{S}\hat \beta_{S} -
  \Psi_{S}\beta^{*}_{S} - W - V\right) + \lambda_{n}\hat g_{j} = \bm{0}.
\end{equation}
Therefore,
\begin{eqnarray}
\nonumber
\hat g_{S^c} &=& \frac{1}{\lambda_n} \left\{
\frac{1}{n}\Psi_{S^c}^{\mytop}\Psi_{S} \left(\beta^*_{S} - \hat\beta_{S}\right)
    + \frac{1}{n}\Psi_{S^c}^{\mytop} \left(W + V\right) \right\} \\
&=& \frac{1}{\lambda_n} \left\{
\frac{1}{n}\Psi_{S^c}^{\mytop}\Psi_{S} 
\left(\frac{1}{n} \Psi_S^T \Psi_S\right)^{-1}
\left(\lambda_n \hat g_S - \frac{1}{n} \Psi_S^T W - \frac{1}{n}
  \Psi_S^T V\right)
\right. \\
&& \qquad 
    + \left. \mbox{\LARGE\strut} \frac{1}{n}\Psi_{S^c}^{\mytop} \left(W + V\right) \right\}
\nonumber
\\
&=& \frac{1}{\lambda_n} \left\{
\Sigma_{S^c S} \Sigma_{S S}^{-1} 
\left(\lambda_n \hat g_S - \frac{1}{n} \Psi_S^T W - \frac{1}{n} \Psi_S^T V\right)
    + \frac{1}{n}\Psi_{S^c}^{\mytop} \left(W + V\right) \right\}
\label{eq:g}
\end{eqnarray}
from equation \eqref{eq:bseq}.

We require that 
\begin{equation}
g_j^T \left(\frac{1}{n} \Psi_j^T \Psi_j\right)^{-1} g_j \;\le\; 1
\end{equation}
for all $j\in S^c$.  Since 
\begin{equation}
g_j^T \left(\frac{1}{n} \Psi_j^T \Psi_j\right)^{-1} g_j \;\le\; 
\frac{1}{C_{\min}} \|g_j\|^2
\end{equation}
it suffices to show that $\max_{j\in S^c} \|g_j\| \leq \sqrt{C_{\min}}$.

From \eqref{eq:g}, we see 
that $\hat g_j$ is Gaussian, with mean
\begin{eqnarray}
\mu_j = \E(\hat g_j) = \Sigma_{jS}\Sigma_{SS}^{-1} \left(\hat g_S - 
\frac{1}{\lambda_n} \left(\frac{1}{n}\Psi_S^T V\right)\right)
- \frac{1}{\lambda_n} \left(\frac{1}{n} \Psi_j^T V\right).
\end{eqnarray}
We then obtain the bound
\begin{eqnarray}
\nonumber
\|\mu_j\| &\leq& \left\|\Sigma_{jS}\Sigma_{SS}^{-1}\right\| \left(\|\hat
  g_S\| + \frac{1}{\lambda_n} \left\|\onen \Psi_S^T
    V\right\|\right)
+ \frac{1}{\lambda_n} \left\|\onen\Psi_j^T V\right\| \\
&=& \left\|\Sigma_{jS}\Sigma_{SS}^{-1}\right\| \left(\sqrt{s C_{\max}}
   + \frac{1}{\lambda_n} \left\|\onen \Psi_S^T
    V\right\|\right)
+ \frac{1}{\lambda_n} \left\|\onen\Psi_j^T V\right\|.
\end{eqnarray}
By our earlier calculations, we have that
\begin{eqnarray}
\| \onen \Psi_j^T V \| &\leq& \sqrt{d_n} \, \| \onen \Psi_j^T V \|_\infty 
\;\leq\; \frac{D s}{d_n} \\
\| \onen \Psi_S^T V \| &\leq& \sqrt{sd_n} \, \| \onen \Psi_S^T V \|_\infty 
\;\leq\; \frac{D s^{3/2}}{d_n} 
\end{eqnarray}
Therefore
\begin{eqnarray}
\|\mu_j\| &\leq& 
\left\|\Sigma_{jS}\Sigma_{SS}^{-1}\right\| \left(\sqrt{s C_{\max}}
  + \frac{Ds^{3/2}}{\lambda_n d_n}\right)
+ \frac{Ds}{\lambda_n d_n} .
\end{eqnarray}
Now suppose that
\begin{eqnarray}
\label{eq:condition2}
\left\|\Sigma_{jS}\Sigma_{SS}^{-1}\right\|  &\leq&
\sqrt{\frac{C_{\min}}{C_{\max}}} \,\frac{1-\delta}{\sqrt{s}} 
  \;\; \text{for some $\delta > 0$}  \\
\frac{s}{\lambda_n d_n} &\rightarrow& 0,
\end{eqnarray}
which are conditions \eqref{eq:design3} and \eqref{eq:conda} in the statement of the theorem.
Then 
\begin{eqnarray}
\|\mu_j\| &\leq& \sqrt{C_{\min}} (1-\delta) + \frac{2Ds}{\lambda_n d_n}  
\end{eqnarray}
and in particular $\|\mu_j\|\leq \sqrt{C_{\min}}$ for sufficiently large $n$.
It therefore suffices to show that
\begin{eqnarray}
\label{eq:toshow}
\P\left( \max_{j\in S^c} \| \hat g_j - \mu_j\|_\infty > \frac{\delta}{2\sqrt{d_n}} \right)
\longrightarrow 0
\end{eqnarray}
since this implies that 
\begin{eqnarray}
\|\hat g_j \| &\leq&\|\mu_j\| + \|\hat g_j - \mu_j\| \\
  &\leq& \|\mu_j\| + \sqrt{d_n} \|\hat g_j - \mu_j\|_\infty \\
  &\leq& \sqrt{C_{\min}} (1-\delta) + \frac{\delta}{2} + o(1)
\end{eqnarray}
with probability approaching one.  %\marginpar{\hazard}
To show \eqref{eq:toshow}, we again appeal
to Gaussian comparison results.
Define
\begin{equation}
Z_j \;=\; \Psi_j^T\left(I - \Psi_S (\Psi_S^T \Psi_S)^{-1} \Psi_S^T\right) \frac{W}{n}
\end{equation}
for $j\in S^c$.  Then $Z_j$ are zero mean Gaussian random variables, and we need to show that 
\begin{eqnarray}
\P\left(\max_{j\in S^c} \frac{\|Z_j\|_\infty}{\lambda_n} \geq \frac{\delta}{2\sqrt{d_n}}\right)
\longrightarrow \infty
\end{eqnarray}
A calculation shows that $\E(Z_{jk}^2) \leq \sigma^2/n$.  %\marginpar{\hazard}
 Therefore, we have
by Markov's inequality and Gaussian comparison that
\begin{eqnarray}
\nonumber
\P\left(\max_{j\in S^c} \frac{\|Z_j\|_\infty}{\lambda_n} \geq \frac{\delta}{2\sqrt{d_n}}\right)
&\leq&
\frac{2\sqrt{d_n}}{\delta\lambda_n} \, \E\left(\max_{jk} |Z_{jk}|\right) \\
&\leq& \frac{2\sqrt{d_n}}{\delta\lambda_n} \left( 3 \sqrt{\log((p-s)d_n)} \max_{jk}
  \sqrt{\E\left( Z_{jk}^2\right)}\right) \\
&\leq& \frac{6\sigma}{\delta\lambda_n} \sqrt{\frac{d_n \log((p-s) d_n)}{n}}
\end{eqnarray}
which converges to zero under the condition that
\begin{eqnarray}
\frac{\lambda_n^2 n}{d_n \log ((p-s)d_n)} \longrightarrow \infty.
\end{eqnarray}
This is condition \eqref{eq:condb} in the statement of the theorem.
\end{proof}

\begin{proof}[{\it Proof of} Theorem \ref{thm::persist}]
We begin with some notation.
If
$\M$ is a class of functions
then the $L_\infty$ bracketing number
$N_{[\,]}(\epsilon,\M)$ is defined as the smallest number of pairs
$B=\{(\ell_1,u_1), \ldots, (\ell_k,u_k)\}$ such that
$\left\|u_j-\ell_j\right\|_\infty \leq \epsilon$,
$1\leq j \leq k$, and 
such that
for every $m\in {\cal M}$ there exists
$(\ell,u)\in B$ such that
$\ell \leq m \leq u$.
For the Sobolev space ${\cal T}_j$,
\begin{equation}\label{eq::ent-of-T}
\log N_{[\,]}(\epsilon, {\cal T}_j) \leq K \left(\frac{1}{\epsilon}\right)^{1/2}
\end{equation}
for some $K>0$.
The bracketing integral is defined to be
\begin{equation}
J_{[\,]}(\delta,\M) =\int_0^\delta \sqrt{\log N_{[\,]}(u,\M)} du.
\end{equation}
\relax From Corollary 19.35 of \cite{vaar:1998},
\begin{equation}\label{eq::emp}
\mathbb{E}
\Biggl(\sup_{g\in {\cal M}}  |\hat\mu(g) - \mu(g)|\Biggr) \leq
\frac{C\, J_{[\,]}(\left\|F\right\|_\infty,\M)}{\sqrt{n}}
\end{equation}
for some $C>0$, where
$F(x) = \sup_{g\in {\cal M}} |g(x)|$,
$\mu(g) = \mathbb{E}(g(X))$ and
$\hat\mu(g) = n^{-1}\sum_{i=1}^n g(X_i)$.

Set $Z\equiv (Z_0,\ldots, Z_p)=(Y,X_1,\ldots, X_p)$ and
note that
\begin{equation}
R(\beta,g) = \sum_{j=0}^p\sum_{k=0}^p \beta_j\beta_k\mathbb{E}( g_j (Z_j) g_k(Z_k) )
\end{equation}
where we define $g_0(z_0) = z_0$
and $\beta_0 = -1$.
Also define
\begin{equation}
\hat{R}(\beta,g) = 
\frac{1}{n}\sum_{i=1}^n \sum_{j=0}^p\sum_{k=0}^p
 \beta_j\beta_k g_j (Z_{ij}) g_k(Z_{ik}).
\end{equation}
Hence $\hat{m}_n$
is the minimizer
of 
$\hat{R}(\beta,g)$ subject to the constraint
$\sum_j\beta_jg_j(x_j) \in {\cal M}_n(L_n)$
and $g_j\in {\cal T}_j$.
For all $(\beta,g)$,
\begin{equation}
|\hat{R}(\beta,g) - R(\beta,g)| \leq \left\|\beta\right\|_1^2 \ 
\max_{jk}  \sup_{g_j\in \S_j,g_k\in \S_k}|\hat\mu_{jk}(g) - \mu_{jk}(g)|
\end{equation}
where
$\hat\mu_{jk}(g) = n^{-1}\sum_{i=1}^n \sum_{jk} g_j (Z_{ij}) g_k(Z_{ik})$
and
$\mu_{jk}(g) = \mathbb{E}( g_j (Z_j) g_k(Z_k) )$.
\relax From (\ref{eq::ent-of-T})
it follows that
\begin{equation}
\log N_{[\,]}(\epsilon, {\cal M}_n) \leq 2\log p_n + K \left(\frac{1}{\epsilon}\right)^{1/2}.
\end{equation}
Hence,
$J_{[\,]}(C,\M_n) = O(\sqrt{\log p_n})$ and
it follows from (\ref{eq::emp}) and Markov's inequality
that
\begin{equation}
\max_{jk}\sup_{g_j\in \S_j,g_k\in \S_k}
 | \hat\mu_{jk}(g) - \mu_{jk}(g)| =
O_P\left(\sqrt{\frac{\log p_n}{n}}\right) = 
O_P\left(\frac{1}{n^{(1-\xi)/2}}\right).
\end{equation}
We conclude that
\begin{equation}
\sup_{g\in {\cal M}}|\hat{R}(g) - R(g)| = O_P\left(\frac{L_n^2}{n^{(1-\xi)/2}}\right).
\end{equation}
Therefore,
\begin{eqnarray*}
R(m^*) &\leq &
R(\hat{m}_n) 
\leq  \hat{R}(\hat{m}_n) + O_P\left(\frac{L_n^2}{n^{(1-\xi)/2}}\right)\\
& \leq &
\hat{R}(m^*) + O_P\left(\frac{L_n^2}{n^{(1-\xi)/2}}\right) \leq
R(m^*) + O_P\left(\frac{L_n^2}{n^{(1-\xi)/2}}\right)
\end{eqnarray*}
and the conclusion follows.
\end{proof}

\section{Acknowledgements}
\label{sec:acks}

This research was supported in part by NSF grant CCF-0625879 and
a Siebel Scholarship to PR.  A preliminary version of part of this work appears in
\cite{Ravikumar:07}.

%\bibliography{revision}
\bibliography{spam}

\end{document}